





\font\rm=cmr10 \rm

\font\bf=cmb10
\font\Rm=cmr9 at 11pt
\rm
\font\it=cmsl9 at 10pt
\font\sc=cmr7 
 at 7pt
\def\Sc #1{{\sc \uppercase{#1}}}
\font\Rrm=cmr17 at 16pt
   \font\Rm=cmr12 at 11.5pt

\long\def\Pf{\par\noindent {\it Proof.} }
\def\({\left(}
\def\){\right)}
\def\st{such that }
\def\qed{\hfill$\bullet$\vskip 4pt}
\def\quotes#1{{\lq\lq #1\rq\rq}}
\def\brcs#1{\left\{ #1\right\}}

\def\iso{\cong}
\def\wrt{with respect to }
\def\:{\,:}
\def\sp{sp\,}

\def\rk{\text{rank\,}}
\def\det{\text{det\,}}

\def\ker{\text{ker\,}}
\def\lker{\text{{$\ell$}\/ker\,}}
\def\Aut{\text{Aut}\,}

\def\sp{{\text{spec\,}}}

\def\J{{\Cal J}}

\def\C{\text{\bf C}}

\def\M{{\text M}}

\def\Im{{\text{Im}\,}}
\def\Exp{\text{Exp\,}}
\def\I{\text{I\,}}

\def\R{\text{\bf R}}
\def\N{\text{\bf N}}
\def\Z{\text{\bf Z}}
\def\Q{\text{\bf Q}}

\def\Arrow #1;#2.{#1\:#2 \to }

\def\Set#1#2{\brcs{#1 \left|\vphantom{#1 #2} \right.#2}}

\def\Oh#1{{\pmb O}\(#1\)}

\def\oh#1{{\pmb o}\(#1\)}

\def\Rrr#1,#2{{\Cal J}_{#1,#2}}
\def\slfrac#1#2{{\raise -.07 ex\hbox{$^{#1}$}}\!/\raise .35 ex \hbox{${}_{#2}$}}
\def\ssf #1/#2{\slfrac {#1}{#2}}

\def\pd #1,#2.{\frac {\partial #1}{\partial #2}}

   \long\def\Lem
#1.#2\par{\vskip4pt{\baselineskip=13pt\font\it=cmsl12 at
11.5pt\Rm
   \noindent {\rm \uppercase{#1}} #2\vskip3pt

   }} 

\long\def\Proclaim #1.#2 \endproclaim{\vskip4pt{\baselineskip=13pt\font\it=cmsl12 at
11.5pt\Rm
   \noindent {\rm \uppercase{#1}} #2\vskip3pt

   }} 

\long\def\remark #1\endremark{\vskip 2pt \noindent {\it Remark\/} #1\par}

\long\def\Sectionhead #1.#2:\par #3{\vskip 4pt \noindent {\bf #1 #2}vskip 2pt\noindent\nospace #3}

\long\def\Title #1\par {\noindent{\Rrm #1}\vskip 9pt}

 \long\def\SubT #1.{\noindent {\it #1\/} } 
 
 \long\def\SecT
#1\par{\vskip 3pt \noindent {\bf #1}\vglue1pt
   \noindent}

\long\def\subtitle #1.{\vskip 2pt \noindent {\it #1}}

\long\def\Rmk#1\par{\vskip 1pt \noindent {\it
Remark.} #1\vskip2pt}

\long\def\Abstract #1\par{{\leftskip= 3 true cm \rightskip = 3 true cm \font\it=cmsl10 \font\rm=cmr10 \baselineskip = 10pt
\parindent=.35 true cm\rm\noindent 
{\it Abstract} #1\vskip 8pt

}}

\long\def\Author #1 \par{\noindent{\it #1}}

\input amstex
 
\scrollmode\NoBlackBoxes
\magnification=1100

\font\rm=cmr10 \rm

\font\bf=cmb10
\font\Rm=cmr9 at 11pt
\rm
\font\it=cmsl9 at 10pt
\font\sc=cmr7 at 7pt 
 
\def\Sc #1{{\sc \uppercase{#1}}}
\font\Rrm=cmr17 at 16pt
   \font\Rm=cmr12 at 11.5pt

\def\Sc #1{{\sc \uppercase{#1}}}

\long\def\Pf{\par\noindent {\it Proof.} }
\def\({\left(}
\def\){\right)}
\def\st{such that }
\def\qed{\hfill$\bullet$\vskip 4pt}
\def\quotes#1{{\lq\lq #1\rq\rq}}
\def\brcs#1{\left\{ #1\right\}}
\def\Set#1#2{\brcs{#1 \left|\vphantom{#1 #2} \right.#2}}

\def\C{\text{\bf C}}

\def\wrt{with respect to }
\long\def\Lemma #1. #2\par{\noindent {\Sc  {#1.}} {\Rm
#2}\vskip 2pt}
\def\Arrow #1;#2.{#1\:#2 \to }

\def\Oh#1{{\pmb O}\(#1\)}
\def\oh#1{{\pmb o}\(#1\)}
\def\R{\text{\bf R}}
\def\N{\text{\bf N}}
\def\Z{\text{\bf Z}}
\def\Q{\text{\bf Q}}
 
\def\slfrac#1#2{{\raise -.07 ex\hbox{$^{#1}$}}\!/\raise .35 ex \hbox{${}_{#2}$}}
\def\ssf #1/#2{\slfrac {#1}{#2}}


\font\Rrm=cmr17 at 16pt
   \font\Rm=cmr12 at 11.5pt
   \long\def\Lem
#1.#2\par{\write1{#1,
p\,\folio\par}\vskip4pt{\baselineskip=13pt\font\it=cmsl12 \Rm
   \noindent {\rm \uppercase{#1}} #2\vskip3pt

   }} 

   \long\def\Title #1\par {\noindent{\Rrm #1}\vskip 9pt}
   \long\def\SubT #1.{\noindent {\it #1\/} }
   \long\def\SecT #1\par{\vskip 4pt \noindent {\bf #1}\vglue1pt
   \noindent}


\def\oneone{1.1}
\def\onetwo{1.2}
\def\onethr{1.3}
\def\onefou{1.5}
\def\onefiv{1.6}
\def\onesix{1.4}
\def\onesev{1.6}

\def\onenin{1.5}
\def\oneten{1.9}
\def\oneele{1.10}
\def\onetwe{1.11}

\def\twoone{2.1}
\def\twotwo{2.2}
\def\twothr{2.3}

\def\throne{3.1}
\def\thrtwo{3.2}
\def\thrthr{3.3}
\def\thrfou{3.4}
\def\thrsix{1.4}

\def\fouone{4.1}
\def\foutwo{4.2}
\def\fouthr{4.3}
\def\foufou{4.4}

\def\fivone{5.1}

\def\sixone{6.1}

\def\sevone{7.1}
\def\sevtwo{7.2}
\def\sevthr{7.3}
\def\sevfou{7.4}

\input diagrams

\long\def\Keywords #1 \par{{}\plainfootnote{}{\noindent Keywords \& phrases: #1}}
\long\def\amsmos #1 \par{{}\plainfootnote{}{AMS (MOS) 2010 classification: #1}}
\Keywords dimension group, trace, Hermite equivalence, Smith normal form, permutation, counting formulas, natural density, Dirichlet convolution, completely monotone, totient

\amsmos  06F20 20B25 13P15 13C05 11R45 15A21 46L80 52A07 46A55 52B20


\def\flo #1{\lfloor #1 \rfloor}
\def\PP{{\Cal P}}
\def\diag{\text{{\rm diag}}\,}

\def\Q{\text{\bf Q}}

\def\paren #1{\/{\rm(}#1\/{\rm)}}
\def\TO{\text{TO\,}}
\def\gl{\text{GL}}

\let\Ct=\ct %
\def\Mn #1{\text{M}_{#1}}
\def\NS{{\Cal N\!\Cal S}}
\def\lcm{\text{lcm}}
\def\op{^{\text{op}}}
\def\Lt #1{[\![ #1 ]\!]}
\let\lt=\Lt %
\let\List=\Lt

\def\Exp{\text{Exp\,}}
\def\Leg#1,#2.{\(\frac{#1}{{#2}}\)}
\def\Ip #1,#2.{\langle\!\langle #1,#2\rangle\! \rangle}
\def\Ipd #1,#2.{\Ip #1,#2._d}
\let\dIp=\Ipd 

\let\hat=\widehat

\let\cont=\Ct 

\def\Aff{\text{Aff\,}}
\NoBlackBoxes

\def\oneone{1.1}
\def\onetwo{1.2}
\def\onethr{1.3}
\def\onefou{1.4}
\def\onefiv{1.5}
\def\onesix{1.6}
\def\onesev{1.8}

\def\onenin{1.17}
\def\oneten{1.5}
\def\oneele{1.6}
\def\onetwe{1.14}

\def\noneone{1.7}
\def\nonetwo{1.9}
\def\nonethr{1.10}
\def\nonefou{1.11}
\def\nonefiv{1.12}
\def\nonesix{1.13}
\def\nonesev{1.15}
\def\noneeig{3.10}
\def\nonenin{3.11}
\def\noneten{3.12}

\def\twoone{2.1}
\def\twotwo{2.2}
\def\twothr{3.8}

\def\throne{3.1}
\def\thrtwo{3.2}
\def\thrthr{3.3}
\def\thrfou{3.6}
\def\thrfiv{3.4}
\def\thrsix{3.5}

\def\nthrone{3.2}
\def\nthrtwo{3.3}
\def\nthrthr{3.4}
\def\nthrfou{3.5}
\def\nthrfiv{3.6}
\def\nthrsix{3.7}
\def\nthrsev{3.8}
\def\nthreig{3.9}

\def\fouone{7.1}
\def\foutwo{7.2}
\def\fouthr{7.3}
\def\foufou{7.4}

\def\nfouone{4.1}
\def\nfoutwo{4.2}
\def\nfouthr{4.3}
\def\nfoufou{4.4}
\def\nfoufiv{4.6}
\def\nfousix{4.7}
\def\nfousev{4.8}
\def\nfoueig{4.5}

\def\fivone{8.1}

\def\nfivone{5.1}
\def\nfivtwo{5.2}
\def\nfivthr{5.3}

\def\sixone{9.1}
\def\sixtwo{9.2}
\def\sixthr{9.3}

\def\sevone{10.1}
\def\sevtwo{10.2}
\def\sevthr{10.3}
\def\sevfou{10.4}

\def\eigone{11.1}
\def\eigtwo{11.2}
\def\eigthr{11.3}
\def\eigfou{11.4}

\def\appone{B.1}
\def\apptwo{B.2}

\def\nappone{A.1}
\def\napptwo{A.2}
\def\nappthr{A.3}
\def\nappfou{A.4}
\def\nappfiv{A.5}
\def\nappsix{A.6}
\def\nappsev{A.7}

\def\appbone{D.5}
\def\appbtwo{D.6}
\def\appbthr{D.7}
\def\appbfou{D.8}
\def\appbfiv{D.9}
\def\appbsix{D.10}
\def\appbsev{D.11}
\def\appbeig{D.12}
\def\appbnin{C.7}

\def\appcone{C.1}
\def\appctwo{C.2}
\def\appcthr{C.3}
\def\appcfou{C.4}
\def\appcfiv{C.5}
\def\appcsix{C.6}
\def\appcsev{D.1}
\def\appceig{D.2}
\def\appcnin{D.3}
\def\appcten{D.4}

\def\End{\text{End\,}}

\def\I{\text{I\,}}

\def\flo #1{\lfloor #1 \rfloor}

\let\iso=\cong

\def\tripnorm #1xxx{\left\|\hglue-.2ex\left|#1\right|\hglue-.2ex\right\|}

\def\Exp{\text{Exp}\,}

\def\Spec{\text{Spec\,}}

\def\flo #1{\lfloor #1 \rfloor}
\def\PP{{\Cal P}}

\def\diag{\text{diag}\,}
\let\id=\I

\Title Invariants for permutation-Hermite equivalence and critical dimension groups%
\plainfootnote{$^0$}{\hglue -.5em {\rm This replaces an earlier version on Ar$\chi$iv, under the title {\it Invariants for critical dimension groups and permutation-Hermite equivalence.}}}

\Abstract Motivated by classification, up to order isomorphism, of   dense subgroups of Euclidean space that are free of minimal rank, we obtain apparently new invariants for an equivalence relation (intermediate between Hermite and Smith) on integer matrices. These then participate in the classification of the dense subgroups.{\par}The same equivalence relation has appeared before, in the classification of lattice simplices. We discuss this equivalence relation (called {\it permutation-Hermite\/}), obtain fairly fine invariants for it, and have density results, and some formulas counting the numbers of equivalence classes for fixed determinant.

\noindent  {\it David Handelman}%
\plainfootnote{$^1$}{\hglue -.5em Supported in part by a Discovery grant from NSERC.}

\SecT Outline

Attempts at classification of particular families of dense subgroups of $\R^n$ as partially ordered ({\it simple dimension\/}) groups lead to two directed sets  of invariants (in the form of finite sets of finite abelian groups, with maps between them) for an equivalence relation on integer matrices.  These turn  up occasionally in the study of lattice polytopes and commutative codes, among other places. The development in our case was classification of the dimension groups first, and then that of integer matrices; for expository reasons, we present the latter first.

Let $B$ and $B'$ be rectangular integer $m \times n$ matrices. We say $B$ is {\it permutation-Hermite equivalent\/} (or {\it PHermite-equivalent,} or {\it PH-equivalent\/}) to $B'$ if there exist $U \in \gl(m,\Z)$ and a permutation matrix $P$  of size $n$ \st $UBP = B'$. Classification of matrices up to PH-equivalence is the same as classification of subgroups of (a fixed copy of) $\Z^{1 \times n}$ as partially ordered subgroups of $\Z^n$ (with the inherited ordering)---the row space of $B$, $r(B)$, is the subgroup, and the order automorphisms of $\Z^{1\times n}$ are implemented by the permutation matrices (acting on the right). With this in mind, we can even define an equivalence relation on matrices $B \in \Z^{m \times n}$ and $B' \in \Z^{m' \times n}$, if we allow the additional operation of deleting a row of zeros any time it appears in the course of row reduction.

For an important subclass of matrices (suggested by the dimension group problem), we construct two families of  invariants that are surprisingly effective. For example,   the Smith normal form (SNF) is an
invariant, but a relatively crude one;  these new invariants easily distinguish matrices with the same SNF in many cases. They also yield information about the matrices themselves, for example, whether the matrix is PH-equivalent to a matrix of the form
$C:= \(\smallmatrix \I_{n-1} & a \\ 0 & d\\ \endsmallmatrix\)$---that is, an identity matrix of size $n-1$, a column $a$, and $d = |\det B|$. When this happens, the cokernel is  cyclic, but the converse fails. The latter forms are particularly amenable to complete classification for PH-equivalence.  We also construct numerous examples with the expected unusual properties.

The motivation came from classification of dense subgroups, $G$, of $\R^n$ that are free of rank $n+1$, viewed as  partially ordered abelian groups, the ordering obtained by restricting the strict ordering on $\R^n$ to $G$; that is, nonzero $g \in G$ is in the positive cone, $G^+$ iff each coordinate is a (strictly) positive real number. This defines (together with the embedding into $\R^n$, which we often suppress in notation) a {\it critical\/} (dimension) group. Equivalently, we can define a critical group to be a simple dimension group that is free of rank $n+1$, and has exactly $n$ pure traces (any affine representation, $G \to \Aff S(G,u)$, for some order unit $u$ will yield the desired dense embedding in $\R^n \iso \Aff S(G,u)$; different order units yield isomorphisms among the images).

Let $e_i$ denote the standard basis elements of $\R^n = \R^{1\times n}$. A class   of critical groups, known as  {\it basic\/} critical groups, consists of those of the form,  $G= \langle e_1, e_2, \dots, e_n; \sum \alpha_i e_i \rangle$, where $\alpha_i$ are real numbers \st that $\brcs{1, \alpha_1, \dots, \alpha_n}$ is linearly independent over the rationals (this is equivalent to density of $G $ in $\R^n$). Basic critical groups are a useful source of examples, as in [BeH]. They admit a characterization among critical  groups in terms of their structure as simple dimension groups, via the pure traces.

For each subset of the pure trace space $\Omega \subset \partial_e S(G,u)$ \st $|\Omega| = n-1$ (that is, $\Omega$ misses exactly one of the pure traces), define $\ker \Omega = \cap_{\tau \in \Omega} \ker \tau$. For any critical dimension group, the rank of $\ker \Omega$ will either be one or zero. We can thus write the kernel as $x_{\Omega}\Z$ where $x_{\Omega}$ is unique \wrt $\sigma(x_{\Omega}) \geq 0$ where $\sigma$ is the pure trace not in $\Omega$). Now form $E(G) :=\sum_{\Omega} x_{\Omega}\Z\subset G$. Then $G$ is basic iff $G/E(G) \iso \Z$; when this occurs, all sets of pure traces are ugly (in the sense of [BeH]).

However, the converse of the latter statement is not correct, but yields a larger family of critical groups. We say a critical group is {\it almost basic,} if it can be written in the form (that is, up to order isomorphism) $G = \langle f_1, f_2, \dots, f_n; (\alpha_1, \dots, \alpha_n) \rangle \subset \R^n$ where  $f_i \in \Z^n$, the set $\brcs{f_1, \dots, f_n}$ is real linearly independent, and $\brcs{1, \alpha_1, \dots, \alpha_n}$ is rationally linearly independent: these are necessary and sufficient for $G$ to be dense in $\R^n$. Then $G$ is almost basic iff the torsion-free rank of $G/E(G)$ is one, and this is equivalent to all sets of pure traces being ugly.

Of course, $G/E(G)$ itself is an invariant of order isomorphism. In the case of almost basic critical groups, we can restrict to the span of the integer rows, and in doing so, not only do we obtain an invariant for integer matrices, but the invariant boils down to PH-equivalence. Moreover, for each subset $\Omega \subset \partial_e S(G,u)$ (this time, we allow arbitrary subsets, not just those of cosize one), we may form the quotient pre-ordered abelian group $G/\ker \Omega$ (in general, the quotient of a partially ordered abelian group by a subgroup that is not an order ideal---$G$ is simple, so it has no proper order ideals---can only be pre-ordered, and does not inherit many properties from the original).

When the set $\Omega$ is ugly
(for example, if $G$ is almost basic), $G_{\Omega} = G/\ker \Omega$ is itself a critical group \wrt the real vector space $\R^{\Omega}$. Thus we can also look at $G_{\Omega}/E(G_{\Omega})$. This gives rise to an onto map from the torsion part of $G/E(G) $ to that of $G_{\Omega}/E(G_{\Omega})$. If we now assume that $G$ is almost basic, we see that the torsion  lives entirely in the integer part of the row space. This implies that it is a PH-invariant for the integer part (this requires an innocuous extra assumption on the integer part).

 Of course, we give a direct proof (avoiding dimension groups) that the resulting family of abelian groups and maps between them (the torsion parts of $G_{\Omega}/E(G_{\Omega}))$ as $\Omega$ varies over the direct set consisting of the subset of a finite set) is a PH-invariant.

The quotient maps are obtained by removing columns (those not in $\Omega$), and recalculating the invariant (or the torsion part) without using the irrational row.
This turns out to be surprisingly easy, and also leads to a second family of PH-invariants (also indexed by subsets of $\brcs{1,2,\dots,n}$), corresponding to a dual operation.


A {\it list\/} of objects is an unordered tuple (equivalently, a set with multiplicities recorded, sometimes known as a {\it multiset\/}). To distinguish between sets, ordered tuples, and lists, we use the notation $\Lt {a_1, a_2, a_3}$ for lists. (There does not appear to be a standard notation  for  this.)

\SecT Introduction

Let $G \subset \R^{1\times n}$ be a finitely generated subgroup of
$\R^{1\times n}$ (or $\R^n$ for short, if there is no ambiguity). We can
associate to $G$ a lot of matrices as follows. Pick a $\Z$-basis, $F:=
\brcs{f_1, \dots, f_m}$ for $G$, and let $B_F \in \R^{m \times n}$ be the
matrix whose $j$th row is $f_j$. Obviously, the row space of $B_F$,
$r(B_F)$, is $G$, still viewed as a subgroup of $\R^n$. We can apply any
element of $\gl(m,\Z)$ on the left to $B_F$, and the row space is
unchanged. So the inclusions $G\subset
\R^n$ are classified (merely as a subgroup of $\R^n$) by the orbits of
$\gl(m,\Z)$ (acting from the left) on $\R^{m\times n}$.

Now suppose we let $G$ inherit the usual topology from $\R^n$, and assume
that the image of $G$ is dense. Suppose $G'$ is another group with the
same properties (free of the same rank, a dense subgroup of $\R^n$, etc),
and we want to decide whether $G$ and $G'$ are isomorphic as topological
subgroups of $\R^n$. Any such isomorphism, by definition, must extend to a
continuous, hence vector space,  automorphism of $\R^n$, and these
are given by the right action of $\gl(n,\R)$. Thus the classification of
(dense) $G \subset \R^n$ up to topological isomorphism is given by orbits of
$\gl(m,\Z) \times \gl(n,\Z)$ acting on a subset of $\R^{m\times n}$
(corresponding to those matrices whose row space is dense in $\R^n$) in
the obvious manner.

Finally, suppose we also impose the strict ordering on $\R^n$, making it
into a simple dimension group, and by restriction, give a dense subgroup
$G$ the structure of a partially ordered abelian group. By [EHS], it is
also a simple dimension group, and every simple dimension group with no
infinitesimals and exactly $n$ pure traces arises in this manner. Now we
wish to determine the order-isomorphism class of such simple dimension
groups. Every order-isomorphism $G \to G'$ (both embedded in $\R^n$ as
dense subgroups and with the inherited strict ordering) will extend to an
order-automorphism of $\R^n$ [H].
The order-automorphisms of the latter are given exactly by the weighted
permutation matrices all of whose nonzero entries are positive: that is,
they factor as $\Delta P$ where $\Delta$ is a positive diagonal matrix and $P$
is a permutation matrix. Let $P(n,\R)^+$ denote the group of such weighted
permutation matrices. Here the classification of $G$ (now viewed as simple
dimension groups with ordering inherited from $\R^n$) is given by the
orbits of $\gl(m,\Z) \times P(n,\R)^+$ on the subset of $\R^{m \times n}$
consisting of the matrices whose row space is dense.

We are specifically interested in the partially ordered case, with $m =
n+1$; that is, $G$ is free of rank $n+1$, and the embedding into $\R^n$
which determines the ordering and also the topology (the ordering
determines the topology in any case) has dense image; these are called
{\it critical\ \  \paren{dimension} groups.}

This is strongly reminiscent of Hermite equivalence of (integer) matrices,
and Smith normal form. If we let $G \subset \Z^n$ (this requires $m \leq
n$), the classification of the subgroups of $\Z^n$ is just the orbit space
of $\Z^{m\times n}$
under the action of $\gl(m,\Z)$ (acting on the left), and this gives rise
to Hermite equivalence. If instead we want to classify the subgroups of
$\Z^n$ up to isomorphism as subgroups of fixed $\Z^n$, we note that the
automorphism group of $\Z^n$ is $\gl(n,\Z)$ (acting on the right), so we
are looking at the classification of matrices under the action of
$\gl(m,\Z) \times \gl(n,\Z)$; this gives rise to Smith equivalence, and the set of elementary
divisors is a complete invariant.

The analogue of the third relation arises when we view the fixed $\Z^n$ as
a partially ordered group, with the coordinatewise ordering, called {\it
simplicial\/}. Subgroups inherit the partial ordering (but are themselves
almost never simplicial), and we classify them up to order isomorphism. If
the subgroup has full rank, such an order-isomorphism to another one
(necessarily of the same rank) extends uniquely to an order isomorphism of
$\Z^n$. These are given precisely by permutation matrices. We arrive at an
equivalence relation that frequently turns up (e.g., [R, R2, ALTPP, TSCS]), but has no name. So we give it one, at least restricted to square matrices.

Two matrices $B$ and $B'$ in $\Mn n \Z$
are {\it PHermite-equivalent\/} (or {\it PH-equivalent\/} for short) if
there exist $U \in \gl (n,\Z)$ and a permutation matrix $P$ \st $UB =
B'P$. (We could of course place the $P$ to the right of $B$.)

We will see that for a large class of critical groups, the classification
problem includes within it a PH-equivalence class question. We will
develop invariants for PH-equivalence on a subclass of $\Mn n \Z$
(appropriate for critical groups), much finer than the usual elementary
divisors. We also obtain (natural) density results for matrices that have a particularly tractable equivalent form; it turns out that for $n \geq 6$, more than 80\% have this property, converging to $(\zeta(2)\zeta(3)/\zeta(6))/\zeta(2) \zeta(3)\zeta(4) \cdots \sim .845$ as $n \to \infty$ (the expression is the quotient of two moderately well-known constants, the Landau totient and $\prod_{n=2}^{\infty}\zeta(n)$.

Critical (simple dimension) groups have been a source
of interesting examples in dimension groups, e.g., [EHS], [H], and
particularly in [BeH], concerning properties of traces (good, ugly, bad).
These can be used to
characterize classes of critical dimension groups.

Let $G$ be an abelian group, free of rank $n+1$, which is embedded as a
dense subgroup of $\R^n$. This embedding imposes both a topology (the
relative one, inherited from $\R^n$), and a partial ordering, inherited
from the strict ordering on $\R^n$ (thus an element $v$ in $\R^n$ is in
the positive cone iff either $v$ is zero, or if each of its components is
strictly positive). The latter ordering makes the group into a simple
dimension group, whose pure traces are precisely the coordinate functions
(from $\R^n$). In the latter case, the ordering induces a metric, which
yields the same topology as the inherited one.

If $G$ is a simple dimension group, free of rank $n+1$, with exactly $n$
pure traces, then it is critical dimension group. These
are precisely the partially ordered groups described in the previous
paragraph, via any affine representation. If we view $G$ merely as a
topological group (free of rank $n+1$, embedded as a dense subgroup of
$\R^n$), with topology inherited from $\R^n$, we call it {\it
topologically critical.}

In the case that $n = 1$, critical subgroups of $\R$ are of the form $\Z +
r \Z \subset \R$, up to order isomorphism, and is well known that $\Z + r
\Z \iso \Z + r'\Z$ as either topological groups or ordered groups if and
only if $r$ is in the PSL$(2,\Z)$-orbit of $r'$, where PSL$(2,\Z)$ acts by
fractional linear transformations [ES]. However, the situation when $n\geq
2$ is much more complicated.

A special class of critical dimension groups, called {\it basic\/} in
[BeH], is relatively easy to classify. Let $\brcs{e_i}$ be the standard
basis of $\Z^n \subset \R^n$, and let $\alpha= (\alpha_1, \dots, \alpha_n)
\in \R^n$ be such that the set $\brcs{1, \alpha_1, \dots, \alpha_n}$ is
linearly independent over the rationals. Set $G$ to be the subgroup of
$\R^n$ generated by $\brcs{e_i}_{i=1}^n \cup \brcs{\sum \alpha_i e_i}$.
This is automatically dense in $\R^n$, and as an ordered group is
critical. We call a critical group {\it basic\/} if it is order-isomorphic
to $G$ for some choice of $\alpha$ (the rational linear independence is
necessary and sufficient for $G$ to be dense).

All critical groups of rank two (that is, $n=1$) are automatically basic,
but this fails drastically when $n > 1$, as we will see. However, if we
fix $n$, and consider classification of basic critical groups of rank
$n+1$, then the role of PSL$(2, \Z)$ is performed by the much more
elementary group, the semidirect product $\Z^n \times_{\Theta}(S_n \times
\brcs{\pm1})$
(the action of the symmetric group and $\pm 1$ is the obvious one).

Basic critical groups are easily characterized in terms of ugly sets of
pure traces, with an extra condition. This suggests a potentially larger
class of critical groups, characterized entirely in terms of ugly sets of
pure traces. These are given by the following construction. Let $A$ be a
rank $n$ subgroup of $\Z^n$, and let $G$ be the subgroup of $\R^n$
generated by $A$ and $\alpha$ (same $\alpha$ as above); this will
automatically be critical, and we call a critical {\it almost basic\/} if
it is order isomorphic to such a choice of $A$ and $\alpha$.

Almost basic critical groups admit a classification analogous to that for
basic ones, but with an additional feature; after making a preliminary
modification to $A$, the additional feature boils down to PH-equivalence.

Restricting to the relevant class of matrices $B$ (for almost basic
critical groups), we develop invariants (finer than elementary
divisors/invariant factors). These are motivated by and apply back to
almost basic critical groups, and correspond to subsets of the pure trace
space. The invariants consist of a family of finite abelian groups, which
are usually easy to calculate.

There are four appendices. The first deals with a general duality for some sets of rectangular matrices over arebitrary rings (related to the examples of section 6). The second is joint work with my colleague Damien Roy, concerning  a truncated form of the reciprocal of the Euler  function, related to the density arguments in section 7. The third shows that the obvious lower bound for the number of PH-equivalence classes of matrices with determinant $d$ is asymptotically correct, with error bounds, at least when $d$ is square-free. The fourth appendix has exact formulas for PH-equivalence classes, with special attention to those with $1$-block size $n-1$, when $ n =3$.

A subset $\brcs{g_{i}}$ of a torsion-free abelian group $A$ is {\it rationally linearly independent\/} (or {\it linearly independent over $\Q$}) if whenever $\brcs{n(i)}$ is a collection of integers with $n(i) = 0$ for all but finitely many $i$, then $\sum n(i)g_i = 0$ implies $n(i) = 0$ for all $i$. This is equivalent to the usual linear independence over the rationals of the set $\brcs{g_i}$ as a subset of the divisible hull of $A$, that is, $A \otimes_{\Z} \Q$, a vector space over the rationals.

\SecT Statement of results

Section 1 contains the definitions of {\it terminal forms\/} (based on a result, [TSCS, Theorem 4.1] on commutative codes) and the prototype  invariant(s), together with their elementary properties, and short exact sequences relating them. The second section  describes the (pseudo-)action of the permutation group $S_{n+1}$ on matrices whose $1$-block size. Section 3 introduces two families of invariants, and gives examples to show how fine these are; it also includes more short exact sequences relating them. Section 4 contains more results and conjectures for matrices PH-equivalent to a matrix with $1$-block size $n-1$. Section\,5 deals with the (rare) phenomenon of matrices PH-conjugate to their duals. And section\,6  discusses the duality conjecture, and some positive results for classes of matrices. 

Section 7 gives a density result for matrices with this last property, at least $.8$ for $n \geq 6$ and converging up to $(\zeta(2)\zeta(3)/\zeta(6))\cdot 1/(\zeta(2)\zeta(3)\zeta(4) \cdots) \sim .845$ as $ n \to \infty$.

Sections 8--11 deal with critical groups, that is, dense subgroups of $\R^n$ that are free of rank $n+1$, equipped (except in section 5) with the strict ordering, making them into simple dimension groups. Section 5 contains a topological classification theorem, which for $ n \geq 3$ corresponds to the classification of a totally ordered subgroups of $\R$. {\it Basic\/} critical dimension groups [BeH] are characterized in section 6, within the class of critical dimension groups, by means of the invariant which led to the development in sections 1--6.

Almost basic critical dimension groups are introduced in section 10, and the principal result is that the classification of these reduces to PH-equivalence of integer matrices associated to them.
When $n=1$, this is partly given by the action of PSL$(2,\Z)$; however, when $ n\geq 2$, the corresponding group is much smaller, a semi-direct product of $S_n \times \brcs{\pm1}$ acting on $\Z^{n+1}$. Section 11 is a result on almost critical basic dimension groups that amounts to showing that the whole family of PH-invariants yields their counterparts for these dimension groups.

Appendix A contains a general duality argument for natural orbit spaces, used in section\,6. Appendix B (joint with Damien Roy) is a short argument showing that the appropriate truncations of a form of the reciprocal of the Euler function yield a better than expected order of convergence. This is used in section 8. Appendix C suggests an asymptotic formula for the number of PH-equivalence classes of fixed determinant and size, and proves it when the determinant is square-free. Appendix D
 contains exact counting results on the numbers of PH-equivalence classes  of size three matrices and fixed determinants, and also the numbers of PH-equivalence classes that contain a $1$-block size two matrix.

\SecT Contents

\item{1} Permutation-Hermite equivalence; first invariants
\item{2} PH-equivalence for some terminal forms
\item{3} Finer invariants
\item{4} Size $n-1$ $1$-block terminal forms
\item{5} Dual-compatibility and dual-conjugacy
\item{6} Duality?
\item{7} Densities for PH-equivalence to $1$-block size $n-1$
\item{8} Topological isomorphism for topologically critical groups
\item{9} Basic critical dimension groups
\item{10} Isomorphisms between almost basic critical groups
\item{11} Unperforation of quotients
\item{Appendix A} General duality
\item{Appendix B} A truncated reciprocal formula (joint with D Roy)
\item{Appendix C} Counting PH-equivalence classes
\item{Appendix D} Counting PH-equivalence classes in size $3$

\SecT 1 Permutation-Hermite equivalence

Let $B$ and $C$ be $n\times n$ integer matrices ($B,C \in \Mn n \Z$). We consider two very well known, and a lesser-known, equivalence relation  between $B$ and $C$.

 The matrices $B$ and $C$ are {\it Hermite equivalent\/} if there exists $U$ in $\gl (n,\Z)$ \st $B = UC$ (this is more frequently defined on the right, rather than the left, but we will use this form here). In other words, $B$ and $C$ are obtainable from each other other by $\Z$-elementary row operations (that is, permutations, multiplication of a row by $-1$, and adding a row to another). Normal forms have been well-studied (for example, see the Wikipedia article, {\it http\/{\rm:}//en.wikipedia.org/wiki/Hermite\_normal\_form\/}).

Matrices   $B$ and $C$ are {\it Smith  equivalent\/} if there exist $U$ and $V$ in $\gl (n,\Z)$ \st $B = UCV$.  Normal forms are even more well known, and correspond to invariant factors; they are used to classify finite abelian groups.

Matrices  $B$ and $C$ are {\it permutation Hermite-equivalent\/} (or {\it PHermite-equivalent\/} or {\it PH-equivalent\/}) if there exists $U$ in $\gl (n,\Z)$ and a permutation matrix $P$ \st $B = UCP$. In order words, $B$ and $C$ are obtainable from each other other by $\Z$-elementary row operations (that is, permutations, multiplication of a row by $-1$, and adding a row to another), together with column permutations.

PH-equivalence classifies subgroups of a fixed copy of $\Z^n$ up to order-automorphism of the latter (when equipped with the simplicial, that is, coordinatewise, ordering); to see this, given the square matrix $B$, let $r(B)$ denote its row space, viewed as a subgroup of $\Z^n$. Left multiplication by elements of $\gl(n,\Z)$ has no effect on the row space---only the generating set for $r(B) $ is changed---and column permutations implement the order-automorphisms of $\Z^{1 \times n}$ when the latter is given the usual coordinatewise partial ordering. It is helpful to permit the matrices $B$ to be $m\times n$ with $m \geq n$; then elementary row operations are now implemented by elements of $\gl(m,\Z)$. These do not change the row space, and it useful to add another  operation: if at some point during a sequence of row and allowed column operations, a row becomes identically zero, then we delete it (and thus reduce the size). This obviously has no effect on the row space, and will be useful in the development of our invariants.

This section deals with an initial pair of invariants (one involving the dual of a matrix) and some of their properties.

Reduced forms for PH-equivalence have been obtained ([TSCS]; a special case is quoted as Theorem\, \oneone\ below), but normal forms have not, as far as I could tell. (Informally,  {\it reduced forms\/} for an equivalence relation constitute a  useful collection of elements which contains representatives of each equivalence class;  {\it normal forms\/} constitute a collection  containing exactly one representative of each class.)

We say a sequence, vector, list, or set of integers, $v$, has {\it content\/} $c$, denoted $\Ct (v) = c$, if $c$ is the greatest common divisor of the nonzero entries of $v$  (and if all the entries are zero, then $\Ct (v) = 0$). We say $v$ is {\it unimodular\/} (not to be confused with {\it unimodal\/}) if $\Ct (v) =1$.

We will restrict ourselves to the following class of matrices in $\Mn n \Z$. Define $B \in \Mn n \Z$ to be {\it weakly nonsingular\/} if the following two conditions apply:
\item{(a)} $\rk B = n$
\item{(b)} every column of $B$ is unimodular.

If $C$ is any element of $\Mn n \Z$ with full rank, then there is a factorization $C = B D$ where $B$ is weakly nonsingular and $D$ is diagonal with positive integer entries thereon.

Let $\NS_n$ (or simply $\NS$ when $n$ is understood) denote the collection of weakly nonsingular   $n \times n$ (integer) matrices. If $U \in \gl (n,\Z)$ and $w$ is any member of $\Z^{1 \times n}$, then $\Ct (Uw) = \Ct (w)$. Permutation of the columns of matrix simply permutes the contents of the columns. It follows that $\NS$ is preserved under PHermite-equivalence.

Given $B \in \NS$, there is a pseudo-algorithm that can be applied to reduce it to a more tractable form. First, apply the usual algorithm to obtain a Hermite normal form: since the content of the first column is one, there exists $U_1 \in \gl(n,\Z)$ \st the first column of $U_1B$ is $e_1 = (1,0, \dots, 0)^T$, the first standard basis element. Delete the first row and column, so that the second column has content possibly exceeding one (it cannot be zero, since the matrix has full rank), and continue in the obvious way, obtaining an upper triangular matrix whose first diagonal entry is $1$, and for which the other diagonal entries are positive integers.

Permute the rows and columns so that all the diagonal ones are grouped together, in a block (it is easy to see how to do this), and now the matrix is in the form
$$\( \matrix \I_s & Y \\ 0 & \Cal D\\
\endmatrix \),
$$
where $\Cal D$ is an upper triangular  matrix of size $n-s$, whose diagonal entries all exceed one. If, in $\Cal D$, the content of any column is one, we may apply the same process to it via row operations, creating an additional standard basis vector via operations on the rows of size $n-s$. By permuting rows and columns, we may enlarge the identity block, and we continue this until there are no more columns of the resulting lower block matrix that are unimodular. (Recall however, that  at every stage of this process, the size $n$ matrix has all of its columns unimodular.) This yields the Hermite normal form; further processing may be  required.

A PH-reduced form is obtained in the following result of [TSCS], for convenience stated here only for full rank matrices.

\Lem Theorem \oneone. [TSCS, Theorem 4.1] Let $B \in \Mn n \Z$ be of full rank. Then there exists a PH-equivalent  upper triangular matrix $C \in \Mn n \Z^+$, \st
\item{(a)} $0 <   C_{ii} \leq C_{i+1,i+1}$ for all $1\leq i < n$;
\item{(b)} $0 \leq C_{i,j}  < C_{jj}$ for all $i < j$;
\item{(c)} $C_{ii} \leq \gcd\Set{C_{sj}}{ i \leq s \leq j}$ for all $ i < j$.

We say $C$ is {\it PH-terminal\/} (or just {\it terminal\/}) if it is in the form described in the theorem. {\it Terminal\/} suggests that there is nothing more that can be done to such matrices to simplify them. The size of the identity matrix that appears in the terminal form is called its $1$-block size. If $B \in \NS_n$, then it has at $1$-block size at least one.

This is described in the cited reference as a normal form, but this is not the usual use of the term---two distinct  matrices $C$ and $C'$ each satisfying the conditions can be PHermite-equivalent. As a trivial example from $\NS$, set
$$
C = \( \matrix 1 & 0 & 1 \\
0 & 1 & 2 \\
0 & 0 & 6
\endmatrix\)  \qquad C' =  \( \matrix 1 & 0 & 2 \\
0 & 1 & 1 \\
0 & 0 & 6
\endmatrix\).
$$
Then $C$ and $C'$ are  conjugate via the transposition $\(\smallmatrix  0 & 1 \\ 1 & 0 \\  \endsmallmatrix \) \oplus (1)$, hence are PHermite equivalent. This type of phenomenon can  be avoided by refining the invariant. For example, we can make the top of the first column to the right of identity block be increasing; if there are ties, we can go to the next truncated column, and break the ties, etc.   However, there is a  less trivial difficulty with terminal matrices.

Applied to an  $\NS$ matrix, the terminal form has an identity block of some size in the upper left corner. If two terminal forms are PH-equivalent, it is natural to ask whether the sizes of the identity blocks are the same. The answer is  no, and we will see that this phenomenon occurs fairly frequently, almost generically (Proposition \twothr).
The equation,
$$
\(\matrix  2 & -1 & -1 \\
3 & -1 & -2 \\
6 & -3 & -4 \\\endmatrix \)
\(\matrix  1 & 1 & 2 \\
0 & 2 & 0  \\
0 & 0 & 3 \\
\endmatrix \) =
\(\matrix  1 & 0 & 2 \\
0 & 1 & 3 \\
0 & 0 & 6 \\
\endmatrix \)
\(\matrix  0 & 0 & 1 \\
0 & 1 & 0  \\
1 & 0 & 0 \\
\endmatrix \) 
$$
is of the form $UC = C' P$ where $C$ and $C'$ belong to $\NS_3$, are in terminal form, $C$ has just one $1$ on the diagonal, $C'$ has two; each has determinant $6$ and $\det P = -1$, so $|\det U| = 1$, and thus $U \in \gl (3,\Z)$. So $C$ and $C'$
 are PH-equivalent but with different block sizes for $1$.

Hermite normal forms of matrices in $\NS_n$, while themselves in $\NS_n$, need not be terminal.

We define $\NS_{n,m}$ to be the class of matrices $B \in \NS_n$ which have a terminal form with $1$-block size at least $m$. Obviously, $\NS_{n,n} = \gl (n,\Z)$, and from the definition, $\NS_{n,1} = \NS_n$. The most important of these classes is $\NS_{n,n-1}$.

First, we give a simple example to distinguish the three equivalence relations. Barely any calculation is required.

For each of $i = 0,1, 2, 3, 4$, set
$$
B_i = \(\matrix 1 & i \\ 0 & 5 \\ \endmatrix \).
$$
Then
\item{(i)} If $i\neq 0$,   $B_i$ is in $\NS$ and is in terminal form.
\item{(ii)} Every $2 \times 2$ matrix  with invariant factors $\brcs{1,5}$ is Hermite-equivalent to one of the $B_i$.
\item{(iii)} all five are mutually Hermite-{\it in\/}equivalent.
\item{(iv)} $B_2$ and $B_3$ are PH-equivalent, but there are no other PH-equivalences among these matrices.
\item{(v)} all five are mutually Smith equivalent, that is, their set of invariant factors is $\brcs{1,5}$.\vskip 4pt

An obvious invariant for PH-equivalence of matrices  $B \in \NS_n$ is simply the cokernel, $J(B) = \Z^{1\times n}/\Z^{1 \times n}B$, the Smith invariant. We will often abbreviate this $\Z^{n}/\Z^{n}B$, or $\Z^{ n}/r(B)$ (so that $r(B)$ denotes the subgroup generated by the rows of $B$). This is a  very coarse invariant. 

A second invariant arises from the dual. Let $B \in \NS$ (it need not be in terminal form); label its rows $f_i$. For each $i$, define $x_i$ to be the unique row in $\Z^{1 \times n}$ with the following properties:
\item{(a)} $x_i = m(i)E_i$ where $E_i$ is the $i$th standard basis element of $\Z^{1 \times n}$ and $m(i)$    is a positive integer;
\item{(b)} $x_i \in \sum f_j \Z$;
\item{(c)} whenever $y \in \sum f_j \Z$ and $y = kE_i$ for some $k \in \Z$, then $m(i)$ divides $ |k|$.

To see that each $x_i$ exists, note that $r(B) = \sum f_j \Z \subseteq \Z^{1\times n}$ is just the row space of $B$, hence is of rank $n$, so it hits every nonzero cyclic subgroup of $\Z^{1\times n}$ in a nonzero element; then the usual well-ordering argument works.

Now form $X(B) = \sum x_i \Z  = \oplus x_i\Z$. Then $I(B) = r(B)/ X(B)$ is a finite abelian group (since the rank of $X(B)$ is obviously $n$). The claim is that this is an invariant for PHermite equivalence between matrices in $\NS$.

To see that it really is a PH-invariant (for matrices in $\NS$), suppose that $C$ is another member of $\NS_n$, and $UCP = B$ where $U \in \gl(n,\Z)$ and $P$ is a permutation matrix. The row space of $B$ is unaffected by the left action of $\gl(n,\Z)$, and the {\it list\/} $\Lt{x_i}$ is similarly unaffected by permutation of the columns.

It would be useless if we couldn't compute with it, but it turns out to be rather easy to deal with.

Unless inconvenient, we write $\Z_k$ (for $k$ a positive integer) in place of $\Z/k\Z$. This is not going to cause confusion with the other meaning of $\Z_k$, the $k$-adic completion, as we never use the latter.

The following will be subsumed by more easily obtained results after we have an equivalent form of the construction of $I(B)$.

\Lem Lemma \onetwo. Let $n, d_i, z_i, d >1$  ($i = 2, \dots n$) be positive integers and let  $a_i$  ($i = 1, \dots, n-1$) be nonnegative integers with $a_i < d$ and $\gcd\brcs{d_i, z_i} = 1$. Suppose $B$ and $B'$ are the following $n \times n$ matrices:
$$
B = \(\matrix  1 & 0 & 0 & \dots & 0 & a_1 \\
0 & 1 & 0 & \dots & 0 & a_2 \\
0 & 0 &\ddots &&&\vdots \\
0 & 0 &0 &\dots &1&a_{n-1} \\
0 & 0 &0 &\dots &0&d \\
\endmatrix \)\qquad
B' =\(\matrix  1 & z_2 & z_3 & \dots & z_{n-1} & z_n \\
0 & d_2 & 0 & \dots & 0 & 0 \\
0 & 0 &\ddots &&&\vdots \\
0 & 0 &0 &\dots &d_{n-1}&0 \\
0 & 0 &0 &\dots &0&d_n \\
\endmatrix \)
$$
Then both $B$ and $B'$ are $\NS$ matrices in terminal form. Moreover,
$$
I(B) \iso \bigoplus_i \(\Z\Bigg/\(\frac d{\gcd\brcs{d, a_i}}\)\Z\),
$$
and $I(B')$ is cyclic of order
$$
\lcm \brcs{d_2, d_3,
\dots, d_n}.
$$


\Pf That the matrices have all their columns unimodular is an immediate consequence of the properties ascribed to the coefficients.  Let $f_j$ ($j=1, \dots, n$) be the rows of $B$. Then for $i < n$, $E_i = f_i - (a_i/d) f_n$, so that $x_i = (d/\gcd\brcs{d,a_i})f_i - (a_i/\gcd\brcs{d,a_i})f_n$. In addition, $x_n = f_n$, so that a basis for $X(B)$ is $ \brcs{(d/\gcd\brcs{d,a_i})f_i} \cup \brcs{f_n}$. As $\brcs{f_1, \dots, f_n}$ is a basis for $r(B)$, we have that $I(B) \iso  \bigoplus_i \(\Z \big/\(\frac d{\gcd\brcs{a_i,d}}\)\Z\)$.

Now let $f_j$ be the $j$ row of $B'$, and let $l   = \lcm\brcs{d_i}$. Then
$$\eqalign{
E_1& = f_1 - \sum_{i \geq 2}\frac{z_i}{d_i} f_i\cr
lE_1& = lf_1 - \sum_{i \geq 2}z_i f_i \cr
}$$
If $t> 1$ is a prime dividing $l$ and all of the $z_i$, then it divides at least one of the $d_j$; but this would contradict $\gcd\brcs{d,z_i } = 1$ for all $i$. Hence $lE_1$ is a unimodular element of $\sum f_j \Z$, so that $x_1 = l E_1$. For $i > 2$, $x_i = f_i$. Hence a basis for $\sum_{i=1}^n x_i \Z$ is $\brcs{lf_1, f_2, \dots, f_n}$, and thus $I(B') $ is cyclic of order $l$.\qed

Here are some very simple examples with $n =2$. Define
$$
B_{a,d} = \(\matrix 1 & a\\ 0 & d\\ \endmatrix \)
$$
where $d > 1$; in order to be terminal, we need $\gcd\brcs{a,d} = 1$ and $1 \leq a < d$. By taking determinants, we see that $B_{a,d} $ PH-equivalent to $B_{a',d'}$  entails $d = d'$ (a peculiarity of the $n=2$ case). So let $a'$ be another integer in the interval $1 \leq a' < d$ relatively prime to $d$. Then  $B_{a,d}$ is PH equivalent to $B_{a',d}$ if and only either $a = a'$ or $aa' \equiv 1 \mod d$ (that is, in $\Z/d\Z$, $[a] = [a]^{\pm1}$). The second choice comes from letting $P$ be the nontrivial permutation matrix, and working out the details. Here $I(B_{a,d}) \iso \Z/d\Z$, not very exciting.

Next, consider variations on the earlier example. Set
$$
B = \( \matrix 1 & 0 & b \\ 0 & 1 & c \\  0 & 0 & 6\\
\endmatrix\), \qquad C_1 = \( \matrix 1 &1 & 1 \\ 0 & 2 & 0 \\  0 & 0 & 3\\
\endmatrix\), \qquad C_2 = \( \matrix 1 &1 & 2 \\ 0 & 2 & 0 \\  0 & 0 & 3\\
\endmatrix\).
$$
In order for $B$ to terminal, we require $\gcd\brcs{b,c,6} = 1$ and $0 \leq b, c < 6$; we may assume $ b \leq c$ (by conjugating with the obvious transposition).
Every terminal form of an $\NS$ matrix with diagonal entries $1,2,3$ is PH-equivalent to one of $C_1$ or $C_2$; this is routine.

We will show that the only choices for $B$   which  are PH-equivalent to a terminal form whose $1$-block has size unequal to two (which means it has size one) correspond to $(b,c) = (2,3)$ and $(3,4)$. The former comes from the earlier example, and it is PH-equivalent to $C_2$. A similar computation (which comes from an easy sequence of row reductions) shows that with $(b,c) = (3,4)$ or $(4,3)$, $B$ is PH-equivalent to $C_1$.

There are no other terminal forms of size three with $2,3$ along the diagonal than $C_1$ and $C_2$, since both numbers are   prime.

We have, by the earlier result, $I(B) = \(\Z/(6/\gcd\brcs{6,b})\Z \)\oplus \(\Z/(6/\gcd \brcs{6,c})\Z\)$. Hence if at least one of $b$ or $c$ is relatively prime to $6$, then $I(B)$ is not cyclic, and has $\Z/6\Z$ as a proper quotient.

Now  $I(C_i) \iso \Z_6$ since $6 = \lcm\brcs{2,3}$. Hence if $b$ or $c$ is relatively prime to $6$, $B$ cannot be PH-equivalent to either $C_i$, and in particular, all terminal forms of $B$ have the same $1$-block size, two.

Finally $C_1$ and $C_2$ are not PH-equivalent, since the corresponding $B$ forms are not; this will come from a general result obtained later.
 \qed

\Lem Lemma \onethr. Let $B = \( \smallmatrix \I_r & X \\ 0 & \Cal D \\ \endmatrix\)$ be in terminal form  with $\Cal D$ upper triangular, and whose diagonal entries satisfy $1 < d_{r+1} \leq d_{r+2} \leq d_n$. Set $l = \lcm\brcs{d_i}$.
\item{(a)} If $\Cal D$ is diagonal, then $I(B) $ is a quotient of $(\Z_l)^r$.
\item{(b)} In general, $I(B)$ is a quotient of
$$
(\Z/l\Z)^r \oplus \( \bigoplus_{j=r+1}^{n-1} \Z/\lcm\brcs{d_{j+1}, d_{j+2},\dots , d_n} \Z \).
$$

\Pf For $1 \leq i \leq r$, $E_i = f_i - \sum_{j> 1} (a_{ij}/d_j) f_j$ for some integers $\brcs{a_{ij}}$. Hence $lE_i \in \Cal C(B)$, and thus $lE_i \in X(B)$. Hence $x_i = t_i E_i$ for some positive integer $t_i$ dividing $l$.

\noindent (a) Here $x_i = f_i$ if $i > r$, and thus $X(B) $ is spanned by $\brcs{x_i}_{i\leq r} \cup \brcs{f_i}_{i> r}$; from the form of $t_i E_i$, we have that $X(B)$ is spanned by $\brcs{t_i f_i}_{i \leq r} \cup \brcs{f_i}_{i > r}$. Since $\brcs{f_i}$ is a basis for $r(B)$, it follows that $I(B) \iso \oplus_{i\leq r} (\Z_{t_i})$. This is a quotient of $(\Z_l)^r$ since each $t_i$ divides $l$.

\noindent (b) If $r < i < n$, we can write $E_i = f_i - \sum_{j> i} (a_{ij}/d_j)f_j$. Obviously, $x_n = f_n$. Set $l_i = \lcm \brcs{d_{i+1}, d_{i+2}, \dots, n}$, so that $l_i E_i \in r(B)$ and thus is in $X(B)$. So again we can write $x_i = t_i E_i$ with $t_i$ dividing $l_i$, and we obtain $I(B) $ is a quotient of $(\Z_l)^r \oplus (\bigoplus_{i> r} \Z_{t_i})$, which is a quotient of the desired group.
\qed

The $1$-block size (that is, the size of the identity matrix in the upper left corner) in terminal forms turns out to be significant, particularly if it is $n-1$---when this occurs, PH-equivalence classes can be determined exactly.

\Lem Corollary \onefou. Suppose $B = \( \smallmatrix \I_s& X \\ 0 & \Cal D \\ \endsmallmatrix\)$
is in terminal form, and let $d = \det B$. If $I(B)$ has a quotient which is isomorphic to $(\Z_d)^s$, then all terminal forms PH-equivalent to $B$ must have $1$-block  size at least $s$.

\Pf Suppose $B' = \( \smallmatrix \I_r& X' \\ 0 & \Cal D'\\ \endsmallmatrix\)$ is a PH-equivalent terminal form with $r < s$. In particular, $\det B' = \det B = d$. All the factors that are quotients of  $\Z/\lcm\brcs{d_j', \dots, d_n'}$ for $j> r$, have order at most $\prod d'_j/d'_{r+2}  < d$. But then the preceding says that $I(B')$ has at most $r$ copies of $(\Z_d)^r$ appearing as a factor, a contradiction.\qed

\comment
\Lem Corollary \onefiv. Suppose $B = \( \smallmatrix \I_{n-1}& X \\ 0 & d \\ \endsmallmatrix\)$ is in  terminal form, and suppose each entry of the column $X = (z_i)^T$ ($1 \leq i \leq n-1$) satisfies $\gcd (z_i, d) = 1$.
Then all terminal forms PH-equivalent to $B$ have $1$-block of size $n-1$.
\Pf By Lemma \oneone, $I(B) \iso (\Z_d)^{n-1}$, so the result follows from the preceding corollary.
\qed

\Lem Lemma \onesix. Let $B \in \NS$, and suppose it has a terminal form with non-unital diagonal entries $d_{r+1}, d_{r+2}, \dots, d_n$. Then the exponent of $I(B)$ divides $\lcm \brcs{d_i}$.

\Pf As in the previous arguments, $x_i = m(i)E_i$, where $m(i) $ divides  $\lcm\brcs{d_j}$, and in particular, the coefficient of $f_i$ in $x_i$ divides $m(i)$; but the least common multiple of these coefficients is exactly the exponent of $I(B)$.\qed

\endcomment

It is convenient to introduce the notion of opposite here, in order to put
the invariant(s) in a broader context.

\noindent {\it A dual formulation of the invariant.}
When we construct the $x_i$ in order to determine $I(B)$, we also create a dual of the matrix $B$, call it $B^{op}$, also in $\NS_n$, and for which $I(B) = \Z^{1\times n}/\Z^{1\times n} B^{op}$, that is, $I(B) \iso J(B\op)$.
To see this, we have a unique representation for each $i$, $x_i = \sum_j c_{ij}f_j$ with $c_{ij} \in \Z$; since $x_i$ is not a nontrivial multiple of any element of $\sum f_i\Z$, it follows that the content of $\brcs{c_{ij}}_{j=1}^n$ is one. Hence the matrix $C = (c_{ji})$ (the transpose of what is expected) belongs to $\NS_n$.

Next, we see that if $B'$ is PH-equivalent to $B$, then $C'$ (constructed out of the canonical $x_i'$) is PH-equivalent to $C$. A row operation on $B$ simply multiplies $C^T$ on the right by an element of $\gl(n,\Z)$, hence multiplies its transpose, $C$, on the left by an element of $\gl(n,\Z)$. A column permutation applied to $B$ multiplies the representation of the $x_i$ by a row permutation of the matrix $C^T$, so induces a column permutation of $C$.

So we  call $C$, $B^{op}$.  In general, when $B$ is in terminal form, $B^{op}$ will be far from terminal, requiring both row operations and column permutations to put it into terminal form. If  we think in terms of the row space of $B^{op}$, then it is almost tautological that $I(B) = \Z^n/r(B^{op})$. That being the case, $I(B)$ is determined from the Smith normal form of $B^{op}$. To some extent this explains some of the loss of information in going from the PH-equivalence class of $B$ to $I(B)$.  
Unsurprisingly,  $(B^{op})^{op} = B$. In general, $|\det B| \neq |\det B^{op}|$; this occurs when $|\det B| \neq |I(B)|$, and we have seen an example for which $\det B = 8$, but $I(B) \iso \Z_8 \oplus \Z_2$. From the equations defining $B\op$, we have $(B\op)^T B = \Delta := \diag(m(i), \dots m(i))$, where the $m(i)$ are defined via the $x_i$, that is, $x_i = m(i)E_i$.

Because of potential confusion caused by the notation, we redefine $J(B) =I(B\op) = \Z^{1\times n}/\Z^{1\times n} B$ (determined by the Smith normal form of $B$), and thus $J(B\op) =I(B)= \Z^{1\times n}/\Z^{1\times n} B\op$. We will soon obtain a simpler description for $B\op$.

\Lem Lemma \oneten. Suppose that $B,B' \in \NS_n$ and $\Delta, \Delta'$ are diagonal real matrices with strictly positive entries. If $B \Delta = B' \Delta'$, then $\Delta = \Delta'$ and $B = B'$.

\Pf Since $B$ is invertible in $\Mn n \Q$, we have $B^{-1}B' = \Delta (\Delta')^{-1}$; thus the latter has only rational entries (all of which are nonnegative). We can therefore write $N\Delta (\Delta')^{-1} = \Delta''$ for some positive integer $N$ and $\Delta'' = \diag(d_i)$ diagonal with only positive integer diagonal entries. From $NB^{-1}B' = \Delta''$, we have $B'N= B\Delta''$. Now the $i$th column of $B'N$ has content $N$, and the $i$th column of $B\Delta$ is just  $d_i$ times the $i$th column of $B$, hence has content $d_i$. Thus $d_i = N$ for all $i$, so $(B'-B)N = 0$ and thus $B' = B$. As $B$ is invertible in $\Mn n \R$, $\Delta = \Delta'$.
\qed

The following shows that $B\op$ can be characterized via a more general equation.

\Lem Proposition \oneele. Let $B \in \NS_n$. Then
\item{(a)}  $B\op \in \NS_n$ and $(B\op)^T B = \diag(m(1), \dots, m(n))$;
\item{(b)} if $C \in \NS_n$ and $C^T B $ is diagonal with only nonnegative entries, then $C = B\op$;
\item{(c)} $(B\op)\op = B$, and the lists $\lt{m(i)}$ are the same for $B$ and $B\op$;
\item{(d)} if $B' \in \NS_n$ is PH-equivalent to $B$, then $(B')\op$ is PH-equivalent to $B\op$.

\Pf (a) is noted above.

(b) Write $(B\op)^T B = \Delta_0$. As $B \in \NS_n$, $B^{-1}$ exists (in $\Mn n \Q$),  we can write $C = (\Delta B^{-1})^T = (B^{-1})^T \Delta$, and similarly, $B\op = (B^{-1})^T \Delta_0$. Then
$$
B\op \Delta =  (B^{-1})^T \Delta_0 \Delta = (B^{-1})^T \Delta \Delta_0 = C\Delta_0.
$$
The result now follows from the preceding lemma.

(c) From $(B\op)^TB = \Delta_0$, on transposing, we obtain $B^T B\op = \Delta_0$; as $B \in \NS_n$ implies $B\op \in \NS_n$, we have $B = (B\op)\op$  from (b). It then  follows  from $B^T B\op = \Delta_0$ that the list $\lt{m(i)}$ (the list of diagonal entries of $\Delta_0$) is the same, whether computed \wrt $B$ or \wrt $B\op$.

(d) There exist $U \in \gl(n,\Z)$ and a permutation matrix $P$ \st $ B = U B' P$; then $\Delta_0 = (B\op)^T B = (B\op)^TUB'P$. Pre-multiplying by $P$ and post-multiplying by $P^{-1} $, we have $P\Delta_0 P^{-1} = P(B\op)^T U B'$. Since $B', (P(B\op)^T U)^T \in \NS_n$ and $P\Delta_0 P^{-1}$ is diagonal, by the lemma, we have $(B')\op = (P(B\op)^T U)^T = U^T B\op P^{-1}$ (since $P^{-1} = P^T$), yielding that $(B')\op$ is PH-equivalent to $B\op$.
\qed

This leads to a fast  construction of $B\op$. From the characterization of $B\op$ in \oneele(a,b), finding $B\op$ and $\Delta$ becomes relatively simple. Pick $B \in \NS_n$; form $B^{-1} \in \text{M}_n \Q$. There exists a smallest positive integer $m(i)$ \st $m(i)$ times the $i$th row of $B^{-1}$ consists of integers---and necessarily, the resulting row has content one. Set $\Delta = \diag (m(i))$; since the entries of $\Delta B^{-1}$ are all integers and the content of each row is one, it is immediate that $(\Delta B^{-1})^T \in \NS_n$. Then $B\op = (\Delta B^{-1})^T $.

In [ALTPP], the authors introduced two numbers associated to a matrix $B \in \NS_n$; the first was denoted $I$, which is $|\det B|$; the second was denoted $I^*$,\plainfootnote{$^1$}{Unfortunately I came across this reference after I had established the notation for this paper, so that their $I$ is $|I(B\op)| = |J(B)|$, and their $I^*$ is $|I(B)| = |J(B\op)|$.}%
and is $|\det B\op|$; they also use $B^*$, the dual matrix emanating from lattice polytopes, for what is called here $B\op$. Among other things, they  constructed very useful tables of numbers of isomorphism classes, and explicit generators, which turned out to be particularly helpful for Appendix D.

Let $B$ belong to $\NS_n$. Defining $E_i$, $x_i$, and $m(i)$ as we have
above, there is an obvious short exact sequence,
$$
0 \to \frac{r(B)}{\sum x_i \Z} \to \frac{\Z^{1 \times n}}{\sum x_i \Z} \to
\frac{\Z^{1\times n}}{r(B)} \to 0.
$$
The left term is just $I(B)$, the right is $I(B\op)$, which is determined
by the invariant factors of $B$. The middle term is naturally isomorphic
to $r(B)/(\sum x_i \Z)B$ via $B$; the map sending $w \in \Z^n$ to $wB$
induces a group homomorphism ${\Z^{1 \times n}}/{\sum x_i \Z}\to
r(B)/(\sum x_i \Z)B$, which is clearly onto; it is also one to one, since
$wB = vB$ (with $v \in \sum x_i \Z$) entails $w = v$. In addition, $x_i B =
m(i)f_i$, so that the middle group is just $\oplus \Z_{m(i)}$. So we can rewrite the short exact sequence,
$$
0 \to J(B\op) \to \oplus \Z_{m(i)} \to J(B) \to 0.
$$
Since we may interchange $B$ with $B\op$ (from $(B\op)^T B = \Delta$, we obtain $B^T B\op = \Delta$), we also obtain a short exact sequence $0 \to J(B) \to \oplus \Z_{m(i)} \to J(B\op) \to 0$. This can be re-interpreted more generally. 

\def\gl{\text{GL}}

For a finite group $G$, the {\it exponent\/} of $G$, denoted $\Exp G$, is the smallest positive integer \st the order of every element divides $e$.

\Lem Proposition \noneone. Suppose $B \in \NS_n$. Then $\Exp J(B) = \Exp J(B\op) = \Exp \Z^n/\Z^n\Delta = \lcm \brcs{m(i)}$.

\Pf From the short exact sequence $0 \to J(B) \to \oplus \Z_{m(i)} \to J(B\op) \to 0$, obviously $\Exp J(B)$ and $\Exp J(B\op)$ divide the exponent of the middle term, which is $\lcm \brcs{m(i)}$. Set $d= \Exp J(B)$. This says that $\Z^n d \subseteq \Z^n B$. Applying $B^{-1}$ (which exists in $\Q^{n\times n}$), we have $\Z^n dB^{-1} \subset \Z^n$, whence $C:= dB^{-1}$ is an integral matrix satisfying $CB = d\I$.  Since $(B\op)^T B = \Delta$, we deduce $(B\op)^T = d\Delta C$, so that $C^T = B\op d \Delta^{-1}$ (as matrices with rational entries). 

The $i$th column of $C^T$ is thus $d/m(i)$ times the $i$th row of $B\op$. As each column of $B\op$ has content one, this entails (as $C^T $ has only integer entries) $m(i)$ divides $d$. Hence $\lcm \brcs{m(i)}$ divides $d$; since $d$ divides $\lcm\brcs{m(i)}$, we have $d = \lcm \brcs{m(i)}$.

Since $B^T B\op = \Delta ^T = \Delta$, we can interchange the roles of $B$ and $B\op$, obtaining the final equality.
\qed

Let $d = \Exp J(B)$; then we can regard each of $J(B)$, $J(B\op)$, and $\Z^{1\times n}/\Z^{1\times n} \Delta$ as $\Z_d$-modules. As $\Z_d$ is self-injective, each of them contains a nonzero free submodule as a direct summand; and $J(\Delta):= \oplus \Z_{m(i)}$ contains a free $\Z_d$-module on two generators as a direct summand.

\Lem Lemma \onesev. Suppose that $p$ is a prime and $B$ in $NS_n$ has terminal form
$$
\( \matrix \I_{n-1} & X \\
0 & p^m \\
\endmatrix\)
$$
for some $m\geq 1$.
If $B'$ is a terminal matrix in $NS$ that is PH-equivalent to $B$, then the $1$-block of $B'$ has size $n-1$.

\Pf  If $B$ is PH-equivalent to $B'$ in terminal form with block size less than $n-1$, then the non-one diagonal entries of the latter are powers of $p$, and their product is the determinant, $p^m$. Their $\lcm$ is thus strictly less than $p^m$, and so the exponent of $I(B') \neq p^m$, a contradiction.
\qed

In particular, if $|\det B| $ is a power of a prime and $I(B)$ has exponent equalling $|\det B|$, then every terminal form PH-equivalent to $B$ must have $1$-block of size $n-1$.

The following is completely elementary, and the use of self-injectivity is like cracking a walnut with a hammer.

\Lem Lemma \nonetwo. Let $0 \to A \to B \to C \to 0$ be a short exact sequence of finite abelian groups. If any of the following holds,
\item{(a)} $\Exp B = \Exp A$ and $A$ is cyclic, \quad or
\item{(b)} $\Exp B = \Exp C$ and $C$ is cyclic, \quad or
\item{(c)} $\Exp B$ is square-free.

\noindent then the sequence splits.

\Pf Let $d = \Exp B$; as $\Exp A$ and $\Exp C$ divide $d$, the sequence is a short exact sequence of  $\Z_d$-modules. 

If $\Exp C = \Exp B$ and $C$ is cyclic, then $C$ is free as a $\Z_d$-module, so the sequence splits. If $\Exp A = \Exp B$ and $A$ is cyclic, then $A$ is free and singly generated; since $\Z_d$ is self-injective, $A$ is injective as a $\Z_d$-module, so the sequence splits.

If $\Exp B: = d$ is square-free, then $C$, being a $\Z_d$-module, is projective, hence the sequence splits.§
\qed

A consequence of the method of proof is the following somewhat interesting result.

\Lem Lemma \nonethr. Suppose that $H \subset \Z^n$ of rank $n$, with invariant factors $(f_1, f_2, \dots, f_n)$. Suppose that $x + H$ has order $f_n$ (the exponent of $\Z^n/H$) in $\Z^n /H$. Then $H + x\Z \subset \Z^n$  has invariant factors $(1,f_1, \dots, f_{n-1})$.

\Pf The onto map $\Z^n/H \to \Z^n/(H+x\Z)$ has kernel $(x\Z + H)/H $, which is free as a $\Z_{f_n}$ module, so is a direct summand. Hence we can write $\Z^n /H = (x+H)\Z \oplus D$ for some $\Z_{f_n}$ module $D$. Since $D \iso \Z_{f_n}$ and the sequence  invariant factors is unique, $D$ must have invariant factors $1,f_2, \dots, f_{n-1}$ (delete the last one, and insert a one at the beginning). Obviously $D \iso \Z^n/(H+x\Z)$.
\qed

The following is presumably well-known, but useful. If $G$ is an abelian group, then $t(G)$ denotes its torsion subgroup.

\Lem Lemma \nonefou. Let $A$ be an $r \times n$ integer matrix. Then
$$
t\( \Z^n/\Z^r A\) \iso t\( \Z^r/\Z^n A^T\).
$$

\Pf Let $s = \rk A$; then $s \leq r,n$. The first step is to reduce to the case that $r = s = n$.

To that end, we observe that the row space of $A$, $\Z^r A$ is free of rank $s$; hence there exists $E \in \gl(r,\Z)$ \st $EA \(\smallmatrix A'\\ 0\\ \endmatrix\)$, where $A'$ is $s\times n$. Since $A'$ has rank $s$, there exists $F \in \gl(n,\Z)$ \st $ A'F = \(\smallmatrix A'' \\ 0 \\ \endmatrix\)$, where $A'' $ is $s \times s$. In particular, 
$$
EAF = \( \matrix A'' & 0 \\ 0 & 0 \\ \endmatrix\).
$$
Hence
$$
\Z^n/\Z^r A \iso  \Z^n/\Z^r EAF \iso \Z^s/\Z^sA'' \oplus \Z^{n-s}.
$$

From $F^T A^T E^T = (EAF)^T =  \( \smallmatrix (A'' )^T& 0 \\ 0 & 0 \\ \endsmallmatrix\)$, we similarly obtain  $\Z^n/\Z^r A \iso   \Z^s/\Z^s(A'')^T \oplus \Z^{r-s}$.

So it suffices to show that if $M \in \Z^{s\times s}$ is of rank $s$, then $\Z^s/\Z^s M \iso \Z^s /\Z^s M^T$. But this is straightforward. Let $(f_1,\dots, f_s)$ be the sequence of invariant factors of $M$; then there exist $J,K \in \gl (s,\Z)$ \st $JMK = \diag(f_1, \dots, f_s): = \Delta$. Obviously $K^TM^T J^T = \Delta$, so $M^T$ has the identical sequence of invariant factors. 
\qed

\Lem Lemma \nonefiv. Let $A \in \Z^{r\times n}$, $B \in \Z^{n\times r}$. If $A$ has rank $r$, then there is a short exact sequence,
$$
0 \to \Z^r/\Z^n B \to \Z^n/ \Z^n BA \to \Z^n /\Z^r A \to 0,
$$
the maps induced  by $v \mapsto  vA$ and $v \mapsto v$.

\Pf Since $\Z^n BA \subset \Z^r A$, the map from middle to the right term, $v + \Z^n BA \mapsto v + \Z^r A$ is well-defined, and obviously onto. Its kernel is $\Z^r A/\Z^n BA $. The map $v + \Z^n B \mapsto vA + \Z^n BA$ is clearly well defined, and maps onto the kernel; it suffices to show it is one to one. But $vA \in \Z^n BA$ entails $v = wBA$ for some $w \in \Z^n$, whence $(v-wB)A = 0$. However, $A$ is $r \times n$ and of rank $r$, so right multiplication by $A$ is one to one. Thus $v = wB$, and the map is one to one.
\qed

Weirdly, even under the (strong) hypotheses that $A,B$ are square of the same size and with nonzero determinant, it need not be true that $\Z^n/\Z^n AB \iso \Z^n/\Z^n  BA$ (both of these are torsion). The following is presumably well-known. 

Set $A = \(\smallmatrix 1 & 1 \\ 0 & 2 \\ \endsmallmatrix\)$. Then $AA^T = \(\smallmatrix 2 & 2 \\ 2 & 4 \\ \endsmallmatrix\)$, so $\Z^2 /\Z^2 AA^T \iso \Z_2 \oplus \Z_2$. On the other hand, $A^TA = \(\smallmatrix 1 & 1 \\ 1 & 5 \\ \endsmallmatrix\)$, and thus $\Z^2/\Z^2 A^T A \iso \Z_4$. Label $B = A^T$.

In general, the sequence of torsion subgroups of a short exact sequence is not exact (take $0\to\Z \to \Z \to \Z_2 \to 0$ where the middle map is multiplication by $2$); however, in this case, it is. 

\Lem Corollary \nonesix. Suppose $A \in \Z^{r\times n}$, $B \in \Z^{n \times r}$, and $BA$ has rank $r$. Then there is a short exact sequence,
$$
0 \to \Z^r/\Z^n B  \to t\(\Z^n/\Z^r BA \) \to t\(\Z^n/\Z^r A\) \to 0.
$$

\Pf Since $B$ has rank (at least, hence exactly) $r$, $ \Z^r/\Z^n B$ is finite, hence torsion. Suppose $v + \Z^r A$ is a torsion element of the third term in the short exact sequence in the previous lemma. Then there exists a positive integer $N$ \st $Nv \in \Z^r A$. Since $\Z^r BA \subset \Z^r A$ and they have equal ranks, there exists an integer $M$ \st $M\Z^r A \subset \Z^r BA$. Hence $MN v \in \Z^r BA$. Therefore, all pre-images in the middle term, of $v + \Z^r A$, lie in the torsion subgroup of $\Z^n/\Z^r BA$. In particular, the right map is onto.

Since $\Z^r/\Z^n B$ is finite, it maps to the torsion subgroup of the middle term, and exactness of the original sequence now yields exactness of the sequence of torsion subgroups.
\qed

Now suppose that $B \in \NS_n$, and form $B\op$. From $(B\op)^T B = \Delta = \diag(m(1), \dots, m(n))$, we deduce (setting $B = A$, etc), a short exact sequence
$0 \to \Z^n/\Z^n (B\op)^T \to \oplus \Z/\Z m(i) \to \Z^n/\Z^n B \to 0$. By \nonefou, the first term is isomorphic to $\Z^n/\Z^n B\op$. This yields a short exact sequence,
$$
0 \to J(B\op) \to \oplus \Z/m(i)\Z  \to J(B) \to 0.
$$

We also have $B^T B\op = \Delta$ (by applying the transpose); this permits us to reverse the roles of $B$ and $B\op$, and we obtain another short exact sequence,
$$
0 \to J(B) \to \oplus \Z/m(i)\Z \to J(B\op) \to 0.
$$

By \noneone\ and \nonetwo, if either $J(B) $ or $J(B\op)$ is cyclic, or more generally, if either $J(B)$ or $J(B\op)$ is a free $\Z_d$-module, then the sequence splits; similarly, if the exponent, $d$, of $J(B)$ is square-free, the sequence splits. There are examples to show that neither of these need split. Unfortunately, because we are taking isomorphisms at various points, the extensions themselves need not be PH-invariants. However, splitting (and not splitting) are PH-invariants.

\Lem Example \onetwe. The extension $0 \to J(B) \to \oplus \Z_{m(i)} \to J(B\op) \to 0$ need not  split; in fact, $J(B) \oplus J(B\op)$ need not be isomorphic to $ \oplus \Z_{m(i)} $.

\noindent  To construct an example, suppose that $n = 3$ and $d $ is a power of a prime $p$. If we can find $B \in \NS_3$ \st $\det B = d$, and both $J(B)$ and $J(B\op)$ are not cyclic, then the extension cannot be split. We   note that  $J(B) \oplus J(B\op)$ cannot be generated by $3$ elements, since it is  a $p$-group and has at least four elementary divisors.  But the list $\lt{m_i}$ consists of $n=3$ elements, so that $\oplus \Z_{m(i)}$ has three generators; in particular, $\oplus \Z_{m(i)} \not \iso J(B) \oplus J(B\op)$.

So it suffices to find a matrix $B$ with these properties. For any prime $p$,  set
$$
B = \(\matrix 1 & 1 & 1 \\
0 & p & p^2 \\
0 & 0 & p^3 \\
\endmatrix \) \quad \text{so that}\quad B\op = \(\matrix p^3 & 0 & 0\\
-p^2 & p & 0 \\
p-1 & -1& 1 \\
\endmatrix \).
$$
Since the cokernel of $C:=\(\smallmatrix p & p^2 \\ 0 & p^3\\ \endsmallmatrix\)$ is $\Z_{p^3} \oplus \Z_p$, $J(B) \iso \Z_{p^3} \oplus \Z_p$ (in $\Z^2/\Z^2 C$, there is no element of order $p^4$, but there are elements of order $p^3$; alternatively, subtract $p$ times the first column from the second, to create $\diag(p,p^3)$); similarly, the cokernel of
$\(\smallmatrix p & 0 \\ -p^2 & p^3\\ \endsmallmatrix\)$ is $\Z_{p^3} \oplus \Z_p$, so this is also $J(B\op)$. Hence $J(B) \oplus J(B\op)$ has elementary divisors $\lt{p^3,p^3,p,p}$, and is thus not $3$-generated as an abelian group. In this case, $m(1) = m(3) = p^3$ and $m(2) = p^2$.
\qed

There is another invariant of PH-equivalence, concerning a particularly strong form of splitting. The imbeddings $J(B) \to \Z^{1\times n}/\Z^{1\times n} \Delta$ and $J(B\op) \to \Z^{1\times n}/\Z^{1\times n}\Delta$ (we will sloppily abbreviate $\Z^{1\times n}/\Z^{1\times n}\Delta$ to $J(\Delta)$ from now on) are given by first identifying $J(B) = \Z^{1\times n}/\Z^{1\times n}B$ with $\Z^{1\times n}/\Z^{1\times n}B^T$ (\nonefou\ and \nonefiv), and then with the latter's image in $J(\Delta)$, given by $v + r(B^T) \mapsto vB\op + \Z^{1\times n}\Delta$, and then doing the same with $J(B\op)$.

This gives us two subgroups of $J(\Delta)$, $Y(B) := r(B\op)/r(\Delta) \iso J(B)$ and  $Y(B\op):= r(B)/r(\Delta) \iso J(B\op)$ (note how the ${}\op$ has switched). Denote by $\Arrow\pi_B; J(\Delta).J(B)$ and $\Arrow\pi_{B\op}; J(\Delta).J(B\op)$ the two quotient maps in the short exact sequences. Then we can ask whether the image of $J(B)$ in $J(\Delta)$ (that is, $Y(B)$)  maps under $\pi_{B}$ {\it onto\/} $J(B)$, that is, $\pi_B (Y(B)) = J(B)$. This of course entails that $\pi_B$ splits, but is  stronger than that (there are easy examples wherein $\pi_B$ splits, but this property does not hold). We say that $B$  {\it super-splits\/} when this occurs (it is a two-sided property, as we will see). 

\Lem Lemma \nonesev. Let $B \in \NS_n$. The following are equivalent.
\item{(a)} $B$ is super-splitting;
\item{(b)} $r(B) + r(B\op) = \Z^{1\times n}$;
\item{(c)} $r(B) \cap r(B\op) = r(\Delta)$. 

\Rmk Since (b) and (c) are symmetric under the interchange $B \leftrightarrow B\op$, we deduce that $B$  super-splits iff $B\op$ does. It is also now clear that super-splitting is a PH-invariant (which was not at all clear from the definition, since the latter uses identifications such as that of $J(B)$ with  $\Z^{1\times n}/\Z^{1\times n}B^T$).

\Pf  Obviously, $r(\Delta) \subseteq r(B) \cap r(B\op)$ directly from $(B\op)^T B = B^T B\op = \Delta$, and the map $\pi_B$ is $v + r(\Delta) \mapsto v + r(B)$. 

\noindent (a) implies (b). The map $\pi_B$ sends $Y(B\op) = r(B)/r(\Delta)$ to zero, but is an isomorphism when restricted to $Y(B) = r(B\op) /r(\Delta)$. Hence $(r(B)+ r(B\op))/r(\Delta) $ has cardinality $|J(B)| \cdot |J(B\op)| $, and the latter is $\det \Delta = |\Z^{1\times n}/\Z^{1\times n}\Delta|$. If $r(B) + r(B\op)$ were strictly contained in  $\Z^{1\times n}$, then $J(\Delta)$ would be strictly larger than $|J(B)| \cdot |J(B\op)| $, a contradiction.

\noindent (b) implies (c). We have the  map $J(B) \times J(B\op) \to \Z^{1\times n}/\Z^{1\times n}\Delta$ (given by the identifications of $J(B)$ with $Y(B)$ and $J(B\op)$ with $Y(B\op)$), and this is onto by hypothesis (b). Since the cardinalities are the same, the map is an isomorphism. However, $(r(B\op) + r(B))/(r(B) \cap r(B\op)) \iso r(B\op)/(r(B) \cap r(B\op))   + r(B)/(r(B) \cap r(B\op)) \iso J(B) \times J(B\op)$, so again by cardinality, $r(B) \cap r(B\op) = r(\Delta)$.

\noindent (c) implies (a). Using the standard isomorphisms (as in (b) implies (c)), we have that $r(B\op) + r(B)/r(\Delta)$ is the direct sum, and by cardinality, we obtain (b). Ontoness of $\pi_B$ is then immediate.
\qed

We know that if $J(B)$ or $J(B\op)$ is cyclic, or a free $\Z_d$-module (where $d = \Exp \J(B)$), then both sequences $J(B) \to J(\Delta) \to J(B\op)$ and $J(B\op) \to J(\Delta) \to J(B)$ split. But not all of them super-split. For example, if $B \in \NS_{n,n-1}$ and $B = \( \smallmatrix  I_{n-1} & a \\ 0 & d\\ \endsmallmatrix\)$ where $a = (a_1,\dots, a_{n-1})^T$ is unimodular modulo $d$, then $B$ is super-split iff $1 + \sum a_i^2$ is relatively prime to $d$ (this is obtained by looking at  criterion (b) modulo $d$). Since we can easily solve $1 + \sum a_i^2 \equiv 0 \mod d$ if $n,d \geq 3$ (with all $a_i$ relatively prime to $d$), we have found many examples wherein the sequences split, but $B$ is not super-split.

At one possible opposite extreme is the case that $B$ is Hermite- (not just PHermite-) equivalent to $B\op$. This means that $B = EB\op$ for some $E \in \gl(n,\Z)$, or equivalently, that $r(B) = r(B\op)$. These obviously are not super-split. There are lots of examples; for instance, if $B \in \NS_{n,n-1}$ and is weakly indecomposable, then it follows from the results of section 5 that  $B $ PH-equivalent to $B\op$ implies $B$ is Hermite-equivalent to $B\op$, and necessary and sufficient conditions were given in that section. 

(There are examples, even at size three, of $B \in \NS_n$ with $B$ PH-equivalent to $B\op$ but not being Hermite equivalent; the smallest $d$ for which this occurs seems to be $13^3$.)

There are other possibilities, e.g,, $r(B\op) $ is strictly contained in $r(B)$; examples with $\det B = p$ (a prime) are easy to obtain (since $r(B)$ is a maximal proper subgroup of $\Z^{1\times n}$, there are only three possibilities: either $B$ super-splits ($r(B\op)$ not contained in $r(B)$), $B\op$ is Hermite-equivalent to $B\op$ ($r(B\op) = r(B)$), or there exists noninvertible $F$ \st $B\op = FB$ ($r(B\op)$ is strictly contained in $r(B)$). All three occur.

\SecT 2 PH-equivalence for some terminal forms

Let $n,d> 1$ and consider all the terminal forms with $1$-block size $n-1$ and determinant $d$; that is, matrices of the form $B_a:= \( \smallmatrix \I_{n-1} & a \\ 0 & d \\ \endmatrix\)$, where $a = (a_i)^T$ is in $\Z^{(n-1) \times 1}$ and satisfies $\gcd\brcs{d, a_1, \dots, a_{n-1}} = 1$ and $0 \leq a_i < d$
 for all $i$. Since we can add or subtract multiples of the bottom row to the others at any time in a sequence of PH-equivalences, we may regard the $a_i$ as elements of $\Z_{d}$.

We wish  to describe PH-equivalence for this class of matrices. Since the absolute value of the determinant is a PH-invariant for matrices in $\NS$, we may fix the determinant, denoted $d$; so the problem boils down to the column $a$. There is an obvious action of $S_{n-1}$ on $a$, and this are implemented by left and right multiplication  of $B$ by the corresponding permutation matrix.  Hence at any time, we may assume that the entries of $a$ are, for example, increasing. Alternatively, we can regard $a$ merely as a list, thereby disregarding the action of $S_{n-1}$.

The equivalence relation on $(\Z_d)^{(n-1)\times 1}$  (that is, the columns $a$) induced by PH-equivalence between the corresponding $B_a$ (with $d$ fixed of course) is  more complicated than merely given by permutations.

First, we describe a well-known action of $S_n$ (not $S_{n-1}$) on $A^{n-1}$ where $A$ is a finite abelian group; for convenience, $A$ is written multiplicatively.
The permutation representation of $S_n$ on $A^n$ admits the diagonal $\delta:= \Set{(z,z,\dots,z) }{z \in A}$ as a set of fixed points. Thus there is an action of $S_n$ on the quotient group $A^n/\delta \iso A^{n-1}$. To see just what the resulting action is, pick $y = (a_i) \in A^{n-1}$; lift it to an element of $A^n$ by setting $y' = (y,1)$ (since $A$ is written multiplicatively, $1$ means the identity element). Apply the permutation action of $S_n$ to $y'$.

For $\pi \in S_n$, if $\pi$ fixes the point $\brcs{n}$, then it comes from an element of $S_{n-1}$, so we just define $\pi (y)$ to be the first $n-1$ coordinates of $y'$, the obvious thing. Otherwise, there exists $j < n$ \st $\pi(j) = n$, so that the last coordinate of $\pi(y') $ is $a_j$ and a $1$ appears in the $\pi(n)$-entry. Multiply the vector $\pi(y)$ by $a_j^{-1}$. Now the final entry is $1$, so we can define $\pi(y)$ to be the first $n-1$-coordinates of $a_j^{-1} \pi(y')$. The multiplication operator is equivalent to performing the group action with the diagonal element $(a_j^{-1}, a_j^{-1}, \dots, a_j^{-1})$ to $\pi(y')$, hence is compatible with the quotient action.

(Replacing $n$ by $n+1$ and $A$ by $\Z$---this time viewed additively---this is the Weyl group action of $S_{n+1}$ on the dual of the maximal torus of SU$(n+1)$.)

\def\Aut{\text{Aut}}
Denote this action $\Arrow \Pi_{A,n} ; S_n . \Aut\, A^{n-1}$. Now replace $A$ by $\Z_d^*$, the group of (multiplicatively) invertible elements in $\Z_d$ (so $\phi(d) = |\Z_d^*|$). Suppose that $a \in (\Z_d)^{n-1}$ consists of elements relatively prime to $d$, that is, members of $\Z_d^*$. Then we will see that the PH-equivalence class of $B_a$ consists of a slightly twisted $S_n$-orbit of $a$ under $\Pi_{\Z_d^*,n}$.

However, if some of the entries of $a$ are zero divisors in $\Z_d^*$, then the situation becomes pear-shaped. We may permute the entries so that the first $k$ are invertible, and the rest are zero-divisors. Then we can apply $S_{k+1}
$ to the column of the first $k$, obtaining (for each element of the group) an element $a_j^{-1}$---and instead of multiplying merely the top $k$ entries by $a_j^{-1}$, we multiply all of $a$ by it.

This of course preserves the entries that are zero-divisors in the ring $\Z_d$, whose locations are unmoved. It also preserve the ideals the elements generate, e.g., there are the same number of zeros in the new element as in the original, the same number that are divisible by any prime $p$ that divides $d$ as in the original, etc.

The upshot is that there is no group structure (except when $n=2$) on the equivalence classes, but instead a union of actions of various groups.

When $n = 2$,   $B_a$ is PH-equivalent to $ B_{a'}$ iff $aa' \equiv 1 \mod d$; this is easy, and can be done directly, since we are dealing only with the transposition matrix. For $n > 2$, if each of the  $a_i$ is {\it not\/} relatively prime to  $d$, then the equivalence class is simply the set of permutations of the entries, that is, via the action of $S_{n-1}$---in this case, there is an obvious normal form,  arranged monotonically.

If $a_i$ are all relatively prime to $d$, then the action is given by permutations and a twisted multiplication by each of the $a_i^{-1}$ modulo $d$; it looks like these should generate a larger orbit, but they don't. (So if all $a_i = -1$, it is not equivalent to anything else.) The orbit consists of $\brcs{(-a_1 a_j^{-1}, \dots, a_{j}^{-1}, \dots,  -a_{n-1}a_j^{-1})}$, together with $(a_i)$ itself, and all their permutations.

To verify these claims, suppose $UB_{a}P = B_{a'}$. First, we note that if also $U_1 B_a P = B_{a''}$ (where $B_a$, $B_{a'}$, and $B_{a''}$ are all terminal) with the same $P$, then  $B_{a'} = B_{a''}$. This follows from the equalities in $\Mn n \Q$, $P =B_{a}^{-1} U^{-1} B_{a'} = B_{a}^{-1} U_1^{-1} B_{a''} $, whence $U^{-1} B_{a'}  = U_1^{-1} B_{a''}$, so that $U_1 U^{-1} B_{a'} = B_{a''}$. Set $V = U_1 U^{-1} \in \gl(n,\Z)$. From the form of the $B$s (first $n-1$ columns are standard basis elements), $V = \(\smallmatrix I_{n-1} & X \\ 0 & t \\ \endsmallmatrix\)$; since the lower right entries of both $B$s are $d$, $ t = 1$, and we have $a' + d X = a''$; but this simply means that $a$ and $a' $ are coordinatewise congruent modulo $d$; since we have assumed the entries are in the interval $0\leq a_i'',a_i' < d$, this forces $a' = a''$. 

Thus for each permutation matrix $P$, there is a most one $a'$ for which $UB_a P = B_{a'}$  for some $U \in \gl(n,\Z)$ (and of course, there may be none).

Let $B$ and $B'$ be matrices in $NS_n$, both in terminal form. Suppose there exists a permutation matrix $P$ together with $U$ in $\gl(n,\Z)$ \st $UB = B'P$; then we say {\it $P$ is realizable over $B$} (in other words, there exists $B'$ in terminal form, etc).

Suppose that $B = \(\smallmatrix \I_{n-1} & a \\ 0 & d\\  \endmatrix\)$; its $1$-block is size $n-1$. Let $\pi$ be the permutation corresponding to the right action by $P$ on columns (that is, if $P$ takes the first column to the second, then $\pi(1) = 2$).

 If $\pi(n)= n$, then $B' = PBP^{-1} $ is also in terminal form, since $a$ has been replaced by $Qa$ (a permutation of the entries of $a$)  where $P = Q \oplus 1$. So in this case, all of $S_{n-1}$ is realizable. Morever, if $P'$ is realizable over $B$ and $P = Q \oplus 1$, then $PQ$ is also realizable, so that the realizable permutation matrices form a coset space over $S_{n-1}$. However, this is fairly complicated.

For $a \in \Z^{(n-1) \times 1}$ \st $\Ct \brcs{a,d} = 1$, recall that $B_a = \(\smallmatrix \I_{n-1} & a \\ 0 & d\\  \endmatrix\)$. We will determine precisely the permutation matrices $P$ \st there exist $a'  \in \Z^{(n-1) \times 1}$ \st $B_a'$ satisfies $UB_a = B_{a'}P$ for some $U \in \gl(n,\Z)$. This is not the full realizability problem, since $P$ may be realizable over $B$, but the outcome, $B'$, although in terminal form, need not have its $1$-block of size $n-1$. (We have already seen such an example.)

For an integer $d>1$, $\Z_d^*$ will denote the group of multiplicatively invertible elements in the {\it ring\/} $\Z_d$ (formerly, we just considered the latter as an additive group). If $x$ is an integer relatively prime to $d$, then $x^{-1}$ will denote a representative $y$ \st $xy \equiv 1 \mod d$.

\Lem Proposition \twoone. Let $d> 1$ be an integer. Let $P$ be a permutation matrix of size $n$ with corresponding permutation $\pi$, and $a   \in \Z^{(n-1) \times 1}$ \st $\Ct \brcs{a,d} = 1$.
Then  $P$ is realizable over $B_a$ with $B_{a'} = UB_a P^{-1}$ having $1$-block of size $n-1$ iff either  $\pi(n) = n$ or $a_{\pi(n)}$ is invertible modulo $d$. In the latter case, modulo $d$,
$$
a'_{t} \equiv \cases
- a_{\pi(t)} a_{\pi(n)}^{-1} &\text{if $t \neq  \pi^{-1}(n)$}\\
a^{-1}_{\pi(n)} & \text{if $t = \pi^{-1}(n)$}.\\
\endcases
$$

\Rmk It is important to emphasize that this result describes only PH-equivalence between terminal forms of $NS_n$-matrices, both of which have $1$-block size $n-1$. It says only a limited amount about PH-equivalences between terminal forms only one of which has $1$-block size $n-1$ (essentially, the statement that each permutation matrix $P$ can contribute at most one new terminal form). In particular, if $\gcd\brcs{a_i,d} > 1$ for all $i$, then the only choices for $P$ are those corresponding to $S_{n-1}$---in this context. Where we are allowed to choices for terminal $B'$ that have a smaller $1$-block, we can obtain more realizable $P$.

\Rmk For $n =3$, this type of action of the symmetric group was discussed in [R].

\Pf First, suppose that $UB_a = B_{a'}P$ for some $U \in \gl(n,\Z)$, and $\pi(n) \neq n$. Then the $i$th column of $UB_a$ is $Ue_i$, except when $i = n$, in which case, it is $U\(\smallmatrix a \\ d \endsmallmatrix\)$. On the other hand, the $i$th column of $B_{a'}P$ is the $\pi^{-1}$th column of $B_{a'}$, which is $e_{\pi^{-1}}$, unless $\pi(i) = n$, in which case it is $\(\smallmatrix a' \\ d \endsmallmatrix\)$.

In particular,
$$\eqalign{
Ue_i &= \cases e_{\pi^{-1}(i)}& \text{if $ i \not\in \brcs{n, \pi(n)}$}\\
\(\matrix a' \\ d \endmatrix\) & \text{if $i = \pi(n)$}
\endcases\cr
U\(\matrix a\\ d \endmatrix\) & = e_{\pi(n)}.
}$$
We have that  $n-2$ of the columns of $U$ are standard basis vectors and  the $\pi(n)$th column is $\(\smallmatrix a' \\ d \endsmallmatrix\)$; let $(h_j)^T$ be the $n$th column of $U$. The basic vectors represented in the columns exclude $e_n$ and $e_{\pi(n)}$; hence in the $\pi(n)$ and $n$th rows of $U$, there are at most two nonzero entries, $a'_{\pi(n)}$ and $h_{\pi(n)}$, \& $d$ and $h_n$, respectively.

Now we can apply this to the third equation, and obtain (after sorting through the subscripts and cases),
$$\eqalign{
a_{\pi(t)} + a'_{t}a_{\pi(n)} + h_t d & = 0 \qquad \text{if $t \neq n, \pi^{-1}(n)$}\cr
a'_{\pi^{-1}(n)} a_{\pi(n)} + h_{\pi(n)} d & = 1.\cr
}$$
The second equation says that $a_{\pi(n)}$ is invertible modulo $d$, and $a'_{\pi^{-1}(n)} \equiv a_{\pi(n)} \mod d$. Now that we know that $a_{\pi(n)}$ is invertible modulo $d$, the first equation yields the rest of what we want.

As to the converse, we can almost reconstruct $U$ from the equations; the $a'_i$ are defined up to multiples of $d$ (so we can perform additional elementary row operations if necessary to ensure that $0 \leq a'_i < d$. There is only one additional condition;  $|\det U| =1$ iff $\left|\det \( \smallmatrix a'_{\pi(n)} & h_{\pi(n)}  \\ d & h_n \\\endsmallmatrix\) \right| = 1$, that is, $a'_{\pi(n)}h_n - h_{\pi(n)}d = \pm 1$, which is easily arranged (since $a'_{\pi(n)}$ is invertible modulo $d$).

The case that  $\pi(n) = n$ has already been discussed.
\qed

In particular, the number of $i$ \st $\gcd\brcs{a_i,d} = 1$ is an invariant of this equivalence relation, as is for each prime $p$ dividing $d$ and each $m$, the number of $a_i$ \st $p^m$ divides $a_i$ since up to permutation, we are multiplying the entries by an invertible modulo $d$, except in one place, where an invertible is replaced by another invertible. Both of these are also obtainable from $I(B_a)$ as in Lemma \onetwo\ above. Generically the number of elements in the equivalence class of $B_a$ is
$$(n-1)! \cdot \left| \Set{i} {\gcd\brcs{d,a_i} = 1}\right|,
$$
but it could be less.  Observe that if $a_i = a_j \in \Z_d^*$, on taking a permutation $\pi$ \st $\pi(n) =i$, the corresponding $a'$, is up to permutation (that is, the $S_{n-1}$-action), obtained by multiplying all the entries by $-a_i^{-1}$ and replacing one of the $-1$ terms that result  by $a_i^{-1}$; the same set, up to the $S_{n-1}$ action, will arise from a permutation sending $n \mapsto j$. In this case, different permutations, even modulo $S_{n-1}$ are realizable, but yield the same matrices.

For $n = 2$, of course the only possible action is $a \mapsto a^{-1}$ (modulo $d$).
In particular,
$$
\( \matrix  1 & a \\ 0 & d \\\endmatrix\) \text{ is PH-equivalent to }
\( \matrix  1 & a' \\ 0 & d' \\\endmatrix\)
$$
iff $d = d'$ and either of $a' \equiv a^{\pm1} \mod d$.


It also allows us to conclude that
$$
\(\matrix 1 & 0 & 2 \\ 0 & 1 & 3 \\ 0 & 0 & 6 \\
\endmatrix \) \text{ and }
\(\matrix 1 & 0 & 4 \\ 0 & 1 & 3 \\ 0 & 0 & 6 \\
\endmatrix \)
$$
are not PH-equivalent. As they are respectively equivalent to
$$
\(\matrix 1 & 1 & 2 \\ 0 & 2 & 0 \\ 0 & 0 & 3 \\
\endmatrix \) \text{ and }
\(\matrix 1 & 1 & 1 \\ 0 & 2 & 0 \\ 0 & 0 & 3 \\
\endmatrix \),
$$
the latter two are not PH-equivalent to each other either. All four matrices have $J(B) \iso \Z_6$.

If two matrices $B, B' \in \NS_{n,n-1}$, then there is a relatively efficient  procedure for deciding whether they are PH-equivalent. The determinants must be the same, $d$, and each has a list $\Lt a$, $\Lt {a'}$ (consisting of the integers in the last column, above the $d$). There are only $n$ cosets of $S_n/S_{n-1}$, and we just have to test those for which the corresponding element of $\Lt a$ is relatively prime to $d$ (testing for relative primeness of $a_i$ and $d$ requires at most $\Oh{\ln a_i}$ steps, usually much less), and for each one of those, do the operation described in Proposition \twoone, and check whether the new list is that of $a'$. To make it more efficient, we may rearrange the lists as they appear so they are descending, etc. This amounts to sorting lists of nonnegative integers with a fixed upper bound, $d-1$, on the entries.   An easy algorithm (good if $d \ll n$) is for each $i= 0, 1, \dots, n-1$, decide which of the numbers in $\brcs{0, 1, 2, \dots, d-1}$ $a_i$ is, and keep $d$ running counts. The final counts determine the ordering.

If merely one of them has  $1$-block size $n-1$, then we first test whether $B'$ does as well, by deleting the $i$th column and testing whether the resulting row space is all of the standard copy of $\Z^{n-1}$---one way is to take the $n$ determinants of the submatrices of size $n-1$, and see if their greatest common divisor is one (it would be enough to show their gcd is relatively prime to the determinant of $B'$). If $B'$ is already in terminal (or merely upper triangular) form, this will likely be very fast.
\vskip 4pt

\SecT 3 A family of invariants

We will use {\it lattice\/} in the sense of partially ordered sets with least upper and greatest lower bounds. 

In this section, we introduce and investigate a family of invariants which form a lattice of abelian groups with factor maps between them. 

Fix $n$, and for $1 \leq i \leq n$, let $\Arrow p_i;\Z^{1\times n} .\Z$ be the coordinate maps, and let $S = \brcs{1,2,\dots, n}$. Let $\Omega \subset 2^S$. For $B \in \NS_n$, define $B_{\Omega} \in \NS_{|\Omega|}$ (up to PH-equivalence) as follows. Delete from $B$ the columns whose index is not in $\Omega$ (thus, if $1  \not\in \Omega$, delete the first column of $B$) to create an $n \times |\Omega|$ matrix, each of whose columns has content one. The rank of the resulting matrix is exactly $|\Omega|$, since the set of columns of $B$ was linearly independent to start with. By applying elementary (integer) row operations to $B$ with columns deleted, we can obtain a matrix of the form $\(\smallmatrix C \\ 0 \\ \endsmallmatrix\)$ where $C$ is square of size $\Omega$. Since elementary row operations preserve the content of  columns, $C \in \NS_{|\Omega|}$. All choices for such $C$ are Hermite- (and therefore PHermite-) equivalent (within $\NS_{|\Omega|}$). We choose one, and call it $B_{\Omega}$.

An alternative approach (leading to the same thing) is to consider the PH-equivalence class of $B$ as a means of studying the row space of $B$, $r(B) \subset \Z^n$, up to the restriction of the action of the permutation matrices acting on the right (that is, as column permutations). When we delete the columns not corresponding to elements of $\Omega$ and take the row space of the resulting matrix, and use that to define $r(B_{\Omega})$, without defining $B_{\Omega}$(!).

If $\Omega = S$, then $B_{\Omega} = B$. If $\Omega$ consists of a singleton, then the resulting column, being unimodular, row-reduces to the first (or any) standard basis element of $\Z^{n \times 1}$, and thus $B_{\Omega} = (1)$, the size one identity matrix.

Define, for each $i = 1,2, \dots, n$, the subset $\Omega(i)= S\setminus\brcs{i}$. 

If $\Omega' \subset \Omega$, let $\Arrow P_{\Omega',\Omega}; \Z^{1 \times\Omega}.\Z^{1 \times \Omega'}$, and $\Arrow P_{\Omega}; \Z^{1\times n}.\Z^{ 1 \times \Omega}$ be the obvious projection maps (sometimes we will rewrite the last as $\Z^{\Omega}$). Then $P_{\Omega} (r(B)) = r(B_{\Omega})$ (and similarly for $P_{\Omega', \Omega}$),  thus  inducing the well-defined, onto group homomorphisms $\Arrow p_{\Omega}; J(B).J(B_{\Omega})$ and $\Arrow p_{\Omega', \Omega};J(B_{\Omega}) . J(B_{\Omega'})$. It is routine to verify that the maps are transitive (that is, if $\Omega'' \subset \Omega' \subset \Omega$,  then $p_{\Omega'',\Omega'}\circ p_{\Omega',\Omega} = p_{\Omega'',\Omega}$). In case there is ambiguity about which $B$ they are referring to, we will occasionally use $p_{\Omega}^B$.

\comment
To see this, we first note that we may assume that $\Omega = S$, so $B = B_{\Omega}$. Now $X(B)$ is spanned by $\brcs{x_i}$ where $x_i$ is uniquely determined from $x_i \in \cap_{j\neq i}\ker t_j$ and $t_i(x_i)$ is minimal among elements of $\cap_{j\neq i}\ker t_j$. The corresponding $t$s for $B_{\Omega'}$ are simply $t_k$ with $k \in \Omega$. Since $x_i$ vanishes at all $t_j$s except $t_i$, its image in the row space of $B_{\Omega'}$ either vanishes at all the associated $t$s (so its image is zero
) or, if $i \in \Omega'$, it vanishes at all but one. In that case, its image is an integer multiple of the corresponding $X_i$. Hence $X(B_{\Omega})$ maps to $X(B_{\Omega'})$, and thus we obtain a family of maps $I(B_{\Omega}) \to I(B_{\Omega'})$, one for each pair $(\Omega', \Omega)$ with $\Omega' \subset \Omega$.

Suppose that $\Omega' \subset \Omega$; we thus have an induced onto group homomorphisms $\Arrow p_{\Omega',\Omega}; J(B_{\Omega}).J(B_{\Omega'})$. It is routine to verify that the maps are transitive (that is, if $\Omega'' \subset \Omega' \subset \Omega$,  then $p_{\Omega'',\Omega'}\circ p_{\Omega',\Omega} = p_{\Omega'',\Omega}$).
\endcomment

Now suppose that $B$ and $B'$ belong to $\NS_n$, and there is a PH-equivalence between them. Then we claim this implies that there exists a permutation of $
S$ together with a compatible family of group isomorphisms $J(B_{\Omega}) \to J(B_{\pi \Omega})$. This is  trivial: if we apply an element of $\gl(n,\Z)$, the row space is unchanged, and we obtain the identity maps. If we permute columns, $\pi$ is the corresponding permutation, etc. We thus see that not only is $J(B)$ a PH-equivalence, but so is (for example), the set of maps $J(B) \to J(B_{\Omega})$ where we restrict the $\Omega$s to have the same cardinality.

The lattice of maps and quotient groups $\Arrow p_{\Omega',\Omega}; J(B_{\Omega}).J(B_{\Omega'})$  will be denoted $ {\Cal J}(B)$. This is a fairly strong invariant, as we will see later, but it is also somewhat more difficult to calculate (except in special cases), compared to the list $\List{J(B)_{\Omega}}_{|\Omega| = n-1}$.

\Lem Lemma \throne. The lattice of finite abelian groups and homomomorphisms, $\Cal J(B)$, is a PH-invariant for matrices $B$ in $\NS_n$. If $ k \leq n$, then the list $\Lt{J(B_{\Omega})}_{|\Omega| = k}$ is also a PH-invariant.

When we put $k = n-1$, we obtain 
 the list of  $n$  groups $\Lt{J(B_{\Omega(i)})}$. This contains a lot of information (although generally less than $\Cal J(B)$).

Originally, the intent of developing $J(B_{\Omega})$ was to find a finer invariant that $J(B\op)$: even together with $|\det B\op|$ (another PH-invariant) $J(B)$  does not determine the family $\brcs{J((B\op)_{\Omega})}$, or  $\brcs{J((B\op)_{\Omega(i)})}$. More interestingly, those $\Omega$ for which $J(B_{\Omega})= \brcs{0}$ play a particularly important role. For example, we will see there exists  $\Omega$ of cardinality $r$ \st $J(B_{\Omega}) = \brcs{0}$ iff $B$ is PH-equivalent to a terminal form whose $1$-block size is at least $r$. This is practically tautological, but provides a useful way of constructing interesting examples.

We also have a second family of PH-invariants,
 specifically, $\Cal J(B\op)$. In general $(B_{\Omega})^{op}$ is not PH-equivalent to $(B^{op})_{\Omega}$, nor need they yield isomorphic invariants. So we have to be careful \wrt the notation, that is, construct the opposite, $B^{op}$, first, then the cut-down matrices, $(B^{op})_{\Omega}$. However, I could not decide whether $\(J((B{op})_{\Omega})\)_{\Omega \in 2^S}$ is determined by $\(J(B_{\Omega})\)_{\Omega \in 2^S}$, that is, whether for  $B,B' \in \NS_n$,
\st $\(J(B_{\Omega})\) \iso \(J(B'_{\Omega})\)$ (as a family, that is, $\Cal J(B) \iso \Cal J(B')$ implies $\Cal J(B\op) \iso \Cal J(B'{}\op)$). This will discussed in more detail in section\,6.

Recall that 
$r(B)$ frequently   denote $\Z^{1\times n} B$, re-inforcing the idea that it is the subgroup of
$\Z^{1\times n}$ generated by the rows of $B$.

\Lem Lemma \nthrone. Suppose that $B \in \NS_n$. Then $\ker p_{\Omega}$ is spanned by
$\brcs{E_j + r(B)}_{i \notin \Omega}$.

\Pf Obviously $P_{\Omega})(E_j) $ is zero if $j \notin E$, so that $E_j + r(B)
\in \ker p_{\Omega}$ for all $j \notin \Omega$.

Suppose for $v \in \Z^n$ that $P_{\Omega}(v) \in r(B_{\Omega})$. Then there
exist $a_i \in \Z$ and corresponding rows $r_i$ of $B$ \st $P_{\Omega}(v) - \sum
a_i P_{\Omega}(r_i) = 0$. Thus the only nonzero entries of $w:= v - \sum a_i
r_i$ can only appear in position $j$ where $j \not\in \Omega$. We can thus write
$w = \sum_{\Omega^c} b_j E_j$, and so $v \in r(B) + \sum_{\Omega^c} E_j \Z$.\qed

\Lem Corollary \nthrtwo. If $\Omega' \subset \Omega$, then $\ker p_{\Omega',\Omega}$ is
spanned by $\brcs{E_j + r(B) + \sum_{i \in \Omega} E_i \Z}_{j \in
\Omega\setminus \Omega'}$.

\Pf One inclusion is obvious; for the other, suppose that $p_{\Omega',\Omega} (v
+ \ker p_{\Omega}) = 0$. Then $p_{\Omega'}(v+ r(B)) = 0$. Hence by the
preceding, there exist integers  $a_i$  and $b_j $ ($j \notin \Omega'$) \st $a -
\sum a_i r_i = \sum_{j \in \Omega'{}^c} b_j E_j $. If $j \not\in \Omega$, then
$E_j + r(B) \in \ker p_{\Omega}$; thus the right side decomposes as $\sum_{j \in
\Omega \setminus \Omega'} b_j E_j$ plus an element of $\ker p_{\Omega}$.
\qed

As a consequence, if $\Omega' \subset \Omega$ and $J(B_{\Omega'})$ is generated (as an abelian group, or as a $\Z_d$-module) by $k$ elements, then $J(B_{\Omega})$ has a generating set of cardinality at most $k + |\Omega| - |\Omega'|$.

As usual, $S_n$ will denote the full permutation group on $S = \brcs{1,2,\dots,n}$; sometimes this will be identified with $\Cal P_n$, the group of $n \times n$ permutation matrices. 

\Lem Proposition \nthrthr. Let $B$, $B'$ belong to $\NS_n$. Necessary and sufficient for
${\Cal J}(B) \iso {\Cal J}(B') $ is the following condition:
\item{} there exist $\pi \in S_n$ and an isomorphism $\Arrow \phi; J(B). J(B')$
\st for all $i$, $\phi(\ker p_{\Omega(i)}^B) = \ker p_{\Omega (\pi i)}^{B'}$.

\Pf Necessity is obvious, so let us prove sufficiency. Let $d$ be exponent of
$J(B)$ (which by the isomorphism, is also the exponent of $J(B')$. If $P$ is the
permutation matrix representing $\pi^{-1}$, then we can replace $B'$ by $B'P$,
and thus assume that $\pi $ is the identity.

By the preceding characterization of $\ker p_{\Omega}$, for any $i$, we obtain
$\phi(\langle E_i + r(B)\rangle) = \phi(\ker p_{\Omega(i)}^B) = \ker
p_{\Omega(i)}^{B'} = \langle E_i + r(B')\rangle$. Hence for any proper subset
$\Omega$, $\phi\(\langle E_i + r(B)\rangle_{i \not\in \Omega} \) = \langle E_i +
r(B')\rangle_{i \not\in \Omega} $. By the preceding proposition, $\phi(\ker
p_{\Omega}^B) = \ker p_{\Omega }^{B'} $. Then we define the map $\Arrow \phi_{\Omega};J(B_{\Omega}).J(B'_{\Omega})$ in the obvious way, $v +
r(B) + \sum_{i \notin \Omega}  E_i\Z \mapsto \phi\(v + r(B' )+\sum_{i \notin
\Omega}  E_i\Z \)$; this is well defined by the preceding sentence, and is an isomorphism. Thus the following diagram commutes.
$$
\diagram
J(B)&\rTo^{\phi} & J(B') \\
\dTo^{p_{\Omega}^B} &&\dTo^{p_{\Omega}^{B'}}& \\
J(B_{\Omega})&\rTo^{\phi_{\Omega}} & J(B_{\Omega}') \\
\enddiagram
$$

If $\Omega' \subset \Omega$, then from $P_{\Omega',\Omega} \circ P_{\Omega} =
P_{\Omega'}$, the corresponding diagram with $B_{\Omega} $ replaced by
$B_{\Omega'}$ and $B$ replaced by $B_{\Omega}$ also commutes. Hence $\phi$
induces an isomorphism of lattices of quotient groups, ${\Cal J}(B) \iso \Cal
J(B')$.\qed

\Lem Corollary \nthrfou. Suppose $B,B' \in \NS_n$ and $J(B)$ is cyclic. Sufficient for
${\Cal J}(B) \iso {\Cal J}(B')$ is that $J(B) \iso J(B')$ and there exist $\pi
\in S_n$ \st for all $\Omega$, $|J(B_{\Omega}) | = |J(B'_{\pi\Omega})|$.

\Pf Any quotient of $J(B)$ (and therefore of $J(B')$ is cyclic and therefore
their cardinality determines uniquely their isomorphism class  (so that
$J(B_{\Omega})  \iso J(B'_{\pi\Omega})$ and the kernel (as cyclic groups have at
most one subgroup of given cardinality). Now the preceding proposition applies.
\qed

Particularly useful are the $J(B_{\Omega(i)})$ (recall that $\Omega(i) = S\setminus \brcs{i}$, the subset of $\brcs{1,2,\dots,n}$ missing only $i$). We define $\NS_{n,m}$ to consist of the elements $B \in \NS_n$ which have a terminal form with $1$-block of size at least $m$. Thus $\NS_{n,n-1}$ consists of elements of $\gl(n,\Z)$ (trivially, these have $1$-block size $n$) and those $B \in \NS_{n}$ with a terminal form having $1$-block of size $n-1$.

A matrix $B \in \NS_n$ is {\it decomposable\/} if it is PH-equivalent to a direct sum of matrices in $\NS$, and {\it indecomposable\/} otherwise. it is {\it weakly indecomposable\/} if it is not PH-equivalent to a matrix of the form $1 \oplus C$ where $C \in \M_{n-1}\Z$ (if such  $C$ exists, it is necessarily in $\NS_{n-1}$).

\Lem Lemma \nthrfiv. Suppose that $B \in \NS_n$. 
\item{(a)} $B \in \NS_{n,m}$ iff there exists $\Omega \subset 2^S$ with $|\Omega| = m$ and $J(B_{\Omega}) = 0$.
\item{(b)} $B \in \NS_{n,n-1}$ iff there exists $i$ \st  $J(B_{\Omega(i)}) = 0$.
\item{(c)}  $B$ is weakly indecomposable iff for all $i =1,2, \dots, n$, the kernel of $\Arrow p_{\Omega(i)}; J(B).J(B_{\Omega(i)})$ is not zero. 

\Pf (a) If such an $\Omega$ exists, the set of $\Omega$-truncated rows of $B$ contain a $\Z$-basis for $\Z^{\Omega}$; by rearranging the rows of $B$, we can assume that the top $|\Omega|$ rows of $B$, $(B_{(i)})_{i=1}^{|\Omega|}$, satisfy $(P_{\Omega}(B_{(i)}))$ is a basis for $\Z^{\Omega}$. By permuting the columns, we can also assume that $\Omega = \brcs{1,2, \dots, |\Omega|}$. Then the upper left $|\Omega| \times |\Omega|$ corner of the current $B$ belongs to $\gl(|\Omega|,\Z)$ (since the rows form a $\Z$-basis for $\Z^{\Omega}$). Hence there exists $E \in \gl(n,\Z)$ of the form $E = F \oplus \I_{\Omega^c}$ with $F \in \gl(|\Omega|,\Z)$ \st $EB =\(\smallmatrix I_{|\Omega|} & X \\ Y & Z \\ \endsmallmatrix\)$. The obvious row operations allow us to reduce to the case that $Y = 0$.  Now we can apply the procedure of [TSCS] to put $Z$ itself in terminal form.  There is nothing to prevent additional $1$s from appearing, so when we proceed to fix $X$ (by appling row operations corresponding to the rows of  the new $Z$) so that the $n \times n$ matrix is in terminal form, the identity block size may have become larger.   The resulting matrix is a terminal form with $1$-block size at least $|\Omega|$.

If $B$ has a terminal form $C = \(\smallmatrix I_{|\Omega|} & X \\ 0 & Z\\ \endsmallmatrix \)$, then with $\Omega = \brcs{1,2,\dots, |\Omega|}$, we have $C_{\Omega}$ consists of the first $\Omega$ standard basis elements as columns, hence $J(C_{\Omega}) = 0$. Since $B$ is PH-equivalent to $C$, there exists $\Omega'$ (obtained from a permutation in $S_n$, hence of equal cardinality) \st $J(B_{\Omega'}) = 0$.

\noindent (b) Apply (a) to subsets consisting  of $n-1$ elements.

\noindent (c) If $B $ is PH-equivalent to $1 \oplus C$, then $\ker p_{\Omega(1)}^C$ is spanned by the image of the standard basis element $E_1$; but this already belongs to $r(C)$, so the kernel is zero. Conversely, suppose $\ker p_{\Omega(i)}^B = 0$. Then $E_i \in r(B)$ (by \nthrone\  above); applying the obvious row operations to eliminate all the other nonzero entries in the $i$th column, and then rearranging the columns (moving the $i$th column to the first), we see that the resulting matrix is a direct sum.
\qed

\Lem Examples \nthrsix. Matrices $B \in \NS_3$ ($\NS_4$) \st $J(B)$ is cyclic, but neither $B$ nor $B\op$ has a size two (respectively, three) $1$-block terminal form. 

\noindent (i) Let
$$
B = \( \matrix 1 & 1 & 8 \\ 0 & 2 & 6 \\ 0 & 0 & 15\\ \endmatrix\); \quad B^{-1} = \( \matrix 1 & -\frac 12 & -\frac 13 \\ 0 & -\frac12 & - \frac 15 \\ 0 &  0 & \frac 1{15}\endmatrix\).
$$
As $\det B = 30$ is square-free, $\Z^3/\Z^3 B$ is cyclic. It is easy to check that the list of $J(B_{\Omega(i)})$ is $\List{\Z_5, \Z_3, \Z_2}$, so $B$ has no size two $1$-block terminal forms by \nthrfiv(b). 

From the inverse, we see that the $\List{m(i)} = \List{6,10,15}$ (the smallest positive integer to make the corresponding row integral), and thus
$$
B\op =\( \matrix 6 &0 & 0 \\ -3 & 5 & 0 \\ -2 & -2 & 1\\ \endmatrix\).
$$
Thus $\det B\op = 30$, so again $\Z^3/\Z^3 B\op = \Z_{30}$. It is straightforward to verify $\List{J((B\op)_{\Omega(i)})}$ is $\List{\Z_5, \Z_3, \Z_2}$ (again); since none of them are zero, $B\op$ has no size two $1$-block terminal forms. 

\noindent (ii) A different type of example arises from decomposable matrices. Recall that $B \in \NS_n$ is {\it decomposable\/} if there exists a matrix $B' = A\oplus C \in \NS_n$ that is PH-equivalent to $B$ (from the fact that $B' \in \NS_{n}$, it follows easily that both $A$ and $C$ belong to their corresponding $\NS_j$). 

Set $A = \(\smallmatrix  1 & 1 \\ 0 & 3 \\\endsmallmatrix \)$ and $C = \(\smallmatrix  1 & 1 \\ 0 & 5 \\\endsmallmatrix \)$, and define $B = A\oplus C$. Each of $A$ and $C$ belong to $\NS_2$, so $B \in \NS_4$. Moreover, $B\op = A\op \oplus C\op$, so $\det B = \det B\op = 15$. The latter being square-free, both $J(B)$ and $J(B\op)$ are cyclic. 

However, for any $i$, $J(B_{\Omega(i)})$ has a direct summand which is one of $A$ or $B$ (this is true for any direct sum); in this case, both $A$ and $B$ are not invertible, so $J(A)$ and $J(B)$ are both nonzero. Thus $J(B_{\Omega(i)})$ is not zero for any $i$, so $B \not\in \NS_{4,3}$, and similarly, $B\op \notin \NS_{4,3}$ by \nthrfiv(b).
\qed

Now we want to address the near-ubiquity of matrices some but not all of whose terminal forms have  $1$-block size $n-1$. 
A useful PH-equivalence tool (found in [AALPT]) is that $C$ and $C'$ are PH-equivalent (via the permutation matrix $P$ or its inverse) iff $C'PC^{-1}$ has only integer entries. We use this frequently, without further comment.

\comment

\Lem Proposition \twotwo. Let $B = \(\smallmatrix 1 & b \\ 0 & D\\ \endsmallmatrix\)$ (below) be a terminal member of $\NS_n$, and suppose $d_1 > 1$ with $d:= \prod d_i$ \st $\gcd\brcs{d_i,d_j} = 1$ for all $i \neq j$. Then there exist $a_i$ divisible by $d/d_{n-i+1}$, \st
$$
B:= \( \matrix 1 & b_2 & b_3 & \dots & b_n \\
0 & d_2 & 0 & \dots & 0 \\
0 & 0 & d_3 & \dots & 0\\
\vdots&\vdots&0&\ddots&\vdots\\
0 & 0 & 0 & \dots & d_n\\ \endmatrix\) \text{ is PH-equivalent to } \( \matrix 1 & 0 & 0 & \dots &0 & a_1 \\
0 & 1 & 0 & \dots & 0 & a_2\\
0 & 0 & 1 & \dots & 0 & a_3\\
\vdots&\vdots&0&\ddots&\vdots& \vdots\\
0 & 0 & 0 & \dots & 1 & a_{n-1}\\
0 & 0 & 0 & \dots & 0 & d\\
\endmatrix\),
$$
where the right matrix is in terminal form.

\Rmk Both matrices are in terminal form; the  left has $1$-block size one, and the right has $1$-block size $n-1$.

\Rmk Example \thrthr\ (fourth matrix) provides an example in terminal form wherein the diagonal entries are pairwise relatively prime, but the matrix is not PH-equivalent to any terminal form with $1$-block size $n-1$; the difference lies in the fact that the lower right block matrix is not diagonal, but merely triangular. This shows that the zeros there are essential to the statement of the result.
 
\Pf Deleting the first column of $B$, the row space of the resulting matrix $B'$ is all of the standard copy of $\Z^{n-1}$. Hence there exists $E \in \gl{n,\Z}$ \st $EB' = \(\smallmatrix \I_{n-1} \\ 0 \\ \endsmallmatrix \)$. Thus $C = EB$ is in the form on the right.  The lists, $\List {J(B_{\Omega(i)})}$ and $\List{J(C_{S\setminus\brcs{i}})}$, are equal. The latter is easy to calculate, $\List{\Z_{(d,a_i}; 0}$.

On the other hand, $J(B_{S\setminus \brcs {1}}) = 0$ by xxx above, and it is elementary that $J(B_{S\setminus \brcs{j}}) = \oplus_{k \neq j} \Z_{d_{k}} \iso \Z_{d/d_j}$, so up to a permutation, the result follows.
\qed

************************************************
\endcomment 

 Let $B$ denote the $n \times n$ integer matrix $ \(\smallmatrix \I_{n-1}& a \\ 0 & d \\ \endsmallmatrix\)$, with $d >1$, where $a = (a_1, \dots, a_{n-1})^T \in \Z^{(n-1)\times 1}$, and assume $B$ is in terminal form. Thus $\Ct \brcs{a,d} = 1$ and $0 \leq a_i < d$.

Now let $n-1 > r > 1$ be an integer, and $d_{r+1}, \dots , d_n $ be integers exceeding $1$ \st $d = \prod d_j$. Form the matrix $C =  \(\smallmatrix \I_{r}& X \\ 0 & \diag (d_{r+1}, \dots , d_n) \\ \endsmallmatrix\)$; also assume that $C \in \NS$, so that the content of any column is one. Here $X \in \Z^{r \times (n-r)}$.

Let $P$ be a permutation matrix. We want to establish conditions (in terms of all the variables) so that $UB = CP$ for some $U \in \gl(n,\Z)$. Since $\det B \neq 0$, $B^{-1} $ exists as an element of $\Mn n \Q$, and so existence of such a $U$ implies $ CPB^{-1} \in \Mn n \Z$; but this is also sufficient as $|\det CPB^{-1}| = 1$. As $B^{-1}$ is particularly easy to calculate, the conditions are not difficult to obtain.

Let $\pi$ denote the permutation corresponding to the action of $P$ on the right; that is, if right multiplication sends the $i$th column to the $j$th column, then $\pi(i) = j$. We have (zeros are omitted)
$$
U_0:= CPB^{-1} = \frac 1d \( \matrix & \I_r &&X&&\\
&& d_{r+1} &&& \\
&&& d_{r+2} && \\
&&&&\ddots& \\
&&&&& d_n \\\endmatrix \) P  \( \matrix d&& &&&-a_1\\
&d&  &&&-a_2 \\
&&\ddots&  && \vdots \\     
&&&&d&-a_{n-1} \\
&&&&& 1 \\ \endmatrix \).
$$
Let $S \subset \brcs{1,2,\dots, n}$ be the image of $\brcs{1,2,\dots, r} $ under $\pi$, and $T$ its complement. First, we must have $n \in S$. If not, say $n = \pi(k) $ (with $k> r$), then the $k$th row of $CP$ is just $(0 \ 0 \ 0 \dots \ 0 \  d_k)$. Thus the $kn$ entry of the product is $d_k/d$, which by hypothesis is not an integer. Thus $n \in S$.

We calculate the $kn$ entry of the product for where $k > r$; set $t = \pi(k)$. By the preceding, $t \neq n$. The $k$th row of $CP$ is $d_k E_{t}$ (where $E_i$ are the standard basis elements of $\Z^{1 \times n}$). Thus the $kn$ entry is $-d_k a_t/d$.
We deduce the following
\item{(1)} for all $t\in T$, $d/d_{\pi^{-1}(t)}$ divides $a_t$.

Now we calculate the $ln$ entries of the product for $l \leq r$; the crucial case is $m = \pi^{-1}(n) $. The $m$th row of $CP$ has $1$ as its final entry, zeros in the entries corresponding to $S$, and various $x$s (entries of $X$, too  complicated to establish a notation for) in the entries corresponding to $T$. Then the $mn$ entry of $dU_0$ is $-\sum_{t \in T} x_{m,\pi^{-1}(t)} a_{t} + 1$, so this expression  is divisible by $d$. This yields 
\item{(2)} $\gcd\brcs{\brcs{a_t}_{t \in T}\cup \brcs{d}} = 1$.

It also yields a corresponding condition on the $m$th row of $X$.

For $l \neq m$, the condition is a let-down. In this case, the $l$th row of $CP$ has a one in the $\pi(l) \neq n$ position and various $x$s located in coordinates corresponding to $T$ (which does not include $n$). This  yields that $-a_{\pi(l)} + \sum x_{l,t} a_t $ is divisible by $d$. So we obtain the additional (semi-) condition.

\item{(2$\slfrac12$)} every $a_s$ (for $s \in S \setminus\brcs{n} $ is an additive combination of $\brcs{a_t}_{t \in T}$ modulo $d$.

Suppose $a$ is given, and we want to decide whether $C$ and $P$ exist so that with $CPB^{-1}$ is an integer matrix. Then conditions (1), (2) are necessary, and (2$\slfrac12$) is a consequence of (2). Moreover, conditions (1) and (2) imply something drastic about the $d_i$s, namely that they must be mutually coprime (that is, $\gcd\brcs{d_i,d_j} = 1$ if $i \neq j$).

To see this, from (1), we may write $a_t = h_t (d/d_{\pi^{-1}(t)})$, which we can rewrite as a product of all the $d_i$s with $d_{\pi^{-1}(t)}$  replaced by $h_t$. If $p$ is a prime dividing both $d_i$ and $d_j$ (with $i \neq j$), then it obviously divides all the $a_t$, contradicting (2). We also see that each $h_t$ is relatively prime to $d_{\pi^{-1}(t)}$ (for the same reason). The fact that $B$ is reduced entails $h_t < d_{\pi^{-1}(t)}$ as well, although this does not seem useful.

Hence  $d_i$s are mutually coprime. In particular, if $d$ has exactly $f$ distinct prime divisors, then $n-r \leq f$ (no prime can divide two of the $d_i$s); when $d$ is a power of single prime, this gives an alternative but much more tedious proof of Lemma \onesev, that the $1$-block size is constant on terminal forms PH-equivalent to $B$. This means that if we write $d = \prod p^{m(p)}$ in its prime decomposition, the only factorizations permitted here are those with \st for all $i$, and all $p$ dividing $d$, we must have either $p$ does not divide $d_i$ or $p^{m(p)}$ does, and in the latter case, $p$ cannot divide the other $d_j$s.

Now suppose that $d$ and $a$, the partition $S\dot\cup T$, etc satisfy the necessary conditions (1) and (2) (and their consequences) with corresponding factorization and indexing $d = \prod_{i > r} d_i$. Then we can pick $X$ ($r = |S|$ is already determined) and $P$ so as to construct the corresponding $C$. This is straightforward.
 
As a consequence, we have the following result about non-stability of $1$-block sizes.

\Lem Proposition \twothr.  Let $d $ be a positive integer, and $n> 2$. Suppose  $B$ belonging to $\NS_n$ and with $|\det B | = d$ has the property that every terminal form has $1$-block size $n-1$. Then $d$ is a power of a prime.

In contrast, we have the following sufficient conditions to have a $1$-block of size $n-1$.

\Lem Lemma \nthreig. Let $n \geq 3$. Suppose that $B \in \NS_n$, and let $p,q$ be distinct primes. 
\item{(i)} If  $B$ is PH-indecomposable and $|\det B|  = pq$, then $B \in \NS_{n,n-1}$;
\item{(ii)} If $B$ is PH-indecomposable, $|\det B| = p^r$ where $r\geq 1$, and $J(B)$ is cyclic,  then $B \in \NS_{n,n-1.}$
\item{(iii)}  If $n=3$, $J(B)$ is cyclic, and $|\det B| = p^a q^b$ for some $a,b \in \N$, then $B \in \NS_{3,2}$.

\Rmk  These results do not contradict \twothr, since these say only that at least one terminal form has $1$-block size $n-1$. 

\Pf 
\noindent Suppose $n \geq 4$ and $|\det B| = pq$. We may assume that $B$ is in terminal form; if the form does not already have $1$-block size $n-1$, then its terminal form is $\(\smallmatrix \I_{n-1} & X \\ 0 & \diag (p,q) \\ \endmatrix \)$ (up to possible transposition of the primes). Label the two columns of $X$, $Y$ and $Z$. Each column consists of zeros and numbers between $1$ and $p-1$ (respectively,  between $1$ and $q-1$). Set $S = \Set{j}{Y_j \neq 0} $ and $T = \Set{l}{Z_l \neq 0}$. If $S \cap T $ is empty, then up to a permutation of the indices, the terminal form is a direct sum of matrices, contradicting PH-indecomposability. Hence we may select $k \in S \cap T$. Then $J(B_{\Omega(j)})$ is a quotient of $\Z^2/\langle (p,0),(0,q),(X_k,Y_k) \rangle$, and the relative primeness ($Y_k$ to $q$, $X_k$ to $p$) yields that this is zero. Hence 
$J(B_{\Omega(j)}) = 0$, so $B \in \NS_{n,n-1}$.

\noindent (ii) Write $B = \(\smallmatrix I_{n-k} & X \\ 0 & \Cal D\\ \endsmallmatrix\)$ in terminal form with $\Cal D$ being $k \times k$ upper triangular and having  nontrivial powers of $p$ along its diagonal, the powers appearing in increasing order of size. Suppose $k > 1$ (if $k =1$, then $B$ already has $1$-block size $n-1$). The lower right $2 \times 2$ block is of the form $\(\smallmatrix p^a & x \\ 0 & p^b\\ \endsmallmatrix\)$, where $a \leq b$, and if $x \neq 0$, then $p^a \leq (p^b,x)$. In the latter case, $p^a $ divides $x$. 

Now the image of $E_n+r(B)$ (as an element of $J(B)$) is easily seen (from the upper triangular form) to be of order exactly $p^b$, and thus $z:= p^{b-1}E_n + r(B)$ has order $p$. Similarly $E_{n-1} + (x/p^a )E_{n}+ r(B)$ has order exactly $p^a$ in $J(B)$, and thus $y:= p^{a-1} E_{n-1} + (x/p) E_{n}$ has order $p$ as well.

Since $J(B)$ is cyclic, there is at most one subgroup of each order, and thus there exists $v$ relatively prime to $p$ \st $y - vz \in r(B)$. But this is impossible, as easily follows from the form of $\Cal D$. 

If instead, $x =  0$, then we set  $y = E_{n-1}$, and deduce the same conclusion. 
Hence $k =1$.

\noindent {(iii)} Suppose $n =3$ and $\det B = p^r q^s$. Put $B$ in terminal form; if it does not have a $1$-block of size two, then 
$$B = \(\matrix 1 & a & b \\ 0 & c & d \\ 0 & 0 & f \\\endmatrix\),
$$
 where $(a,p) = 1$, $\cont(b,d,f) = 1$,  $cf = p^r q^s$, and $c \leq (d,f)$. We will show that $B_{\Omega(1)} = (0)$; this implies that $B \in \NS_{3,2}$. Sufficient is that $\cont(ad-bc,cf,af) = 1$. Without loss of generality, we may relabel the primes so that $q|f$. 

If $d=0$, then $J(B) \iso \Z_c \oplus \Z_f$; the latter being cyclic entails that $(c,f) =1$, which forces $c = p^r$, $d = q^s$. Then $(ad-bc,f) = (bp^r,q^s)$; since $p\neq q$, and $(b,q) = 1$, we have $(ad-bc,f) = 1$. Also, $\cont(ad-bc,c,a) = \cont(bp^r,p^r,a) = 1 $, so $\cont(ad-bc,cf,af) = 1$.

If $d \neq 0$, then $(d,f) \geq c > 1$. We may interchange $p$ and $q$ if necessary, and thus assume that $q|(d,f)$. If $q|c$ as well, is routine to see that $J(B)$ cannot be cyclic (since $q$ would divide all the entries of $\(\smallmatrix a & d \\ 0 & f\\\endsmallmatrix\)$). Hence $c = p^t$ for some $1 \leq t \leq r$. Since the content of the third column is one, we must have $(b,q) = 1$. Then $(ad-bc,q) = (bc,q) = 1$. It remains to show that if $p| f$, then $(ad-bc,p) = 1$. 

If  $p$ does not divide $d$, then $(ad-bc,p) = (ad,p) = 1$. 

Hence we may assume that $p|d$. Then $p$ cannot divide $f$, as then it would divide all the entries of  $\(\smallmatrix a & d \\ 0 & f\\\endsmallmatrix\)$, which contradicts cyclicity of $ J(B)$. Thus $f $ can only be a power of $q$, so $c = p^r$ and $f = q^s$. Hence $p$ does not divide $f$, and we are done.
\qed

Example \nthrsix(i) shows that (i) can fail if $|\det B|$ is a product of three distinct primes.
If we try to generalize \nthrsev(i) by assuming  $|\det B| = p^2 q$, $B$ is PH-indecomposable, and  $J(B)$ cyclic, the result fails already at $n=4$:
require $q < p$, and set 
$$
B = \(\matrix 1 & 0 & 0 & 1 \\ 0 & 1 & 1 & p \\
0 & 0 & q & p\\
0 & 0 & 0 & p^2\endmatrix \)
$$
Then $B \in \NS_4$ and  is in terminal form. Since $\gcd(p, q) = 1$, there exists $U \in \gl(2,\Z)$ \st $(q \ \ p)U = (1 \ \ 0)$, and since absolute determinant are preserved, $\(\smallmatrix q &p \\ 0& p^2\endsmallmatrix\)U = \(\smallmatrix 1 &0\\ x & \pm p^2q \endsmallmatrix\)$. It follows that $J(B) $ is cyclic of order $p^2 q$. 

It is also routine to calculate $\Lt{B_{\Omega (i)}} = \Lt{\Z_q,\Z_p, \Z_p, \Z_q}$. Since none of them are zero, $B \not\in \NS_{4,3}$ (\nthrfiv(b)). If $B$ were equivalent to $A \oplus C$ and both $A$ and $C$ have determinants bigger than $1$, then two of the four in the list would have to have direct summands isomorphic to $J(A)$ and the other two would have direct summands isomorphic to $J(B)$. This would force $J(A) \iso \Z_p$ and $J(B) \iso \Z_q$ (or vice versa), entailing $|\det B | = pq$, a contradiction. The only remaining possibility is that $B$ is PH-equivalent to a matrix of the form $1 \oplus C$ where $C \in \NS_{3}$---but this is excluded by the fact that none of the kernels of $p_{\Omega(i)}$ are trivial (\nthrfiv(c)).

Now we can obtain some results about relations between $J(B)$ and $J(B\op)$. 
Let $\Arrow \pi_B; J(\Delta) . J(B)$ and $\Arrow \pi_{B\op}; J(\Delta) . J(B\op)$ be the respective onto maps in the two short exact sequences; these are given by $v + r(\Delta) \mapsto v + r(B)$ and $v + r(\Delta)\mapsto J(B\op)$ respectively. Then   $E_i + r(\Delta)\mapsto E_i + r(B)$ and to $E_i + r(B\op)$ (via $\pi_{B\op}$). It follows from \throne\  that $\pi_B (\ker p^{\Delta}_{\Omega(i)}) = \ker p^{B}_{\Omega(i)}$, and the same with $B\op$ replacing $B$. We claim that $  \ker p^{B}_{\Omega(i)} $, $  \ker p^{B\op}_{\Omega(i)} $,  and $\Z_{m(\Z)} \iso \ker p^{\Delta}_{\Omega(i)}$ are isomorphic to each other.

Since $\pi_B$ restricts to an onto map from $ \ker p^{\Delta}_{\Omega}$ to $\ker p^{B}_{\Omega(i)}$, it suffices to show the map is one to one. If for some positive integer $t$, $t E_i + r(\Delta)$ maps to zero (under $\pi_B$), then $tE_i \in r(B)$, i.e., $tE_i = wB$ for some $w \in \Z^{1\times n}$. From the original definition of $m(i)$, we must have $m(i)$ divides $t$. Hence the restriction of $\pi_B$ is an isomorphism. The same applies with $B\op$ replacing $B$. Thus we have the following.

\Lem Lemma \noneeig. Let $B \in \NS_n$. For each $i$, $\ker p^B_{\Omega(i)} $ and $\ker p^{B\op}_{\Omega(i)} $ are cyclic of order $m(i)$, the isomorphism $ \ker p^{\Delta}_{\Omega(i)} \to \ker p^{B}_{\Omega(i)}$ being induced by the restriction of $\pi_B$.

\Lem Corollary \nonenin. Let $B \in \NS_n$. Then 
\item{(a)} $m(i) = |J(B)|/|J(B_{\Omega(i)})| = |J(B\op)|/|J((B\op)_{\Omega(i)})|$; 
\item{(b)}$|\det B\op| = \frac{|\det B)|^{n-1}}{\prod_{i=1}^n |J(B_{\Omega(i)})|}
$;
\item{(c)} $|J((B\op)_{\Omega(i)})| = \frac{|\det B|^{n-1}}{m(i)\cdot\prod_{i=1}^n |J(B_{\Omega(i)})| }.$

\Pf Part (a) is an immediate consequence of Lemma \noneeig. For (b), we have $\prod_i \(|J(B)|/|J(B_{\Omega(i)})| \) = \prod m(i) $; the latter is $\det \Delta = |J(B)|\cdot |J(B\op)|$, and now we can solve for $|J(B\op)| = |\det B|$. Part (c) is a consequence of (a) and (b).
\qed

Corollary \nonenin(a) entails that $\Delta = \diag(|J(B)|/|J(B_{\Omega(i)})|)$, so is determined by $J(B)$ and $J(B_{\Omega(i)})$Ña small fragment of $\Cal J(B)$. However, it is difficult to see how to obtain the embedding (up to equivalence) $J(B) \to J(\Delta)$ from $\Cal J(B)$.

Parts (b) and (c) imply that $|J(B\op)|$ and the $|J((B\op)_{\Omega(i)})|$ are determined by
$|J(B)|$ and the $|J(B_{\Omega(i)})|$; see the discussion in section 6 concerning the Duality conjecture. 

There is a form stronger than $m(i) = |J(B)|/|J(B_{\Omega(i)})|$. Identify, as usual, $J(B)$ with $r(B\op)/r(\Delta)$. We define $\Delta_{\Omega(i)}$ to be the  {\it square\/} matrix with $i$th row and column deleted, and of course, $J(\Delta_{\Omega(i)})$ is $\oplus_{j \neq i} \Z_{m(j)})$.

\Lem Corollary \noneten. For each $i = 1, 2, \dots,n$ there are short exact sequences $0\to J(B) \to J(\Delta_{\Omega(i)}) \to J((B\op)_{\Omega(i)}) \to 0$ and $0\to J(B\op) \to J(\Delta_{\Omega(i)}) \to J(B_{\Omega(i)})\to 0$.

\Pf Let $\Arrow \pi; J(\Delta) . J(B\op)$ be the quotient map in the original short exact sequence. We have seen that $\pi$ maps $\ker p^{\Delta}_{\Omega(i)}$ isomorphically onto $\ker p^{B\op}_{\Omega(i)}$. This allows us to define $\Arrow {\overline \pi}; J(\Delta_{\Omega(i)}) . J((B\op)_{\Omega(i)})$, via $\pi$; explicitly, $v+ r(\Delta) + E_i \Z \mapsto v + r(B\op) + E_i \Z$.

Now consider $r(B\op)/r(\Delta)$, the kernel of $\pi$. We claim that $r(B\op) \cap E_i \Z \subseteq r(\Delta)$. Pick $wB\op = tE_i + v \Delta$; writing $\Delta = B^T B\op$, we have $(w -vB^T)B\op = tE_i \Z$. From the original definition of $m(i)$, we must have $m(i)$ divides $t$. Since $m(i)E_i \in r(\Delta)$, we have that $r(B)/r(\Delta)$ misses $\ker p^{\Delta}_{\Omega(i)}$. Thus the composed map $r(B\op)/r(\Delta) \to J(\Delta) \to J(\Delta_{\Omega(i)})$ is one to one, and clearly contained in the kernel of $\overline \pi$. Since $|r(B)/r(\Delta)| = |J(B\op)|$ and $|J(B\op)|\cdot |J(B)| = |\det \Delta|/m(i) = \det \Delta_{\Omega(i)}$, the sequence $0 \to r(B)/r(\Delta) \to J(\Delta_{\Omega(i)}) \to J(B\op)_{\Omega(i)}) \to 0$ must be exact.

The other sequence comes from interchanging $B$ with $B\op$. \qed

Let $\Omega \subset 2^S$ and $j \not\in \Omega$. There are natural onto homomorphisms $J(\Delta_{\Omega\setminus \brcs{j}}) \to J((B\op)_{\Omega})$ and $J(\Delta_{\Omega\setminus \brcs{j}}) \to J(B_{\Omega})$; but it is very difficult to relate their kernels to obvious invariants of $B$ and $B\op$ respectively. It is not clear (but likely true) that there are exact sequences $J(B_{\Omega(i)}) \to J( \Delta_{\Omega(i)}) \to J(B\op) $ or $J(B_{\Omega(i)}) \to J( \Delta_{\Omega(i,j)}) \to J(B\op)_{\Omega(j)})$ for all $i \neq j$.

\SecT 4 Size $n$ $1$-block

In this section, we deal with $B \in \NS_n$ \st $B$ has a terminal form with $1$-block size $n$; we write this as $B \in \NS_{n,n-1}$. Computations are relatively tractible, and lead to conjectures for general $B$ in $\NS_n$, that can be proved in our restricted case. We will see in subsequent sections that the density of $\NS_{n,n-1} $ in $\NS_n$ approaches approximately $.845$ as $n \to \infty$ (already at $n =6$, the density exceeds $.8$), so that this special cases covers a large proportion of matrices.

 In addition, computations are also easy in the case that $B\op \in \NS_{n,n-1}$. The density (or even whether it exists) of $\NS_{n,n-1} \cup \NS_{n,n-1}\op$ in $\NS_n$ is not known, but I speculate that it exists and is at least $.99$. 

Perhaps the most important reason for studying this special class is that it is easier to formulate and verify conjectures than in the general case. Results \noneone, \nthrtwo, and \nonenin\ were obtained first for matrices in $\NS_{n,n-1}$, suggesting their validity in general. (Of course, not everything extends in this fashion!)

Let $B = \(\smallmatrix \I_n & a \\ 0 & d \endsmallmatrix\)$ where $a = (a_1, a_2, \dots ,a_{n-1})^T \in \Z^{n-1}$, $d$ is a positive integer exceeding one, $\cont(d,a) =1$, and the entries of $a$ are ordered so that $\gcd(d,a_i)$ are monotone decreasing. We call this a {\it standard form\/} for $C \in \NS_{n,n-1}$ if $C$ is PH-equivalent to $B$ in this form. There can be several standard forms, arising from the column entries being permuted.

Standard forms are terminal, and every $C \in \NS_{n,n-1}$ is PH-equivalent to one in standard form. To arrange the latter, from the definition, there is a matrix of the form $B' = \(\smallmatrix \I_n & a' \\ 0 & d \endsmallmatrix\)$ PH-equivalent to $C$, which is almost in standard form, the only obstruction being that the entries of $a'$ need not be ordered. However, we can conjugate $B'$ by any permutation matrix of the form $Q = P \oplus \brcs{1}$ (where $P$ is a permutation matrix of size $n-1$). Then $Q^{-1} B' Q$ is still in terminal form, but the $a'$ entries have been permuted (according to the permutation induced by $Q$). 

\Lem Proposition \nfouone. Suppose that $B = \(\smallmatrix \I_{n-1}& a \\ 0 & d \\ \endsmallmatrix\) \in \NS_{n,n-1}$ is in standard form.  Then 
$$\eqalign{
B\op = & \(\matrix &&\Cal D &  &  0_{n-1} \\
\\
\frac{-a_1}{(d,a_1)} &\frac{-a_{2}}{(d,a_{2})} &\dots&\frac{-a_{n-1}}{(d,a_{n-1})}&1 \\ 
&&\\
\endmatrix \)\text{ and is PH-equivalent to }\cr 
& \(\matrix 1 & \frac{-a_1}{(d,a_1)} & \dots &\frac{-a_{n-1}}{(d,a_{n-1})}\\ \\
0_{n-1} & &\Cal D \\ \endmatrix \),\cr
}$$
where $\Cal D = \diag\(\frac{d}{(d,a_i)}\) \in \M_{n-1}\Z$. The matrix on the right is in terminal form, and
$m(i) = d/(d,a_i)$. 

\Pf As $B^{-1} =  \(\smallmatrix \I_{n-1}& -a/d \\ 0 & 1/d \\ \endsmallmatrix\)$, we see that $m(i) = d/(d,a_i)$, and $B\op = (\Delta B^{-1})^T$. Conjugating with the obvious cyclic permutation puts it into the indicated form. It is in terminal form, since $(d/(d,a_i))$ is increasing.
\qed

\Lem Lemma \nfoutwo. Let $B = \(\smallmatrix \I_{n-1} & a \\ 0 & d \\ \endsmallmatrix\)$ be in terminal form, with $a = (a_1, \dots, a_{n-1})^T$. Then
\item{(i)} $J(B) \iso \Z_d$ and $J(B\op) \iso \oplus_{i=1}^{n-1} \Z_{d/(a_i, d)}$;
\item{(ii)} if $n \in \Omega$,  then $J(B_{\Omega}) \iso \Z_{\gcd(d,\cont(a_i; i \in \Omega))}$, and if $n\not\in \Omega$, $J(B_{\Omega}) = \brcs{0}$;
\item{(iii)} if $n \in \Omega$, then $J((B\op)_{\Omega}) \iso \oplus_{i \in \Omega\setminus \brcs{n}} \Z_{d/(a_i, d)}$, and  if $n \not\in \Omega$, then $J(B\op)_{\Omega}) $ is isomorphic to any quotient of 
$\oplus_{i \in \Omega\setminus \brcs{n}} \Z_{d/(a_i, d)}$
 by a cyclic subgroup of order $\lcm{\Set{d/(a_i,d)}{i \in \Omega}}$.

\Lem Corollary \nfouthr. Suppose $C \in \NS_{n,n-1}$ and is PH-equivalent to $B$ in standard form with determinant $d$. Then 
\item{(i)} $J(B_{\Omega}) =0$ if $n \not \in \Omega$, and otherwise, $J(B_{\Omega}) \iso \Z_{\cont \(\brcs{a_i}_{i\in \Omega} \cup \brcs{d}\)}$;
\item{(ii)} $J(B\op) \iso \oplus \Z_{d/(d,a_i)}$ and  $J((B\op)_{\Omega}) \iso \oplus_{i \in \Omega} \Z_{d/(d,a_j)}$ if $1  \in \Omega$; $J((B\op)_{\Omega(1)}) \iso J(B\op)/A$ where $A$ is cyclic of order equalling $d$. 

\Pf All the computations easily follow from the forms in the previous result, together with $\Exp J(B\op) = \Exp J(B) =d$. 
\qed

As a consequence of this and \nthrfou, for all $B \in \NS_{n,n-1}$, the lattice $\Cal J(B)$ determines $\Cal J(B\op)$. Whether this is true for all $B \in \NS_{n}$ is unknown (this is the duality conjecture of section 6.

A special case arises when all $J(B_{\Omega(i)}) = 0$. This implies $J(B_{\Omega}) = 0$ for all proper $\Omega \subset S$, and occurs iff $(a_i,d) = 1$ for all $i$. Other special cases will be addressed in section 5. 

\Lem Corollary \nfoufou. Suppose that $B \in \NS_n$ has absolute determinant $d$. Then $J(B_{\Omega(i)}) = 0$ for all $i$ if and only if, in one (or all) of its standard forms, all $a_i$ are relatively prime to $d$. When this occurs,  $J((B\op)_{\Omega}) \iso \Z_d^{|\Omega| -1}$.

\Pf Since one of $J(B_{\Omega(i)})$ is zero, $B \in \NS_{n,n-1}$, and so we can assume $B$ is in standard form. Then (i) of the preceding, with $\Omega = \Omega(i)$ yields $(a_i,d) = 1$ for all $i$. The converse is trivial. The rest follows from \nfouthr(ii).
\qed

Recall what it means for $B \in \NS_n$ to super-split (end of section 1). 

\Lem Lemma \nfoueig. Suppose $B = \(\smallmatrix \I_{n-1} & a \\ 0 & d\\ \endsmallmatrix\) \in \NS_{n,n-1}$ where $a = (a_i)^T \in \Z^{1\times (n-1)}$ 
\item{(i)} Sufficient for $B$ to super-split is that $1 + \sum a_i^2/(d,a_i)$ be relatively prime to $d$.
\item{(ii)} If all $a_i$ are relatively prime to $d$, then the condition in (i) is also necessary for $B$ to super-split.

\Rmk It is probably true that the condition $(d,a_i) = 1$ is unnecessary. 

\Pf Necessary and sufficient for $\Z^{1\times n} B + \Z^{1\times n} B\op = \Z^{1\times n}$ is that the same hold modulo $d$, since $r(\Delta) \subseteq r(B) \cap r(B\op)$. Necessary and sufficient for this to occur is that the set of  all $n \times n$ determinants obtained  from the $2n$ rows of $\(\smallmatrix B\\ B\op\\\endsmallmatrix\)$ has content relatively prime to $d$. If we take the first $n-1$ rows of $B$ and the bottom row of $B\op$, we obtain the matrix
$$
C= \(\matrix \I_{n-1} & a \\ -\frac{a_1}{(d,a_1)}, \dots, -\frac{a_{n-1}}{(d,a_{n-1})}& 1\endmatrix \).
$$
The determinant of this is $1 + \sum a_i^2/(d,a_i)$. This yields (i). 

When $(a_i , d) = 1$ for all $i$, modulo $d$, the only nonzero row of $B\op$ is the bottom one, and it is simply $(-a^T,1)$. Modulo $d$, the bottom row of $B$ is zero---so the only combination of rows to give a nonzero determinant (modulo $d$) consists of the top $n-1$ rows of $B$ with the bottom row of $B\op$, which is  the matrix $C$. Super-splitting thus implies $(\det C,d) = 1$, proving (ii).

It is of interest to give criteria for both $B$ and $B\op$ to belong to $\NS_{n,n-1}$. These will be used when we determine when $\Cal J(B) \iso \Cal J(B\op)$ and the stronger property that  $B $ be PH-equivalent to $B\op$.

\Lem Corollary \nfoufiv. Let $B \in \NS_{n}$. The following are equivalent.
\item{(a)} $J(B)$ and $J(B\op)$ are cyclic;
\item{(b)} $\oplus \Z_{m(i)} \iso \Z_d^2$ for some $d>1$.

\Pf (a) implies (b). Set $d = \Exp J(B) = \Exp J(B\op)$; the groups being cyclic, they are cyclic of order $d$. By \nonesix, the sequence $J(B) \to \oplus \Z_{m(i)} \to J(B\op)$ splits.

\noindent (b) implies (a). Since $\Exp \oplus \Z_{m(i)} = d$, we have $\Exp J(B) = \Exp J(B\op) = d$. From $\left| \oplus \Z_{m(i)}\right| = | J(B)|\cdot |J(B\op)|$; as the exponent of a group is at most the order, we deduce $|J(B)| = |J(B\op)| = d$; since the exponents equal the order, the groups are cyclic. \qed

The following is an obvious consequence of the preceding.

\Lem Corollary \nfousix. 
Suppose that  $B \in \NS_{n,n-1}$, and has determinant of absolute value $d$. The following are equivalent.
\item{(i)} $B\op \in \NS_{n,n-1}$;
\item{(ii)} $J(B\op)$ is cyclic;
\item{(iii)} $|\det B\op| = |\det B|$;
\item{(iv)} $\oplus \Z_{m(i)} \iso \Z_d^2$;
\item{(v)} in (a) standard form, $\gcd(d/(d,a_i), d/(d,a_j)) = 1$ for all $i \neq j$.


\SecT 5 Dual-conjugacy and dual-compatibility

When is ${\Cal J}(B) \iso {\Cal J}(B\op)$ (as lattices of groups), or the stronger condition, $B $ is PH-equivalent to $B\op$? An obvious way to obtain such examples (of the stronger property) is to take $B = C \oplus C\op$ (since, as is evident from \oneele\  or otherwise, $(A \oplus A')\op = A\op \oplus (A')\op$). To avoid such trivial examples, we recall notions of  indecomposability, applied to subgroups of $\Z^{1\times n}$, not just to matrices.

Let $H \subset \Z^n$ be a subgroup of $H$ of full rank, and for which there exists no $m>1$ \st $H \subset m\Z^n$ (this corresponds to the content one condition of all the columns in the corresponding matrix). As usual, let $S = \brcs{1,2,\dots,n}$. We say that $H \subset \Z^n$ is {\it decomposable\/} if there exists a proper subset $T \subset S$ \st $H = H_1 \oplus H_2$ where $H_1 \times \brcs{0} \subset \Z^T \times \brcs{0}$, $H_2 \subset \brcs{0} \times\Z^{S\setminus T} $, and  neither $H_1$ nor $H_2$ is  contained in any $m\Z^n$ for $m >1$. As defined, this is clearly a PH-invariant property. If $H \subset \Z^n$ is not decomposable, then it is {\it indecomposable.}

If we translate this back to square matrices (with $H = r(B)$), then $B \in \NS_n$ is (PH)-indecomposable iff the corresponding subgroup is indecomposable. The same applies to weak indecomposability. 

So we look for PH-indecomposable $B \in \NS_n$ \st either ${\Cal J}(B) \iso {\Cal J}(B\op)$ [$B$ is {\it dual-compatible}] or $B $ is PH-equivalent to $B\op$ [$B$ is {\it dual-conjugate}]. With indecomposability and $1$-block size $n-1$, the first property is fairly drastic; the second property is even more drastic.

We make an obvious comment about the ordered $n$-tuple (not merely the list, with which we have been dealing up to now) $(J(B_{\Omega(i)})_{i=1}^n)$. Suppose $B,C \in \NS_n$, and $B$ is PH-equivalent to $C$. Thus there exists a permutation matrix $P$ and $U \in \gl(n,\Z)$ \st $B = UCP$. The invertible matrix $U$ has no effect on the subgroups spanned by subsets of the columns of $CP$. Let $\pi$ be the permutation induced by $U$, extended in the obvious way to subsets of $S$.  We must have, for all $\Omega$, $J(B_{\Omega}) = J((CP)_{\Omega}) = \Z^{\Omega \pi^{-1}}/\Z^{\Omega \pi^{-1}} C_{\Omega \pi^{-1}}$. If we specialize to $\Omega(i) := S \setminus \brcs{i}$, we have $J(B_{\Omega(i)}) \iso J(C_{\Omega (i\pi^{-1})})$. 

In particular, if $J(B_{\Omega(i)})$ are distinct (meaning, mutually nonisomorphic), then $\pi$ is uniquely determined, and thus $P$ is uniquely determined---so we know exactly which $P$ to use (this can also be extended to the permutation action on ${\Cal J}(B)$, but we never use use this), and thus if $B$ is PH-equivalent to $C$, then the $P$ is uniquely determined (by the $n-1$-tuples of abelian groups $(J(B_{\Omega(i)}))$ and $(J(C_{\Omega(i)}))$). Thus to show PH-equivalence, it is necessary and sufficient that $CPB^{-1}$ have only integer entries. 

In the special case that $P$ must be the identity (that is, $J(B_{\Omega(i)}) \iso J(C_{\Omega(i)}) $ for all $i $ and the $J(B_{\Omega(i)})$ are pairwise nonisomorphic), then the test is merely that $CB^{-1}$ have only integer entries. 
If $B$ is in terminal form with $1$-block size $n-1$, then $B^{-1}$ is especially  simple: if $B = \(\smallmatrix \I_{n-1} & a \\ 0 & d \\ \endsmallmatrix\)$ (where $a = (a_1, \dots, a_{n-1})^T$), then $B^{-1} = \(\smallmatrix \I_{n-1} & -a/d \\ 0 & 1/d \\ \endsmallmatrix\)$.

\Lem Lemma \nfivone. Suppose $B,C \in \NS_n$ and the following conditions hold.
\item{(a)} For all $i=1,2,\dots,n$, $J(B_{\Omega(i)}) \iso J(C_{\Omega(i)})$.
\item{(b)} The $J(B_{\Omega(i)})$ are pairwise nonisomorphic.{\par}
\noindent Then  $B$ is PH-equivalent to $C$ iff  $CB^{-1} \in \M_n \Z$ iff $B$ is Hermite equivalent to $C$.

\Rmk It can happen that (a) holds, but the conclusion does not. There are examples with  $J(B_{\Omega(i)}) \iso \Z_{13}$ for all $i$, as is the case for $J(B\op)_{\Omega(i)})$, and $B$ is PH-equivalent to $B\op$. However $B\op B^{-1}$ is not an integer matrix, so $B$ cannot be Hermite equivalent to $B\op$.

\Pf Suppose $B = ECP$ where $E \in \gl(n,\Z)$ and $P$ is a permutation matrix. Then the permutation corresponding to $P$, call it $\pi$, induces isomorphisms $J(B_{\Omega(i)}) \iso J(C_{\Omega(\pi^{i})})$ directly from the equation. The two conditions (a) and (b) together force $\pi$ to be the identity permutation, hence $P = \I$, and thus $B = EC$.

If $B$ is PH-equivalent to $C$, then $|\det B | = |\det C|$, and thus $\det CB^{-1} = \pm 1$. Thus $CB^{-1} \in \M_n \Z$ entails that $CB^{-1} \in \gl(n,\Z)$.
\qed

\Lem Proposition \nfivtwo. Let $B$ be a weakly indecomposable element of $\NS_{n,n-1}$, and let $|\det B| := d = \prod_{p \in U} p^{m(p)}$ be the prime factorization of the absolute determinant of $B$. Then $B$ is dual-compatible (that is, $\Cal J(B) \iso \Cal J(B)$ as lattices of abelian groups) iff the following holds.
\item{(\dag)}There exists a partition $U = \dot\cup_{i=1}^{n-1} T_i$ with $|T_i| \geq 1$ \st  on  defining $d(i) = \prod_{T_i} p^{m(p)}$, we have $J(B_{\Omega(i)}) \iso \Z_{d/d(i)}$, up to a permutation on the indices.

\Rmk In other words, $B$ is PH-equivalent to 
$$
B'= \(\matrix I_{n-1} & a \\ 0 & d \endmatrix\)
$$
where $a = (a_1, \dots, a_{n-1})^T$ satisfies $d(i) = d/(a_i,d)$. Just observe that $J(B'_{\Omega(i)}) = \Z_{(d,a_i)}$ for $i\leq n-1$; of course, $J(B'_{\Omega(n)}) = (0)$.

\Pf We may assume that $B$ is already in terminal form, and of the form $B = \(\smallmatrix \I_{n-1}& a \\ 0 & d \\ \endsmallmatrix\)$, where $a = (a_1,\dots, a_{n-1})^T$ and $\cont(a_1, \dots, a_{n-1};d) = 1$. We have already seen (\nfouone) that $J(B_{\Omega(i)}) \iso \Z_{(d,a_i)}$ for $1\leq i \leq n-1$, and $J(B_{\Omega(n)}) = (0)$. 

We also have, by \nfouone, that $\det B\op = \prod d/(a_i,d)$ and $J((B\op)_{\Omega(i)}) = \oplus_{j\neq i} \Z_{d/(a_j,d)}$ for $1 \leq i \leq n-1$, and moreover, $J((B\op)_{\Omega(i)}) \iso \oplus \Z_{d/(a_j,d)}/ \Z_{\lcm{d/(a_i,d)}}$. Since $\Cal J(B) \iso \Cal J(B\op)$, we must have $d = \prod d/(a_i,d)$; since $J(B)$ is cylic, so must $J(B\op) $ be; this forces $d/(a_j,d)$ to be pairwise relatively prime. 

Let $U$ be the set of prime divisors of $d$; the fact that $d/(a_j,d)$ are pairwise relatively prime and $d $ is their product forces each $d(i):= d/(a_j,d)$ to be expressible as a product $d(i) = \prod_{T_i} p^{m(p)}$ for some partition $U = \dot\cup T_i$ of $U$. If any of the $T_i$ were empty, we would obtain the corresponding $d_i = 1$, so that $(a_i,d) = d$; that $0 \leq a_i < d$ forces $a_i = 0$. But then $B$ would be, up to a permutation, decomposable as $1 \oplus C$ for some $C \in \NS_{n-1}$, contradicting weak indecomposability. Hence $|T_i| \geq 1$. The rest of necessity is straightforward.

Conversely, suppose that $B$ satisfies the conditions. We may assume it is in terminal form, and the $a_i$ satisfy $d = (a_i,d)d(i)$. Then it is easy to verify that $J(B_{\Omega})$ and $J((B\op)_{\Omega})$ are isomorphic, and the isomorphisms are compatible with the projections, $p_{\Omega,\Omega'}$. 
\qed

\Lem Theorem \nfivthr. Let $B$ be a weakly indecomposable dual-conjugate matrix in $\NS_{n,n-1}$ with $n \geq 3$.  There exists a partition $U = \dot\cup_{i=1}^{n-1} T_i$ with $|T_i| \geq 1$ \st  on  defining $d(i) = \prod_{T_i} p^{m(p)}$ \st $B$ is PH-conjugate to a matrix of the form
$$
B'= \(\matrix I_{n-1} & a \\ 0 & d \endmatrix\)
$$
where $a = (a_1, \dots, a_{n-1})^T$ satisfies $d_i = d/(a_i,d)$, and on writing $v_i = a_i/(a_i,d)$ (where $(v_i,d_i) = 1$), we have
$$
\sum v_i^2 \prod_{p \in T_i^c} p^{m(p)} \equiv -1 \mod d.
$$
Conversely, any such $B'$ is dual-conjugate, and $B$ is Hermite-equivalent to $B\op$.

\Rmk  The congruence condition can be rewritten  in much simpler form, suitable for computing with. It imposes a strong condition on the possible determinants $d$ for which such matrices exist. 

\Pf  By \nfivtwo, we can assume $B$ is already in the form described therein; with $v_i$ defined as $a_i/(a_i,d)$, we have $a_i = v_i\prod_{T_i^c} p^{m(p)}$. We see that $J(B_{\Omega(i)}) = \Z_{(a_i,d)}$ for $1 \leq i \leq n-1$, and $J(B_{\Omega(n)}) = 0$. The set $\brcs{\Z_{(a_i,d)}} \cup \brcs{(0)}$ consists of $n$ distince elements (none of the $a_i$ can be relatively prime to $T$ since the partition is nontrivial).

Now consider $B\op$; this is given in \nfouone, and we have (if $i \neq n$) $J((B\op)_{\Omega(i)})$ is $\oplus_{j\neq i} \Z_{d_j}$; as the $d_j$ are pairwise relatively prime, $J(B\op)_{\Omega(i)}) \iso \Z_{d/\prod_{j\neq i} d_j} \iso \Z_{(a(i),d)} \neq (0)$. Now assume that $B$ is dual-conjugate. Then $B\op \in \NS_{n,n-1}$, so at least one of the collection $J(B\op)_{\Omega(j)})$ must equal zero; hence $j = n$. In particular, we have $J(B_{\Omega(i)}) \iso J((B\op)_{\Omega(i)}) $ for {\it all\/} $i = 1,2, \dots, n$. 

From $B$ PH-equivalent to $B\op$, lemma \nfivone\ applies, and thus the only choice for permutation matrix $P$ \st $B = EB\op P$ with $E \in \gl(n,\Z)$ is $P = \I$;  in particular, $B$ is Hermite-equivalent to $B\op$. Hence $B\op B^{-1}$ is an integer matrix. But $B^{-1} = \(\smallmatrix I_{n-1} & -a/d \\ 0 & 1/d\\ \endsmallmatrix \)$, and the computation of $B\op B^{-1}$ is particularly easy: the constraint that all the entries be integers is exactly the sum of squares condition, resulting from the $(n,n)$ entry of the product. 

The converse is completely straightforward.
\qed

In the case of $n=2$ (in the proof, we used $n \geq 3$ in order to obtain that none of the $a_i$ could be zero---equivalently, divisible by $d$---so that the  $J(B_{\Omega(i)})$ are distinct), the condition on $B = \(\smallmatrix 1 & a\\ 0 & d \\ \endsmallmatrix\)$ (this time, $a$ is just an integer), that $a^2 \equiv -1 \mod d$. Such an $a$ will exist iff all odd primes dividing $d$ are congruent to one modulo $4$, and $4$ does not divide $d$. For larger $n$, the situation is much more complicated.

For example, to obtain $B \in \NS_{n,n-1}$ that is weakly indecomposable, dual-conjugate, of determinant $d$, the partition condition requires that $d$ have  least $n-1$ distinct prime divisors (so that a nontrivial partition of $U$ is possible). If $d$ is a square (and has at least $n-1$ distinct prime divisors), the pair $(d,n)$ can be realized iff $d$ is odd, and every prime divisor is congruent to one modulo four. But if $d/2$ is a square,  the condition is  more complicated: there  should exist an odd prime \st $\Leg -2,p. =1$, and  all primes $q$ with $\Leg -2,q. \neq 1$ must be congruent to $1$ modulo $4$.  And if $d$ is $2^k$ times a square for some $k \geq 2$, the pair cannot be realized at all.

The equation in (\nfivthr) can be rewritten in much simpler form. Write $d = 2^{m(2)} \prod p^{m(p)}$. The Chinese remainder theorem implies solvability of the equation is equivalent to $v_i^2 \prod_{p \in T_i^c} p^{m(p)} \equiv -1 \mod q^{m(q)}$ for every prime $q$ in $T_i$, for all $i$; that is, the negative product  is a square modulo $q^{m(q)}$. For odd $q$, we can replace $q^{m(q)}$ by $q$, and if $q = 2$, by $8$ if $m(2)\geq 3$, and we can delete the condition if $m(q) =2$. So the conditions for the existence of solutions to the equation boil down (after deleting the even powers of primes) to 
$$
\eqalign{
\Leg -1,q. \prod_{p\in T_i^c; \text{ odd } m(p) } \Leg p,q. = 1&\text{ for all odd $q \in T_i$, for all $i$,} \cr
-\prod_{p\in T_i^c; \text{ odd } m(p) } p\quad  \text{ is a square modulo $2^{m(2)}$}&\text{ if $2 \in T_i$} \cr
}$$ 
If $m(2) = 0$ or $1$, the last condition is vacuous; if $m(2)\geq 3$, the term $2^{m(2)}$ can be replaced by $8$. 

These are fairly drastic conditions on the possible determinants. For example, if $d$ is an odd square, then all the prime divisors of $d$ must satisfy $\Leg -1,p. = 1$, that is, $p\equiv 1\pmod 4$, and if $d$ is an even square, then there are no partitions possible, that is, there does not exist $B \in \NS_{n,n-1} $ \st $\det B = \pm d$ and $B$ is PH-equivalent to $B\op$. In some other situations, some partitions will work and others won't. 

Things change if we ask merely for indecomposable dual-conjugate matrices in $\NS_{n}$ (note the switch to indecomposable: for any $B$ in $\NS_ n$, $B\oplus B\op$ is trivially an dual-conjugate member of $\NS_{2n}$). For example, we can realize $(p^3,3)$ for any odd prime $p$ (not all of these are Hermite-equivalent to their opposite, unlike the situation in Theorem \nfivthr), while $(8,3)$ is not realizable---but $(8,4)$ is. Since in these examples $U$ consists of a single prime, the situation is obviously quite different when we drop the requirement on the $1$-block size.

\SecT 6 Duality?

Now we refer to definitions and results in Appendix A. We discuss what we have called the duality conjecture (briefly, $\Cal J(B) \iso \Cal J(B')$ implies $\Cal J(B\op) \iso \Cal J((B')\op)$), and prove it for a class of matrices.

Let $n > k$, let $d$ be an integer, and let $X$ be an $(n-k) \times k$ matrix with integer entries between
$0$ and $p-1$ inclusive; assume that the content of each column of $X$ is relatively prime to $d$. Form the matrices (as $X$ varies), 
$$
B(X) = \(\matrix I_{n-k} & X \\ 0 & dI_k \\\endmatrix \).
$$
Each of these is in terminal form. Obviously $J(B(X)) \iso \Z_d^{k}$, and $\Exp J(B(X)) = \Z_d$. So we can regard $J(B(X))_{\Omega}$ as $\Z_d$-modules. It is easy to calculate  $J(B(X)_{\Omega(i)})$. If every row of $X$ has content relatively prime to $d$ (a reasonable assumption), then $J(B(X))_{\Omega(i)} \iso \Z_d^{k-1}$ for every $i$. If we further
assume that every $j \times j$ submatrix of $X$ has nonzero determinant  which
is relatively prime to $d$ (in particular, all entries of $X$ are
units modulo $d$), then $J(B(X)_{\Omega}) \iso \Z_p^{|\Omega|-1}$. This is equivalent to the matrix
$ \(\smallmatrix X \\ I_k\endsmallmatrix \)$ belonging to $F_{{n\choose k}}(n,k)$, if we regard
the entries of $X$ as belonging to $\Z_d$ (see the comment between Propositions \nappfiv\ and \nappsix\ in Appendix\,A).

\def\Aut{\text{Aut\,}}
\def\End{\text{End\,}}


In these case, the lists $\Lt {J(B(X))_{\Omega}}_{|\Omega| = j}$ are thus useless for distinguishing $\Cal J (B(X))$ from $\Cal J(B(X'))$. However, we can say when $\Cal J (B(X)) \iso\Cal J(B(X'))$, by appealing to the orbit spaces under the actions of $W(n) \times \gl(k,\Z_d)$, as discussed in Appendix\,A.

Let us make a minimal assumption on the $X$s: they have no zero rows (modulo $d$). This is equivalent to $B(X)$ being weakly indecomposable, and entails that $\ker p_{\Omega(i)}$ is nonzero (Lemma \nthrfiv).

For each $B(X)$, define an explicit isomophism $J(B(X)) \to \Z_p^k = \Z_p
\times \Z_p \times \dots \times \Z_p$ (with $k$ copies of $\Z_p$), sending (for $j = 1,
\dots, k$) $E_{n-k+j} + r(B(X))$ to $e_j := (0,\dots,0,1,0,\dots,0)$, the $1$
appearing in the $j$th position. Since the kernel of $ J(B(X)) \to
J(B(X)_{\Omega(i)})$ is $E_i + r(B(X))$ (here $i \in \brcs{1,2,\dots, n}$), we
can identify the kernels with the following subgroups of $\Z_p^k$: for $i \leq
n-k$, $\langle -r_i (X)\rangle$ (the $i$th row of $X$, viewed as an element of
$\Z_p^k)$, since $E_i + (0^{n-k}, r_i(X)) \in r(B(X))$, and for $i > n-k$,
$e_{i-(n-k)}$.

Putting the generators into an $n \times (n-k)$ column with entries in $\Z_p$,
we obtain the matrix $M(X):= \(\smallmatrix -X \\ I_k \\ \endsmallmatrix \)$
(although the minus sign plays no role in terms of subgroups, it does play a
role when we work out the corresponding $\phi$). Obviously $M(X) \in F(n,k)$.

Now $\Aut J(B(X)) = \gl(k,\Z_d)$, and we have a natural action of $W(n) \times \gl(k, \Z_d)$ on $F(n,k)$. If two points, say $M(X)$ and $M(X')$ are in the same orbit, then there exists an automorphism $\Arrow \psi; J(B(X)).J(B(X'))$ \st the subsgroups match, that is, there exists a permutation $\pi$ \st $\psi(\langle E_i + r(B(X))\rangle) = \langle E_{\pi i} + r(B(X))\rangle $. This is exactly the condition, $\psi (\ker p^{B(X)}_{\Omega(i)}) = \ker p^{B(X')}_{\Omega(\pi i)}$ discussed in \nthrthr. Hence $M(X)$ being in the same orbit as $M(X')$ implies that $\Cal J(B(X)) \iso J(B(X'))$. The converse is straightforward.

For example, suppose $d = p$ a prime, $n = 5$, and $k=2$. Form $F_{10}(5,2)$ over $\Z_p$. The condition that $M(X) $ belongs to this is simply that every entry of $X$ is relatively prime to $p$, all pairs of rows of $X$ are linearly independent modulo $p$ (that is, the three $2\times 2$ determinants are invertible modulo $p$). When $p = 5$ or $7$, there is only a single orbit (that is, $W(5) \times \gl(2,\Z_p)$ acts transitively on $F_{10}(5,2)$), whereas when $p > 7$, the action is not transitive. In the former case, $M(X),M(X') \in F_{10}(5,2)$ implies $\Cal J(B(X)) \iso \Cal J(B(X'))$; but in the latter ($p > 7$), we can choose $M(X), M(X')$ in different orbits, and then $\Cal J(B(X)) \not \iso \Cal J(B(X'))$.

From the earlier comments, $M(X) \in F_{10}(5,2)$ implies $J(B(X)_{\Omega}) \iso \Z_p^{k+|\Omega|-n}$ for all $\Omega$, i.e., $J(B(X)_{\Omega})$ depends only on $\Omega$. If we put $p =11$, the presence of more than two orbits yields an example of two matrices, $B$ and $B'$ with the property that $J(B)_{\Omega} \iso J(B')_{\Omega}$ for all $\Omega$, but $\Cal J(B) \not \iso \Cal J(B')$.

\vskip4pt\noindent {\it Duality.} We state the duality conjecture.

\Lem Duality Conjecture. Suppose $B, B' \in \NS_n$ and $\Cal J(B) \iso \Cal J(B')$. Then $\Cal J(B\op) \iso \Cal  J((B')\op)$.

This is known if $J(B)$ is cyclic, or if either of $B$ or $B\op$ belongs to $\NS_{n,n-1}$ (Corollary \nthrfou). The conjecture is also true when both $B = B(X)$ and $B' = B(X')$ above, as we will show. We will put the conjecture in the form of a possible generalization of the dualities established in Appendix A.  There is also a stronger form. 

\Lem Constructive Duality Conjecture. Determine $\Cal J(B\op)$ from $\Cal J(B)$. 

A small step in this direction appears in \nonenin: $\Delta$, $\det B\op$, and $|J((B\op)_{\Omega(i)})|$  are determined from $\det B$ and the ${|J(B_{\Omega(i)})|}$.

Suppose that $B = B(X)$ and $B' = B(X')$. Then $B^{-1} = \(\smallmatrix \I_{n-k} & -X/d \\ 0 & \I_k/d \\\endsmallmatrix\)$. If we assume that all rows of $X$ have content relatively prime to $d$, then $\Delta = d\I$, and thus $B\op = \(\smallmatrix d\I_{n-k} & 0 \\ -X^T & \I_k \\\endsmallmatrix\)$. It is not in terminal form, but this does not matter. We wish to verify the duality conjecture for a subclass of these matrices. 

It is easy to check that $J(B\op) \iso \Z_d^{n-k}$ (this also follows from the short exact sequence, $0\to J(B) \to \Z_d^{n} \to J(B\op) \to 0$; here $I = \Z^n/d\Z^n$ since $\Delta = d\I_n$), so that its automorphism group is $\gl(n-k,\Z_d)$. Identifying the kernels of $p^{B\op}_{\Omega(i)}$ with the rows of $X^T$, we form the analogue of $M(X)$, that is, $N(X) = \(\smallmatrix \I_{n-k}\\ X^T \endsmallmatrix\)$.

Then $M(X)^T N(X) = 0$, and it easily follows from Appendix A that the map $\phi$ therein sends $[M(X)] \to [ N(X)]$. We do the same thing for $M(X')$ and $N(X')$, and then we have the sequence of implications (from the main result of the Appendix),
$$\eqalign{
\Cal J(B) \iso \Cal J(B') \implies &M(X), M(X') \text{ are in the same orbit} \implies \cr 
N(X), N(X') \text{ are in the same orbit} \implies &\Cal J(B\op) \iso \Cal J((B')\op).\cr
}$$
So the duality conjecture is true for matrices in this class. \qed

The duality conjecture can be rephrased so that it vaguely resembles the results in Appendix\,A.
Let $B \in \NS_n$ and define $d= \Exp J(B) = \Exp J(B')$ (\noneone); we view $J(B)$ and $J(B\op)$ as $\Z_d$-modules. Let $E = \End J(B)$ and $E^o = \End J(B\op)$. Then $\Aut J(B)$ is just the group of units of $E$, and $\Aut J(B\op)$ is the group of units of $E^o$.

The centres of $E$ and $E^o$ are both $\Z_d$ (this is true for the endomorphism ring of any finite abelian group with exponent $d$). Pick a representative for a generator of each of $\ker p^B_{\Omega(i)}$, and form them into a column of size $n$, that is, an element of $J(B)^n$. Let $\Cal P_n$ be the group of permutations of $n$-element sets, and defined $\Cal D_n$ to be the diagonal matrices with entries from $\Z_d^{\times}$, and define $W(n)$ to be $\Cal  P_n \Cal D_n$. Then we view $J(B)^n$ as a set with the obvious $W(n) \times \Aut J(B)$ action. We do the same with $B\op$. Then the duality conjecture boils down to a bijection between these orbit spaces.

In the $B(X)$ examples, the corresponding rings $E$ and $E^o$ are just $\M_k \Z_d$ and $\M_{n-k} \Z_d$, and in particular, they are Morita equivalent; moreover, $J(B)$ and $J(B\op)$ are free $\Z_d$-modules. In general, $E$ and $E^o$ are not Morita equivalent and neither $J(B)$ nor $J(B\op)$ need be free $\Z_d$-modules.

In addition, the condition on the elements of the column, that they generate a cyclic subgroup corresponding to a $\ker p_{\Omega(i)}$ is somewhat restrictive. For example, if $n =3$ and $J(B)\iso \Z_{p^2} \oplus \Z_p$ (lots of such examples exist), then $J(B_{\Omega(i)})$ must be cyclic (by \nthrone\ and \nthrtwo). Hence we must rule out $p\Z_{p^2} \oplus 0$ as a subgroup appearing as $\ker p_{\Omega(i)}$, hence $(p,0)$ cannot appear as an entry in the column.

Perhaps the key feature of the $B(X)$ matrices is that $J(B(X))$ and $J(B(X)\op)$ are free $\Z_d$-modules, and thus $E$ is Morita equivalent to $E'$. In addition, their ranks add up to exactly the right number, in order that the duality of Appendix A can be applied; this is a consequence of  $\Delta = d\I$. 

\noindent {\it A bilinear function.} The identification of $J(B)$ and $J(B\op)$ with subgroups of $J(\Delta):= \Z^{1\times n}/r(\Delta)$ leads to a bilinear function, potentially useful for the duality conjecture. 

Given $B \in \NS_n$, $B\op$, and $\Delta$, as usual, let $d = \Exp J(B) = \lcm\brcs{m(i)}$ where $\Delta = \diag (d_i)$. Then $d\Delta^{-1}$ is an integer matrix, and we define the bilinear function, $\Z^{1\times n} \times \Z^{1\times n} \to \Z$ given by $\Ip v,w. = vd\Delta^{-1}w^T$. This clearly induces a faithful bilinear function $\Z^{1\times n}/\Z^{1\times n}\Delta \times\Z^{1\times n}/\Z^{1\times n}\Delta \to \Z_d$, denoted $\Ipd \overline v,\overline w. $, which is $vd\Delta^{-1}w^T$ modulo $d$, and the overlines indicate equivalence classes modulo $r(\Delta)$.

Recall from the discussion at the end of section 1, the two subgroups of $J(\Delta)$, $Y(B) := r(B\op)/r(\Delta) \iso J(B)$ and  $Y(B\op):= r(B)/r(\Delta) \iso J(B\op)$, which are the images of $J(B)$ and $J(B\op)$ in the two short exact sequences discussed therein.

Now we note that  $Y(B)$ and $Y(B\op)$ are dual (even if they are equal, as could well be the case, e.g., if $B$ is Hermite-equivalent to $B\op$ (examples appear in Theorem \nfivthr). Specifically, if $\dIp \overline x \overline B,\overline y. = 0$ for all $x$, then $\overline y \in Y(B\op) = r(B)/r(\Delta)$. To see this, we have $xB d\Delta^{-1} y^T \in d\Z$ for all $x$; then $B\Delta^{-1}y^T $ has only integer coefficients.  Replacing $\Delta^{-1}$ by $B^{-1}((B\op)^T)^{-1}$, we see that $((B\op)^T)^{-1}y^T $ has only integer coefficients; applying the transpose, it follows that $y(B\op)^{-1}$ has integer coefficients, and thus $y \in r(B\op)$. The reverse inclusion is trivial. So the dual of $Y(B\op)$ is $Y(B)$ \wrt this bilinear function on $J(\Delta)$, and vice versa.  

In particular, we can realize elements of $J(\Delta)$ as $\Z_d$-module homomorphisms  $J(\Delta) \to \Z_d$, with those that kill $Y(B\op)$ coming from elements in $Y(B)$ (and again, vice versa).

So now the duality conjecture can be translated to this setting. Pick a set of $n$ elements of $J(B)$ (or better, $Y(B)$), each generating the cyclic subgroup which is the kernel of $p^B_{\Omega(i)}$, and form them into a column, $M$, that is, an element of $Y(B)^{n\times 1}$. On the right, $\Aut J(B)$ acts, and on the left, $\Cal P_n \Cal D_n$, where $\Cal D_n$ consists of diagonal matrices with entries in $\Z_d^{\times}$. We do the same with $J(B\op)$. Each entry of $M$ can be viewed as a module homomorphism $J(\Delta) \to \Z_d$, so we can view $M$, essentially the transpose, as a $\Z_d$-module homomorphism $\Arrow \tilde M; J(\Delta)^n. \Z_d^n$. Then the kernel should correspond to the analogous matrix made out of $B\op$, rather than $B$, as in the arguments in Appendix A. But it is not clear how to proceed.

The identifications of  $J(B)$ with $r(B\op)/r(\Delta)$ and the corresponding one interchanging $B$ with $B\op$ are particularly well-behaved \wrt applying $p_{\Omega}$. We can create $\Delta_{\Omega}$, obtained by deleting all the rows and columns indexed by an integer not in $\Omega$, and it is easy to check that $((B\op)_{\Omega})^T B_{\Omega} = \Delta_{\Omega}$ (let $c_i (\cdot)$ denote the $i$th column; then $(B\op)^T B = \Delta$ simply means $c_i (B\op)^T c_j(B) = m(i) \delta_{ij}$, and the columns of $B_{\Omega}$ and $B$ are identical  unless they are completely eliminated. The results in \noneeig--\noneten\ suggest that more can be done along these lines.

\SecT 7 Densities for PH-equivalence to $1$-block size $n-1$

\def\TF{{\Cal T\!\Cal F}}
Here we give estimates for the likelihood that a matrix $B \in \NS_n$ has a terminal form with $1$-block  of size at least $n-1$. Although we give an explicit formula, valid for each $n$, it is difficult to compute with; however, it   converges (as $n \to \infty$) to a product of two known constants, the Landau totient, and the reciprocal of $\prod_2^{\infty} \zeta(k)$,
$$
\frac{\zeta(2)\cdot \zeta(3)}{\zeta(6)}  \cdot \frac 1{\zeta(2)\zeta(3)\zeta(4) \dots} \sim .845
$$
This is almost double the likelihood that $B \in \Mn n\Z$ has a Hermite normal form with at least $n-1$ ones [MRW]. The methods  derive from that  reference, with a few added twists.

First, we obtain  an upper bound. Suppose that $B \in \Mn n \Z$. Then $B \in \NS_n$ iff modulo every prime, each column is not zero. That by itself together with usual notion of natural density (see [MRW] for very clear explanations) says that the likelihood that $B$ is in $\NS_n$ is $1/\zeta(n)^n = 1 - n2^{-n} - \Oh{n3^{-n}}$, which goes to one  quickly.

Now suppose that $B \in \Mn n \Z$ is PH-equivalent to a terminal form having $1$-block size at least $n-1$. Then for every prime $p$, the matrix $B + p\Mn n \Z \in \Mn n \Z_p$ has rank at least $n-1$. The converse fails---examples are ubiquitous. Let $\TF_n$ denote the collection of matrices in $\Mn n \Z$ PH-equivalent to a matrix with at least $n-1$ ones in its terminal form (for large $n$, $\TF_n \cap \NS_n$ is practically the same as $\TF_n$, so we do not require members of the latter collection to be in $\NS_n$). This does give an upper bound for the natural density (assuming it exists) of $\TF_n$.

In fact, we can do a bit  better. For fixed $n$ and for every prime $p$, let $\Arrow \pi_p; \Mn n \Z. \Mn n \Z_p$ be the usual modulo $p$ onto homomorphism. We define a property for $n \times n$ matrices in terms of its reduction modulo every prime. We say that a matrix $B \in \Mn n \Z$ is of {\it deficiency\/} at most  $s$ if for every prime $p$, the image, $\pi_p (B)$ has rank at least $n-s$. For fixed $n$, the collection of these has a natural density, and if $ n \geq (s+1)^2$, it is
$$
\frac 1{\psi_{(s+1)^2 + 2} \cdot \zeta((s+1)^2) \cdot \zeta((s+1)^2 + 1) \cdot \dots \cdot \zeta(n)},
$$
where $\psi_{(s+1)^2 + 2} $ is defined as $\prod_p f(1/p) $ where $f$ is  a function (given explicitly below) with the property that $f(z) = 1 - z^{-(s+1)^2 + 2} + \Oh{z^{-(s+1)^2 + 3}}$ (except for small primes, the product is more or less $\zeta((s+1)^2+2)$). At $s =1$ (so for $n \geq 4$), the outcome is at least $.845$, at $s =2$, it is bigger than $.99$, at $s = 3$, it is at least $.9999$, and each addition of  one to $s$ results in the difference from one approximately squaring. 

The case of $s=1$ gives the upper bound.

However, when we look at the original problem, density of $\TF_n$, the situation is more complicated, and the best we can do is to use the inclusion-exclusion principle to obtain a formula, which is difficult to evaluate, except for small or large $n$.

Throughout this section, we refer to {\it natural density\/} of families of integer matrices, although most of the effort is spent on counting matrices modulo primes, and multiplying the results over all the primes. The problem is then to relate the relatively easily obtained infinite product expressions to the usual or somewhat stronger notion of natural density, as discussed, for example, in [MRW, Ma].

The methods of [op\,cit] can be used to justify the expression natural density, and we will outline what has to be done, at various points.

\noindent {\it Upper bound.} Fix integers $s,n$ with $ n > (s+1)^2 +1$ and let $p$ be a prime. The normalized number of matrices in $\Mn n\Z_p$ of rank at least $n-s$  (that is, divided by the cardinality of $\Mn n \Z_p$, which is $p^{n^2}$) is given by Landsberg's theorem [L] (quoted in Appendix A) as $\left.\(\prod_{i=1}^{n} (1-z^i)\) (1 + \sum_{1 \leq j \leq s} c_j(z))\right|_{z=1/p}$ where
$$
c_j (z) = \frac{z^{j^2}(1-z^n)(1-z^{n-1})\dots (1-z^{n-j+1})}{(1-z)^2(1-z^2)^2\dots (1-z^{j})^2},
$$
although for some computations we could take the simplified (and slightly less accurate)
$$
c_j \sim \frac{z^{j^2}}{(1-z)^2(1-z^2)^2\dots (1-z^{j})^2}.
$$
By  Proposition \apptwo, the Maclaurin series of $a_s = \(\prod_{i=1}^{(s+1)^2 -1} (1-z^i)\) (1 + \sum_{1 \leq j \leq s} c_j(z))$ (or $c_j$ replaced by its simpler form) expands as $1 - z^{(s+1)^2 + 2} + {\text{higher order terms}}$. Then the normalized number of matrices of rank at least $n-s$ in $\Mn n \Z_p)$ is
$$
n_{s,p} := a_s \cdot \left.\prod_{i = (s+1)^2}^n (1-z^i)\right|_{z=1/p}
$$

 Form the infinite product $\psi_{n,s} = \prod_{p} a_s (1/p)$ (this converges---very fast---since $(s+1)^2 + 2 \geq 2$). Then $\prod_p n_p$ is $\psi_{n,s} /\(\zeta((s+1)^2)\cdot \zeta((s+1)^2 + 1)\cdot \dots \cdot \zeta (n)\)$. For very large $n$, $\psi_{n,s}$ is extremely close to $1$ (just as $\zeta(n)$ is). So as $n \to \infty$, the limiting value is
$$
\frac 1{ \prod_{j \geq (s+1)^2} \zeta (j)}.
\tag 1$$

The case of interest occurs when $s =1$, and  an easy computation reveals that $a_1 = 1 - z^6$ (exactly!). Hence
$$\eqalign{
\prod_p n_{1,p} & = \frac 1{\zeta(6)\cdot \prod_{j=4}^n \zeta (j)}\cr
& = \frac{\zeta(2) \zeta(3)}{\zeta(6)} \cdot \frac 1{\ \prod_{j=2}^n \zeta (j)}.
}\tag2$$
The left factor is Landau's totient constant (On-line Encyclopedia of Integer Sequences [oeis] A082695); about $1.94\dots$; the right factor, for large $n$, is about $.436$ [MRW] (with extremely fast convergence in $n$), so the product is about $.845$ or so. As $n$ increases, the value decreases.

When $s=2$, the limiting value in (1) is in excess of $.99$, and when $s=3$, the limiting value exceeds $.9999$ (with the distance from $1$ approximately squaring with each addition of $1$ to $s$). 

To check that the expressions $\prod_p n_p$, (1), and (2) really do represent natural densities (that is, the number of $B \in \Mn n \Z$ with all entries in $[-N,N)$ \st for every prime $p$, the rank of $\pi_p(B) \in \Mn n \Z$ is at least $n-s$, divided by $(2N)^{n^2}$, tends as $N \to \infty$ to the corresponding expression), we note that the method of [MRW] works almost verbatim. Specifically, the Chinese remainder theorem argument in the proof of [MRW, Lemma 3] applies here, as does the argument of [MRW, Lemma 4]. This is made easier by the fact that we are defining the property of matrices in terms of properties modulo every prime. In contrast, when we deal with $\TF_n \cap \NS_n$, there does not appear to be simple characterization of the set by properties modulo $p$.

In particular, if $ n \geq 6$,  the density of matrices   $M \in \Mn n \Z$ with the property that for every prime $p$, the rank of $\pi_p (M)$ is at least $n-1$ is given by the expression in (2), and is at least the limiting value as $n \to \infty$. This gives an upper bound for the (upper) density of matrices \st $M \in \TF_n \cap \NS_n$.

\noindent {\it Counting $\TF_n \cap \NS_n$.}
First, we count the number of matrices $b \in \Mn n \Z_p$ the leftmost $n-1$ columns form a linearly independent set, and the last column is not zero. (If this happens modulo $p$ for every prime $p$, then the original matrix belongs to $\NS_n \cap \TF_n$.) This is almost the same as a special case of  [M; Corollary 7].

There are $N_p = (p^n -1)(p^n - p)\cdots (p-1)$ full rank matrices. If the last column is dependent on the preceding $n-1$ columns and they form a linearly independent set, then we can write it $c_n = \sum_{i < n} a_i c_i$; since we have required that the last column be not zero, we must also have $(a_i) \neq (0,0,\dots,0)$, and every such choice will do. The number of $(n-1) \times n$ matrices of full rank is just $N_p/(p^n - p^{n-1})$. Thus the total number of matrices whose  set of leftmost  $n-1$ columns is not zero  and whose $n$th column is not zero is
$$\eqalign{
N_p \cdot \( 1 + \frac{p^{n-1} - 1}{p^n - p^{n-1}}\) &=  p^{n^2} (1-1/p)(1-1/p^2) \cdots (1-1/p^n)\( 1 + \frac{1-1/p^{n-1}}{p(1-1/p)} \)\cr
& =  p^{n^2} (1-1/p)(1-1/p^2) \cdots (1-1/p^n)\frac{1 - 1/p^n}{1-1/p} \cr
& = p^{n^2} (1-1/p^2) \cdots (1- 1/p^{n-1})(1-p^n)^2. \cr
}$$

This yields that the natural density (see below) of $B \in \NS_n$ \st removing the last column yields a matrix with full row space (equivalently, the Hermite normal form of $B$ is $\( \smallmatrix I_{n-1}  & a \\ 0 & d \\ \endsmallmatrix\) $) is
$$
 \frac 1{\zeta(2)\cdot \zeta (3) \cdot \dots \zeta(n-1)\cdot \zeta(n)^2} . \tag *
$$
This differs from the natural density of matrices with Hermite normal form with at least $n-1$ ones [Ma, Corollary 7] only by the extra factor of $1/\zeta(n)$, which appeared because we insisted that the last column be nonzero (in order to ensure that it came from a matrix in $\NS_n$).

As in all of these computations, the $1-1/p$ factor that appears in $N_p/p^{n^2}$ has conveniently been wiped out, thereby removing the singularity that would have arisen from $\zeta(1)$. If $\Phi$ is a subset of $\brcs{1,2, \dots, n}$, let $D_{\Phi}$ be the set of matrices in $\NS_n$ \st for {\it every\/} $j \in \Phi$, the $\gcd$ of the $(n-1)\times (n-1)$ determinants of the matrix with the $j$th column deleted is one. Clearly, if $|\Phi| = |\Phi'|$ and $D_{\Phi}$ has a natural density, then so does $D_{\Phi'}$ and their natural densities are equal.

That this number is the  natural density for this problem is practically immediate from the special case of [Ma, Corollary 7] with $d_1 = d_2 = \dots = d_{n-1} = 1$ in the notation there---the only (slight) difference is that we have insisted here the the final column be unimodular, so nonzero modulo every prime. This resulted in the extra factor of $\zeta(n)$.

We have just shown that if $|\Phi| = 1$, then $D_{\Phi}$ has a natural density, given by the number  in (*). Now $\cup D_{\Phi}$, where $\Phi$ ranges over all one-element sets, is precisely the set of $ B \in \NS_n$ \st $B$ is PH-equivalent to a terminal form with at least $n-1$ ones.

The inclusion-exclusion formula now can be used. We will obtain a density  for every $D_{\Phi}$. At various points, it will be convenient to use a variable $z$ which will be evaluated at $z = 1/p$ for $p$ prime.

Say $|\Phi| = s> 1$; then we may assume that $\Phi = \brcs{n, n-1, \dots, n-s+1}$, that is, corresponding to the final $s$ columns. Again, if we restrict to invertible matrices, there are $N_p$; otherwise, the first $n-1$ columns  constitute a linearly
independent set, and we can write $c_n = \sum_{i < n} a_i c_i$. Only this time, we also require that if $i \in \Phi$, then $a_i \neq 0$ (this occurs iff the $i$th column can be expressed as a linear combination of all the other columns; it also guarantees all the columns are nonzero). Hence the number of choices for the $(a_i)$ is $p^{n- |\Phi|} (p-1)^{|\Phi|-1} = p^{n-s}(p-1)^{s-1}$. Hence the normalized number of such matrices is
$$\eqalign{
\frac{N_p}{p^{n^2}}\(1 + \frac{p^{n-s}(p-1)^{s-1}}{p^n - p^{n-1}}\)& =
(1-1/p)\dots (1-1/p^n) \(1 + \frac{(1-1/p)^{s-2}}p \); \text{ setting $z = 1/p$,}\cr
& = (1-z)(1-z^2) \dots (1-z^n) \(1 + z (1-z)^{s-2}\)\cr
}$$
Denote by $f_s$ the polynomial (now in the variable $z$) $(1-z)  \(1 + z (1-z)^{s-2}\) $; this is  $(1-z)(1 + z - (s-2)z^2 + \dots)$, so $f_s = 1 - (s-1)z^2 + \Oh{z^3}$. This permits us to define a function (which it turns out is entire),
$$
F(s) :=  \prod_p f_s(1/p) = \prod_p \(1 - \frac{p^{s-1} - (p-1)^{s-1}}{p^s} \).
$$
Provided the (now, complex) $s$ is such that for every prime $p$, $p^{s-1} - (p-1)^{s-1} \neq p^s$ (this simplifies), it is easy to check that $F$ is analytic on a neighbourhood of $s$, and a routine verification assures us that at any of the trivial zeros, $t$, $\lim_{s\to t} F_{s}/(s-t)$ exists and is not zero, hence $F$ is also analytic on neighbourhoods of the zeros; so $F$ is entire. Its zeros are precisely the set, $\Set{s \in \C}{\exists \text{ prime }p \text{ \st } p^s = p^{s-1} - (p-1)^{s-2}}$; this  can be rewritten as
$$\brcs{1 + \frac{(2k+1)\pi i + \ln (p-1)} {\ln \frac p{p-1}} }_{ p \in \Spec \Z,\  k \in \Z}
$$
The reciprocals of the moduli of the zeros is thus absolutely summable along any infinite strip of the form $|\Im z|< N$.

The values of $F$ at various integers are interesting, and will play a role in what follows.
$$\eqalign{
F(0) & = \prod_p \(1 + \frac1{p(p-1)} \); \quad \text{this is $\zeta(2)\zeta(3)/\zeta(6) \sim 1.94$, the Landau totient constant, again}\cr
F(1) & = 1\cr
F(2) & = \prod_p \(1 - \frac1{p^2} \) = \frac 1{\zeta (2)}\cr
F(3) & = \prod_p  \(1 - \frac{2p-1}{p^3} \);  \quad \text{the {\it carefree\/} constant, $\sim .426$ [M]}\cr
}$$
The values at the other integers (both positive and negative) have likely appeared before, but I couldn't locate them in the huge literature on constants.
The density of $D_{\Phi}$ (when $|\Phi| = s > 1$) is thus
$$
\frac {F(s)}{  \zeta(2)\dots \zeta(n)}.
$$

Once again, we may use the methods of [Ma, section 4] to justify the  natural density. With this, we also see that the inclusion-exclusion principle applies  (first to  subsets of $\TF_n \cap \NS_n$ inside $[-N,N)^{n^2}$ and their translations, then letting $N\to \infty$).

For  $s=2 $, the density of $D_{\Phi} $ is $ 1/\zeta(2)^2 \zeta(3) \dots \zeta(n)$.
The inclusion-exclusion principle   reveals that the density of matrices in $\NS_n$ PH-equivalent to a terminal form with $1$-block size at least $n-1$  is
$$
\frac{\frac n{\zeta(n) }- \frac{{n \choose 2}}{\zeta(2)}
+ \sum_{j= 3}^n {(-1)^{j-1}{n\choose j}} F(j)}{\zeta(2)\zeta(3)\dots \zeta(n)}.
\tag{**}
$$
The leading term does not involve $F(1)$, as we would have expected; however, for large $n$, $1/\zeta(n)$ is practically $1 = F(1)$; and we have substituted $F(2) = 1/\zeta(2)$.
Now we have to estimate this. The denominator converges extremely rapidly, and has been calculated as around $.44$ for large (and not so large) $n$ [Ma]. Also, $\brcs{F(j)}_{j\in \N}$ forms a decreasing, log convex  sequence, as easily follows by taking the logarithmic derivative of $F$. The logarithmic derivative, $F'/F$, is analytic except at the zeros of $F$, and is given by 
$$
\sum_p \frac{\ln (1- 1/p)}{\(\frac{p}{p-1}\)^{s-1} (p-1) + 1}.
$$
This converges uniformly on compact subsets of $|\Im s | < \pi/\ln 2$. Viewed as a real function (that is, restricting $s$ to be real), each summand is the negative of a completely monotone function and $F$ is nonnegative on $\R$, so that $F$ is {\it logarithmically completely monotone\/} (meaning that $F > 0$ and $-F'/F$ is completely monotone) which implies $F$ is completely monotone.

With single-digit accuracy, I managed to approximate (with pencil and paper) the values of the expression in (**) for $n = 3,4,5,6$; they are respectively, $.55, .6, .7, .8$. The last is surprisingly close to the upper bound computed from (2) above, which is  $(\zeta(2)\zeta(3)/\zeta (6) )\cdot 1/\zeta(2)\zeta(3)\dots \sim.845 $. This suggests that the numerator of (**) tends to $\zeta(2)\zeta(3)/\zeta (6) $; in other words, that  the upper bound be approximately achieved. We will prove  this  after putting it in a more recognizable form.

Let us rewrite the numerator, substituting innocuously (when $n$ is large) $F(1) = 1$ for $1/\zeta(n)$ and  $F(2) = 1/\zeta(2)$; then, subtracting the expression from $F(0) = \zeta(2)\zeta(3)/\zeta(6)$, we obtain
$$
D(n):= F(0) - nF(1) + {n \choose 2}F(2) -  \dots + (-1)^n F(n) = \sum_{i=0}^n (-1)^n {n \choose i} F(i).
$$

We will show
$$
\lim_{n\to \infty} D(n) = 0.
$$
This is equivalent to the numerator in (**) converging (in $n$) to $F(0) = \zeta(2)\zeta(3)/\zeta(6)$.

A function $\Arrow f; \R.\R$ is {\it completely monotone\/} if $(-1)^n f^{(n)}(r) \geq 0$ for all $n \in \Z^+$ and $r \in \R$ (here $f^{(n)}$ is the $n$th derivative); it is {\it logarithmically completely monotone\/} if $f(r) > 0$ for all $r$ and $\ln f$ is completely monotone. It is known that logarithmically completely monotone functions are completely monotone.

Let $\Delta$ denote the usual difference operator, acting on functions on $\Z$ or $\R$, that is, $\Delta f(k) = f(k+1) - f(k)$. If $\Arrow f; \Z.\R$ satisfies $(-1)^n\Delta^n f(k) \geq 0$ for all $n \in \Z^+$ and $k \in \Z$, then we say that $f$ is {\it completely monotone.}

It is routine that $D(n) = (-1)^n \Delta^n F(0)$; so it is enough to show that $(-1)^n \Delta^n F (0) \to 0$, which turns out to be completely elementary. Consider $d_{n} (k-1) = d_{n}(k) + d_{n+1}(k-1)$; iterating this, we quickly see that since all $d_n(m) \geq 0$, we have $d_{n} (k-1)\geq j d_{n+j}(k)$. As $d_{n+j}(k) \geq 0$, this forces $d_{n+j}(k) = \Oh{1/j}$; in particular,  $d_n (k) \to 0$ as $n \to \infty$.

Now suppose that $\Arrow f;\R.\R$ is completely monotone; then it is routine to see that $f|\Z$ (or any other discrete subgroup) is completely monotone (in the sense of functions on $\Z$). By the higher order mean value theorem, given $r \in \R$, and $ n \in \Z^+$, there exists $\xi \in [r, r+n]$ \st $\Delta^{n} f(r) = f^{(n)}(\xi)$; setting $r = k \in \Z$, the sign of $\Delta^{n} f(k)$  is the same as the sign of $f^{(n)}(\xi)$ at some real number, and we are done.

The following is elementary, and presumably standard.

\Lem Proposition \fouone. Suppose that $\Arrow f; \Z.\R$ is completely monotone. Then for all $k \in \Z$
$$
\lim_{N\to \infty} \sum_{j=0}^N (-1)^j \Delta^j f (k) \quad \text{exists and equals } f(k-1).
$$

\Rmk Formally, this means that $\I + \sum_{j=1}^{\infty} (-1)^n \Delta^n = (\I + \Delta)^{-1}$ (as would be expected from the power series expansion) when applied to completely monotone functions (and therefore to the vector space they span).

\Pf Apply $ \I + \Delta$ to the expression on the left of the display; this yields $(\I + (-1)^{N+1}\Delta^{N+1}) f(k) = f(k) + d_{N+1} (k) \to f(k)$. On the other hand, $(\I + \Delta)f(k-1)   = f(k)$.

Set $g_N(l):= \sum_{i}^N d_i(l)$. Then $(\I + \Delta )g_N (l) = f(l) + d_{N+1}(l)$, but also $(\I + \Delta )g_N(l) = g_N(l+1)$. Setting $l = k-1$, we have $g(k) = f(k-1) + d_{N+1}(k-1)$; this says $| g_N(k) - f(k-1)| \leq d_{N+1}(k-1)$, which goes to zero as $ N\to \infty$.
\qed

\Lem Proposition \foutwo.  The restriction of $F$ to $\R$ is logarithmically completely monotone.

\Pf With $\ln F$ given above, we note that $F|\R$ is strictly positive, and the logarithmic derivative $F'/F = (\ln F)'$ is a locally convergent (on compact subsets of the strip $|\Im z | < \pi/\ln 2$) sum of terms each of which is the negative of a completely monotone function.
\qed

\Lem Corollary \fouthr. The natural density of matrices in $\TF_n \cap \NS_n$ increases upwards (as $n \to \infty$)
to
$$
\frac{\zeta(2) \zeta(3)}{\zeta(6) } \cdot \frac 1{\zeta(2)\cdot \zeta (3) \cdot \zeta (4) \cdots} \sim .845.
$$

\Rmk In fact, it also follows from the last two propositions that if $T(n)$ is the (strong) natural density of $\TF_n \cap \NS_n$, then $\brcs{T(n)} $ is increasing, and if $\epsilon (n)$ is the difference between the limit and $T(n)$, then $\sum \epsilon(n) < \infty$. So convergence is somewhat faster than expected.

\noindent {\it Motivation.} Why the emphasis on $1$-block size $n-1$ (for PH-equivalence classes of matrices in $\NS_n$)? For one thing, if $B$ and $B'$ are in terminal form with $1$-block size $n$, we can easily decide (from Proposition \twoone) whether they are PH-equivalent (and the procedure can be made very fast).

For another, the condition that $B \in \NS_n$ have a terminal form with $1$-block size $n$, for $n \geq 6$, has density at least $.8$, tending in $n$ to $.845\dots$---meaning five out of six random matrices should have such a terminal form.

If we consider $1$-block size at least $n-2$ instead, the upper bound is then in excess of $.99$; so if the upper bound is achieved (as $n \to \infty$), then for sufficiently large $n$, over $99\%$ of random integer matrices will have a terminal form with $1$-block size at least $n-2$. This suggests that it might be worthwhile obtaining the analogue of Proposition \twoone\ for $n-2$, describing the equivalence classes containing terminal form of this type).

For the classification, it would be reasonable to determine the likelihood that at least {\it one\/} of $B$ and $B\op$  be PH-equivalent to a terminal form with $1$-block size $n-1$. The simplest possible form of inclusion-exclusion would yield a likelihood of $2a - b$ where $a$ is the likelihood that $B$ have a terminal form with $1$-block size $n-1$ (about $.845$ as just calculated above), and $b$ is the likelihood that both $B$ and $B\op$ have such a terminal form. Computing $b$ appears to be difficult ($b \neq a^2$; the properties are not independent). Towards this, the characterizations for  $ J(B)$ and $ J(B\op)$ to both belong to $\NS_{n,n-1}$ (Corollary \nfoufiv) might be useful.

\comment 
If $p$ is a prime and $m \in \N$, then $v_p (m) = r$ means that $p^r$ maximally divides $m$; this is the usual valuation.

\Lem Lemma \foufou. Suppose that $B \in \NS_n$ is PH-equivalent to the terminal form
$$
\(\matrix  \I_{n-1} & a \\ 0 & d \\ \endmatrix\)
$$
where $a = (a_i)^T \in \Z^{(n-1) \times 1}$ and $\cont \brcs{a,d} = 1$. Then $B\op$ is PH-equivalent to a terminal form with $1$-block size $n-1$ iff for every  prime $p$ dividing $d$, there exists at most one (hence exactly one) nonzero $a_i$ \st $v_p(a_i) < v_p (d)$, and for this $i$, $v_p(a_i) = 0$.

\Pf If $B$ is in the displayed terminal form, we may reorder the entries of the column $a$ (by conjugating with a permutation matrix) so that if there are any $i$ \st $a_i = 0$, they form an initial segment; let $l$ (possibly zero) denote the number.
It is a routine computation that
$$\eqalign{
B\op &= \(\matrix \I_l & \\ &\frac{d}{(d,a_{l+1})} &&&& \\
&& \frac{d}{(d,a_{l+2})} &&& \\
&&& \ddots && \\
&&& &\frac{d}{(d,a_{n-1})}& \\
0^l&-\frac{a_{l+1}}{(d,a_{l+1})}& -\frac{a_{l+2}}{(d,a_{l+2})}& \dots& -\frac{a_{n-1}}{(d,a_{n-11})}& 1
\endmatrix \) \cr
&\sim \(\matrix \I_{l} & \\
 0^l&1&-\frac{a_{l+1}}{(d,a_{l+1})}& -\frac{a_{l+2}}{(d,a_{l+2})}& \dots& -\frac{a_{n-1}}{(d,a_{n-11})}
\\
&&\frac{d}{(d,a_{l+1})} &&& \\
&&& \frac{d}{(d,a_{l+2})} && \\
&&&& \ddots & \\
&&&& &\frac{d}{(d,a_{n-1})} \\
\endmatrix \): = B'.
}$$
The second matrix is obtained by cyclically permuting the rows below the identity block, then correspondingly permuting the columns; if we order the $a$s so that $\gcd\brcs{d,a_i}$ are decreasing for $l+1 \leq i \leq n-1$, then the second matrix is in terminal form, although not a particularly nice one. If $l = n-2$, its $1$-block size is already $n-1$, and there is nothing to do. Otherwise, $l < n-2$.

Now $B'$ has a terminal form with a $1$-block size $n-1$ iff there exists $j$ \st $I(B_{\Omega(j)}) = \brcs{0}$; in other words, iff there exists a column, which after its deletion, the resulting matrix has $1$ as the $\gcd$ of the  determinants of the square size $n-1$  matrices. Let $s_i = d/\gcd\brcs{d,a_i}$ and $u_i = -a_i/\gcd\brcs{d,a_i}$. For all $i$, $\gcd\brcs{s_i,u_i} = 1$

Deleting any of the first $l+1$ columns results in the $\gcd$ being the product of the diagonal entries, $\prod s_i$. Removing the $i$th column with $l+2 \leq i \leq n-1 $ results in a $\gcd$ of $t_i:=s_i^{-1}\prod_j s_j $. However, deleting the $l+1$st column results in the determinants, $\Lt{\prod s_i, s_1^{-1}u_1\prod s_i, s_2^{-1}u_2 \prod s_i, \dots, s_n^{-1}u_n \prod s_i}$. It is easy to check (since $(s_i,u_i) = 1$ for all relevant $i$) that this list has the same greatest common divisor as $\Lt{{s_j^{-1}\prod s_i}}_{j=l+1}^n = \Lt {t_j}$.

Hence $B'$ is PH-equivalent to a terminal form with $1$-block size $n-1$ iff $\gcd \brcs{t_j} = 1$. Let $p$ be a prime; if $v_p (a_i) < v_p (d)$ for two values of $i$ (with $a_i \neq 0$), then $p$ divides all the $t_j$; the converse is clear.
\qed

The condition in the lemma seems rather restrictive, but it is not clear how to convert this into a likelihood estimate for both $B$ and $ B\op$ having $1$-block size $n-1$.
\endcomment

\SecT 8 Topological isomorphism for topologically critical groups

In this section, we state some well-known and not-so-well known results about topologically critical groups; see also [H]. Suppose $G \to V$ and $H \to W$ are group homomorphisms from abelian groups to ordered real Banach spaces. We say $\Arrow f; G.H$ is {\it continuous\/} if there exists continuous and linear $\Arrow F; V.W$ whose restriction to $G$ is $f$ (typically, the images of $G$ and $H$ will be dense in their respective Banach spaces; in this case, continuity is equivalent to the usual notion \wrt the relative topologies on $G$ and $H$).

A subgroup $G$ of $\R^n$ is {\it topologically critical of rank $n+1$} if it is free of rank $n+1$ and dense. Any subgroup of lesser rank of a topologically critical group is discrete. In this section (only), when  we regard $g \in G$ as an element of $\R^n$, we denote it $\hat g$.  Associated to a topologically critical group is an isomorphism class of  rank $n+1$ subgroups of $\R$, $\TO(G)$, defined as follows. Select any ordered $\Z$-basis for $G$, $(g_i)_{i=1}^{n+1}$. Since $\brcs{g_i}_{i=1}^n$ generates a discrete subgroup, it is a real basis for $\R^n$; hence we can write $\hat g_{n+1} = \sum \alpha_i \hat g_i$. It is easy to check that $\brcs{1, \alpha_1,\dots, \alpha_n}$ is rationally linearly independent, and so we may form the subgroup of $\R$, $\Z + \sum \alpha_i\Z$, of rank $n+1$. Every topologically critical subgroup of $\R^n$ is topologically isomorphic to the group generated by $\brcs{e_i; \sum e_j \alpha_j}$ (where $e_i$ are the standard basis elements of $\R^n$) by this construction (for example, see [H]).

Topologically critical groups have an interesting property: every subgroup is either dense (those of full rank) or discrete (those of lesser rank).
 
Let $\TO(G)$ denote the isomorphism class of the inclusion $\Z + \sum \alpha_i \Z \subset \R$, that is, \wrt continuous maps. Alternatively, we may view the group as a totally ordered group (the ordering inherited from $\R$), and use order-preserving group isomorphisms between $
G = \Z + \sum \alpha_i \Z $; the resulting equivalence classes are the same, since in this case, any continuous map is either order-preserving or its negative is.
 
 \Lem Lemma \fivone. Suppose $G$ and $H$ are topologically critical groups \st $\TO (G) \iso \TO (H)$. Then $H$ and $G$ are continuously isomorphic.

\Pf
Suppose $\brcs{\alpha_i}_{i=0}^{n} $ and $\brcs{\beta_i}_{i=1}^{n+1}$ are
subsets of $\R$ that are linearly independent over the rationals, and
$\alpha_{n+1} = 1 = -\beta_{n+1}$, and moreover, $\sum \alpha_i \Z = \sum
\beta_i \Z$ (as subgroups of $\R$).
Let $G = \langle e_i;e_{n+1}:=  \sum_{i=1}^n {\alpha_i e_i} \rangle$ be
the (dense) subgroup of $\R^n$, where $\brcs{e_i}_{i=1}^n$ is the standard
basis for $\R^n$.
Then there exist $\brcs{h_i}_{i=1}^n$ \st $G = \sum h_i \Z$ and $h_{n+1} =
\sum_{i=1}^n \beta_i h_i$.
 
For each $i=1, 2, \dots, n$, there exist integers $a_{i,t}$ ($t = 0, 1,
\dots, n$) \st $\alpha_i = \sum_{t = 1}^{n+1} \beta_i a_{it}$. Complete
$(a_{it})$ to an $(n+1)\times (n+1)$ matrix $
A$ by defining $a_{n+1,t} = \delta_{n+1,t}$ (so the bottom row is $(0, 0,
\dots, 0,1)$.
 
Set $g_i = e_i$ (to avoid confusion between the standard bases) for $i =
1, 2, \dots, n+1$. Define for each $j = 1,2, \dots, n+1$,
$$
h_j = \sum_{i=1} a_{i,j} g_i
$$
(so here we are using $A^T$). Obviously, $h_j \in G$. We first show that
$h_{n+1} = \sum_{t=1}^{n+1} \beta_t h_t$.
On one hand,
$$\eqalign{
h_{n+1} & =  \sum_{i=1}^{n+1} a_{i,n+1}g_i \cr
& = \sum_{i=1}^{n} (a_{i,n+1}+ a_{n+1,n+1}\alpha_i)g_i;\qquad \text{ on the
other hand, }\cr
\sum_{t=1}^n \beta_t h_t & = \sum_{t=1}^n \beta_t \sum_{i=1}^{n+1}
a_{i,t}g_i \cr
& = \sum_{t=1}^n \beta_t  \( \sum_{i=1}^n a_{i,t} g_i + \alpha_i \beta_t
a_{n+1,t}\)\cr
& = \sum_{i=1}^n g_i \cdot \(\sum_{t=1}^n a_{i,t} \beta_t + 0\)\cr
& = \sum_{i=1}^n g_i\cdot  ( \alpha_i  - \beta_{n+1}a_{i,n+1}).\cr
}$$
Since $\beta_{n+1} = -1$, we are done.
 
As $\sum \alpha_i\Z = \sum \beta_i \Z$, we can find the inverse map
(both are free abelian groups of rank $n+1$ to $A$; this takes the $h_j$
to $g_j$, and it follows immediately that $\sum h_j \Z = \sum g_j \Z$, and
the rank condition guarantees that the sums are direct.
\qed

\SecT 9 Basic  critical dimension groups

A {\it dimension group\/} is a direct limit of simplicial (partially ordered abelian) groups; see [G], the standard reference for partially ordered abelian groups, for far more information than can be given here. By [Gr], [EHS], a partially ordered abelian group $G$ is a dimension group iff it is {\it unperforated\/} (for $n \in \N$ and $ g \in G$, $ng \geq 0$ entails $g \geq 0$) and satisfies {\it Riesz interpolation\/} (for $a_i, b_j \in G$ with $i,j \in \brcs{1,2}$ with $a_i \leq b_j$ for all $i,j$, there exists $c \in G$ \st $a_i \leq c \leq b_j$ for all $i,j$). All partially ordered groups will be abelian.

An {\it order unit\/} of a partially ordered group $G$ is an element  $u \in G^+$ \st for all $g \in G$, there exists $n\in \N$ \st $-nu \leq g \leq nu$. A partially ordered abelian group is {\it simple\/} if every nonzero element of $G^+$ is an order unit. A {\it trace\/} (or {\it state\/}) of $G$ is a nonzero positive real-valued group homomorphism; it is {\it normalized\/} at the order unit $u$ if its value thereat is $1$. The collection of normalized traces, denoted $S(G,u)$ and equipped with the point-open (weak) topology, is a compact convex subset of a Banach space. The {\it value group\/} of a trace $\tau$ is simply $\tau(G)$, its set of values.

The real vector space consisting of convex-linear continuous ({\it affine}) real-valued functions $\Arrow f; S(G,u). \R$ is denoted $\Aff S(G,u)$. It is a Banach space \wrt the supremum norm.  There is a natural order preserving group homomorphism, the {\it affine representation\/} (\wrt $u$), $\Arrow \hat\ ; (G,u). \Aff S(G,u) $ given by $g\mapsto \hat g$, where $\hat g(\tau) = \tau(g)$ for $\tau \in S(G,u)$. This imposes a pseudo-norm topology on $G$, which is a norm if the affine representation is one to one.

When $G$ is a dimension group, $S(G,u)$ is a Choquet simplex. When $G$ is also simple, there is a complete characterization available, the affine representation $G \to \Aff S(G,u)$ (\wrt any, or equivalently all, choices of order unit $u$) has dense range, and $G^+\setminus\brcs{0}$ consists of $\Set{g \in G}{\hat g \text{ is strictly positive}}$. The converse is also true.

A trace is {\it pure\/} (or {\it extremal\/}) if it is not a proper convex-linear combination of other traces. The {\it extremal boundary\/} (of $S(G,u)$), denoted $\partial_e S(G,u)$, consists of the pure normalized traces. When $S(G,u)$ is finite-dimensional, it is a simplex in the usual sense (as a compact convex subset of Euclidean space), and in that case, $\Aff S(G,u)$ can be identified with $\R^n$ for some integer $n$,  the standard basis elements identified with the pure traces  (possibly with normalization). The {\it strict ordering\/} on $\R^n$ or $\Aff S(G,u)$
is the partial ordering whose positive cone consists of the strictly positive functions.

A consequence is that if $G$ is a simple dimension group with finitely many, say $n$, pure traces and the kernel of the affine representation is zero, then $G$ is order isomorphic to a dense subgroup of $\R^n$ equipped with the strict ordering. The pure traces are just (up to renormalization) the coordinate maps.

We say a simple dimension group $G$ is {\it critical\/} if it is free of rank $n+1$ and has $n$ pure traces. By the preceding, this means it can be identified with a dense subgroup of $\R^n$, and since the partial ordering determines the topology (here the affine representation is automatically one to one), it is also topologically critical.

We are interested in classification of critical groups. It turns out that there is a class of them whose classification incorporates PH-equivalence.

A critical  group is called {\it basic\/} if it is order isomorphic to a dense subgroup of $\R^n$ (equipped with the strict ordering) with generators $\brcs{e_1, \dots, e_n; \sum \alpha_i e_i}$, where $e_i$ are the standard basis elements, and $\alpha_i$ are real numbers. For  a subgroup so generated, density is equivalent to the set $\brcs{1, \alpha_1, \dots, \alpha_n}$ being rationally linearly independent. We will give a characterization that avoids such a specific realization, referring only to internal properties.

Critical, and especially basic critical groups, are a useful source of examples. For example, in [BeH], we translated Akin's notion of {\it good measure\/} on a Cantor set to dimension groups, and we were able use these to illustrate various properties of good and non-good traces. Following [BeH], we say that a trace $\tau$ on a  dimension group $G$ is
 {\it good\/} if for all $b \in G^+$ and $a \in G$ \st $0 < \tau(a) < \tau(b)$, there exists $a' \in G^+$ \st $a' \leq b$ and $\tau(a') = \tau(a)$. For simple dimension groups, this is equivalent to a much simpler criterion (in  context),
that the image of $\ker \tau$ in the affine representation of $G$ be norm-dense in $\tau^{\vdash}: = \Set{h \in \Aff S(G,u)}{h(\tau) = 0}$.

This lead to the definition of ugly for a trace on a dimension group; $\tau$ is {\it ugly\/} if $\ker \tau $ has discrete image in $\Aff S(G,u)$ and the trace $\tau \otimes 1_{\Q}$ on $G \otimes \Q$ is good.

For sets of traces, there are corresponding definitions, which become rather complicated---but if $S(G,u)$ is finite-dimensional, and $\Omega \subset \partial_e S(G,u)$, the relevant ones for this article reduce to the following:

\item{(i)} $\Omega$ is {\it good\/} if whenever $b \in G^+$ and $a\in G$ satisfy $0< \tau(a) < \tau (b)$ for all $\tau \in \Omega$, then there exists $a' \in G^+$ \st $a - a' \in \ker \Omega:= \cap_{\tau \in \Omega} \ker \tau$ and $a' \leq b$
\item{(ii)} $\Omega$ is {\it ugly\/} if the image of $\ker \Omega$ is discrete in $\Aff S(G,u)$ and the extension of $\Omega$ to a set of traces on $G \otimes \Q$ is good.

\noindent These are not equivalent to the definitions in general; the restriction to $\Omega \subset \partial_e S(G,u)$ allowed considerable  simplification. Among other things, these correspond to faces in $S(G,u)$.
For critical groups in general and any nonempty family of traces, $\ker \Omega$, being a subgroup of rank at most $n-1$, is automatically discrete. So the definition of ugly simplifies further.

 Necessarily, when $G$ is a basic critical group, for all pure traces $\tau$, $\rk \tau (G) = 2$, and this forces all the pure traces to be ugly. Conversely, the pure trace $\tau$ is ugly if $\rk \tau(G) = 2$. There are examples (for every $n \geq 2$, that is, rank at least $3$) of  critical groups  all of whose pure traces are ugly, and even with the additional property that $\brcs{\tau_i(G)}$ are mutually order isomorphic as real subgroups, that are not basic (or even a modest extension, to be defined later, almost basic).

Let $r$ be a real number that is neither rational, quadratic, nor cubic over the rationals; that is, the set $\brcs{1,r,r^2,r^3}$ is linearly independent over the rationals.
Let $G$ be the subgroup of $\R^3$ spanned by $\brcs{E_1:=(1,1,1),E_2:= (1,1,r),E_3:= (1,r,0), E_4 = (r,0,0)}$. The set of four $3\times 3$ determinants of the spanning set is rationally linearly independent. Hence $G$ is dense in $\R^3$, and thus with the strict ordering, is a critical group (of rank three).

The  pure traces on $G$ are the three coordinate maps, denoted $\tau_i$. Then we see that $\tau_1 (G) = \Z + r\Z = \tau_2(G) = \tau_3(G)$, free of rank two. In all three  cases, the kernel is free of rank two, and since the affine representation is one to one, and since the kernels are discrete subgroups, the corresponding pure traces are ugly.
However, as we will see later, $G$ is not basic.

This leads to a class of non-basic critical groups free of rank $n+1$ \st all $\tau_i (G)$ are equal and rank two (hence all the pure traces are ugly). Pick $r$ \st $\brcs{1,r, \dots, r^{n}}$ is rationally linearly independent (that is, either $r$ is transcendental or its algebraic degree is at least $n+1$). Define elements of $\R^n$
$$\eqalign{
F_n & = \(\matrix r & 0 & 0 &\dots & 0  & 0 & 0 \endmatrix\)\cr
F_{n-1} & = \(\matrix 1 & r & 0 &\dots & 0  & 0 & 0 \endmatrix\)\cr
F_{n-2} & = \(\matrix 1 & 1 & r &\dots & 0  & 0 & 0 \endmatrix\)\cr
\vdots &\ \ \ \  \matrix & &  & & &\ddots &  &  &  &\endmatrix \cr
F_{1} & = \(\matrix 1 & 1& 1 &\dots & 1  & 1 & r \endmatrix\) \cr
F_0 & = \(\matrix 1 & 1& 1 &\dots & 1  & 1 & 1 \endmatrix\) \cr
}$$
That is, $F_i$ has $i-1$ zeros (for $i \geq 1$), immediately preceded  by $r$, which in turn is immediately preceded by enough ones to fill up the row. Let $M_i$ be the $n \times n$ matrix obtained by deleting $F_i$, and throwing together the rest of the $F_j$s.
Then $ \det M_0 = r^{n}$ and $|\det M_{1}| = r^{n-1}$ as is easily seen from the lower triangular forms. For $i >1$, $M_i$ is a block lower triangular matrix, and it is straightforward to check that $\det M_i = r^{n-i}(1-r)^{i - 1}$. (At one point, multiply the matrix $ rN^T + \I + N  + N^2  + \dots  $ by $\I - N$, creating an upper triangular matrix. See the lemma below.)
Next we claim that the set $\brcs{r^n, r^{n-1}, r^{n-2}(1-r), \dots , r(1-r)^{n-2}, (1-r)^{n-1}}$ spans $\sum_{i=0}^n r^i \Q$, which is easily checked by induction. Hence the set is rationally linearly independent.

Thus $G \equiv G(n,r)$ is a critical group of rank $n+1$, so with the strict ordering inherited from $\R^n$ is a simple dimension group with $n$ pure traces, the latter arising as the coordinate functions. Their value groups, that is the ranges of the pure traces, are all equal to the rank two group, $\Z + r\Z$. In particular, their kernels are necessarily of rank $n-1$ and discrete (the latter from being a critical group), and it easily follows that they are all ugly. We will soon show that if $n>2$, then $G(n,r)$ is not basic (or even satisfy a more general property, almost basic).

We have $G (n,r)\subset (\Q + r\Q)^n$ of rank $n+1$ and $G$ is dense in $\R^n$; we have assumed $r$ does not satisfy a rational equation of degree $n$ or less.

\Lem Lemma \sixone. Let $N$ be the lower triangular $k \times k$ matrix with $1$s in the $(j+1,j)$ entries and zeros every where else. Let $r$ be any number, and set $Q =r N^T + \I + N + N^2 + \dots $. Then $\det Q = (1-r)^{k-1}$.

\Pf Multiply $Q$ from the left by $\I - N$ (which has determinant $1$); the outcome is $\I -rNN^T + rN^T$. Now $NN^T$ is just the identity matrix less the first $1$, so that
$(\I - N)Q$ is upper triangular, with diagonal entries $(1, 1-r, 1-r, \dots, 1-r$. Hence $\det Q = (1-r)^{k-1}$.
\qed

Basic critical groups admit rather strong properties. The first is that every proper subset of the pure trace space is ugly. For a simple dimension group $(G,u)$ with one to one affine representation and finite-dimensional $S(G,u)$, and  $\Omega \subset \partial_e S(G,u)$, the definition of ugliness of $\Omega$ simplifies to  (i) $\ker \Omega:= \cap_{\tau \in \Omega} \ker \tau$ is discrete, and (b) $\ker \Omega \otimes \Q$ is dense in $\Omega^{\perp} = \Set{h \in \Aff S(G,u)}{h|\Omega \equiv 0}$ ($\Omega$ can be replaced by the face it spans).

When $(G,u)$ is critical of rank $n+1$, and $\Omega  \subset \partial_e S(G,u)$, then it is fairly easy to decide whether $\Omega $ is ugly. First, every subgroup of rank $n$ or less is automatically discrete, hence any $\Z$-linearly independent subset is real linearly independent. Second, if $\Omega  \subseteq \partial_e S(G,u)$, then $\Omega ^{\perp}$ has (real) dimension exactly $n - |\Omega |$ (the  set of  pure traces is a dual basis for $\Aff S(G,u)$). The following is then immediate. Note that although the definition involves a choice of order unit, the criterion does not. In other words, it does not matter at which order unit $u$ we choose to normalize the traces.

\Lem Lemma \sixtwo. Let $(G,u)$ be a critical group of rank $n$, and let $\Omega $ be a proper set of pure traces. Then $\Omega $ is ugly iff $\rk \ker \Omega  = n-|\Omega |$.

It is trivial that if $G$ is basic, then the criterion is satisfied for every proper subset $\Omega $ of $\partial_e S(G,u)$. However, there exist non-basic but critical groups which also have the property that for every proper $\Omega  \subset \partial_e S(G,u)$, $\Omega $ is ugly. In this case, there is a finite obstruction to being basic.

In the examples above, $r$ is a real number that satisfies no nonconstant rational polynomial of degree $n$ or less, and we formed the group $G(n,r) \subset \R^n$.
These are critical dimension groups with the interesting property that for all pure traces $\tau$, $\tau(G) $ are equal to each other. Equality of the value groups is not an invariant (since by changing the order unit, we change the value groups), except in the case that we are looking at invariants for $(G,u)$, that is, where $u$ is specified. However, what is an invariant is that all $\tau (G) $ be order-isomorphic as subgroups of the reals as $\tau$ varies over the pure traces.

Moreover, in these examples, we have that $\rk \tau(G) = 2$, so that $\rk \ker \tau = n-1$; thus all pure traces are ugly, just as in the case of basic critical groups. However, if $n \geq 3$, $\rk (\ker \tau_1 \cap \ker \tau_n) = n-3 \neq n-2$; specifically, a $\Z$-basis for the intersection is $\brcs{F_n-F_2, F_{n-1}-F_2, \dots, F_3-F_2}$). Hence there exists a two-element subset of the pure trace space that is not ugly, so that if $n\geq 3$, these critical groups are not basic.

We analyze potential isomorphisms of critical groups of rank $n+1$ as follows. Begin with any ordered $\Z$-basis, $\brcs{v_1, v_2, \dots, v_n, v_{n+1}}$, which we regard as elements of  $\R^{1\times n}$, that is, rows of real numbers. We construct an $(n+1) \times n$ real matrix $A$ by letting its $i$th row be $v_i$.

Applying any element of GL$(n+1, \Z)$ to $A$ (from the left) just changes the $\Z$-basis, hence leaves the group they generate the same.

As in the earlier sections, let $P(n,\R)^+$ denote the group weighted permutation matrices of size $n$ with only positive weights---that is, the set of products $P\Delta$ where $P$ is a permutation matrix, and $\Delta$ is a diagonal matrix with only strictly positive real entries along the diagonal. The group of order-automorphisms of $\R^{1\times n}$ \wrt either the strict or the usual ordering is just $P(n,\R)^+$, and since any order isomorphism between critical groups (necessarily of the same rank) extends uniquely to an order automorphism of $\R^{1 \times n}$ (after identifying the two sets of pure traces), we have that the order isomorphisms between critical groups are determined by right actions of $P(n,\R)^+$.

So we can act on $A$ from the left by GL$(n+1,\Z)$ and from the right by $P(n,\R)^+$. In particular, we can permute rows, we can permute columns, perform elementary row operations (over the integers), and multiply columns by  positive real scalars. If after a sequence of such actions, we arrive at a matrix $A'$ where the the top $n\times n$ part is just the identity, then the critical dimension group is basic.

We illustrate this with a simple example, the case $n = 2$ of $G(n,r)$. Here $r$ is a real number that is not quadratic or rational. Let $G = \langle (r,0), (1,r), (1,1)\rangle \subset \R^2$. We have the following series of transformations,
$$
\( \matrix 1 &1 \\ 1 & r \\  r & 0 \\
\endmatrix\) \mapsto \( \matrix 1 &1 \\ 0 & r-1 \\  r & 0 \\
\endmatrix\) \mapsto \( \matrix r & 0 \\ 0 & r-1 \\  1 & 1 \\
\endmatrix\)  \mapsto \( \matrix 1 & 0 \\ 0 &1 \\  \frac 1r & \frac 1{r-1} \\
\endmatrix\).
$$
Thus $G$ is basic (since $\brcs{1,1/r, 1/(r-1)}$ is linearly independent over $\Q$ iff $\brcs{1,r,r^2}$ is).
 It also satisfies  the property that $\tau(G)$ are mutually isomorphic as $\tau$ varies over the pure trace space.

Suppose $A$ is partitioned as $\(\smallmatrix B \\ \alpha \\ \endsmallmatrix\)$, where $B$ is $n \times n$ (so $\alpha = (\alpha_1, \dots, \alpha_n)$ is just the bottom row), and now assume that $B$ is a rank $n$ matrix (necessary for it to yield a critical group anyway) with only integer entries. Some of the time (but not always), we restrict the actions of GL$(n+1,\Z)$ to be those of GL$(n,\Z) \times \brcs{1}$, that is, performing only elementary row operations not affecting the bottom row, $B$. Necessary and sufficient for the row space of $A$ to be a critical dimension group is that the set $\brcs{1, \alpha_1, \dots, \alpha_n}$ be rationally linearly independent.

Since multiplying on the right by weighted positive diagonal matrices preserves order isomorphism, we may assume that each column of $B$ is unimodular (of course, the corresponding entry of $\alpha$ is multiplied by a rational at the same time). Hence we may assume that $B \in \NS_n$.

Every $U\in \gl(n,\Z)$ and permutation matrix $P$ yields an order isomorphism of the dimension group (by extending $U$ to $C = U \oplus (1)$), so we may assume that $B$ is in terminal form.

In particular, if the terminal form is simply the identity (of size $n$), then $G$ is basic. More generally, let $G'$
be the subgroup of $\R^n$ generated by the  rows of the current matrix, renamed $A = \(\smallmatrix B \\ \alpha \\ \endsmallmatrix\)$; as we have observed, this is order isomorphic to $G$. The pure traces are still the coordinate functions, $\tau_i$. It is easy to check that $\tau_i (G') = \Z + \alpha_i \Z$, and the latter being of rank two implies that all pure traces are ugly. But more is true. If we manipulate further using GL$(n,\Q)$ (that is rational elementary row operations), we can reduce $B$ to the identity matrix. This means that $G' \otimes \Q$ is order isomorphic to $G_0 \otimes \Q$ for some basic critical group $G_0$. It follows immediately that every proper subset of the pure trace space of $G'$ is ugly.

 We investigate the converse. For any critical dimension group with pure trace space  $\partial_e S(G,u) = \brcs{\tau_i}$, set $J_i = \ker \Omega(i)$. It is easy to see that either $J_i = \brcs{0}$ or $\rk J_i = 1$. In the latter case, pick a generator $x_i$ for $J_i$ (we only have two choices, $\pm x_i$). Now form $E\equiv E(G):= \sum x_i \Z$ where the $i$ varies over those \st $J_i$ is not zero. The $x_i$ are the same as those in the original construction of the invariant for the integer part of $G$.

\def\tor{\text{Tor}}

The latter ensures that the isomorphism $G \to G'$ induces a group isomorphism $E(G) \to E(G')$, and thus yields an isomorphism $G/E(G) \to G'/E(G)$. In particular, the torsion parts are respectively isomorphic. We claim that this induces an isomorphism $\tor(G/E(G)) \to J(B\op)$. We are not done yet, since $G = (\sum f_j \Z ) \oplus \alpha \Z$ as abelian groups.

It suffices to show that $  r(B)/X(B)$ (a subgroup of $G/E(G)$) is exactly the torsion part of $G/E(G)$ (and similarly with $C$ replacing $B$). Since the former is torsion, we have inclusion. Now suppose that $g + E(G)$ is a torsion element in $G/E(G)$. There thus exists $n > 0$ \st $ng \in E(G)$, in particular, we can write $ng$ as an integer combination of elements of $x_i$, so that $ng \in \sum f_j \Z$ (as the $x_i \in \sum f_j \Z$). On the other hand, since $\brcs{f_j} \cup \brcs{\alpha}$ is a $\Z$-basis for $G$, we may write $g $ uniquely as $\sum t_j f_j + m\alpha$, so that $ng = \sum nt_j f_j + nm \alpha$; since $ng \in \sum f_j \Z$, we deduce $nm\alpha $ is in the span of $f_j$, which of course is impossible unless $nm =0$, that is, $m = 0$. So $g \in \sum f_j\Z$, and thus $g + E(G) \in r(B)/X(B)$. Of course, the same works with $C$ replacing $B$.

First, $\sum x_i \Z = \oplus x_i \Z$ (routine). Next, $E$ and $G/E$ are invariants for order isomorphism; that is, any order isomorphism  between critical dimension groups $G_1 \to G_2$ maps $E(G_1)$ isomorphically (as abelian groups, of course) onto $E(G_2)$, so that the induced map on their cokernels $G_1/E(G_1) \to G_2/E(G_2)$ is also an isomorphism.

When $G$ is basic, $G/E \iso \Z$, as is obvious from its matrix $A$ representing it. When every proper subset of the pure trace space is ugly, then the torsion-free rank of $G/E$ is one, but it may have torsion elements. If not every proper subset is ugly, then the torsion-free rank of $G/E$ must exceed one, and there can also be torsion. 
The following is practically tautological.

\Lem Lemma \sixthr. Let $G$ be a critical dimension group. Then $G$ is basic iff $G/E(G) \iso \Z$.

\Pf One way is trivial. Suppose $G/E \iso \Z$. Then $G \to G/E$ splits, and thus we may find $y \in G$ \st $E \oplus y\Z = G$. We can write $E = \oplus x_i \Z$, and since the rank of $E$ is $n$, there are $n$ of the $x_i$. Now each $x_i$ vanishes at all the traces except $\tau_i$; by replacing $x_i $ by $-x_i$ if necessary, we can also assume that $\tau_i (x_i) > 0$. Set $u = \sum x_i$, so that $\tau_i (u) = \tau_i(x_i) > 0$ for all $i$. Thus $u$ is an order unit. Now renormalize the traces \wrt $u$, that is, $\tau_i$ is replaced by $\sigma_i:= \tau_i/\tau_i(x_i)$. Then $\sigma_i(x_j) = \delta_{ij}$ (Kronecker delta), and in the affine representation \wrt $u$, each $x_j$ simply maps to the $j$th standard basis element. Now $y$ (or more accurately $\hat y$) is a real linear combination of $x_i $, say $\hat y= \sum \alpha_i \hat x_i$. As $G$ has dense range, it easily follows that $\brcs{1, \alpha_1, \dots , \alpha_n}$ is rationally linearly independent, and we have exhibited an order-isomorphic copy of $G$ as a basic critical group.
\qed

In the examples we just computed, we see that the torsion-free part is rank one (also follows from the fact that all proper sets of pure traces are ugly). The torsion part is determined by the elementary divisors in the final form. Here is a simple example. Set $f_1 = (1,1)$, $f_2 = (0,2)$, $f_3 = (\alpha, \beta)$ where $\brcs{1,\alpha,\beta}$ is linearly independent over the rationals, and set $G = \langle f_1, f_2 ,f_3 \rangle = \oplus f_i\Z $. The matrix $A$ is already reduced as far as it can be (if we insist that the  top $2\times 2$ matrix has only integer entries),
$$
A = \( \matrix  1 & 1 \\ 0 & 2 \\ \alpha & \beta \\
\endmatrix
\).
$$
Then $\ker \tau_1 = f_2 \Z$, so we set $x_1 = f_2$; $\ker \tau_2 =( 2 f_1 - f_2)\Z$, so we set $x_2 = 2f_1 -f_2$. But $\langle x_1, x_2 \rangle = \langle 2f_1, f_2 \rangle$, so $G/E \iso \Z \oplus \Z_2$; in particular, this dimension group is not basic. (It is the presence of the $1$ in the $(1,2)$ entry, that ensures that we obtain $2$-torsion; if $f_1 = (1,0)$ instead, then the group would be basic, since we could  divide the second column by $2$).

Now let $n = 3$, and define $f_i$ to be the four rows of the matrix
$$
 \( \matrix  1 &0 & 11 \\ 0 & 1 & 2 \\ 0 & 0 & 12 \\ \alpha & \beta & \gamma \\
\endmatrix
\),
$$
where $\brcs{1, \alpha , \beta , \gamma}$ is rationally linearly independent. Then $x_1 = 12f_1 - 11f_3$ (up to sign), $x_2 = 6f_2 - f_3$, and $x_3 = f_3$. Then the torsion subgroup of $G/E$, that is, $J(B\op)$, is isomorphic to  $\Z_{12} \oplus  \Z_6$, which has $72$ elements, not the expected $12 = 1 \times 2 \times 6$.

We will see (next section) that the invariant really boils down to PH-equivalence, together with an action on the bottom row.

When $n=2$, we saw an example of a basic critical group \st $\tau(G)$ are all isomorphic as $\tau$ varies over all (two) pure traces. When $n> 2$, the corresponding construction $G(n,r)$, does not yield a basic critical group, but we can still construct basic ones with this property.

Let $r$ be a positive real number that satisfies no nontrivial integer polynomial of degree $n$ or less. Then the set $\brcs{1,r, r/(1+r), r/{1+2r}, \dots, r/(1+(n-1)r)}$ is rationally linearly independent. This is an easy exercise, which becomes trivial if we assume $r$ is transcendental. Hence there is a basic critical group whose last row is
$(r, 1/(1+r), \dots, 1/(1+(n-1)r))$. The respective value groups of the pure traces are
$\Z + r\Z, \Z + (r/(1+jr)\Z$ ($1 \leq j \leq n-1$). But these are all isomorphic (multiply $\Z + (1/(1+jr)\Z$  by $1 + jr$; this is an order isomorphism to $(1+jr)\Z + r\Z = \Z + r\Z$).

\SecT 10 Isomorphisms between almost basic critical groups

A critical group of rank $n+1$ is {\it almost basic\/} if it is order isomorphic to a dimension group $G$ given by the matrix $\(\smallmatrix B \\ \alpha \endsmallmatrix \)$ where $B \in \Mn n\Z$; necessarily (in order to have dense image in $\R^n$), $\rk B = n$ and $\brcs{1, \alpha_1, \dots, \alpha_n}$ is rationally linearly independent. As above, we may assume that all the columns of $B$ are unimodular, that is, $B \in \NS_n$. We will show that two almost basic groups (with corresponding $(B, \alpha)$ and $(B', \alpha')$ are order isomorphic iff $B = UB'P$ (with $U \in \gl(n,\Z)$ and $P$ a permutation matrix, i.e., $B$ is PH-equivalent to $B'$) and one of $ \alpha \pm \alpha'P \in r(B)$. We also obtain an internal characterization of almost basic among critical groups, independent of how it is realized, that is, every subset of $\partial_e S(G,u)$ is ugly.

Suppose $r$ and $s$ are irrational real numbers. Then the critical groups of rank $2$ ($n=1$), $\Z + r\Z $ and $\Z + s\Z$ with orderings inherited from the reals, are order-isomorphic iff $r$ is in the PGL($2,\Z$)-orbit of $s$, that this, there exist integers $a,b,c,d$ \st $|ad-bc| = 1$ and $r = (as + b)/(cs + d)$ [ES]. This easily follows from $(as + b)\Z + (cs + d)\Z = \Z + s\Z$ when $\( \smallmatrix a & b \\ c & d \\ \endsmallmatrix\) \in \text{GL}(2,\Z)$. For $n >1$ and basic critical groups, perplexingly, the role of PGL$(2,\Z)$ is replaced by the semi-direct product $\Z^n \times_{\pi \times \rho} (S_n \times \brcs{\pm1})$ where $S_n$ is the symmetric group. This is  abelian by finite, rather different from PGL$(2,\Z)$. A similar, but somewhat more restrictive description for isomorphism classes of almost basic groups, follows from the same result.

\noindent {\it Notation for the statement of the theorem.} Let $B \in \Mn n \Z$ be of rank $n$. Suppose  $ \alpha = (\alpha_1, \dots, \alpha_n) \in \R^{1 \times n}$ is such that $\brcs{1, \alpha_1, \dots, \alpha_n}$  is rationally linearly independent. Form the augmented matrix $\Cal B = \(\smallmatrix B \\ \alpha \\  \endsmallmatrix \)  \in \R^{(n+1) \times n}$. Set $G_{B,\alpha}$ to be the subgroup of $\R^{1\times n}$ generated by the rows of $\Cal B$. Then $G_{B, \alpha}$ is a critical dimension group of rank $n+1$.
If the content of $i$th column of $B$ is $\delta_i \in \Q$, then applying $\Delta^{-1}$ on the right, where $\Delta = \diag (\delta_1, \dots, \delta_n)$, we see that $B'':= B\Delta^{-1}$ is still an integer matrix, but now in $\NS_n$, and $G_{B,\alpha} \iso G_{B'', \alpha \Delta^{-1}}$ as partially ordered abelian groups.
Hence (at a cost of multiplying the entries of $\alpha$ by various fractions of the form $1/k$), we may assume that $B$ is already in $\NS$.

\Lem Theorem \sevone. Let $G_{B,\alpha}$ and $G_{B', \alpha'}$ be almost basic critical groups, where $B, B' \in \NS_n$. If they are order isomorphic, then there exists $C \in \gl(n+1,\Z)$ and $\Delta P \in P(n,\R)^+$ (with $P$ a permutation matrix) \st $C  \Cal B \Delta P =  \Cal B'$. Moreover,
\item{(i)} In the $n,1$ partition of $C = \(\smallmatrix U & c \\ r & t \\ \endsmallmatrix\)$, $c = (0,0,\dots,0)^T \in \Z^{(n-1) \times 1}$, $U \in \gl(n,\Z)$, and $ t\in \brcs{\pm 1}$.
\item{(ii)} $\Delta = \I$ and $UBP = B'$.
\item{(iii)} $\alpha' $ belongs to  one of $\pm \alpha P + r(B)$.{\par}
\noindent In particular,  $G_{B,\alpha} \iso G_{B', \alpha'}$ iff (ii) and (iii) hold.

\Rmk Condition (iii) says that one of $\alpha' \pm \alpha P$ belongs to the row space of $B$.

\Pf  First, suppose that $B$ and $B'$ are in $\NS$, $\alpha$ is given (so that $G_{B,\alpha}$ is a critical group), and $B$ is PH-equivalent to $B'$. Then it is elementary that $G_{B,\alpha} \iso G_{B', \alpha^{\pi}}$ (as partially ordered groups), where $\pi$ effects a permutation of the entries. To see this, suppose $U B  P = B'$ where $U \in \gl(n,\Z)$ and $P$ is a permutation matrix. Let $C = U \oplus 1$. Then
$$
C \( \matrix  B\\ \alpha \\
\endmatrix\)P =   \( \matrix  UBP\\ \alpha P \\
\endmatrix\) =  \( \matrix  B'\\ \alpha P \\
\endmatrix\),
$$
and of course, $\alpha P$ is just a permutation of $\alpha$. By our usual construction, this yields an order isomorphism $G_{B,\alpha} \to G_{B',\alpha P}$.

Thus given full rank $B \in \Mn n \Z$ and $\alpha $ \st $\brcs{1} \cup \brcs{\alpha_i}$ is rationally linearly independent, there exists a terminal $B'' \in \NS_n$ \st $G_{B, \alpha} $ is order isomorphic to $G_{B', \alpha'}$ (where $\alpha$ is obtained from $\alpha$ by applying some weighted permutation to the latter).

Hence we may  suppose that $\alpha$ and $\alpha'$ are given (and satisfy the usual rational linear independence condition), $B$ and $B'$ are terminal forms in $\NS_n$, and there is an order isomorphism $G_{B,\alpha} \to G_{B', \beta}$. We will show (i--iii) hold.

The isomorphism entails there exist $C \in \gl(n+1,\Z)$ and   a weighted permutation matrix with  positive real entries (here factored as diagonal times permutation), $\Delta P$, \st $C \Cal B  \Delta P = \Cal B'$. Partitioning the matrices as we did before and writing $B = \(\smallmatrix\I_s & X \\ 0 & \Cal D \\ \endsmallmatrix\)$ and
$B = \(\smallmatrix\I_{s'} & X' \\ 0 & \Cal D' \\ \endsmallmatrix\)$, in terminal form (thus $\Cal D$ is upper triangular with positive increasing entries along the diagonal, none of the them $1$, etc)
$$
\(\matrix U & c \\ r & t\\
\endmatrix \) \(\matrix {\(\matrix\I_s & X \\ 0 & \Cal D \\ \endmatrix\)} \\ \alpha\\
\endmatrix \) \Delta P = \(\matrix {\(\matrix \I_{s'} & X' \\ 0 & \Cal D' \\ \endmatrix\)} \\ \alpha'\\
\endmatrix \).
$$
Our objective is to show that the column $c = (c_1, \dots, c_n)^T$ is zero, and we achieve this by exploiting the numerous zeros in the matrices. Then it is elementary that $\Delta$ must be the identity and $U B P = B'$, and moreover, $|\det U| =1$ is immediate.

From the equation,
$$\eqalign{
\(\matrix U & c \\ r & t \\ \endmatrix \) \( \matrix B \\ \alpha
\endmatrix \) \Delta & = \( \matrix B'P^{-1} \\ \alpha P^{-1}\\ \endmatrix
\),\qquad \text{we obtain,} \cr
(UB + c \alpha) \Delta & = B'P^{-1} \cr
rB + t \alpha& = \alpha' P^{-1}.\cr
}$$
One of the columns of $B'P^{-1}$, say the $h$th, is the first standard
column basis element. Hence for all $i$,
$$
\((UB)_{ih} + c_i \alpha_h\) \delta_h = \cases 1 & \text{if $i =1$}\\
0 & \text{if $i > 1$.}
\endcases
$$
Hence if $i > 1$, $(UB)_{ih} + c_i \alpha_h = 0$. As the first term and
$c_i$ are integers, and $\brcs{1, \alpha_h}$ is rationally linearly
independent, we deduce $c_i = 0$ (and $(UB)_{ih} = 0$). Assume $c_1 \neq
0$; we will obtain a contradiction.
 
Write $U = \brcs{\gamma_{ij}}$. As the first column of $B$ is the first
standard basis element, we have $(UB)_{i1} = \gamma_{i1}$. Thus
$(\gamma_{i1} + c_i \alpha_1)\delta_1 \in \Z$ (as these are the entries of
a column of $B'P^{-1}$). Hence for $i > 1$, $\gamma_{i1} \in
\delta_1^{-1}\Z$ (as the corresponding $c_i$ are zero). If for some $i >
1$, $\gamma_{i1} \neq 0$, then $\delta_1$ is rational. From $\gamma_{11} +
c_1 \alpha_1 \in \delta_1^{-1}\Z$ together with rational linear
independence of $\brcs{1, \alpha_1}$, we deduce $c_1 = 0$, a
contradiction. Hence $\gamma_{i1} = 0$ for all $i > 1$.
 
Now consider the second column of $UB$; as $B$ is upper triangular,
$(CB)_{i2} = \gamma_{i1}B_{12} + \gamma_{22}B_{22}$, and $B_{22} \neq 0$.
Hence
$$
\gamma_{i1} B_{12} + \gamma_{i2} B_{22} + c_i \alpha_2 \in \frac 1
{\delta_2}\Z.
$$
If $i > 1$, this simplifies to $\gamma_{i2} B_{22} \in
{\delta^{-1}_2}\Z$; thus if $\gamma_{i2}\neq 0$ for some $i > 1$, then
$\delta_2 $ is rational. Hence, $\gamma_{11}B_{12} + \gamma_{12}B_{22} +
c_1 \alpha_2 \in \delta_2^{-1}\Z \subset \Q$. As $c_1 \neq 0$, rational
linear independence of $\brcs{1, \alpha_2}$ is impossible, a
contradiction.
Hence $\gamma_{i2} = 0$ for $i > 1$.
 
Thus the first, second, and $n+1$st columns of the matrix $C \in
\gl(n+1,\Z)$ are
$$
\(\matrix \gamma_{11}\\ 0 \\ 0 \\ \vdots \\ 0 \\ r_1 \\
\endmatrix\), \(\matrix \gamma_{12}\\ 0 \\ 0 \\ \vdots \\ 0 \\ r_2 \\
\endmatrix\), \(\matrix c_1\\ 0 \\ 0 \\ \vdots \\ 0 \\ t \\
\endmatrix\).
$$
These generate a subgroup of rank only two, so that $\rk C < n+1$. This
final contradiction shows that $c_1 = 0$, and thus $c $ is zero.
 
Thus $C = \( \matrix U & 0 \\ r & t\\ \endsmallmatrix \)$, and so $1 = |\det
C| = |t \det U|$. Thus $t = \pm 1$ and $U \in \gl (n,\Z)$, and of course
the equations simplify to
$UB\Delta = B'P^{-1}$. Since $B$ and $B'$ are invertible in $\Mn n \Q$,
this forces $\delta_i \in \Q^+$ for all $i$. In particular, there exists
an integer $N$ \st $N\Delta$ is an integer matrix (with positive entries).
 
As $B$ and $B'$ have all their columns unimodular, so do $UB$ (as $U \in
\gl (n,\Z)$) and $B'P^{-1}$. Thus the content of the $i$th column of
$NUB\Delta$ is $N\delta_i$ while that of $NB'P^{-1}$ is just $N$. Hence
$N\delta_i = N$, so $\Delta = \I$. Thus (finally)
$$\eqalign{
UBP& = B' \cr
rB \pm \alpha P & = \alpha'.
\cr
}$$
This yields (i--iii), and the final statement is a consequence of this and the remarks early in the proof.
\qed

\Lem Proposition \sevtwo. Let $n > 1$, $\alpha = (\alpha_1, \dots,\alpha_n), \alpha' = (\alpha'_1, \dots, \alpha'_n) \in \R^n$ be such that both $\brcs{\alpha_i}_{i=1}^n \cup \brcs{1}$ and $\brcs{\alpha'_i}_{i=1}^n \cup \brcs{1}$ are linearly independent over the rationals. The basic critical dimension groups $G_{\alpha}$ and $G_{\alpha'}$
(generated by $\brcs{e_1, \dots, e_n, \alpha}$ and $\brcs{e_1, \dots, e_n; \alpha'}$ respectively)  are order isomorphic iff $\alpha'$ is in the orbit of $\alpha$ under the action of $\Z^n \times_{\Pi \times \rho}(S_n \times \Z_2)$.

\Pf Here $B = B' = \I$, so the criterion of the theorem simplifies to $\alpha'
\pm \alpha P \in \Z^{1\times n}$.
\qed

\Lem Corollary \sevthr. Almost basic critical groups of rank at least three admit no nontrivial
order-automorphisms.

\Pf From $C \Cal B \Delta P = \Cal B$ (the order-automorphisms on a
critical group automatically extend to order-automorphisms of the closure,
$\R^n$, hence must be given by weighted permutation matrices), the
preceding yields $\pm \alpha + rB = \alpha P $. If $\pi$ is the
permutation induced by $P$, and $\pi(i) = j \neq i$ for some $i$ and $j$,
then
$\alpha_j \pm \alpha_i \in \Z$; but this contradicts the rational linear
independence of $\brcs{1, \alpha_i, \alpha_j}$. Hence $P$ is the identity.
Thus by the preceding $B = UBP = UB$; as $B$ is of full rank, this forces
$U $ to be the identity, and thus the only automorphism is the identity.
\qed

This contrasts with the critical groups discussed in [H]; those arise from integral orders in totally real fields with one real embedding discarded, and are classified by their ideal class structure. In those cases, there are plenty of order automorphisms, arising from some of the units in the number field.

\def\tfrk{\text{tf\,}\rk}

For an abelian group $J$, the torsion-free rank of $J$, that is, the rank of $J$ modulo its torsion subgroup, is denoted $\tfrk J$.

\Lem Proposition \sevfou. Let $G$ be a critical group of rank $n +1$, with $n > 1$. Let $u$ be any order unit for it. The following are equivalent.
\item{(a)} the torsion-free rank of $G/E(G)$ is one;
\item{(b)} for all $\sigma \in \partial_e S(G,u)$, the intersection $\cap_{\tau \in \partial_e S(G,u) \setminus \brcs{\sigma}} \ker \tau$ is nonzero;
\item{(c)} there exists a basic critical group $G'$ \st $G \otimes \Q$ is order-isomorphic to $G' \otimes \Q$;
\item{(d)} every proper subset of $\partial_e S(G,u)$ is ugly;
\item{(e)} $G$ is almost basic.

\Pf (a) iff (b): Property (b) is equivalent to $\rk E(G) = n$ (since the sum $\sum x_i \Z$ is direct), which is equivalent to $\tfrk G/E(G) = 1$.

\noindent (b) implies (c). For each $i$, there exists $x_i \in G$ unique \wrt the properties $\cap_{\tau_j \in \partial_e S(G,u)\setminus \brcs{\tau_i}} \ker \tau_j = x_i \Z$ and $\tau_i (x_i) > 0$. Then $E(G):= \sum x_i \Z$ is free of rank $n$ and there exists $y \in G$ \st $G_0 = E(G) \oplus y\Z$ is of finite index in $G$. As $G$ is dense in $\R^n$, so is $G_0$.

As a subgroup of $G$ of rank less than $n+1$, $E(G)$ is discrete; being of rank $n$, any $\Z$-basis for it is also an $\R$-basis for $\R^n$. Hence there exist $\alpha_i \in \R$ \st $y = \sum \alpha_i x_i$. Density of $G_0$ in $\R^n$ entails that $\brcs{1, \alpha_i,\dots, \alpha_n}$ be rationally linearly independent, and  $G_0$ is a dimension group \wrt the strict ordering, which  obviously agrees with the relative ordering inherited from $G$, and its pure traces are the restrictions of $\tau_i$, which we will also call $\tau_i$.

Set $u = \sum x_i$, so that $\tau_i (u) = \tau_i (x_i) > 0$ for all $i$. Thus $u$ is an order unit in both $G_0$ and $G$. Normalize the traces of $G_0$ \wrt $u$---the pure traces are now $\tilde \tau_i$ given by $\tilde \tau_i (g) = \tau_i(g)/\tau_i (u)$. The normalized traces now satisfy $\tilde \tau_i (x_j) = \delta_{ij}$ (Kronecker delta). Hence the embedding $(G_0, u) \to \Aff S(G_0,u)$ realizes $G_0$ as a basic critical group.

Since $G_0$ is of finite index in $G$, $G_0 \otimes \Q = G \otimes \Q$.

\noindent  (c) implies (d). For any trace $\tau$ on any dimension group $G$, $\ker \tau \otimes 1_{\Q} = (\ker \tau ) \otimes \Q$. Thus $\rk \ker \tau = \rk \ker (\tau \otimes 1_{\Q})$. As $G$ is critical, every subgroup of less rank than $n+1$ is discrete, and the result follows.

\noindent (d) implies (e). For a pure trace $\tau_i$, let  $\Omega(i)$ be the complement of $\brcs{\tau_i}$ in $\partial_e S(G,u)$. As $\Omega(i)$ is ugly, $\cap_{\tau \in \Omega(i)} \ker \tau$ is not zero, and being discrete and spanning $\Omega(i)^{\perp}$ over the reals, it must be rank one. As it is a subgroup of a free group, it is free, so it equals $x_i \Z$ for some $x_i \in G$, and we may assume $\tau_i(x_i) $ is positive. Now we are in a position to use the method of proof in (b) implies (c), coming up with a basic critical group $G_0$ of finite index in $G$. There thus exists an integer $N$ \st $NG \subseteq G_0$, and $NG$ is obviously order isomorphic to $G$, while $G_0 \subset \Z^n$. Any subgroup of a free group is free, so we can find the desired basis.

\noindent (e) implies (a). Trivial.
\qed

\SecT  11 Unperforation of quotients

In this section, we want to ensure that the quotient pre-ordered groups of almost basic critical groups by kernels of subsets of $\partial_e S(G,u)$ are themselves almost basic; the crucial property is that these quotients are unperforated. We will prove the following. This construction is what motivated the $\Cal J(B)$ invariant of PH-equivalence.

\Lem Proposition \eigone. Let $G$ be an almost basic critical group of rank $n+1$. Let $\Omega \subset \partial_e S(G,u)$, and define $L =\ker \Omega:= \cap_{\tau \in \Omega} \ker \tau$. Then $G/\ker \Omega$, equipped with the quotient ordering, is an almost basic critical group with pure trace space $\Omega$.

 This boils down to showing the quotient is unperforated, something that is obvious for basic critical groups (and the quotients are themselves basic critical groups), but  not so obvious for almost basic ones. This provides an alternative path to the definition of $I(B_{\Omega})$ as the torsion subgroup of $G_{\Omega}/E(G_{\Omega})$ where $ G_{\Omega} = G/\ker \Omega$ (for $\Omega \subset \partial_e S(G,u)$, the latter identified with $\brcs{1, \dots, n}$).

The following is a slight improvement on [BeH, Appendix B, Propositions 1 and 2], not covered by  the results there.

\Lem Lemma \eigtwo. Let $(G,u)$ be a simple unperforated group with order unit, and let $L$ be a convex subgroup of $G$ \st $G/L$ is torsion-free. Suppose that the closure of the image of $L$, $\overline L$, in $\Aff S(G,u)$, contains a  subgroup of the form $D + P$, where $D$ is a rational vector space and $P$ is generated by nonnegative elements of $ \Aff S(G,u)$, and $\overline L/(D + P)$ is torsion. Then equipped with the quotient ordering, $G/L$ is unperforated.

\Rmk For example, if $G$ is basic, say with generators $\brcs{e_i; \sum \alpha_j e_j}$, and $L = \ker S$ (where $S \subset \brcs{\tau_i}$), then $L$ is generated by $\brcs{e_i}_{\Set{i}{\tau_i \not\in S}}$. Each $e_i$ has image in $\Aff S(G,u) = \R^n$ as $e_i$ itself, which is nonnegative in the affine function space (of course, the $e_i$ is not in the positive cone of $G$, since the ordering is the strict one). By [BeH, Appendix B], the quotient is nicely behaved. {\par}
The lemma above removes the density condition on $D  + P$ (that it be dense in $\overline L$) [op\.cit.], and replaces it with a  different requirement. This is particularly useful when $L$ is already discrete, hence closed in the affine representation; then $D = 0$, but $P$ need only be a subgroup; this will automatically be closed, so that $L $ need not equal $P$. But the lemma here says that sufficient for unperforation is that $\rk P = \rk L < \infty$, which is easy to verify for almost basic critical groups.

\Pf The convexity condition (which in the simple case boils down to $L \cap G^+ = \brcs{0}$) is sufficient to guarantee that the quotient pre-ordering is a partial ordering, that is, an element that is both positive and negative must be zero.

If $kg + L = L$, then torsion-freeness of the quotient entails $g \in L$. Hence we may assume that $kg + L \in (G/L)^+ \setminus\brcs{0}$.

Suppose $g \in G$ and $k$ is a positive integer \st $k g + L \in G^+ \setminus \brcs{0}$. We may thus find $x \in L$ \st $kg + x $ is an order unit. Let $\epsilon = \inf_{\sigma \in S(G,u)} \sigma(kg + x) =  \inf_{\sigma \in S(G,u)} \widehat{(kg + x)} (\sigma) > 0 $. There exists a positive integer $N$ \st $N\hat x$ is in the norm closure of $D + P$. Select an integer $M$ to be determined (as a function of $k$ and $N$).

We may find $d \in D$ and $p \in P$ \st $\|N\hat x - d - p \| < \epsilon/M$. There exists (from the definition of $D$), $d' \in D$ \st $d = Nk d'$; so $\|N\hat x - Nkd' - p \| < \epsilon/M$. We may write $p = p_1 - p_2$ where $p_i \geq 0$ and $p_i \in P$ (in particular, $p_i \in \overline L$).

There exists $f \in L$ \st $ \| \hat f - d'\| < \epsilon/M$ and $q_i \in L$ \st $\| \hat q_i - p_i \| < \epsilon/M$. In particular $\hat q_i \geq -\epsilon \pmb 1/M$  as functions on $S(G,u)$.

Set
$$
z = Nkg + Nkf  +Nk q_1 = Nk(g + f + q_1).
$$
If we can show $z \in G^+\!\!$, then as $G$ itself is unperforated, it would follow that $g + f + q_1 \in G^+$, and so $g + L$ is in the positive cone of the quotient. So it suffices to show $z \in G^+$.

We have
$$\eqalign{
z - N(kg + x) & = Nkf + Nkq_1 - (Nx - Nkf - q_1 + q_2) -Nkf - q_1 + q_2 \cr
& = Nkq_1 + q_2 - (Nx - Nkf - q_1 + q_2); \quad\text{evaluating at $\sigma \in S(G,u)$,} \cr
\sigma(z) & \geq N\sigma(kg + x) +Nk \sigma(q_1) + \sigma(q_2) - \| N\hat x - Nk\hat f -\hat q_1 + \hat q_2\|\cr
& \geq N\epsilon -\frac{Nk \epsilon}M - \frac{\epsilon}M - \| N \hat x - Nk d' -p\| - Nk \| \hat f - d'\| - \|\hat q_1 - p_1 \| -  \|\hat q_2 - p_2 \| \cr
& \geq \epsilon \(N -\frac{Nk +1}M  - \frac 1M - \frac{Nk\epsilon}{M} - \frac{2}M\) \cr
& = \epsilon  \(N -\frac{2Nk + 4}M\).
}$$
If we select $M > 2k + 4/N$ (e.g., $M = 2k + 6$), $\hat z $ is strictly positive, so that $z$ is an order unit of $G$, and we are done.
\qed

\Lem Corollary \eigthr. Suppose $G$ is a simple dimension group with an order unit $u$. Let $L$ be a convex subgroup of $G$ \st $G/L$ is torsion-free and the image of $L$ in $\Aff S(G,u)$ is discrete. Sufficient for $G/L$ to be a simple dimension group (\wrt the quotient ordering) is that there exist a subgroup $L_0$ of $L$ \st $L/L_0$ is torsion, and the image of $L_0$ is generated by its nonnegative elements (\wrt the usual ordering on $\Aff S(G,u)$). In this case, the trace space of $(G/L, u+L)$ is a closed face of $S(G,u)$.

\Pf Since the image of $L$ is discrete, its image is already closed in $\Aff S(G,u)$; the hypothesis ensures that $\hat L_0 = P$ satisfies $\hat L/P$ is torsion, so the preceding applies with $D = 0$. Hence $G/L$ is unperforated. Simplicity is automatic.

Since $L/L_0$ is torsion, if $\tau \in S(G,u)$ kills $L_0$, it automatically kills $L$. Hence $L^{\vdash} = L_0^{\vdash}$. Let $P^+ = P \cap \Aff S(G,u)^+$ (the latter with the usual, not the strict ordering), so that $P = P^+ - P^+$. Since $L_0$ maps to $P$, and $P^{\vdash} = (P^+)^{\vdash}$, we have that $L^{\vdash}$ is a (closed) face, call it $F$, of $S(G,u)$. In particular, $F$ is a Choquet simplex.

Let $\phi$ be a trace of $G/L$; then $\phi$ induces a trace of $G$, with kernel containing $L$. Thus $\phi \in L^{\vdash} = F$. Conversely any element of $F$ kills $L$ and thus induces a trace on $G/L$. Hence the map $S(G/L,u + L) \to F$ is  an affine  bijection; it is obviously continuous, so by compactness of $S(G/L,u+L)$, it is an affine homeomorphism.

Select an element $h \in \Aff F$; this lifts to an element $j \in \Aff S(G,u)$. Given $\epsilon$, there exists $g \in G$ \st $\|\hat g - j \| < \epsilon$, that is, $\sup_{\sigma \in S(G,u)} |\sigma(g) - j(\sigma)| < \epsilon$. This implies $\sup_{\sigma \in F}  |\sigma(g) - j(\sigma)| < \epsilon$, and together with the affine homeomorphism, this forces the image of $G/L$ to be dense in $\Aff F$, hence in its affine representation (\wrt $u+L$). As $G/L$ is unperforated and simple, its ordering must be the strict one inherited from affine functions on a Choquet simplex, and thus $G/L$ is a dimension group.
\qed

\Lem Corollary \eigfou. If $G$ is an almost basic critical group and $\Omega \subset \partial_e S(G,u)$, then $G/\ker \Omega$ is a simple dimension group whose pure trace space is $\Omega$.

\Pf Let $F$ be the face spanned by $\Omega$ (since $\Aff S(G,u)$ is a finite dimensional simplex, it is simply the convex hull of $\Omega$). By Proposition \sevfou(c), there exists a basic critical group $G_0$ of finite index in $G$ (whose relative ordering agrees with its usual one). Then $\ker \Omega \cap G_0$ is generated  by elements with nonnegative image in $\Aff S(G,u)$, and this is of finite index in $\ker \Omega$. By the result above, $G/\ker S$ is a simple dimension group, and its pure trace space is just the set of extreme points of $F$, which is $\Omega$.
\qed

\noindent{\it Connections to PH-equivalence.}  This was the starting point for the development of $(J(B\op_{\Omega}))_{\Omega \subset S}$, the directed family of PH-invariants; when $G$ is generated by the row space of $B$ and $\alpha$, then the torsion subgroup of $G_{\Omega}/E(G_{\Omega})$ is just $I(B_{\Omega})$.

For almost basic critical groups, $G_{B,\alpha}$, $G_{B', \alpha'}$ with $\(\smallmatrix B \\ \alpha \\ \endsmallmatrix\), \(\smallmatrix B' \\ \alpha' \\ \endsmallmatrix\) \in \R^{n \times (n+1)}$ \st $B,B' \in \Mn n \Z$ and $\det B, \det B' \neq 0$, we immediately reduce to the case that $B,B' \in \NS_n$, by factoring out a positive diagonal matrix, as in section\,7. Then by Theorem \sevone, $G_{B,\alpha}$ is order-isomorphic to $G_{B', \alpha'}$ iff $B$ is PH-equivalent to $B'$ and the permutation involved in the PH-equivalence sends to $\alpha$ to $\alpha'$.

For example, if
$$
\Cal B_{\alpha} := \(\matrix 1 & 0 & 15 \\ 0 & 1 & 2 \\ 0 & 0 & 30 \\ \alpha_1& \alpha_2 &\alpha_3 \\ \endmatrix \) \text{ and } \Cal B'_{\alpha'} := \(\matrix 1 & 0 & 5 \\ 0 & 1 & 6 \\ 0 & 0 & 30 \\   \alpha'_1& \alpha'_2&\alpha'_3 \\ \endmatrix \),
$$
given $\alpha = (\alpha_1, \alpha_2, \alpha_3)$ (with $\brcs {1, \alpha_1, \alpha_2, \alpha_3}$ rationally linearly independent), there is no choice of  $\alpha' = (\alpha_1, \alpha_2, \alpha_3)$ \st $G_{B, \alpha} $ is order isomorphic to $G_{B', \alpha'}$, since from Example \thrthr\ (first two matrices), $B$ and $B'$ are not PH-equivalent.

If we set $B = B'$ to be the leftmost example in Example\,\thrthr\ (the $3\times 3$ integer matrix in $\Cal B_{\alpha}$ above), and let $\alpha = (\sqrt 2, \sqrt 3, \sqrt 5)$ and $\alpha' = (\sqrt 3, \sqrt 2, \sqrt 5)$, then even though the integer matrix parts are the same, the resulting critical dimension groups are not order isomorphic, because the permutation $\pi = (12)$ and its corresponding matrix $P = \(\smallmatrix 0 & 1 \\ 1 & 0 \\ \endsmallmatrix \) \oplus (1)$ does not fix $B$, as follows from Proposition \twoone, with no invertible elements in the column (modulo $d = 30$).

\def\lArrow #1;#2.#3{{}_R #1 \leftarrow {}_R #2\: #3}\def\lker{{}_L \ker }

\SecT Appendix A. A general duality result

This appendix gives fairly general duality results about orbits under natural actions, that appear in section\,6. Here $R$ will be a not necessarily commutative ring (of course with $1$), but not much additional effort is required to prove the corresponding results over noncommutative rings. 

Let $R$ be a ring, and $n >k$ be positive integers. A matrix with entries from a ring will be called {\it invertible\/} if it is square and two-sided invertibleÑsometimes we add {\it two-sided,} for emphasis. We follow  [C] in saying  an $n \times k$ matrix $M$ is {
\it completable\/} if there exists an $n \times (n-k)$ matrix $W$ \st the $n
\times n$ matrix $\(W \ \ M \)$ is invertible.
Invertibility of this matrix is equivalent to the columns constituting an $R$-basis for $R_R^{n
\times 1}$ (as a right $R$-module).

If instead, $n < k$, then we say $M$ is {\it completable\/} if there exists a $(k-n) \times k$ matrix $W$ \st the $k \times k$ matrix $\(\smallmatrix W \\ M \\ \endsmallmatrix \)$ is invertible.
If $n = k$, then completable is simply invertible. 

These notions date back to the origins of K-theory.

The ring of $n \times n$ matrices will be denoted $\M_n R$, but the set of non-square rectangular matrices with $k$ rows and $n$ columns will be denoted $R^{k \times n}$. We denote by $\gl(k,R)$ (or simply $\gl (k)$ if no ambiguity will result) the group of invertible matrices in $\M_n R$. The  group of invertible elements of $R$ will be denoted $R^{\times}$. 

The next two lemmas are obvious.

\Lem Lemma \nappone. Let $M \in R^{n \times k}$ with $n > k \geq 1$. The following are equivalent.
\item{(a)} $M$ is completable;
\item{(b)} the set of columns of $M$ can be enlarged to a basis of size $n$ of $R^{n\times 1}$ as a right $R$-module;
\item{(c)} there exists a right $R$-submodule $L$ of $R^{n\times 1}$ that is free on $n-k$ generators \st $L \oplus M(R^{k\times 1}) = R^{n\times 1}$.

\Lem Lemma \napptwo. Let $M \in R^{n \times k}$ with $k > n \geq 1$. The following are equivalent.
\item{(a)} $M$ is completable;
\item{(b)} the set of rows of $M$ can be enlarged to a basis of size $k$ of ${}_R R^{k\times 1}$ as a left $R$-module;
\item{(c)} there exists a left $R$-submodule $L$ of $R^{1\times k}$ that is free on $k-n$ generators \st $L \oplus(R^{1\times n})M = R^{1\times k}$.

\Lem Lemma \nappthr. Let  $A,B,C$ be respectively in  $\M_{n-k} R$, $R^{(n-k) \times k}$, and $\M_k R$ with $C \in \gl (k,R)$, and set $E: = \(\smallmatrix A & B \\ 0 & C \endsmallmatrix  \) \in \M_n R$. Then $A \in \gl(n-k,R)$ iff $E \in \gl(n,R)$.

\Pf Suppose that $E$ is invertible. Let $U$ be the inverse of $E$; partitioned in the same way as $E$, we can write $U = \(\smallmatrix Q & S \\ T & V \\ \endsmallmatrix\)$, and  we have 
$$\eqalign{
\(\matrix \I_{n-k} & 0 \\ 0 & \I_k \\ \endmatrix \) & = UE = \(\matrix QA & QB + SC \\ TA & TV + VC \\ \endmatrix \)\cr
& = EU = \(\matrix AQ + BT & AQ + BT \\ CT & CV\\ \endmatrix \).\cr
}$$
From $CT = 0$ (second matrix, lower left), and invertibility of $C$ (two-sided invertibility, that is), we deduce $T = 0$ (of the appropriate dimensions). From the upper left corners of each matrix, we then have $AQ = QA = \I_{n-k}$, so $A \in \gl (n-k,R)$.

Conversely, suppose that $A$ is invertible. Multiply $E$ on the left by $A^{-1} \oplus C^{-1} \in \gl(n,R)$, resulting in  $\(\smallmatrix \I_{n-k}  & A^{-1}B\\ 0 &\ I_{k} \endsmallmatrix\)$; this has  inverse $\(\smallmatrix \I_{n-k}  &- A^{-1}B\\ 0 &\ I_{k} \endsmallmatrix\)$, so $E$ is a product of two invertible matrices, hence is invertible. \qed

A ring $R$ is {\it stably finite\/} if for all matrix rings $\M_n R$ and elements $x \in \M_n R$, $x$ is right invertible implies $x$ is invertible. This is a two-sided condition, and is equivalent to onto right module homomorphisms $\Arrow x; R^{n\times 1}. R^{n\times 1}$ always being    isomorphisms, or equivalently, if the left module homomorphism $\lArrow R^{1\times n}; R^{1 \times n}. x$ is onto, then $x$ is an isomorphism. 

The ring $R$ has the property that {\it stably free modules are free\/} (SFF) if whenever $R^n \iso R^k \oplus Q$ as right $R$-modules, then $Q$ is free. This property is also right/left symmetric. This property also harkens back to the origins of K-theory,  e.g., what was formerly Serre's conjecture, now known as the Quillen-Suslin theorem.

Finally, $R$ has {\it invariant basis number\/} (IBN) if $R^m \iso R^n$ as right $R$-modules implies $m=n$. This is again left/right symmetric, and it is easy to check that stably finiteness implies IBN. The reverse implication is well known not to be trueÑe.g., the ring generated by the unilateral shift and its transpose (defined on $l^2(\Z^+)$) has the IBN property but is not stably finite. However, SFF and IBN together imply stable finiteness. ({\it Proof\/}: Suppose  that  the right $R$-module homomorphism $\Arrow x;R^n . R^n$ is onto. It splits since the image is free;  this yields  $R^n \iso R^n \oplus Q$ for the  module $Q = \ker x$. SFF implies that $Q$ is free, and IBN entails that it must be free on zero generators, hence zero. So $x$ is an isomorphism.)

For the computations in section\,6, we only deal with rings of the form $R = \Z_d$. All of these, and $\Z$ itself, satisfy both SFF and stable finiteness (the latter being trivial, since the rings are commutative). 

Kaplansky [K, p\,498] had a limited definition of {\it Hermite rings.} TY Lam [La, p\,26] defines Hermite to mean  a ring satisfying SFF. Cohn [C, 0.4] refers to a ring satisfying SFF and IBN as an Hermite   ring (Charles Hermite was French, so the initial {\it H\/} is   pronounced as a stop, requiring {\it an\/}, not {\it a\/}; this practise is adopted in [K] and many   subsequent papers).   To avoid confusion, particularly in view of the main subject of this paper, we will not use this term at all, nor IBN. 

Sometimes, to emphasize the chirality of a module (left or right) over the ring $R$, we place a subscripted $R$ beside it: thus ${}_R Q$ means $Q$ considered as a left module, and $Q_R$ means as a right $R$-module. 

\Lem Lemma \nappfou. Let $R$ be a stably finite ring satisfying SFF. Suppose $n > k$ are positive integers. Let $M \in R^{n\times k}$. The following are equivalent.
\item{(a)} The set of columns of $M$ is a right $R$-module basis for a free direct summand of $R_R^{n\times 1}$.
\item{(b)} There exists $W \in R^{n \times (n-k)}$ \st $U:= \( W \ \ M \) \in \gl(n,R)$.
\item{(c)} The left $R$-submodule of ${}_R R^{1 \times n}$ spanned by the rows of $M$ is all of ${}_R R^{1 \times n}$.

\Rmk It is  not sufficient (in (c)) that the row space be free on $k$ generators---$R = \Z$ yields an example.

\Rmk The right versus left hypotheses are important (unless $R$ is commutative), as shown in Example\, \nappfiv.

\Pf (a) implies (b). Let $V$ be the (right) submodule of $R^{n\times 1}$ spanned by the columns of $V$; by hypothesis, it is free with $k$ generators, and there exists a  submodule $X$ of $R^{n\times 1}$ \st $R^{n\times 1} = X \oplus V$. Since $V$ is free, SFF implies that $X$ is free, and stable finiteness guarantees that it is free on $n-k$ generators. Label them $w_1, w_2, \dots, w_{n-k}$, and let $W$ be the resulting $n \times (n-k)$ matrix (whose $i$th column is $w_i$). From the direct sum decomposition, the set of columns of $U =  \( W \ \ M \) $ is a basis for $R^{n\times 1}$. Thus the map $\Arrow U; R^{n\times 1}.R^{n\times 1}$ is a homomorphism of right $R$-modules which is onto. Stable finiteness now yields that $U$ is (two-sided) invertible. 

\noindent (b) implies (c). Given $U =  \( W \ \ M \) $ invertible, define the map of left modules, ${}_R R^{1\times n} \leftarrow {}_R R^{1\times n}:U$ (given by right multiplication of course, hence the weird notation). This is an isomorphism (since $U$ is invertible). Now consider the images of $M$ and $W$ separately. Identifying ${}_R R^{1\times (n-k)}$ and ${}_R R^{1\times k}$ with respectively the submodules of ${}_R R^{1\times n}$ having zeros in the bottom $k$ positions and zeros in the top $n-k$ positions, we have that the respective ranges satisfy $({}_R R^{1\times n})W \subset {}_R R^{1\times (n-k)}$ and $({}_R R^{1\times n})M \subset {}_R R^{1\times (k)}$. Since $U$ (as a homomorphism of left modules) is onto, the sum of the two ranges is all of $({}_R R^{1\times n})$, and since the respective images have zero intersection, it follows that $({}_R R^{1\times n})M = {}_R R^{1\times (k)}$.  But this is precisely condition (c). 

\noindent (c) implies (b). Again, view $M$ as a left module homomorphism ${}_R R^{1\times k} \leftarrow {}_R R^{1\times n}:M$. Hypothesis (c) says that $M$ is onto, so splits. Hence there exists $V \in R^{k \times n}$ \st $VM$ is the identity on ${}_R R^{1\times k} $ (one of the peculiarities of  left module homomorphisms is that they compose in the correct order, unlike what we're used to with right modules). Moreover, we have a direct sum decomposition ${}_L\ker M \oplus Z = {}_R R^{1\times n}$, where ${}_L\ker M $ is the kernel of $M$ as a left module homomorphism, and $Z$ is obtained from the splitting, and moreover, left multiplication by $M$ induces an isomorphism $ {}_R R^{1 \times k} \leftarrow Z$. 

As in the proof of (a) implies (b), SFF and stable finiteness imply that ${}_L\ker M$ is free on $n-k$ generators. Identifying ${}_R R^{1\times k}$ with the submodule of ${}_R R^{1\times n}$ consisting of elements whose leftmost $n-k$ entries are zero (and similarly ${}_R R^{1\times (n-k)}$ with the obvious complementary submodule), we define $W \in R^{n \times (n-k)}$  acting (on the right of course) as an isomorphism ${}_R R^{1 \times (n-k)} \leftarrow {}_L\ker M $, and zero on $Z$. 

Now define the $n\times n$ matrix $U = \( W \ \ M\)$, and observe that is range (as a left module homomorphism) is all of ${}_R R^{1 \times n}$. As $R$ is stably finite, this implies that $U$ is invertible.

\noindent (b) implies (a).  Let $Z$ be the inverse of $U = (W \ \ M)$. Invertibility entails that the  set of columns of $U$ is a basis for $R^{n\times 1}$. Let $V_1$ be the right $R$-module span of the set of columns of $W$, and $V_2$ the span of the set of columns of $M$.

Then $R^{n\times 1} = V_1 + V_2$; from the fact that the set of all columns is a basis, we have $R^{n\times 1} = V_1 \oplus V_2$. Since any subset of a basis is itself a basis for the submodule it generates, we have that $V_2$ is free, necessarily on the $k$ generators arising from the columns of $M$. 
\qed

Suppose that $n > k$. Define the following subsets of rectangular matrices over $R$.
$$\eqalign{
C(n,k) & = \Set{M \in R^{n \times k}}{M \text{ is completable}}\cr
C'(n,k) & = \Set{M \in R^{k \times n}}{M \text{ is completable}}\cr
F(n,k) & = \Set{M \in R^{n \times k}}{M \text{ contains a set of $k$ rows whose corresponding matrix is invertible}}\cr
F'(n,k) & = \Set{M \in R^{k \times n}}{M \text{ contains a set of $k$ columns whose corresponding matrix is invertible}}\cr
F_i (n,k) & = \Set{M \in R^{n \times k}}{M \text{ contains exactly $i$ sets of $k$ rows whose corresponding matrix is invertible}} \cr&\qquad \text{for each $i = 1,2, \dots, {n \choose k}$}\cr
F_i' (n,k) & = \Set{M \in R^{k \times n}}{M \text{ contains exactly $i$ sets of $k$ columns whose corresponding matrix is invertible}}\cr&\qquad \text{for each $i = 1,2, \dots, {n \choose k}$}.\cr
}$$

If $R$ is a local ring (in the not necessarily commutative setting, this means that $R$ modulo its Jacobson radical is a division ring), then $C(n,k) = F(n,k)$. It is probably true that $C(n,k) = F(n,k)$ for {\it some\/} $n > k$ implies that $R$ is local. 

It can happen that some of the $F_i(n,k)$ are empty when $i$ is at or near the maximum. For example, if $R $ is the finite field $\Z_p$ with $p$  prime, then with $n = 5$ and $k=2$, ${5 \choose 2} = 10$, and $F_{10}(5,2)$ is empty for $p = 2$ or $3$, while $F_9(5,2)$ is empty if $p = 2$. However if $R$ is an infinite field, then $F_{{n\choose k}}(n,k)$ is generic (and if the field is $\R$ or $\C$, its complement in $F(n,k)$ is just a union of lower dimensional varieties). 

Set $\Cal D_n$ to be the group of invertible diagonal $n \times n$ matrices with
entries from $R^{\times}$ (the group of invertible elements of $R$), and let $\Cal P_n$ denote the group of $n \times n$
permutation matrices (regarded as elements of GL$(n,R)$). The group they
generate (consisting of weighted permutation matrices), will be denoted $W(n)$.
Then $W(n)$ acts from the left on  each of the sets $C(n,k), F(n,k), F_i (n,k)$ and $\gl (k,R)$ acts on the right,
yielding a $W(n) \times \gl(k,R)$ action; similarly, $\gl(n,R)$ acts from the left on $C'(n,k), F'(n,k), F'_i (n,k)$, yielding an action of $\gl(k,R) \times W(n)$ on each of these.

We see that $F(n,k) = \dot\cup F_i(n,k)$ ($i=1,2,\dots, {n\choose k}$), and it is not difficult to see that $F(n,k) \subset C(n,k)$ (below). We denote their orbit spaces by replacing the roman capital letters by their script forms, e.g., $\Cal C(n,k)$, $\Cal C'(n,k)$, etc.

We will obtain what amounts to duality by showing that there exists a natural bijection $\Cal F(n,k) \to \Cal F'(n, n-k)$, which restricts to dualities $\Cal F_i(n,k) \to \Cal F'_i(n, n-k)$ for each $i =1, \dots, {n \choose k}$ (conveniently, ${n \choose k} = {n \choose {n-k}}$). The corresponding groups implementing the actions are $W(n) \times \gl (k,R)$ and $\gl(n-k,R) \times W(n)$. When $R$ is commutative (or more generally admits an anti-automorphism), there are bijections $\Cal F'_i(n, n-k) \to \Cal F_i(n, n-k)$ (determined by composing the anti-automorphism with transpose), yielding bijections $\Cal F_i(n,k) \to \Cal F_i(n, n-k)$; if the anti-automorphism is either the identity or involutive, these are dualities.

For $C(n,k)$, at the moment, the situation requires an additional hypothesis: that $R$ be stably finite and satisfy SFF (even if $R$ is commutative). Then there is a duality $\Cal C(n,k) \to \Cal C'(n,n-k)$ that extends the duality $\Cal F(n,k) \to \Cal F'(n, n-k)$, and the same comments about the presence of an anti-automorphisms apply.  

My colleagues, Kirill Zaynullin and Damien Roy, pointed out that if  $R = F$ is a field, then  one of the dualities, $\Cal F(n,k) \to \Cal F(n,n-k)$, is implied by  the usual Grassmannian duality, $\text{Gr}(k,n) \to \text{Gr}(n-k,n)$. To see this, pick $M \in F^{n \times k}$; its columns form a basis for a $k$-dimensional subspace of $F^{n\times 1}$; the right action by $\gl(k,F)$ (given by elementary column operations) has no effect on the subspace, and the duality takes the transpose and looks at its kernel.

If $R$ is commutative, we will see that the dualities $\Cal F(n,k) \to \Cal F(n,n-k)$ and $\Cal F_i(n,k) \to \Cal F_i(n,n-k)$ are  implemented by $[M] \to [N]$ where the columns of $N$ constitute a right basis for the kernel of the map $\Arrow M^T; R^{n \times 1}.R^{k\times 1}$. This is particularly simple, as $M^T N = 0$ iff $N^T M = 0$. However, if $R$ is not commutative, the transpose does not do what we expect. For example, at one point in the argument of the commutative case, we would use the obvious fact that if $g \in \M_k R$ is invertible, then so is $g^T$. This is no longer the case in the noncommutative situation; in fact, it is almost never true. Example \nappfiv\ illustrates this. This is a minor modification of an  example in [GKKL], the main result of which is that a ring $R$ in which the transpose of invertibles is always  invertible, must be commutative modulo its Jacobson radical.

\Lem Example \nappfiv. [GKKL] Suppose $R$ is a ring and there exist elements $x,y \in R$ \st $xy-yx$ is invertible (in $R$). Then there exists $g \in \gl(2,R)$ \st $g^T$ is a two-sided zero divisor in $\M_2 R$.

\noindent Set $g = \( \smallmatrix 1 & x \\ y & xy \\ \endsmallmatrix \)$. Elementary column operations (using {\it right\/} multiplications) reveal that $g$ is invertible, in fact even a  product of elementary matrices; explicitly,  
$$g^{-1} = \( \matrix 1+xe^{-1}y &\ \ \ - x e^{-1}\\ -e^{-1}y & e^{-1} \\ \endmatrix \) = 
\( \matrix 1 &- x \\ 0 & 1 \\ \endmatrix \)
\( \matrix 1& 0 \\ 0 & e^{-1} \\ \endmatrix \)
\( \matrix 1 &0 \\ -y & 1 \\ \endmatrix \),
$$ 
where $e = xy-yx$. However, $g^T = \( \smallmatrix 1 & y \\ x & xy \\ \endsmallmatrix \)$ kills the column $\( \smallmatrix -y \\ 1\\ \endsmallmatrix \)$ and the row $(-x \ \ 1)$, so $a := \( \smallmatrix -y \\ 1\\ \endsmallmatrix \)(-x \ \ 1) = \( \smallmatrix yx & -y \\-x &  1\\ \endsmallmatrix \)$ satisfies $ag^T = g^T a = 0$. 

Any division ring which is not commutative admits such a pair $x,y$, as does any ring containing a full set of matrix units (e.g., $x = \( \smallmatrix 0 & 1 \\ 0& 0 \\ \endsmallmatrix\)$ and $y = \( \smallmatrix 0& 0 \\ 1 & 0 \\ \endsmallmatrix \)$, and similar examples for larger sets of matrix units). \qed

In particular, we must avoid the temptation to use the transpose. We avoid it  by sometimes considering matrices as left $R$-module homomorphisms, acting on the right. (An alternative is to use the opposite ring of $R$; however, I found this extremely confusing.) We use the weird but logical notation, $\lArrow Q_2 ; Q_1. x$, to denote the left module homomorphism $x$ from $Q_1$ to $Q_2$, the subscripted $R$s on the left of the names of the modules emphasizing the fact that they are left modules. 

First, we show that $F(n,k) \subset C(n,k)$. Pick $M \in F(n,k)$; there exists a permutation $P \in \Cal P_n$ \st the bottom $k$ rows of $PM$ constitutes an invertible matrix, $g \in \gl(k,R)$. Then $PMg^{-1} = \( \smallmatrix X \\ I_k \\ \endsmallmatrix \)$. Set $W_0 = \( \smallmatrix \I_{n-k}\\ 0 \\ \endsmallmatrix \)$ (where the zero matrix is size $k \times (n-k)$). Then $h: = \( W_0 \ \ P Mg^{-1}\)$ is invertible (by the lemma about upper triangular matrices iff), and set $W = P^{-1}W_0$, and $U = \( W \ \ M \)$. Then $PU(\I_{n-k} \oplus g^{-1}) = h$, so $U = P^{-1} h (\I_{n-k} \oplus g)$ is a product of invertibles, hence is invertible. Thus $M \in C(n,k)$. The same argument ({\it not\/} using the transpose) works to show $F'(n,k) \subset C'(n,k)$.

To construct the map on equivalence classes, pick $M \in C(n,k)$, and view $M$ as a homomorphism of left modules, $\lArrow R^{1 \times k}; R^{1\times n}. M $ (it acts by right multiplication, of course). Denote its kernel, $\lker M$. By \nappthr(a implies c), the image of $M$ is all of ${}_R R^{1 \times k}$, so the map splits; in particular, $\lker M$ is a direct summand, and there exists a submodule $Q$ \st $R^{1 \times n} = \lker M \oplus Q $, and the restriction to $Q$ of right multiplication by $M$ is an isomorphism $ R^{1 \times k} \leftarrow Q$. In particular, $Q$ is free on $k$ generators as a left module. Now we make the assumption,
\item{(*)} $\lker M$ is free (as a left $R$-module) on $n-k$ generators.

(We will see later that this applies if either $M \in F(n,k)$ or $R$ satisfies SFF and stable finiteness.) Pick a basis with $n-k$ generators, $\brcs{r_1, r_2, \dots, r_{n-k}} \subset R^{1 \times n}$. Let $N$ be the $(n-k) \times n$ matrix whose $i$th row is $r_i$. We will show that the assignment $M \to N$ (which is highly dependent on choices of bases) yields a well-defined map $[M] \mapsto [N]$ on the orbit spaces.

First, we claim that $N \in C'(n,n-k)$. But this is easy: let $W$ be the $k \times n$ matrix consisting of a $k$-element basis for $Q$. Then the rows of $V:= \( \smallmatrix W \\ N \\ \endsmallmatrix \)$ constitute an $n$ element basis for ${}_R R^{1 \times n}$, and it is immediate that $V$ is invertible.

Next, if we choose a different $n-k$-element basis for $\lker M$ and corresponding $N'$, there exists $g \in \gl(n-k,R)$ \st $gN = N'$ (express each element of one basis as a {\it left\/} linear combination of the the elements of other basis; the matrices of coefficients are mutually inverse).

If we replace $M$ by $wMh$ where $w \in W(n)$ and $h \in \gl(k)$, we observe that $\lker (wMh) = \lker (wM) =(\lker M) w^{-1}$. Thus if we choose a basis for $\lker wM$, $(s_i)$, and form the matrix $N$ (whose rows are the $s_i$, then $N := N w$ will have rows constituting a basis for $\lker M$ (recall that $W(n)$ acts on the right of rows, multiplying by a weighted permutation matrix). Thus if $M$ belongs to the $W(n) \times \gl (k)$-orbit of $M_0$, then any choice for $N$ (that is, whose rows form a basis) will belong to the $\gl(n-k) \times W(n)$ orbit of (any choice of) $N_0$, and $N \in C'(n,n-k)$.

All this was under the assumption (*). Now assume that $M \in F(n,k)$; we claim that (*) holds. By assumption, there exists $P \in \Cal P_n$ and $g \in \gl(k,R)$ \st $M_0:= PMg^{-1} = \(\smallmatrix X \\ I_k \\ \endsmallmatrix\)$. It is an easy computation (essentially column-reduced echelon form, but with noncommuting entries) that $\lker M_0$ is spanned as a left $R$-module by the rows,
$$
\brcs{\(-e_j;X_{j,n-k+ 1}, X_{j, n-k + 2},\dots, X_{j,n} \)}_{j=1}^{n-k}\qquad\text{$e_j$ is the $j$th standard basis element of ${}_R R^{1 \times (n-k)}$},
$$
and the module is clearly free on these generators, and a direct summand (with complementary basis $\brcs{E_s}_{s=n-k+1}^n$ where $\brcs{E_i}$ is the standard basis for ${}_R R^{1 \times n}$. In particular, (*) holds.

But we have more: $N \in F'(n,n-k)$: if we choose the displayed basis, then $N = \( \smallmatrix -\I_{n-k} \\ X \\ \endsmallmatrix\) $, which clearly belongs to $F'(n,n-k)$.

Thus the assignment, $\phi: [M] \to [N]$, sends $\Cal F(n,k) \to \Cal F'(n-k)$. If $M \in F_{i} (n,k)$, for each of the $i$ sets of $k$ rows yielding an invertible matrix, there exists a distinct permutation matrix $P$ \st $PM$
has that particular set of $k$ rows moved to the bottom (of course different permutations can yield the same subset of $i$ rows), and we easily see from the preceding explicit form that for each, there is a corresponding subset of $n-k$ columns in $N$. This yields a one to one map from the sets of $k$ rows of $M$ constituting an invertible matrix to the sets of $n-k$ columns of $N$ yielding an invertible matrix. To show it is a bijection, we work in reverse.

Begin with $N \in C'(n,n-k)$, and consider the map of right modules, $\Arrow N; R^{n\times 1}. R^{(n-k) \times 1}$. By \nappfou(a), it is onto, so the map splits, and we a direct sum decomposition of right modules, $\ker N \oplus Q = R^{1\times n}$, with the restriction of $N$ to $Q$ being an isomorphism with $R^{1\times (n-k)}$. If $\ker N$ is free on $k$ generators, then we send $N$ to any matrix whose set of $k$ columns is a basis for $\ker N$ (in parallel with what we did before). Now suppose $N \in F(n,n-k)$; by applying a permutation matrix and pre-multiplying it by an element of $\gl(n-k)$, we may assume $hNP = \( \smallmatrix -\I_{n-k} \\ X \\ \endsmallmatrix\) $. Then $\ker hNP = \ker NP = P^{-1} \ker N$, and we see that the kernel is spanned by the rows of $ \(\smallmatrix X \\ I_k \\ \endsmallmatrix\)$. It follows that if $N$ came from an $M \in F(n,k)$, then we recover $M$ up to the action of $W(n) \times \gl(k)$.

In particular, this yields an inverse map $[N] \mapsto [M]$, and we obtain that $\phi$ is a bijection $\Cal F(n,k) \to \Cal F'(n,n-k)$; moreover, parallel arguments show that the inverse map sends $\Cal F'_i(n,n-k) \to \cup_{j\geq i} \Cal F_j(n,k)$. Combined with the previous, we must have $\phi(\Cal F_i(n,k)) = \Cal F'_i(n,n-k)$.

For the action of $W(n) \times \gl(k,R)$ on $F(n,k)$ or $C(n,k)$, we define the {\it
stabilizer\/} of a point, $M$, to be the subgroup of $\Cal P_n$ consisting of
$$\Set{P \in \Cal P_n}{PDMg = M \text{ for some $D \in \Cal D_n$ and $g \in \gl
(k)$}}
$$
 (and similarly \wrt the actions on $F'(n,k)$ and $C'(n,k)$). It is easy to check that this is a subgroup of $\Cal P_n$ (which we typically identify with $S_n$), with the usual properties of stabilizers, e.g.,
if $M$ and $M'$ are in the same $W(n) \times \gl (k)$ orbit, then the stabilizer
of $M$ is isomorphic to that of $M'$ via an inner automorphism of $\Cal P_n$ (of
course, the word {\it inner\/} is only significant if $n=6$).

Thus far, we have most of the following.

\Lem Proposition \nappsix. Let $R$ be any ring and $ n> k \geq 1$. Then the map $\phi$ induces bijections $\Cal F_i(n,k)) \to \Cal F'_i(n,n-k)$ for each of $i = 1, 2, \dots, {n \choose k}$. Moreover, it induces   isomorphisms on the stabilizers. 

\Pf We have already seen that $\phi$ is a bijection. Select $M$ in $F(n,k)$,
and construct  $N \in F^{n\times(n-k)}$ \st $M^T N $ is zero, as in the definition of
$\phi$. Suppose $P $ belongs to the stabilizer of $M$; then $PDMg = M$ for
$(D,g) \in \Cal D \times \gl(k)$. Thus $M^T( D^T)^{-1} P N =  0$. Since the
columns of $N$ constitute a basis for $\ker M^T$, it follows there exists $h \in
\gl (n-k)$ \st $(D^T)^{-1}PN h = N$. Since we can write $(D^T)^{-1}P = P(P^{-1}
(D^T)^{-1} P)$ and the second factor is diagonal, we have that $P $ belongs to
the stabilizer of $N$. It follows by interchanging $M$ and $N$, that their
stabilizers are equal.\qed

It is clear from the definitions that an $n \times k$ matrix of the form $U(X) := \(\smallmatrix X\\ \I_k \\ \endsmallmatrix \)$ belongs to $F(n,k)$, and moreover, every orbit in $F(n,k)$ contains a matrix of this form. More interestingly, such a matrix $U(X)$ belongs to $F_m (n,k)$, where $m = {n\choose k}$, if and only if all  $s \times s$ submatrices of $X$ are invertible, for all $s = 1, 2, \dots,k$. In particular, all the entries of $X$ have to be invertible (corresponding to $s = 1$), and if $k > 1$, then all matrices of the form
$$
\(\matrix x_{i(1),j(1)} &  x_{i(1),j(2)} \\
x_{i(2),j(1)} &  x_{i(2),j(2)} \\
\endmatrix \)
$$
($s= 2$) have to be invertible as well. There are thus $\sum_{j=1}^k {n-k \choose j}{k \choose j} = {n \choose j} -1$ conditions.

\Lem Proposition \nappsev. Let $R$ be a stably finite ring satisfying SFF. Then the map $\phi$ induces a bijection $\Cal C(n,k) \to \Cal C'(n,n-k)$.

\Pf First, we show that (*) holds; that is, if $M \in C (n,k)$, then $\lker M$ is free on $n-k$ generators. As in the proof above, we have that ${}_R R^{1\times n} = \lker M \oplus Q$ where $Q$ is a free left $R$-module (because the restriction of right multiplication by $M$ to $Q$ is an isomorphism with ${}_R R^{1\times k}$). By SFF, $\lker M$ is free, and by stable finiteness, it can only be free on $n-k$ generators. Hence (*) holds.

Thus the corresponding $N$ exists, and we showed above that $N \in C'(n,n-k)$. In particular, $\phi$ is a well-defined map $\Cal C(n,k) \to \Cal C'(n,n-k)$. Now we can work in reverse to show it is a bijection. Pick $Z \in C'(n,n-k)$; write $\ker Z \oplus Q = R^{n\times 1}$ as right $R$-modules, with $Z$ restricted to $Q$ being an isomorphism. Then $Q$ is free, so SFF implies $\ker N$ is free, and stable finiteness yields freeness on exactly $k$ generators; pick such a basis. Define $M$ to be the matrix whose $i$th column is the $i$th basis element. It is easy to check that if we construct the corresponding $N$, the $R$-module span of its rows will be that of $N$, and both being bases, they are bases for the same submodule, hence there exists $g \in \gl(n-k)$ \st $gN = N$. This shows that the map $\Cal C' (n,n-k) \to \Cal C (n,k)$ is the inverse to $\phi$.
\qed

Suppose that $R$ is commutative, or more generally, admits an anti-automorphism $\psi$ (for commutative rings, we can take the identity; if $R$ is a *-ring, we can take *). Then there is a (relatively) natural map (depending on $\psi$) $\Cal C'(n,n-k) \to \Cal C(n,n-k)$ (and corresponding maps on the $F$s): send $M$ to $\psi(M)^{T}$ (defining $\psi$ on matrices entrywise). So composing this with $\phi$, we obtain bijection $\Cal F_i(n,k) \to \Cal F_i (n,n-k)$ (and if $R$ satisfies SFF and is stably finite, on the corresponding $\Cal C$s).

In particular, if $R$ is commutative, and we take the identity as $\psi$, the columns of $N$ are given by a basis for $\ker M^T$, where we regard $\Arrow M^T; R^{n\times 1}. R^{k \times 1}$ as a homomorphism of right modules. Then $M^T N = 0$ and thus $N^T M = 0$, and we quickly see that the map $\phi$ really does behave like a duality. Similarly, if we assume that $\psi$ is involutive (as in *-rings), then there is a natural duality $\Cal F(n,k) \to \Cal F'(n,n-k) \to \Cal F(n,k)$ implemented via $\ker (M^*)^T$. 

 If, however, $\psi$ is just an anti-automorphism, there are still mutually inverse  bijections $\Cal F(n,k) \to \Cal F'(n,n-k)$ and vice versa, but they are not really the same (the first is implemented by $\ker \psi(M)^T$, the reverse by $\ker \psi^{-1}(M^T)$). 

Suppose that  $k = 2$, $n=5$, and $F = \Z_p$ for a prime $p$. We
note that the orbit space of $F_{10}(5,2)$ can be interpreted as the orbit space of the collection of $5$ (distinct)-element subsets of $1$-dimensional projective space over $\Z_p$, $\PP_1$, acted upon by P$\gl(2,\Z_p)$. Since
P$\gl(2,\Z_p)$ acts $3$-transitively on $\PP_1$, it is easy to check that it acts
transitively on
$F_{10}(5,2)$ if $p = 5$ or $7$ (but not for any larger prime---the order of $F_{10}(5,2)$
does not divide that of P$\gl (2,\Z_p)$ if $p > 7$). In particular, the result
above shows that P$\gl(3,\Z_7)$ acts transitively on $F_{10}(5,3)$---something that
is routine to check with a computer (it suffices to show that the stabilizer of
one or any point has at most six elements), but extremely tedious to check by
hand.

\SecT Appendix B. A truncated reciprocal formula

\noindent {\it David Handelman \& Damien Roy}\vskip 7pt

\noindent Fix a prime $p$. The following goes back to 1893.

\Lem Theorem \appone. [L] The number of rank $n-s$ matrices in $\gl(n,\Z_p)$ is
$$
C_s:= \frac{\((p^n - 1)\dots (p^n - p^{n-s -1})\)^2}{(p^{n-s} - 1)\cdot (p^{n-s}-p)\dots (p^{n-s}-p^{n-s-1})}.
$$

Now we can prove the result of this section. The limiting case of this is the identity   [HW, 19.7], due to Euler. However, we cannot obtain the result below simply by truncation, since there is a bonus of an extra bit in the exponent of the error term.

\Lem Proposition \apptwo. Let $n,s$ be positive integers, with $n > (s+1)^2+ 1$, and let $z$ be a variable. Then as functions analytic on the open unit disk, we have
$$\eqalign{
\(\prod_{i=1}^{(s+1)^2 - 1} (1-z^i)\)&\(1 + \sum_{j=1}^{s} \frac{z^{j^2}}{(1-z)^2 (1-z^2)^2 \cdots (1-z^{s})^2}\) \qquad \text{and}\cr
\(\prod_{i=1}^{(s+1)^2 - 1} (1-z^i)\)&\(1 + \sum_{j=1}^{s} \frac{z^{j^2}(1-z^n)(1-z^{n-1})\dots (1-z^{n-j+1})}{(1-z)^2 (1-z^2)^2 \dots (1-z^{s})^2}\)
}$$
are polynomials, and their Maclaurin expansions are
$$1 - z^{(s+1)^2 + 2} + \text{ higher order terms}.$$

\Pf Since $(s+1)^2 - 1 \geq 2s$,  it follows that all the denominators of the right factor are eliminated by the left (count the multiplicities of the various  roots of unity that are zeros of the denominators, and do the same for the first $2s$ terms in the product on the left). Hence the polynomial assertion is verified.

With $N_p = \prod_{i=0}^{n-1}(p^n - p^i)$ being the number of invertible matrices, we have,
$$\eqalign{
C_s &= \frac{N_p}{\prod_{i=0}^{n-1}(p^n - p^i)} \cdot  \frac{\((p^n - 1)\dots (p^n - p^{n-s -1})\)^2}{(p^{n-s} - 1)\cdot (p^{n-s}-p)\dots (p^{n-s}-p^{n-s-1})}\cr
& = N_p \cdot  \frac{(p^n - 1)\dots (p^n - p^{n-s -1})}{\((p^n-p^{n-s})(p^n - p^{n-s+1})\dots (p^n- p^{n-1})\)\cdot \((p^{n-s} - 1)  (p^{n-s}-p)\dots (p^{n-s}-p^{n-s-1})\)}\cr
& = p^{n^2}\cdot \frac{N_p}{p^{n^2}}\cdot  \frac{p^{(n-s)(n-s-1)/2}(p^{n} -1)(p^{n-1} -1) \dots (p^{s+1} - 1)}{p^{n(n-1)/2}\((p^s - 1)(p^{s-1}-1 ) \dots (p-1) \)\cdot \( (p^{n-s} -1) (p^{n-s-1} - 1)\dots (p-1)\)} \cr
& = p^{n^2}\cdot \frac{N_p}{p^{n^2}}\cdot  \frac{p^{(n-s)(n-s-1)/2}(p^{n} -1)(p^{n-1} -1) \dots (p^{n-s+1} - 1)}{p^{n(n-1)/2}\((p^s - 1)(p^{s-1}-1 ) \dots (p-1) \)^2}; \text{ divide by $p^{n^2}$ and set $z = 1/p$;} \cr
\frac{C_s}{p^{n^2}}& = \frac{N_p}{p^{n^2}}\cdot\frac{z^{s^2}(1-z^n)(1-z^{n-1})\dots (1-z^{n-s+1})}{(1-z)^2(1-z^2)^2\dots (1-z^{s})^2}:= \(\prod_{i=1}^n (1-z^i)\)\cdot c_s (z).\cr
}$$

Set $c_0 = 1$, and let $m(z) = \prod_{i=1}^n (1-z^i)$. We see that each $c_s(z)$ is analytic on the unit disk; moreover, for each prime $p$, $m(1/p) \sum_{j=0}^n  c_j (1/p) = 1$, since the unnormalized forms count the total number of matrices; this equality is also true at $z = 0$. Since each of the factors is analytic on the open disk, and the product agrees with the constant function $1$ on a limit point ($\{ 0,1/2,1/3, \dots \}$), it follows that the product, $m \cdot (\sum_{j=0}^n c_j)$ is $1$ on the unit disk. We  use this to determine some Maclaurin coefficients.  Each $c_i$ is expressed as
$$
\frac {z^{i^2}}{(1-z)^2 \dots (1-z^i)^2} \times (1-z^n)(1-z^{n-1})\dots (1-z^{n-s+1}).
$$
When we expand this in its Maclaurin expansion, we see that $c_i  = z^{i^2} +2 z^{i^2 + 1} + {\text{terms of higher order}}$. Now suppose that $s^2 \leq n$, and consider the truncated sum, $\sum_{i=0}^{s} c_i$. The missing terms are of the form $m_p \cdot c_t$ where $t > s$. It follows immediately that the smallest degree term in the Maclaurin expansion of what is missing is $z^{(s+1)^2} + 2 z^{(s+1)^2 + 1}$. Thus $E_s:= \sum_{i=0}^n c_i - \sum_{i> s} c_i= 1- z^{(s+1)^2} - 2 z^{(s+1)^2 + 1} + \text{terms of higher order}$.

Now truncate $m$ at $(s+1)^2-1$, that is, $m_s = \prod_{i\leq (s+1)^2-1} (1-z^i)$. Then
$$\eqalign{
m - m_s &= m_s\cdot  ((1-z^{(s+1)^2})(1-z^{(s+1)^2  +1})\dots -1) \cr
&= m_s\cdot (-z^{(s+1)^2 } - z^{(s+1)^2 + 1} - z^{(s+1)^2+2} + \dots)\cr
& = -z^{(s+1)^2 }(1 + z + z^2 + \dots)\cdot (1-z) (1-z^2)(1-z^3)\dots\cr
& = -z^{(s+1)^2 }(1- z^2 + \dots); \qquad\text{ so }\cr
m_s & = m + z^{(s+1)^2 }(1- z^2 + \dots).\cr
}$$

Now we have
$$\eqalign{
m_s\cdot E_s  & = m + z^{(s+1)^2 }(1- z^2 + \dots) \cdot \(\sum_{i=0}^n c_i - \sum_{i> s} c_i\)\cr
& = 1 - m \sum_{i> s} c_i + \(z^{(s+1)^2}(1-z^2 + \dots)\)E_s\cr
&= 1 -\((1-z) (1-z^2) \dots \cdot (c_{s+1} + \dots) \)+ \(z^{(s+1)^2}(1-z^2 + \dots)\)(1+ z + 2z^2 + \dots) \cr
&= 1 - (1-z -z^2 + z^3 + \dots)z^{(s+1)^2}(1 + 2z + 5z^2 + \dots) + z^{(s+1)^2}(1+z +z^2 + \dots)\cr
& = 1 - z^{(s+1)^2}\( (1 + z + 2z^2 + \dots) - (1 + z + z^2 + \dots) \)\cr
& = 1 - z^{(s+1)^2 + 2} + \dots.
}$$

This is exactly the desired assertion for the more complicated product. For the less complicated (first) product, from $n-j+1 + j^2 > (s+1)^2 + 2$ (this is equivalent to $n + 1 > (s+1)^2 + 2$), the extra terms in the numerator of the right hand term do not contribute to any Maclaurin series terms of degree less than or equal $(s+2)^2 + 2$, so the first product has the same Maclaurin expansion up to  that degree.
\qed

 The simpler expression (the first one) does not involve $n$ and  product behaves as $1 - z^{(s+1)^2 + 2} (1 + \Oh z)$ without reference to $n$. If we let $s\to \infty$, then the left function converges uniformly on compact subsets of the open unit disk to the Euler function, and since the latter has no zeros, it follows that the infinite sum on the right also converges uniformly on compact subsets, so is also analytic on the open disk; necessarily, the limit is the reciprocal of the Euler function, giving yet another proof of the identity [HW, 19.7]. For all values of  $s$ that we could calculate with, the coefficients of the higher order terms oscillate in a particularly interesting way, and the maximum increases as $s$ does, according to {\it Maple.}

 \long\def\Rf[#1] #2, #3. #4\par%
{\vskip 2pt \itemitem{[#1]} #2, {\it #3,} #4\par\vskip2pt}

\SecT Appendix C. Counting PH-equivalence classes

In [ALTPP], the authors compiled tables of PH-equivalence isomorphism
types, based on (what amounts to) $d = |\det B|$ and $|\det B\op |$ for
$n=3$ and $4$. Using Proposition \twoone, one can obtain explicit
formulas for the the numbers of equivalence classes that contain a
terminal form with $1$-block size $n-1$ of determinant $d$, and subdivide it
according to the possible values of $|\det B\op|$. Aside from the
complicated nature of the expressions, these only deal with $1$-block size
$n-1$.

In this appendix, we will see that the lower bound obtained for the number of PH-classes of $C \in \NS_n$ of determinant $d$  obtained in Lemma \appcone,
$$
F_n (d): = \frac{\phi* J_2 * \dots* J_{n-1}}{n!}
$$
is asymptotically (in $d$) correct (with an estimate of a factor $1+d^{-1}$), when $d$ is square-free. We do this by showing that the vast majority of the $S_n$-orbits on \quotes{weakly terminal} matrices (defined below) of determinant $d$ are of full size, that is, $n!$, via estimates (and in some cases, exact formulas) for numbers of matrices fixed by an arbitrary permutation.

A matrix $C$ is called {\it weakly terminal,} if it is in Hermite normal form and belongs to $\NS$; in particular, it is upper triangular, and its $(1,1)$ entry is $1$.  It is quite easy to count the weakly terminal matrices of given size and determinant.

Let $C$ be a weakly terminal matrix of size $n$, let $\pi \in S_n$ be a permutation, and let $P \equiv P_{\pi}$ be the permutation matrix right multiplication by which implements $\pi$ as a permutation of the set of columns. There exists $U \equiv U_P \in \gl(n,\Z)$ \st $C_P:=  UCP$ is in Hermite normal form, and in fact, given $P$, $C_P$ is unique. If $C'$ and $C''$ are weakly  terminal matrices \st $C'$ \st $C' = U C P$ for some permutation matrix $P$ and $U \in \gl(n,\Z)$,  and $C'' = U'' CP$ (same permutation matrix), then   $C'' = U'' (U^{-1} C' P^{-1})P = U'' U^{-1}C'$, so that $C''$ is Hermite equivalent to $C'$---but both are in Hermite normal form, so must be equal.

Since the property of being in $\NS$ is preserved by Hermite equivalence, it follows that $\brcs{C_P}_{P \in S_n} $ is a finite set consisting of weakly terminal elements, and an orbit, under the action of $S_n$. Moreover,  this orbit must contain at least one terminal matrix (since every matrix is PH-equivalent to a terminal matrix, and terminal implies weakly terminal). Thus the orbits of the form $\brcs{C_P}_{P \in S_n}$ (with $C$ varying over weakly terminal matrices) are in bijection with the PH-equivalence classes of $ C \in \NS_n$.

In particular, for fixed weakly terminal $C$ (weakly terminal is required, since otherwise $C_P$ is not necessarily uniquely determined) is an $S_n$-space. The difficulty in counting arguments is that the orbits need not all be full, that is, there will be some fixed points---$C_P = C$ for some nontrivial permutation matrix $P$.

One case in which the orbit will be full (of cardinality $n!$) occurs when  $ J(C_{\Omega(i)})$ are distinct. The obvious action of $P$ (acting on the columns) implements a permutation of the $n$-tuple, $(J(\Omega(1), J(\Omega(2)), \dots, J(\Omega(n)))$ (the subsequent left action by the unimodular matrix does not affect the order of these groups). If  $J(C_{\Omega(i)})$ are distinct, this action is just the permutation representation of $S_n$ on a set with $n$ elements. It follows that the orbit of the action $C \mapsto C_P$ is full.

For $\pi \in S_n$, let $P \equiv P_{\pi}$ denote the permutation matrix right multiplication by which implements $\pi$ as a column permutation. Then a weakly terminal matrix $C \in \NS_n$ is fixed by $\pi$ (or $P \equiv P_{\pi}$) iff $CPC^{-1}$ has only integer entries. For a subgroup $H$ of $S_n$ and a positive integer $d$, Let $\Cal Z (H)(d)$ denote the set of all weakly terminal matrices of determinant $d$ that are fixed by all elements of $H$. When $H$ is the cyclic group generated by $\pi$, we use the notation $\Cal Z (\pi)(d)$. The cardinality of $\Cal Z(\pi)(d)$ is denoted $\Cal S(\pi)(d)$ (and similarly for subgroups $H$). The function $d \mapsto \Cal S(\pi)(d)$ is multiplicative (in the number-theoretic sense) for all $\pi$.

There is an obvious procedure for counting  PH-equivalence classes. First, we count all the weakly terminal matrices of fixed determinant $d$ (done in Lemma \appcone). Then we count the number of weakly terminal matrices whose orbits are not full, and subtract them off, keeping track of the number of PH-equivalence classes they constitute, and apply Burnside's lemma. When $n =3$, it is barely possible to do this, but for larger sizes, obtaining the exact number seems horrible. (In fact, when $n=3$, we obtain a convenient subdivision into various cases with $1$-block size two, and the remainder; this goes most smoothly when $d$ is square-free.)

However, for $n$ arbitrary and $d$  prime (and thus for $d$ square-free), we can obtain relatively explicit formulas for $\Cal S(\pi)(d)$ for every $\pi \in S_n$; since the matrices with orbit size less than $n!$ must be in $\Cal Z(\pi)$ for some non-identity $\pi \in S_n$, we can easily obtain an upper bound for the number of PH-equivalence classes. This will verify the conjecture below when $d$ is square-free.

Recall the definition of the $k$th Jordan totient, $J_k (n) = n^k\prod_{p|n} (1- p^{-k})$. Then $J_1 = \phi$, $J_k$ is multiplicative (in the number-theoretic sense), and $J_k(d)$ counts the number of content one columns of size $k+1$ with  $d$ in the bottom entry, and all the other entries belonging to $\brcs{0,1,\dots, d-1}$. Recall that for multiplicative functions $f$ and $g$, $f*g$, the convolution, is defined by $(f*g)(t)= \sum_{x|d} f(x) g(d/x)$,  and is multiplicative; moreover, $f*g = g*f$.  There is an identity for constructing $J_k$, namely if $\xi_k$ is the multiplicative function $n \mapsto n^k \phi(n)$, then $\xi_{k-1} * \xi_{k-2}* \dots* \xi_1 * \phi = J_k$. (The Dirichlet series for the function on the left telescopes.)

\Lem Lemma \appcone.  Let $d$ be a positive integer. For $n > 1$, the number of weakly terminal $n \times n$ matrices of determinant $d$ is
$$
\(\phi * J_2 * \cdots * J_{n-1}\)(d).
$$

\Pf This simple proof is by induction on $n$. If $n = 2$, we are counting the matrices $\( \smallmatrix 1 & a \\ 0 & d \\ \endsmallmatrix\)$ where $0 \leq a < d$ and $(a,d) = 1$---the number of choices for $a$ is obviously $\phi(d)$.

For $n > 2$, given a weakly terminal matrix $C$ say with $(n,n)$ entry $x$ (which divides $d$, as the matrix is upper triangular), deleting the last row and column, yields a weakly terminal matrix of size $n-1$, and with determinant $d/x$; moreover the $n$th column has content one. Conversely, given a weakly terminal matrix of size $n-1$ and a content one column of size $n$, we created a weakly terminal matrix of size $n$ by attaching the column, and embroidering $n-1$ zeros on the bottom, and of course the determinant multiples.

If $H_j(t)$ denotes the number of weakly terminal matrices of size $j$ and determinant $t$, we thus have $H_{n-1}(t) = \(\phi * J_2 * \cdots * J_{n-2}\)(t) $   by the induction hypothesis, and
$$\eqalign{
H_n (d)  &= \sum_{x|d} J_{n-1} (x) H_{n-1}(d/x) \cr
& = (J_{n-1} * H(n-1)) (d) = (H_{n-1}*J_{n-1} ) (d) \cr
& = \(\phi * J_2 * \cdots * J_{n-1}\)(d), \cr
}$$
completing the induction.
\qed

Set $F(n,d) =\(\phi * J_2 * \cdots * J_{n-1}\)(d)$. It follows immediately that $F(n,d)/n!$ is a lower bound for the number of PH-equivalence classes of matrices in $\NS_n$ of determinant $\pm d$.

\Lem Conjecture. For $d$ a positive integer, the number of PH-equivalence classes  of matrices in $\NS_n$ having determinant $\pm d$ is
$$
\frac{(\phi*J_2 * \dots * J_{n-1})(d)}{n!} \cdot \( 1 + {\frac {{n \choose 2}}d}(1+\oh{1})\).
$$

One way to proceed, and even obtain a slightly sharper result is as follows. Fix $n$, then  $d$, and a permutation $\pi \in S_n$, and its corresponding matrix $P$ (right multiplication by which implements $\pi$ as a column permutation. We wish to obtain an asymptotic estimate  for the number of weakly terminal $n \times n$ matrices $C$ of determinant $d$ invariant under the  action of $P$, that is, $CPC^{-1} \in \Mn n \Z$.

Let $K(\pi)$ denote the number of cycles in the decomposition of the permutation $\pi$ associated to $P$; fixed points of course are $1$-cycles, so are counted. Then $K(\pi)$ is just the co-rank of the matrix $\I - P$, that is, $n  = \rk (\I - P) + K(\pi)$, as it simply counts the algebraic and geometric multiplicities (they are the same for permutation matrices) of $1$ as an eigenvalue of $P$.

Motivated by the  counting arguments above, the following is likely to be true.

\Lem Specific Conjecture. Let $\pi \in S_n$ be a non-transposition. Then  for all $\epsilon > 0$,
$$
\frac{\Cal S(\pi)(d)}{F(n,d)} = \oh{d^{K(\pi) -n + \epsilon}}.
$$

Without   $\epsilon$, the specific conjecture fails (in
Appendix C, we will see that when $n =3$ and $\pi$ is a $3$-cycle, then $\Cal S(\pi)(d)/\Cal S(\id)(d)$ is infinitely greater than $d^{-2}$).

If the specific conjecture were  true,  the conjecture preceding it would follow (as we will see when we discuss $\Cal S(\pi)$ when $\pi$ is a transposition). Of course, it would be sufficient to prove this when $d$ is restricted to powers of primes.

\comment
What makes the specific conjecture  plausible is an elementary observation. The weakly terminal matrix $C$ is $P$-invariant iff $CPC^{-1}$ has only integer entries, which in turn is equivalent to $C(\I -P)C^{-1} = \I - CPC^{-1}$ having only integer entries. The matrix $C(\I - P)$ has rank $n - K(P)$, and we expect that the smaller the rank, the fewer conditions will be imposed, hence the greater the number of solutions, but that the conditions imposed will reduce the exponent of $d$ appearing in the number of solutions.

\endcomment

Presumably, this is part of a theory of an arithmetic version of varieties, corresponding to subvarieties having measure zero when imbedded in a variety.

We will show that the original conjecture is correct when    limited  to square-free $d$.

 \Lem Lemma \appctwo. Let $H $ be a subgroup of $S_n$, and $\pi \in S_n$. Then for all $d> 0$, 
$$
\Cal S (H)(d) = \Cal S (\pi H \pi^{-1})(d).
$$

\Pf Fix $d$, and let $Q$ be the permutation matrix representing $\pi$. Select a permutation matrix $P$ that corresponds to an element of $H$, and suppose weakly terminal  $C$ is fixed under the action of $P$, that is, $CPC^{-1} \in \Mn n \Z$. There exists $U\equiv U_{C,Q} \in \gl(n,\Z)$ \st $UCQ^{-1}$ is in Hermite normal form; since both left multiplication by elements of $\gl(n,\Z)$ and right multiplication by permutation matrices preserve $\NS_n$, $UCQ^{-1}$ is itself weakly terminal, and  of the same determinant as $C$ (it is of the same absolute determinant, but being in Hermite normal form, the determinant is positive). In addition, $U_{C,Q}$ is unique (\wrt the property that $UCQ^{-1}$ is in Hermite normal form), as $\det C = d \neq 0$.

Next, we observe that
$$
UCQ^{-1}(QPQ^{-1})Q C^{-1}U^{-1} = U(CPC^{-1})U^{-1} \in \Mn n \Z.
$$
Since this is true for every $P$ corresponding to an element of $H$, we have a set map $\Cal Z(H) (d) \to \Cal Z(\pi H\pi^{-1}) (d) $ given by $C \mapsto U_{C,Q} C Q^{-1}$. Since $C$ is itself weakly terminal, it follows  that $U_{UCQ^{-1}, Q^{-1}}$ must be $U_{C,                                 Q}^{-1}$, so the corresponding map $ \Cal Z (\pi H \pi^{-1})(d)) \to \Cal Z (H)(d)$ is the inverse of the original.
This shows that $C \mapsto U_{C,Q} C Q^{-1}$ is a bijection.
\qed

Define for each positive integer $k$, $\Arrow \Cal N_k; \N. \N$ via
$$
\Cal N_k (d) = \left|\Set{z \in \Z_d}{z^k =1} \right|.
$$
The Chinese remainder theorem implies that for each $k$, the function $\Cal N_k$ is multiplicative.
The following is routine, and follows from $\Z_{p^m}^*$ being cyclic of order $p^{m-1}(p-1)$ when $p$ is odd, and isomorphic to $\Z_{2^{m-2}} \times \Z_2$ when $p =2$ (with the interesting convention that $\Z_{2^{-1}} \times \Z_2$ is the trivial group).

\Lem Lemma \appcthr. For a prime $p$,
$$
\Cal N_k (p^m) = p^{\min \brcs{v_p(k),m-1}} \cdot \cases \gcd\brcs{p-1,k} & \text{if $p$ is odd, or $p^m=2$}\\
2 & \text{if $p=2$ and $m \geq 2$. }
\endcases
$$

If $k$ is an odd prime, then $\ln \Cal N_k (d) \leq |\Set{p|d}{p \equiv 1 \mod k}| \cdot \ln k$, and in general, $\Cal N_k (d) \leq 2k^{w(d)+1}$, although the latter is almost always a gross overestimate.

Let $\pi \in S_n$, and let $i \in \brcs{1,2,\dots, n}$ ($n$ will be fixed). Define the {\it orbit of $i$ \wrt $\pi$,} $\Cal O_{\pi}(i)$ (or $\Cal O (i)$ if $\pi$ is understood), to be $\brcs{\pi^k(i)}_{k \in \Z}$.

For $\pi \in S_n$, let   $P_{\pi}$ (or $P$ if $\pi$ is understood)  be the corresponding permutation matrix that implements the action of $\pi$ on the columns of $n \times n$ matrices by right multiplication. If $C$ is weakly terminal, then $C$ is fixed by the action of $\pi$ (or $P$) if $CPC^{-1}$ has only integer entries: explicitly, $CP^{-1}$ is put in Hermite normal form by a matrix $U \in \gl(n,\Z)$, that is, $UCP^{-1}$ is in Hermite normal form; then $UCP^{-1} = C$ iff $U  = CPC^{-1}$, which is equivalent (since the determinant of the right side is plus or minus one) to $CPC^{-1}$ having only integer entries.

We will determine $\Cal S(\pi)(d)$ exactly,  when $d$ is square-free. We obtain a formula involving some of the orbits of $\pi$ and their cardinalities, relating to Jordan totients. It is explicit enough that we can easily verify the specific conjecture for square-free $d$.

As before fix $n$ and fix $d$ as well. Let $2 \leq j \leq n$, and let $u = (a_1, a_2, a_3, \dots, a_{j-1}, d, 0,0,\dots,0)^T \in \Z^{n}$ \st $0 \leq a_i < d$ and $\gcd \brcs{d, a_1, \dots, a_{j-1}} = 1$. Let $C \equiv C(j,u)$ be the  weakly terminal matrix whose
$j$th column is $u$, and whose $i$th column for $i \neq j$ is the standard basis elements $E_i \in \Z^n$. In other words, $C-\I$ has exactly one nonzero column, and it is $u - E_j$. Note that $C$ so constructed is automatically weakly terminal.

Define
$$
V_j \equiv V_j (d)= \Set{C \in \Mn n \Z}{\det C = d; \  \text{$C$ is weakly terminal; the only nonzero column of $C-\I$ is the $j$th}.}
$$
The definition forces the $j$th column to be of the form  $u$ as given above.

Given $\pi \in S_n$, we will determine the number of matrices $C \in V_j$ \st $C$ is invariant under $\pi$, that is, for which $CP_{\pi}C^{-1} \in \Mn n \Z$. If $d = p$, a prime, then every weakly terminal $C$ of determinant $d$ is in $V_j$ for some $j$, so we obtain $\Cal S(\pi)(p)$ as a sum over $j =2, \dots, n$ of these numbers. This  yields a formula for $\Cal S(\pi)(d)$ when $d$ is square-free. The formula is sufficiently explicit to determine asymptotic behaviour (that is, for large, square-free $d$).

Begin with $C \equiv C(j,u) \in V_j$, and $P \equiv P_{\pi}$. The $i$th column of $CP$ is given by
$$
(CP)_i = \cases E_{\pi^{-1}(i)} & \text{if $i \neq \pi(j)$}\\
u = \sum_{l \leq j-1} a_l E_l + dE_j & \text{if $i = \pi (j)$.}\\
\endcases  
$$
Thus the entries are given by
$$
(CP)_{i,m} = \cases 1 & \text{if $i \neq \pi(j)$ and $i = \pi(m)$}\\
0 & \text{if $i \neq \pi(j)$ and $i \neq \pi(m)$}\\
a_m & \text{if $i = \pi(j)$ and $m < j$}\\
d & \text{if $i = \pi(j)$ and $m = j$}\\
0 & \text{if $i = \pi(j) $ and $m > j$.}\\
\endcases
$$
Extend the definition of $a_i$, so that $a_j = d$ and $a_l = 0$ if $l > j$. We will usually write $\pi k$ rather than $\pi(k)$ (unless ambiguity may result) from now on. We can thus write the $m$th row of $CP$, $(CP)^{(m)}$ as
$$
(CP)^{(m)} = a_m e_{\pi j} + \cases e_{\pi m} & \text{if $m \neq j$} \\
0 & \text{if $m = j$.}
\endcases
$$
(We are using the convention that $E_i$ represent the standard basic columns, while $e_i$ represent the standard basic rows, so that $e_k E_l$ is the matrix product whose outcome is $\delta_{kl}$.)

Now $C^{-1}$ is calculated easily by factoring $C$ into a product of a diagonal matrix and a unipotent. The outcome is that all the columns of $C^{-1}$ except the $j$th are just the standard basic columns, and the $j$th column is $d^{-1}(-a_1, \dots, -a_{j-1}, 1, 0, 0, \dots, 0)^T$. In particular, $(C^{-1})_j = d^{-1}\(E_j - \sum_{i \leq j-1}a_i E_i \)$.

We see immediately that $CPC^{-1}$ has only integer entries iff its $j$th column does. We calculate the entries of the $j$th column.
$$\eqalign{
(CPC^{-1})_{m,j} & = (CP)^{(m)} (C^{-1})_j \cr
& = \cases
 \frac{-a_m a_{\pi j}}{d} + \cases -\frac{a_{\pi m}}d & \text{if $m, \pi m \neq j$}\\
0 & \text{if $m = j$}\\
\frac{1}d & \text{if $ m \neq j$ and $\pi m = j$.}\\
\endcases &
\text{if $\pi j \neq j$}\\
\frac{a_m}d + \cases  \frac{-a_{\pi m}}d & \text{if $m \neq j$} \\
0 & \text{if $m = j$.} \\
\endcases
 & \text{if $\pi j = j$.}\\ \endcases\cr
}$$

Now we count the $\pi$-invariant matrices in $V_j$. First suppose that $\pi (j) = j$. Then the conditions for all the entries to be integers boil down to $a_{\pi m} \equiv a_m \mod d$ for all $m \neq j$. Hence, if $m \neq j$,  then $a_m \neq 0$ implies $a_{i} \neq 0$ for all $i \in \Cal O(m)$. Thus $a_m \neq 0$ entails $\Cal O (m) \subseteq \brcs{1,2,\dots, j-1}$. Conversely, if $\Cal O(m) \subseteq \brcs{1,2, \dots, j-1} $, we can put any value in we like for $a_m$, and the same value for $a_i$ as $i$ varies over   the orbit of $m$. The only constraint is that the resulting column $u$ must have content one. The number of such columns is thus exactly $J_{s(j)} (d)$ (the Jordan totient) where $s(j)$ is the number of orbits that are contained in $\brcs{1,2, \dots, j-1}$.

So if $\pi(j) = j$, the number of matrices in $V_j (d)$ that are fixed by the action of $\pi$ is exactly $J_{s(j)}(d)$. (If $s(j) =1$, $J_1 = \phi$, the usual totient; if $s(j)=0$, the outcome is zero.)

Now suppose that $\pi (j) \neq j$. First, we consider conditions arising from the coefficients corresponding to $\Cal O (j)$, the orbit of $j$ itself. Suppose the orbit of $j$ has $k > 1$ elements, so that if $m = \pi j$, then $m, \pi m, \dots , \pi^{k-2} m$ are distinct from each other and $j$, and all but the last one satisfies $\pi s \neq j$ (if $k = 2$, then all we have is $\brcs{m}$).

For $k>3$, we deduce $a_m^2 \equiv - a_{\pi m}$, and then $a_{\pi m}a_m  \equiv - a_{\pi^2 m}$, until we reach $a_{\pi^{k-3}m}a_m  \equiv - a_{\pi^{k-2} m} $, and finally, $a_{\pi^{k-2}m}  \equiv -1$. We can rewrite these as functions of $a_m$, obtaining $a_{\pi m} \equiv  - a_m^2$, $a_{\pi^2 m} \equiv a_m^3$, and for $r \leq k-2$, $a_{\pi^{r}m} \equiv (-a_m)^{r+1}$, and finally $(-a_m)^k \equiv 1 \mod d$. So we have $\Cal N_k (d)$ choices for $a_m$, and every other $a_i$ for $i \in \Cal O(j)\setminus \brcs{j}$ is determined by the choice of $a_m$.

For $k=2$ and $k =3$, the same result applies (and is easily checked); the indexing was confusing.

Now we come to $a_m$ for $m \not\in \Cal O (j)$. Then the equations become $a_{\pi m} \equiv -a_m a_{\pi j}$, $a_{\pi^2 m} \equiv a_m a_{\pi j}^2$, and in general $a_{\pi^r m} = a_m (-a_{\pi j})^r$ (this is true for all $r$). Thus the choice of $a_{\pi j}$ (which has to be a $k$th root of unity in $\Z_d$) and the choice of $a_m$ will determine the rest of the $a_i$ for $i \in \Cal O (m)$. However, there are constraints on the choice of $a_m$ if $k(m) := |\Cal O (m)|$ is not divisible by $k = |\Cal O (j)|$. Write $k(m) = ck + f$ with $c \geq 0$ and $0 \leq f < k$. Then $a_m \equiv  a_{\pi^{k(m)}m}  \equiv a_m (-a_{\pi j})^{k(m)} \equiv a_{\pi^f m} \equiv a_m(-a_{\pi j})^f$. Hence $a_m (1 - (-a_{\pi j})^j) \equiv 0 \mod d$.

Set $z = -a_{\pi j}$, so that $z^k \equiv 1 \mod d$. The restriction, that  $a_m (1 - z^f) \equiv 0$, is trivial if $z^f \equiv 1$. At this point, for simplicity, we assume that $d$ is a prime. In that case, $1-z^f$ is a zero divisor iff $z^f \equiv 1$, and  $a_m$ can be anything; otherwise, $a_m = 0$ is forced. Moreover, if $z^f \equiv 1$, then the remaining $a_i$, determined by $a_{\pi^r m} = a_m (-a_{\pi j})^r$, are consistent with the conditions for invariance.  Hence there are exactly $\gcd\brcs{f,p-1} = \gcd \brcs{k(m),k,p-1}$ selections for $a_{\pi j}$ for which we obtain $p$ choices for $a_m$, and for all the rest ($\Cal N_k (p) - \gcd \brcs{k(m),k,p-1}$), there is exactly one choice, $a_i = 0$ for all $i$ in the orbit of $m$. If $k$ divides $k(m)$, then the latter does not occur (as $f = 0$).

Now we can count the number of matrices in $V_j (p)$ fixed by $\pi$,  for $p$ prime.

\noindent  (a) If $\pi j = j$, there are $J_{s(j)}(p)$, where $s(j)$ is the number of $\pi$-orbits in $\brcs{1,2,\dots, j-1}$.

\noindent (b) Suppose $\pi j \neq j$. If $\Cal O(j) $ is not contained in $\brcs{1,2,\dots, j}$, then there are zero choices. Assuming $\Cal O (j) \subseteq \brcs{1,2,\dots, j}$ (that is, $j = \max \Cal O(j)$), select $z$ in $\Z_p^*$ with order dividing $|\Cal O(j) |$, and set $a_{\pi j} = -z$. For each of the $s(j)$ orbits in $\brcs{1,2,\dots, j-1}$, we select $m $ in the orbit, and either set $a_m$ to zero (if $z^{|\Cal O (m)|} \neq 1$) or let it be arbitrary (if $z^{|\Cal O (m)|} = 1$), and define the $a_i$ for other $i \in \Cal O(m)$ according to the formulas. The constraint that the column must have content one is automatically satisfied, since $a_{\pi j} \equiv-z$ is relatively prime to $p$.

For $z$ fixed, the number of choices (with $a_{\pi j} \equiv -z$) is thus (provided $\Cal O (j) \subset \brcs{1,2,\dots, j}$)
$$
p^{\left|\Set{\Cal O(m)}{\Cal O(m) \subset \brcs{1,2,\dots, j} \text{ and } z^{|\Cal O(m)|}\equiv 1 \mod p}\right|}.
$$
Now we sum this over all choices for $z$, of which there are $\gcd\brcs{|\Cal O(j)|,p-1}$ (the number of $k$th roots of unity in $\Z_p^*$).

Finally, we observe that if $d$ is  prime, then the set of weakly terminal matrices in $\Mn n \Z$ of determinant $d$ is simply $\dot\cup_{j=2}^n V_j (p)$, since a matrix in Hermite normal form with prime determinant can only have one column that is not the corresponding standard basis element. This leads to the following expression. Recall that $s(j) \equiv s(j,\pi)$ is the number of $\pi$-orbits that are contained in $\brcs{1,2,\dots,j-1}$

\Lem Proposition \appcfou. Let $\pi \in S_n$. Then for a prime $p$,
$$\eqalign{
\Cal S(\pi)(p) &= \sum_{\Set{2 \leq j \leq n}{\pi j = j}} J_{s(j, \pi)}(p) \cr
&\quad + \sum_{\Set{2 \leq j \leq n}{\pi j \neq j \text{ and } \Cal O_{\pi}(j) \subseteq \brcs{1,2,\dots,j}}} \sum_{\Set{z \in \Z_p^*}{z^{|\Cal O_{\pi}(j)|} = 1}} p^{\left|\Set{\Cal O(m)}{\Cal O(m) \subset \brcs{1,2,\dots, j} \text{ and } z^{|\Cal O(m)|} \equiv 1 \mod p}\right|}.
}$$

If $\pi$ is a transposition, then by Lemma \appctwo, we may assume that $\pi = (12)$. In that case, there are $n-2$ fixed points, and $s(j, \pi) = j-2$ for $j \geq 3$. The second sum reduces to the case that $j=2$, and there are exactly two solutions to $z^2 \equiv 1 \mod p$ if $p$ is odd ($z \equiv \pm 1$), and just one if $p =2$. So we obtain
$$
\sum_{l=1}^{n-2} J_l(p) + \cases 2 & \text{if $p $ is odd}\\
1 & \text{if $p =2$}.\\
\endcases
$$
The left sum is  $p^{n-2} + p^{n-3} + \dots + p - (n-2)  = (p^{n-1} -1)/(p-1) -n + 1$; perhaps unsurprisingly, this is $\phi_1 * J_2 * \dots * J_{n-2}(p)$. It is easy to see than any nonidentity permutation other than a  transposition will have  leading term   at most $p^{n-3}$.

If $\pi$ is a cycle of order $n$, then the count is hardly anything, just $\Cal N_n (p) = \gcd\brcs{n,p-1}$.

Now we make some crude estimates for the number of PH-equivalence classes of determinant $\pm d$ matrices in $\NS_n$, denoted $\Cal P\Cal H (n,d)$, when $d$ is square-free. We see that $\cup_{\pi \neq \I} \Cal Z(\pi)(d)$ consists of the weakly terminal matrices (of determinant $d$) whose orbit size is strictly less than $n!$. Let $F(n,d) = (\phi* J_2 * \dots * J_{n-1})(d)$, the number of weakly terminal matrices of size $n$ and determinant $d$. Let $T$ be the set of transpositions in $S_n$, together with the identity element. By Burnside's lemma (actually the lemma that is not Burnside's) and Lemma \appctwo,
$$\eqalign{
\Cal P \Cal H (n,d) &= \frac{F(n,d) + \sum_{\pi \neq \I} \Cal S(\pi)(d)}{n!} \cr
& = \frac{F(n,d) + {n \choose 2}\Cal S(12)(d) +  \sum_{\pi \in S_n\setminus T} \Cal S(\pi)(d)}{n!}
}$$
We know that if $\pi$ is any of the ${n \choose 2}$ transpositions, then $\Cal S(\pi)(p) /F(n,p) \leq 1/p$; hence for $d$ square-free, $\Cal S(\pi)(d) /F(n,d) \leq 1/d$, and if $\pi$ is not a transposition, then $\Cal S(\pi)(p) /F(n,p) \leq \Oh{1/p^{2-\epsilon}}$ (as $d\to \infty$) for all $\epsilon >0$, hence $\Cal S(\pi)(d) /F(n,d) \leq 1/d^{2-\epsilon}$. Thus
$$
\frac{F(n,d)}{n!} \( 1 + \frac{{n \choose 2}}d \( 1- \oh{1}\)\)\leq \Cal P\Cal H(n,d) \leq \frac{F(n,d)}{n!} \(1 + \frac{{n \choose 2}}d +  \frac{n!}{d^{2-\epsilon}} \)
$$
for all $\epsilon > 0$. This is not effective until $d\gg  n!$, but it does yield the conjecture (for square-free $d$).

\Lem Corollary \appcfiv. If $ n\geq 3$ is fixed and $d$ is square-free, then the number of PH-equivalence classes of matrices in $\NS_n$ of determinant $\pm d$ is given by
$$
\Cal P\Cal H(n,d) = \frac{(\phi* J_2 * \dots * J_{n-1})(d)}{n!} \( 1 + \frac {{n \choose 2}}d \(1+ \Oh{\frac 1{d^{1-\epsilon}}}\)\)
$$

For general $d$ (not assumed to be square-free), we can obtain  results for transpositions.
For each  integer $k > 1$, define the function $\Arrow \Cal M_k; \N.\N $ via
$$
\Cal M_k (d) = \left| \Set{(a,t_1,t_2, \dots,t_{k-1} )\in \Z_d^k}{a^2 = 1 \text{ and for all $i$, $t_i (a+1) = 0$} }\right|.
$$
By the Chinese remainder theorem, each $\Cal M_k$ is multiplicative.

We also define the multiplicative function $\Arrow \Cal P_k (d); \N.\Z^+$, via
$$
\Cal P_k (d) = \cases J_k(\sqrt d) & \text{if $d$ is a square}\\
0 & \text{if $d$ is not a square.}\\
\endcases
$$
Informally, $\Cal P_k (d) = \chi_2 (d)\cdot J_k (\sqrt d)$, where $\chi_2$ is the indicator function of the set of square integers.

\Lem Lemma \appcsix. For $p$ a prime,
$$
\Cal M_k (p^m) = \cases p^{m(k-1)} +1 & \text{if $p$ is odd} \\
2^{k-1} & \text{if $p^m = 2$}\\
2^{2(k-1)} + 2^{k-1} & \text{if $p^m = 4$}\\
2^{m(k-1)} + 2^{(m-1))(k-1)} + 2^k & \text{if $8$ divides $p^m$.}
\endcases
$$

\Pf If $a = -1$, then there are $p^{m(k-1)}$ choices for $(t_1, \dots, t_{k-1})$. If $a = 1$ and $p$ is odd, as $2$ is invertible in $\Z_d$, we must have $t_i = 0$, that is, just one solution. When $p$ is odd, the solutions to $a^2 = 1$ are exactly $\pm 1$, hence we have a total of $p^{m(k-1)} + 1 $ solutions.

If $p = 2$ and $a = 1$, the number of  $t$ \st $2t \equiv 0 \mod 2^m$ is $2$; hence this case accounts for $2^{k-1}$ solutions. When $m \geq 3$, there are two other roots of  $a^2 = 1$, $2^{m-1} \pm 1$. When $a = 2^{m-1} -1$, the equations reduce to $2^{m-1}t_i \equiv 0 \mod 2^m$, so there are $2^{(m-1)(k-1)}$ solutions. When $a = 2^{m-1} +1$, the equations become $2(1+ 2^{m-2})t_i \equiv 0 \mod 2^m$, and as the middle factor is invertible, there are just $2^{k-1}$ solutions.

When $p^m = 2$, $a = 1$ is the only square root of $1$, so there are $2^{k-1}$ (preceding paragraph, first line) solutions, as indicated in the statement of the result. If $p^m = 4$, then there are two roots of $1$, $\pm 1$; we have $2^{2(k-1)}$
solutions from $a = -1$ (first line of first paragraph) plus $2^{k-1}$ solutions arising from $a = 1$ (first line of second paragraph).

Finally if $p^m = 2^m$ with $m \geq 3$, we have $2^{m(k-1)} + 2^{k-1} + 2^{(m-1)(k-1)} + 2^{k-1}$ solutions, arising respectively from $a = -1,1, 2^{m-1}-1, 2^{m-1}+ 1$. \qed

Recall that we have abbreviated $(\phi * J_2* \dots *J_{n-1} )(d) $ to $F(n,d)$.

\Lem Lemma \appbnin. Let $\pi \in S_n$ be a transposition, with $ n > 2$. For $d$ a positive integer, the number of $n \times n$ weakly  terminal matrices of determinant $d$ that are invariant under $\pi$ is
$$
\Cal S(\pi)(d) = \(\phi *J_2 * \cdots * J_{n-3}* \Cal P_{n-2} * \Cal M_{n-1} \)(d).
$$

\Pf As $d \mapsto S(\pi(d))$ is multiplicative, we may assume $d = p^m$. By Lemma \appctwo, we may assume that $\pi$ interchanges $n-1$ and $n$. Fix $(k,l)$ with $k+l \leq m$, and let $C $ be a weakly terminal matrix whose last two diagonal entries are respectively $p^k$ and $p^l$, respectively. Denote the entries above the diagonal
by $a_{i,j}$ as usual; for convenience, denote $a_{n-1,n} = a$. Let $P$ be the permutation matrix corresponding to $\pi$, that is, right multiplication by it implements the interchange of the last two columns. Then $\I - P = 0_{n-2} \oplus \(\smallmatrix 1 & -1 \\ -1 & 1 \\ \endsmallmatrix\)$

The condition that $C$ is $\pi$-invariant is equivalent to $CPC^{-1} \in \Mn n \Z$, equivalently, $C(\I - P) C^{-1} \in \Mn n \Z$. To calculate $(\I - P) C^{-1}$, we need only calculate the bottom $2 \times 2$ block of $C^{-1}$, which is found in a matter of seconds to be $p^{-(k+l)}\(\smallmatrix p^l & -a \\ 0 & p^k\endsmallmatrix \)$.
Thus, letting $C_0$ be the upper $n-2$ square block of $C$, 
$$
C(\I- P)C^{-1} = \frac{1}{p^{k+l}}\( \matrix C_0 & a_{1,n-1} & a_{1,n}\\
& a_{2,n-1} & a_{2,n} \\
& \vdots  & \vdots \\
& a_{n-2,n-1} & a_{n-2,n}\\
0 & p^k & a \\
0 & 0 & p^l \\
\endmatrix\) \( 0_{n-2} \oplus  \(\matrix p^l & -(a+ p^k) \\ -p^l & a + p^k\\ \endmatrix \)\).
$$
This multiplies easily (we can ignore $C_0$), and we deduce necessary and sufficient conditions for all the entries to be integers:
\item{(i)} for all $1 \leq j \leq n-2$, $a_{j,n} \equiv a_{j,n-1} \mod p^k$;
\item{(ii)} for all $1 \leq j \leq n-2$, $(a_{j,n} - a_{j,n-1})(a+ p^k) \equiv 0 \mod p^{k+l}$;
\item{(iii)} $a  \equiv 0 \mod p^k$
\item{(iv)} $p^{2l} \equiv 0 \mod p^{k+l}$
\item{(v)} $p^{2k} \equiv a^2 \mod p^{k+l}$

From (iv),  we must have $k \leq l$. We may write $a = p^k a_0$ by (iii), with $a_0 < p^{l-k}$ (as $C$ is weakly terminal). Then (v) yields $a_0^2 \equiv 1 \mod p^{l-k}$. By (i), we may write $a_{j,n} = a_{j,n-1} + t_j p^{k}$, for some $t_j < p^{l-k}$. Then (ii) translates to $t_j (1+a) \equiv 0 \mod p^{l-k}$. Conversely, if $l > k$, given $a_0^2 = 1$ and $t_j$ satisfying $t_j (1 + a) = 0$, then for each choice of $(a_{1,n-1}, a_{2,n-1}, \dots , a_{n-2,n-1})^T$, we obtain a fixed point of $\pi$. If $l = k > 0$, then $a = 0$, and $t_i = 0$ for all $i$, so we obtain exactly one solution for each $(a_{1,n-1}, a_{2,n-1}, \dots , a_{n-2,n})^T$. Finally, if $l =k = 0$, there is only once choice.

 The arbitrary weakly terminal matrix in the upper block, $C_0$, is of determinant $p^{m-k -l}$ and size $n-2$; thus there are $F(n-2,p^{m-k -l})$ choices for $C_0$. For the $(n-1)$st column, there are no constraints on the entries (assuming $k \leq l$), so there are $J_{n-2}(p^k)$ choices (since the column has to be unimodular). Finally, once the $(n-1)$st column entries are determined, we have, by the previous paragraph, $\Cal M_{n-1}(p^{l-k})$ choices. Hence the number of weakly terminal matrices is (on setting $t = l-k$)
$$\eqalign{
\Cal S(\pi)(p^m) & = \sum_{k+l \leq m, \ k \leq l} F(n-2, p^{m-k-l})J_{n-2}(p^k) \Cal M_{n-1}(p^{l-k})\cr
& = \sum_{2k + t \leq m} F(n-2, p^{m - 2k -t}) J_{n-2}(p^k) \Cal M_{n-1}(p^t) \cr
& = \sum_{K  + t \leq m} F(n-2, p^{m - K - t})P_{n-2} (p^{K}) \Cal M_{n-1}(p^t)\cr
& = \(\(\phi * J_2 * \cdots * J_{n-3}\)* P_{n-2}* \Cal M_{n-1}\)(p^m)\cr
}$$
The third line is obtained from the second line via the observation that if $K$ is odd, then $P_{n-2}(p^K) = 0$.
\qed

If $n =3$, the result is $P_{1}* \Cal M_2$, and if $n=4$, the result is $\phi*P_{2}* \Cal M_3$ ($J_1 = \phi$ and $J_0$ is the constant function). So for $\pi$ a transposition,  $\Cal S(\id) - \Cal S(\pi) = (\phi * J_2 * \cdots * J_{n-3}) *(J_{n-2}* J_{n-1} - P_{n-2}* \Cal M_{n-1})$, which is   sufficient to show $\Cal S(\pi)(d)/F(n,d) = \Oh{1/d}$. So the conjecture (for general $d$) would be true if the specific conjecture were true (as it almost certainly is).

\SecT Appendix D: counting PH-equivalence classes of size 3

Here we obtain exact counts for various situations involving the PH-equivalence classes when the matrix size  is $3$, without assuming the determinants are square-free (as always, we are assuming the matrices are in $\NS$). For example, those equivalence classes that contain a terminal matrix with $1$-block size two can be subdivided into three interesting subcases, and we can count each. As a result, we show that in terms of PH-equivalence classes, those of fixed determinant with a $1$-block size two matrix are generically swamped by those not containing one, even when we restrict to square-free determinant (generally, for determinant $d$, the larger $\sum_{p|d} 1/p$ is, the smaller is the ratio of $1$-block size two equivalence classes to the rest).

For $n = 3$, again by Burnside's lemma, the number of PH-equivalence classes is
$$
\frac{\phi * J_2 (d) + 3\Cal S(23)(d) +2\Cal S(132)(d)}6. 
\tag1$$

If $ m > 1$ and $p$ is a prime,
$$
(\phi* J_2)(p^m) =(p^{m-2}(p+1)^2 + 1)p^{m-1} (p-1) = p^{2m} (1 + \frac 1p - \frac 1{p^2} - \frac 1{p^{3}}  - \frac{1}{p^{m+1}} + \frac{1}{p^m})
$$
At $m = 1$, the outcome is simply $p-1 + p^2 - 1 =( p-1)(p+2)$. Hence as a function of $d$, it is a bit less than  $ d^2\prod_{p|d} (p+1) (\prod_{p|d} (1-1/p^2))$. At $2^m$, the outcome is asymptotically $9\cdot 2^{2m-3}(1- \Oh{2^{-m}})$

Now $\phi*J_2 (p^m) = p^{2m} + p^{2m-1} + \dots$, so $\phi* J_2 (d) = d^2 \prod_{p|d} (1 + 1/p + 1/p^2 + \dots)$. We will find that $\Cal S(12)(p^m)  = p^m + p^{m-1} - \dots$, except for $p =2$, when it begins $3\cdot 2^m/4$ rather than $2^m$, so $\Cal S(12)(d) \leq d \prod_{p|d} (1+1/p + 1/p^2 \pm \dots)$ and $\Cal S(132)(d)  = \oh{d^{\epsilon}}$ for all $\epsilon >0$; both of these will result from exact expressions.

There are relatively straightforward asymptotic estimates: for example,
with fixed $n$, the number of equivalence classes of terminal forms with
$1$-block size $n-1$ is bounded below by 
$$
\frac{\max\brcs{(d- \phi(d))^{n-1}, \phi(d)^{n-1}}}{n!}.
$$

However, there are some cases (with  $n=3$), wherein  the formulas
become quite simple. If
$d$ is a prime or a product of two distinct primes,  automatically all terminal forms have $1$-block size
$n-1$. More generally, we obtain exact numbers of PH-equivalence classes for fixed absolute determinant $d$ when  $n =3$.

Let $\Arrow w, w', w''; \N.\C$ be defined, respectively, by $w(d)$ is the number of distinct prime divisors of $d$, $w'(d)$ is the number of distinct prime divisors of $d$ that are congruent to $1$ modulo $3$, and $w'' (d)$ is $1$ if $9$ divides $d$, otherwise it is zero (so $w''$ is the indicator function of $9\N$). Each of them is additive (in the number-theoretic sense), so each of $3^w$, $3^{w'}$, and $3^{w''}$ is multiplicative.

We also define $\Arrow \Cal M_2, \Cal M; \N.\C$ by setting, for $d = \prod_{p|d} p^{m(p)}$,
$$\eqalign{
\Cal M_2 (d) &= \cases 1 & \text{if $m(2) = 0$} \\
2 & \text{if $m(2) = 1$}\\
6 & \text{if $m(2) = 2$} \\
3\cdot 2^{m-2} + 4 & \text{if $m(2) \geq 3$}.\\
\endcases \cr
\Cal M (d)&= \Cal M_2(d) \prod_{\text{odd }p|d} (p^{m(p)} + 1)
}$$
Obviously, $\Cal M$ and $\Cal M_2$ are multiplicative, but not completely multiplicative.

Recall  that $\Cal N_3(d)$ denotes the number of solutions to the polynomial $X^3 - 1 = 0$ in $\Z_d$. On replacing $X$ by $-X$, we see that $\Cal N_3$ also counts the solutions to $X^3 =-1$. By the Chinese remainder theorem, the function $\Cal N_3$ is multiplicative. The following is completely elementary.

\Lem Lemma \appcsev. $\Cal N_3 = 3^{w' + w''}$.

\Pf Both sides are multiplicative, so it suffices to show equality for $d = p^m$, with $p$ prime. The set of solutions to $X^3-1 =0$ is a subgroup of $\Z_d^*$ of exponent three or $1$. If $p \equiv 2\mod 3$, $|\Z_d^*| = \phi(p^m) = p^{m-1}(p-1)$ is relatively prime to $3$, so the solution is unique. In this case, $\Cal N_3(p^m) = 1 = 3^{w'(d) + w''(d)}$. If $p \equiv 1 \mod 3$, $\Z_d^*$ is cyclic of order $p^{m-1}(p-1)$; the latter is divisible by $3$, and as the group is cyclic, it has a unique subgroup of order three. Hence $\Cal N_3(p^m) = 3 = 3^{w'(d) + w''(d)}$.

Finally, if $p =3$, with $m =1$, $\Z_3^*$ is order two, so $\Cal N_3(3) = 1 = 3^{w'(3) + w''(3)}$; if $m \geq 2$, then $\Z_{d}$ is cyclic of order $2\cdot 3^{m-1}$, hence has a unique subgroup of order $3$, and thus $\Cal N_3(d) = 3 = 3^{w'(d) + w''(d)}$.
\qed

\Lem Lemma \appceig. Let $p$ be a prime, and $m$ a positive integer. The number of solutions $(y,k)$ to the equations $Y^2 -1 = 0$ and $(Y+1)K = 0$ in $\Z_{p^m}$ is
$$
\Cal M(p^m) = \cases p^m +1 & \text{if $p$ is odd}\\
\cases 2 & \text{if $p^m = 2$}\cr
6 & \text{if $p^m = 4$}\cr
3\cdot 2^{m-1} + 4 & \text{if $p = 2$ and $ m\geq 3$}\cr
\endcases & {}\\
\endcases
$$

\Pf Since $\Z_{p^m}/p^{m-1}\Z \iso \Z_p$ and the latter is embedded in the former, if $y^2 = 1$, then we can write $y = w + p^ts$ for some some $w \in \brcs{0,1,2,\dots, p-1}$, $1 \leq t \leq m-1$ (so if $m=1$, the second summand disappears), and $(p,s) = 1$ with $1 \leq s \leq p-1$, or $y = w$. In the field $\Z_p$, the only solutions are $\pm1$, so $w = \pm 1$. Thus $1 = y^2  = 1 + p^t s (\pm 2 + p^t s^2)$.
As $(p,s) = 1$, we must have $p^m$ divides $p^t (\pm 2 + p^t s^2)$.

If $p$ is odd, then $\pm 2 + p^t s^2$ is invertible in $\Z_{p^m}$, which forces $y = \pm 1$. When $y = -1$, we can set $k$ to be any element of $\Z_{p^m}$, so we obtain $p^m$ choices, $(-1,k)$. When $y = 1$, $y+1 = 2$ is invertible modulo $p$ and thus modulo $p^m$, and so the only choice is $(1,0)$. Hence there are $p^m + 1$ solutions.

If $p = 2$, and $m =1$, then obviously $y = 1$ and then $k$ can be anything, i.e., we obtain two solutions, $\brcs{(1,0), (1,1)}$. If $m=2$, there are two square roots of unity in $\Z_4$, $\pm 1$ (or $\brcs{1,3}$ if you prefer); if $y = -1$, we obtain the four solutions $(-1, k)$, while if $y = 1$, there are only two, $\brcs{(1,0), (1,2)}$. Thus when $m=2$, there are a total of $6$ solutions.

If $p = 2$ and $m \geq 3$, there are now four square roots of $1$, $y = \pm 1 + 2^{m-1}u$ where $u \in \brcs{0,1}$, as follows easily from $2^m$ dividing $2^t (\pm 2 + p^t s^2)$. If $y = -1$, we have the $2^m$ solutions $\brcs{(-1,k)}$; if $y = 2^{m-1} -1$, then $y + 1 = 2^{m-1}$, so we obtain $2^{m-1}$ solutions, $\brcs{(2^{m-1} -1, 2j)}_{0 \leq j < 2^{m-1}}$. If $y = 1 + 2^{m-1}u$, then $y+1 = 2(1 + 2^{m-2}u)$; as the second factor is a unit (since $m \geq 3$), it follows that in order that $k(y+1) = 0$, we must have $2^{m-1}$ divides $k$. Hence in both cases, there are only two solutions.

Thus if $8 $ divides $p^m$, we must have $2^m + 2^{m-1} + 4$ solutions in total.
\qed

By the Chinese remainder theorem, the number of solutions  $(k,y) \in (\Z_d)$ of the equations $Y^2 =1$ and $(Y+1)K = 0$ is $\Cal M(d)$.

Now we determine $\Cal S(23)(p^m)$  and $\Cal S(132)(p^m)$.
The generic weakly terminal matrix is given by
$$
C = \(\matrix 1 & a & b \\ 0 & e & gy \\ 0 & 0 & gx \endmatrix \), \text{ and its inverse is } C^{-1} = \frac 1{egx} \(\matrix egx & -a & ayg-be \\ 0 & gx & -gy \\ 0 & 0 & e \endmatrix \) \in \Mn 3 \Q,
$$
where all of $\brcs{a,b,e,g,y}$ are nonnegative and $ a < e$, $y <x$, $b < gx$, and $\gcd\brcs{a,e} = \gcd\brcs{b,g}=  \gcd\brcs{x,y} = 1$ (by convention, $\gcd\brcs{0,m} = m$). When we have a permutation acting on $C$, it also acts on the triple $ (J(C_{\Omega(i)}))_{i=1}^3$ by permuting according to its action on the columns. Since the three invariants for the generic matrix are, in order (that is, deleting the first column, then deleting the second column), $(\Z_{(\delta,gx)},\Z_g, \Z_e)$ where $\delta = agy -be$  (the determinant of the upper right block), if for example $\pi = (23)$ or $(132)$ and $C$ is invariant under the action of $\pi$ (meaning that $CP_{\pi}C^{-1} \in \Mn 3 \Z$), then we must have $e = g$.

Define the multiplicative functions, $\chi_2$ and $\Cal P$,
$$\eqalign{
\chi_2  &\text{ is the indicator function of the set of square integers and}\cr
\Cal P(d) &= \chi_2 (d)\cdot \phi(\sqrt d). \cr
}$$
\Lem Lemma \appcnin. For $d$ a positive integer, the number of weakly terminal $3 \times 3$ matrices of determinant $d$ that are invariant under a transposition $\pi \in S_3$ is
$$
\Cal S(\pi)(d) = (\Cal P * \Cal M)(d).
$$

\Pf This is Lemma \appbnin.

\qed

This isn't useful unless we can describe the resulting convolution product. The formula below when $p=2$ is obtained by considering a number of special cases, and then summing geometric series; it did not seem worthwhile to transcribe the tedious argument.

\Lem Lemma \appcten. (a) If $p$ is an odd prime and $\pi$ is a transposition, then
$$
\Cal S(\pi)(p^m) = (\Cal P * \Cal M)(p^m) = p^m + p^{m-1} + 1 + \cases p^{m/2} - p^{m/2 -1} & \text{if $m$ is even}\\
-p^{(m-1)/2} & \text{if $m$ is odd.}
\endcases
$$
\noindent (b) For $ p =2$,
$$
\Cal S(\pi)(2^m) = (\Cal P * \Cal M)(2^m)  = \cases 2 & \text{if $m=1$}\\
 7  & \text{if $m = 2$}\\
12    & \text{if $m = 3$}\\
2^m + 2^{m-3} + 2^{m/2} -1    & \text{if $m \geq 4$ and is even} \\
2^m + 2^{m-3} + 2^{(m-1)/2}  + 2^{(m-3)/2}&\text{if $m \geq 5$ and is odd.}\\
\endcases
$$

\Pf When $p$ is odd and $s> 0$, $\Cal M(p^s) = p^s +1$ and $\Cal M(1) =1$. Taking into account the latter, we have
$$\eqalign{
\Cal S(\pi)(p^m) & = \sum_{0 \leq n \leq  m/2} \phi(p^n) p^{m-2n} +\sum_{0 \leq n < m/2}\phi(p^{n}) \cr
& = p^m + 1+ (p-1) \sum_{1 \leq n \leq m/2} p^{m-n-1} + p^{\flo{(m-1)/2}}\cr
& = p^m + 1+  \cases
(p-1)\frac{p^{m-1} - p^{m/2 -1}} {p-1}    +   p^{m/2}  & \text{if $m$ is even}\\
  (p-1)\frac{p^{m-1} - p^{(m-1)/2}} {p-1}                                      & \text{if $m$ is odd.}\\
\endcases \cr
& = p^m + p^{m-1} + 1 + \cases p^{m/2} - p^{m/2 -1} & \text{if $m$ is even}\\
-p^{(m-1)/2} & \text{if $m$ is odd.}
\endcases \cr
}$$

When $p =2$, the computation is more complicated because of the definition of $\Cal M(2^r)$. Fortunately, there is still massive cancellation, and after a battle keeping track of the limits of summation, we obtain the result in the statement of the lemma.
\qed

In particular, the number of invariant $C$ is $d \cdot \prod_{p|d, \text{ $p$ odd}} (1+1/p \pm \dots) \cdot \alpha (v_2(d))$ where $\alpha$ is the function obtained in the last lemma, divided by $2^{v_2(d)}$ (for $m \geq 4$, $\alpha(m) = 1 + 1/8 + \dots$; the dependence on the exponent, $m$, is tiny if $m$ is large).

Now to deal with $\Cal Z(132)(d)$, the set of weakly terminal matrices invariant under the permutation matrix corresponding to $(123)$ or $(132)$. This is fairly horrible, but is not as bad as it could be. It is marginally better to use $(132)$, rather than $(123)$ (the groups they generate are the same, but the computations are a bit less tedious in the former case).

\noindent {\it $\Cal S(\pi)(d)$ with $\pi = (132)$.} From the column action, we have $e = g = (\delta,gx)$ (the last equality,
in the presence of the first, is equivalent to $(\delta,x) = 1$), so $d =
e^2 x$. The equations are then
$$
-a^2 + b \equiv 0 \mod e; \qquad ay - b - y^2 \equiv 0 \mod x; \qquad a^2
y - ab - by + 1 \equiv 0 \mod ex.
$$
Rewrite the third one as $(a^2 - b)y + 1 - ab \equiv 0 \mod ex$. Taking
this modulo $e$, we obtain $ab \equiv 1 \mod e$, which in combination with
the first, yields $a^3 \equiv 1 \mod e$ (and also $b^3 \equiv 1 \mod e$).
Write $b = a^2 + ke$, where $k$ is defined modulo $x$. Plugging this into
the second and third equations yields
$$
-key + {1-ab} \equiv 0 \mod ex; \qquad  y^2 - ay + a^2 + ke \equiv 0 \mod x.
$$
The former yields $-key - a^3 -ake +1  \equiv 0 \mod ex$, so $ke(y+a)
\equiv 1 - a^3 \mod ex$. Multiplying the second displayed equation by $y +
a$ yields $ke(y+a) \equiv - (y^3 + a^3)\mod x$, whence $y^3 \equiv -1 \mod
x$. Since $e|(ab-1)$, we also have $-ky \equiv (1-ab)/e \mod x$. In
particular, if $x \neq 1$, then $k$ (and thus $b$) is uniquely determined
by $y$ modulo $x$.
 
We recall from Lemma \appcsev\  that $X^3 \pm 1 =0$ each has three distinct solutions in
$\Z_{p^m}$ iff $p \equiv 1 \mod 3$ or $9 |p^m$, and otherwise each has
one.
 
Now set $d = p^m$, $e = x^n$, and $x = p^r$  with $2n + r = m$.
First, suppose that $n = 0$, so $r = m$, and the equations boil down to $a
= 0$, $y^3 \equiv -1 \mod p^m$, $b \equiv y^2 \mod p^m$ (so $b$ is
determined by $y$), and $by \equiv -1 \mod x$, but the last is a
consequence of the second last.
 
If $p \equiv 2\mod 3$ or $p^m = 3$, then $y^3 \equiv -1 $ entails $ y
\equiv -1\mod p^m$. Hence $b\equiv 1\mod p^m$, so there is exactly one
solution for $(a,b,y) = (0,1,-1)$. If $p\equiv 1 \mod 3$ or $9|p^m$, there are three choices for
$y$, and thus a total of three choices for $(a,b,y)$ when $n =0$.
 
Now suppose that $n > 0$. If $r=0,$ then $m = 2n$, and the only conditions
imposed are $a^3 \equiv 1\mod p^n$, $y = 0$, and $b \equiv a^2\mod p^n =
ex$. Hence we obtain three solutions for $(a,y,b)$ if $p \equiv 1\mod 3$
or $9|p^m$ (since $b$ is determined by $a$), and $1$ otherwise.
 
Finally suppose that $n,r> 0$, so that $1 \leq n < m/2$. Here $y$ is
defined modulo $p^r = x$ and $b$ is defined modulo $p^{n+r} = ex$. We have
$a^3 \equiv 1 \mod p^n$, and we can write $b = a^2 + kp^n$ (where $k$ is
defined modulo $p^r$). We also have $y^3 \equiv -1\mod p^r$, that is,
$(-y)^3 \equiv \mod p^r$.
 
If $p \equiv 2 \mod 3$, then $a = 1$ (defined modulo $p^n$)
and $y = -1$ (defined modulo $p^r$), and thus $kp^n \equiv 1 + 1 + 1 \mod
p^r$. This forces (since both $n$ and $r$ are positive), $p = 3$ a contradiction, so that in this case, there are no solutions.

If $p^m =3$,
then $n + r = 1$, contradicting $n,r> 0$.
 
If $p \equiv 1 \mod 3$ or $9|p^m$, there are three choices for $a$, and
also for $y$. However, they are not independent of each other. Modulo $p$,
either $y+a$ is $0$ (which corresponds to taking the same cube root of
$\pm 1$) or invertible. But if $p|(y+a)$, as in the previous paragraph, we
obtain $kp^n\equiv -(y^2 -ay + a^2) \mod p^r \equiv -3a^2 + pX \mod p^r$.
This yields a contradiction, unless $p = 3$---and in that case, we have $m
\geq 2$, so either $r = 1$ (in which case $k$ has three values), or $n =
1$ and $r>1$, and in that case $k$ is uniquely determined.
 
Finally, if $p \equiv 1 \mod 3$ or $9|p^m$ and $y+a \neq 0 \mod p$, then
there are six choices for $(a,y)$, namely so that $y+a$ is invertible
modulo $p$, hence modulo any power of $p$, and for each of these, the
equation $k(y+a) \equiv (1-a^3)/e \mod p^r$ uniquely determines $k$.
 
Now we count all these possibilities. Let $H(t)$ be $1$ if $t$ is odd, and $2$ if
$t$ is even. If $p \equiv 2\mod 3$,  there are
zero solutions for $n,r > 0$, giving us a total of $1$ solution (arising
from $n = 0$) plus an additional $1$ iff $m$ is even. So the formula is
$H(m(p))$.
 
If $p \equiv 1 \mod 3$, we obtain $3$ solutions from the case $n = 0$ plus
an additional $3$ if $m$ is even ($r = 0$), plus $\sum_{1 \leq n < m/2} 6
= 6\flo{(m-1)/2}$. Here the formula is $6\flo{(m(p)-1)/2} + 3 H(m(p))$.
This is $3(m(p)-1) + 6 = 3(m(p)+1)$ if $m(p)$ is even, and $3(m(p)-1) + 3
= 3m(p)$ if $m(p)$ is odd, so we can rewrite the expression as $3(m(p)
+H(m(p))-1)$.

Now we look at the totals for the various situations. We recall  $r = m-2n$, so that $r = 0$ entails $m$ is even and $r =1$ entails $m$ is odd. If $p \equiv 2 \mod 3$, then
$$
\Cal S(132)(p^m) = 1 + 0 + H(m)-1 = H(m).
$$
If $p \equiv 1 \mod 3$, then
$$\eqalign{
\Cal S(132)(p^m) &= 3 + \sum_{1 \leq n < m/2} 6 + 3(H(m)-1)\cr
& = 3 + 6\flo{\frac{m-1}2} + 3(H(m)-1) \cr
& = \cases {3m + 3} & \text{if $m$ is even}\\
3m & \text{if $m$ is odd.}\\
\endcases
}$$
If $p^m =3$, then
$$
\Cal S(132)(3) = 1 .
$$
Finally, if $p =3$ and $m \geq 2$,
$$\eqalign{
\Cal S(132)(3^m) &= 3 + \sum_{1 \leq n < m/2} 6 + 3(H(m)-1) +1\cr
& = \cases {3m + 4} & \text{if $m$ is even}\\
3m +1& \text{if $m$ is odd.}\\
\endcases
}$$

We can combine these in one gigantic formula,
$$
\Cal S(132)(d) = 2^{\left|\Set{p|d}{p \equiv 2 \mod 3; \, m(p) \text{ even}}\right|}\cdot 3^{w'(d)} \cdot \text{\hglue -.75 cm}\prod_{{p|d, \ p \equiv 1 \mod 3}} (m(p) + H(m(p))-1)\cdot \cases 1 & \text{if $m(3) \leq 1$}\\ 3m(3) +4 & \text{even $m(3)>0$}\\ 3m(3)+1 &\text{odd $m(3)>1$}
\endcases
$$
 
It is not necessary for  the counting formula, but a similar computation (much easier than the others) reveals that the number of $S_3$-invariant weakly terminal matrices of determinant $d$ is
$$
\Cal S(S_3)(d) = 2^{w(d)}\cdot \cases 1 & \text{if $m(3) = 0$}\\
\frac{m(3)+1}2 & \text{if $m(3) > 0$.}\\
\endcases
$$
 In particular, $\Cal Z (123)(d) = \Cal Z (S_3)(d)$ iff $d$ is a square all of whose prime divisors are congruent to $1$ modulo $3$, and in that case, their cardinality is $2^{w(d)}$.

Equation (1) at the beginning of this section now yields the number of PH-equivalence classes of $ C \in \NS_3$ with determinant $\pm d$:
$$
\Cal P \Cal H(3,d):= \frac{(\phi * J_2)(d) + 3 (\Cal P * \Cal M)(d) + 2\Cal S(132)(d)}6.
$$
The last term is too complicated to expand compactly, but it is given explicitly above. When $d$ is square-free, the formula simplifies considerably, and we will discuss this later.

We see from Lemma \appcten\ and the formula for $\phi * J_2$ that $(1-1/p^2)\Cal S(23)(p^m) \leq (\phi*J_2 )(p^m)/p^m$, so $\Cal S(23)(d)/ ( \phi*J_2 (d) \leq \zeta(2)/d$. And similarly, $\Cal S(132) (d)= \oh{d^{-2+\epsilon}}\cdot \phi*J_2 (d)$. The number of PH-equivalence classes, $\Cal P \Cal H(3,d)$  thus satisfies
$$
1 + \frac{3}d  \leq \frac{\Cal P \Cal H (3,d) }{(\phi * J_2)(d)/6 }\leq 1 + \frac{3\zeta(2)}d + \oh{\frac 1{d^{2-\epsilon}}}
$$
for all $\epsilon > 0$. This of course is close to the  Conjecture of Appendix B, when $n = 3$.   The little oh term may be unnecessary.

\noindent {\it $1$-block size two PH-equivalence classes.} We now investigate the number of PH-equivalence classes of fixed absolute determinant that contain a $1$-block size two weakly terminal (and thus terminal) matrix, so that we can compare them with the total number of PH-equivalence classes.
This time, the set of matrices that we are looking at are not invariant under the action of $S_3$, so somewhat different methods are used.

So fix   $d>1$, and consider the PH-equivalence classes having a $1$-block size two terminal form. We perform operations within the ring $\Z_d$. The third
column of one of these terminal forms is
$$
\(\matrix a_1 \\ a_2 \\ d \\ \endmatrix \),
$$
where the the ideal generated by $\brcs{a_1,a_2}$ is $\Z_d$ (when we regard $a_i$ as elements of $\Z_d$), and $0 \leq a_i <
d$ (when we regard $a_i$ as integers).


Now we count the number of  of PH-equivalence classes of matrices  $B \in \NS_3$ with absolute determinant $d$, having a terminal form with $1$-block size two.

Recall the multiplicative function $\Arrow w'; \N.\Z^+$;  $w'(d)$ is the number of distinct prime divisors of $d$ that are congruent to $1$ modulo $3$.

\noindent {\bf Case 1:}  {\it $J(B\op) \iso \Z_d^2$.} In this case, by Corollary \onefou\ (even without the hypothesis that $B$ has $1$-block size $2$), $B$ has a terminal form,
$$
\( \matrix  1 & 0 & a_1 \\
0 & 1 & a_2 \\
0 & 0 & d\\ \endmatrix\),
$$
where $\gcd\brcs{a_1,d} = \gcd\brcs{a_2,d} = 1$. We now view the entries of the truncated column $(a_1, a_2)^T$ as elements of $\Z_d^*$. The equivalence class of such a truncated column, renamed $(x,y)^T$, is given by
$$
\brcs{\(\matrix x \\ y \\ \endmatrix\),\(\matrix y \\ x \\
\endmatrix\),\(\matrix x^{-1} \\ -x^{-1}y \\ \endmatrix\),\(\matrix -x^{-1}y
\\ x^{-1} \\ \endmatrix\),\(\matrix -xy^{-1} \\ y^{-1} \\
\endmatrix\),\(\matrix y^{-1} \\ -xy^{-1} \\ \endmatrix\)}.
$$
Most of these equivalence classes have size six, but some have size $1$, $2$, or $3$. We count the latter, and then obtain a fairly simple formula for the number of equivalence classes.

\noindent {\bf 1a.} {\it Equivalence class  size $1$.} There is only one element with this property, explicitly $(-1,-1)^T$.

\noindent {\bf 1b.} {\it Equivalence class size $3$.} An inspection of the six possible elements in the equivalence class reveals that the only such with exactly three elements are those of the form,
$$
\brcs{\(\matrix x\\ x\\  \endmatrix \),\(\matrix -1\\ x^{-1}\\  \endmatrix \),\(\matrix x^{-1}\\ -1\\  \endmatrix \)},
$$
provided $x \neq -1$. There are thus $\phi(d) -1$ equivalence classes here, covering $3\phi(d) - 3$ elements.

\noindent {\bf 1c.} {\it Equivalence classes of size $2$.} These are of the form
$$
\brcs{\(\matrix \alpha\\ \beta \\  \endmatrix \),\(\matrix \beta\\ \alpha\\  \endmatrix \)},
$$
where $\alpha^3 = -1$, $\beta = -\alpha^2$, and $\alpha \neq -1$. To count the number of choices for $\alpha$ (and $\beta$), we first observe that if $p$ is a prime exceeding $3$, then the equation $z^3 = -1$ has a solution other than $-1$ in $\Z_p$ iff $-3$ is a square modulo $p$, equivalently iff $p \equiv 1 \mod 3$, and when this occurs, the solutions are distinct. It is easy to verify that these properties hold for any power of $p$ as well.

If $p = 3$, then $-1$ is the only solution to $z^3 = -1$ modulo $3$, but modulo any higher power of $3$, there are exactly $3$ distinct solutions: modulo $p^m$, the solutions are $\brcs{p^{m-1} -1, 2p^{m-1}-1, -1}$, including $-1$.

If $p=2$, then there is only one solution to $z^3 = -1$ modulo $2^m$.

Write $d = 3^{m(3)}\cdot \prod_{p\in P} p^{m(p)}\prod_{q\in Q} q^{m(q)}$ where $P$ is the set of primes congruent to one modulo three, and $Q$ is the set of primes (including $2$) congruent to two modulo three. By the Chinese remainder theorem, the number of solutions (including $-1$) to the equation $z^3= -1$ is thus
$3^{|P|}\cdot 3^a$ where $a = 0$ if $m(3) \leq 1$ and otherwise equals $1$.  After discarding the solution $x = -1$ (which is  in 1a), the number of columns covered is
$3^{w'(d) + a} -1$, accounting for half as many equivalence classes.

The remaining columns (out of the original $\phi(d)^2$) have six-element equivalence classes. Hence the total number of equivalence classes is
$$
\frac{\phi(d)^2 - 1 - (3\phi(d)-3) - (3^{w'(d)+ a}-1)}6 + 1 + (\phi(d)-1) + \frac{3^{w'(d)+ a}-1}{2} = \frac{\phi(d)^2 + 3\phi(d) + 2\cdot 3^{w'(d)+a} }6,
$$
where $a = 0$ if $m(3) \leq 1$, and $1$ otherwise.

This is worth stating as a result on PH-equivalence classes.

\Lem Proposition \appbone. (Case 1) For $d$ a fixed positive integer, the number of PH-equivalence classes of matrices $B \in NS_3$ with $|\det B| = d$ and $|J(B\op)| = d^2$ is
$$
 \frac{\phi(d)^2 + 3\phi(d) + 2\cdot 3^{w'(d)+w''(d)} }6.
$$

\Pf The only thing we have to note is that if $|J(B\op)| = |\det B|^2$
for $B \in \NS_3$, then by Lemma \onethr\ and Corollary \onefou, $B$ has a terminal form of the type discussed in case 1 above (it also follows that $J(B\op) \iso (\Z_d)^2$).
\qed

\noindent {\bf Case 2:} {\it Exactly one of $\brcs{a_1, a_2}$ is invertible modulo $d$.} In this case, $a_1 \not \equiv a_2 \mod d$.  By \twoone, the equivalence classes are then of the form,
$$
\brcs{\( \matrix x \\ y \\ \endmatrix\), \( \matrix x^{-1} \\ -x^{-1}y \\ \endmatrix\),\( \matrix y \\ x \\ \endmatrix\),\( \matrix -x^{-1}y \\ x^{-1} \\ \endmatrix\)
}
$$
where $x \in \Z_d^*$ and $y \not\in \Z_d^*$. The only possible equivalence classes with fewer than four elements are those with two, and this occurs iff $x = x^{-1}$ and $y = - x^{-1}y$; this reduces to $x^2 = 1$ and $(1+x) y =0$. By Lemma~\appceig, the number of choices for $(x,y)$ is $\Cal M(d)$.

Now the only case in which $y$ can be a unit occurs when $x = -1$, and in that case $y$ can be anything.  So to obtain the number of solutions in which $y$ is a nonunit, we simply subtract $\phi(d)$ from ${\Cal M}(d)$. The number of solutions to $x^2 =1$ and $(1+x)y = 0$ for $(x,y) \in \Z_d^* \times (\Z_d \setminus \Z_d^*)$ is thus
$$
N_2(d):=  \Cal M (d)-\phi(d)
$$
All of the other possible $2\phi(d)\cdot (d-\phi(d))$ columns have four-element equivalence classes; hence the total number of equivalence classes for case 2 is
$$
\frac{2\phi(d)\cdot (d- \phi(d)) - 2N_2(d)}4 + {N_2(d)} = \frac{\phi(d)\cdot (d- \phi(d)) + N_2(d)}2.
$$

\Lem Proposition \appbtwo. (Case 2) The number of PH-equivalence classes corresponding to case 2 is
$$
\frac{\phi(d)\cdot (d- \phi(d) -1) + \Cal M(d)}2.
$$

Case two corresponds to all the situations in which $J(B\op) $ contains a proper direct summand isomorphic to $\Z_d$ but $|J(B\op)| < d^2$. The remainder are covered by case\,3.

\noindent {\bf Case 3:} {\it Both $a_1$ and $a_2$ are nonunits in $\Z_d$.} This can be restated as $J(B_{\Omega(i)})$ is not zero for exactly two choices of $i$. Since $B\in \NS_3$, we also have to have $\gcd\brcs{a_1,a_2, d} = 1$, equivalently, that in $\Z_d$, the ideal generated by $\brcs{a_1, a_2}$ is the improper one.

So let $(x,y)^T$ correspond to such a truncated column; for most of this, we regard them as integers (rather than elements of $\Z_d$), each with $\gcd\brcs{x,d}, \gcd\brcs{y,d} > 1$, and of course, $1 \leq x , y \leq d-1$ (we cannot have $x = 0$, $y =0$, or $x = y$, since $\gcd\brcs{x,y,d} = 1$). All the equivalence classes here consist of exactly two elements (the column and its flip), so it is simply a matter of counting the number of pairs, and dividing by two.

First, we note that if $w(d) =1$ (that is, $d$
is a power of a single prime), then there are no equivalence classes. So we assume $k:= w(d) \geq 2$, and write $d = \prod_{i=1}^{k}p_i ^{m(i)}$, and $S = \brcs{1,2,\dots, k}$. For a subset $\Omega$ of $S$, write $d_{\Omega} = \prod_{i \in \Omega} p_i^{m(i)}$ and $D_{\Omega} = \prod_{i \in \Omega} p_i$. Thus $d_{\emptyset} = D_{\emptyset} = 1$, $d_S = d$, and $D:= D_S = \prod_{p|d} p$.

For an eligible truncated column $(x,y)^T$, we may write uniquely
$x = D_{\Omega_1}\cdot t_1$, $y = D_{\Omega_2}\cdot t_2$, subject to the following conditions:
\item{(i)} $\Omega_i \neq \emptyset$
\item{(ii)} $\Omega_1 \cap \Omega_2 = \emptyset$
\item{(iii)} for all $p \in \Omega_i^c$, $\gcd\brcs{t_i,p} = 1$.

We see that since $1\leq x,y < d$, we have $1 \leq t_i < d/D_{\Omega_i}$. If we fix the ordered pair $(\Omega_1, \Omega_2)$, then the number of choices for $t_i$ is $$
\frac{d}{D_{\Omega_i}}\prod_{j \in \Omega_i^c} \(1-\frac 1{p_j}\).
$$
Thus the number of eligible truncated columns corresponding to fixed $(\Omega_1, \Omega_2)$ is the product,
$$\eqalign{
\frac{d}{D_{\Omega_1}}\prod_{j \in \Omega_1^c} \(1-\frac 1{p_j}\)\cdot\frac{d}{D_{\Omega_2}}\prod_{j \in \Omega_2^c} \(1-\frac 1{p_j}\) &= \frac{d^2}{D_{\Omega_1 \cup \Omega_2}} \cdot  \prod_{p|d} \(1- \frac 1p\) \cdot \prod_{j \in \Omega_1^c \cap \Omega_2^c} \(1-\frac 1{p_j}\) \cr
& = d\phi(d)\frac{\prod_{j \in \Omega_1^c \cap \Omega_2^c} \(1-\frac 1{p_j}\)}{D_{\Omega_1 \cup \Omega_2}}\cr
& = \phi(d)^2 \cdot \frac 1 {\phi(D_{\Omega_1 \cup \Omega_2})} \cr
}$$
Now let $\Omega$ be a subset of $S$, say with $|\Omega| =s$; the number of ways
of writing it as a disjoint union of $\Omega_1$ and $\Omega_2$ (maintaining the ordering) with neither being the empty set, is zero if $s \leq 1$, and otherwise
$$
\sum_{i=1}^{s-1} {s\choose i} = 2^s - 2.
$$

Define the polynomial $f(x) = \prod_{p|d} (1 + x/(p-1))$.

The total number of truncated columns is thus
$$\eqalign{
\phi(d)^2 \sum_{s=2}^k (2^s -2)\sum_{|\Omega| = s} \frac{1}{\phi(D_{\Omega})} & =
\phi(d )^2\(1+  \sum_{s=0}^k (2^s -2)\sum_{|\Omega| = s} \frac{1}{\phi(D_{\Omega})} \)
\cr
& = \phi(d)^2 \( 1 + f(2) - 2 f(1)\) \cr
&= \phi(d)^2 \( 1 + \prod_{p|d} \(1 + \frac{2}{p-1} \)- 2 \prod_{p|d} \(1 + \frac{ 1}{p-1}\) \)\cr
& = \frac{\phi(d)^2}{\phi(D)} \(\prod_{p|d} (p+1) - 2 \prod_{p|d} p + \prod_{p|d} (p-1)\)\cr
& = d \phi(d)  \(\prod_{p|d} \(1+ \frac 1p\) - 2  + \prod_{p|d} \(1 - \frac 1p\)\)
}$$

The number of equivalence classes for case three is half of this.

\Lem Proposition \appbthr. (Case 3) The number of PH-equivalence classes of $ B\in \NS_3$ with $|\det B | = d$ corresponding to case\,3 is
$$
\frac {d \phi(d)}{2}  \(\prod_{p|d} \(1+ \frac 1p\) - 2  + \prod_{p|d} \(1 - \frac 1p\)\).
$$

When $\sum_{p|d} 1/p$ is large, the two rightmost summands are small compared to $\prod (1+ 1/p)$; in that case, this is asymptotic with (provided we choose $d$s so that $\sum_{p|d} 1/p$ becomes arbitrarily large)
$$
\frac{d^2}2 \prod_{p|d} \(1 - \frac 1{p^2}\).
$$
Given $\epsilon$, there exists $N$ \st $\sum_{p\geq N} 1/p^2 < \epsilon$; hence given $M$, we can find $d \equiv d(\epsilon)$ \st $\sum_{p|d} 1/p^2 < \epsilon$ and $\prod_{p|d} (1 + 1/p) > M$. It follows that the least upper bound for the number of equivalence classes is at least $d^2/2$ (and we can choose square-free $d $ to asymptotically reach this). On the other hand, initially, we only have a choice of $(d-\phi(d))^2/2$ columns, so this is the best possible (and note that $\phi(d) /d \to 0$ for these sequences).

This means that case 3 overwhelms the other two cases (asymptotically) for the appropriate choice of $d$s (with large numbers of prime divisors). On the other hand, with few prime divisors (or simply small $\sum_{p|d}1/p$), cases 1 and 2 together are dominant. With just one prime divisor, case 3 is empty.

An amusing example occurs when $d(j) $ is the product of the first $j$ primes. Then
$$
\lim_{j\to \infty}\frac{\text{number of case 3 PH-equivalence classes for $B \in \NS_3$ with $|\det B| = d(j)$}}{d(j)^2} = \frac 1{2\zeta(2)}.
$$
For case 2 with the same sequence, the number of PH-equivalence classes is asymptotic to $\phi(d)d/2$, which is smaller. With case 1, the number is about $\phi^2(d)/6$,  smaller still. So in the display we could replace \quotes{case 3} by PH-equivalence classes that contain a terminal form with $1$-block size two.

If $B$ is classified in case\,3, then $I(B\op) \iso \Z_d$; however, there are also examples as part of case\,2 with the same property (case 2 examples with $I(B\op) \iso \Z_d$ automatically have the property that $B\op$ also has a terminal form with $1$-block size two; however, not all case 3 classes satisfy this).

There are a couple of  situations in which we can go directly to the number of PH-equivalence classes, without requiring the restriction to those with $1$-block size $n-1$.

\Lem Lemma \appbfou. If $B \in \NS_3$ and $d:= |\det B|$ is either a prime or of the form $pq$ for distinct primes $p$ and $q$, then $B$ is PH-equivalent to a terminal form with $1$-block size $2$.

\Rmk We have seen that the conclusion can fail if $d $ is a product of three distinct primes, in fact, $d = 30 = 2\cdot 3\cdot 5$, and of course, it can also fail if $d = p^2$.

\Pf This is a special case of \nthreig.
\qed

Adding the results from case\,1, case\,2, and case\,3 yields the next result, without referring to the general horrible formula (1).

\Lem Proposition \appbfiv. Suppose the positive integer $d$ is of one of the following forms,
$d = p, 2p, pq$ where $p$ and $q$ are distinct odd primes. Then the number of PH-equivalence classes of $B \in \NS_3$ \st $|\det B | = d$ is
$$\eqalign{
\frac{p^2 + 4p + 1 + 2\cdot 3^{w'(p)}}6 \qquad& \text{if $d = p$}\cr
\frac{2p^2 + 5p- 1 + 3^{w'(d)}}3 \qquad& \text{if $d = 2p$}\cr
\frac{\phi(d) (3d - 2\phi(d) + 3) + 2 \cdot 3^{w'(d)}}6 + d + 1 \qquad & \text{if $d = pq$},\cr
}$$
where  $w'(d)$ is the number of distinct prime divisors of $d$ that are congruent to $1$ modulo $3$.

The number of equivalence classes with $|\det B| = 2p$ is itself divisible by $p$ iff $p \equiv 2 \mod 3$. 

There is one more bit of low-hanging fruit.

\Lem Proposition \appbsix. If $p$ is a prime, then the number of PH-equivalence classes of $B \in \NS_3$ with $|\det B| = p^2:= d$ is given by
the number of PH-equivalence classes for $1$-block size two of determinant $d$ (cases 1 and 2 for $d = p^2$) plus the number of PH-equivalence classes for case~ 1 with $d =p$.
This is
$$
\eqalign{
7 \qquad & \text{if $p = 2$} \cr
\frac{p^4 + p^3 + 2p^2 + p + 1 + 2\cdot 3^{w'(d)}(1 + 3^{w''(d)})}6 \qquad & \text{if $p \neq 2$}.\cr
}$$

\Pf Let $B \in \NS_3$ have determinant $\pm p^2$. Any of its terminal forms has diagonal either $(1,1,p^2)$ or $(1,p,p)$. In the former case, it has a terminal form with $1$-block size two, so is covered by cases 1, 2, and 3; however, for a power of prime, case 3 is empty.

Suppose that the diagonal is   $(1,p,p)$. Then the terminal form must be
$$
B':= \(\matrix 1 & b_1 & b_2 \\
0 & p & 0 \\
0 & 0 & p \\ \endmatrix\),
$$
where  $1 \leq b_i < p $ and $\gcd\brcs{b_1, p} = \gcd\brcs{b_2, p} = 1$
(recall the condition in the terminal form that the diagonal entry in the second row from the bottom must be less than or equal to the  greatest common divisor of the bottom diagonal entry and the entry immediately above; this explains the zero in the $(2,3)$ position). Now for $i = 1,2,3$, each of $I(B_{\Omega(i)})$ is $\Z_p$, a trivial computation. Hence $B'$ (and thus $B$) is not PH-equivalent to a terminal form with $1$-block size two, so these equivalence classes are disjoint from the former case.

However, if we calculate ${B'}\op$, we find that it is PH-equivalent to a $1$-block size two terminal form, with determinant $p$, corresponding to case 1 of the latter class:
$$
{B'}\op =  \(\matrix  p & 0 & 0 \\
-b_1 &1 & 0 \\
-b_2 & 0 & 1 \\ \endmatrix  \) \sim  \(\matrix 
1 & 0 & -b_2 \\
0 & 1 & -b_2 \\
0 & 0 & p \\ \endmatrix  \).
$$
(The PH-equivalence was implemented by conjugation with the permutation matrix that transposes $1$ and $3$.) Thus ${}\op$ implements a bijection between the current matrices and the matrices covered by case 1 for $d = p$, and of course, this bijection preserves PH-equivalence classes. Hence the number of equivalence classes arising from terminal forms with diagonal $(1,p,p)$ is the same as the number from case 1 of the equivalence classes with $d=p$.
\qed

The function $w''$ is nonzero only when $p =3$; in that case, the outcome is $138/6 = 23$, which of course agrees with the entry for $I = 9$ in [ALPPT]. For $p > 3$, the expression simplifies (?) to
$$
\frac{p^4 + p^3 + 2p^2 + p + 1 + 4\cdot 3^{w'(d)}}6.
$$
I was relieved to find that for $p =5$ ($w'(d) = 0$), and $p=7$ ($w'(d) = 1$), this yields  $135$ and $477$ respectively, agreeing with the table entries for $I = 25$ and $49$.

Table 1 of [ALTPP] was particularly useful in checking examples in order  to see whether  the formulas were very likely correct! With other values of $d$ than those covered in \appbfou, there will be PH-equivalence classes that contain no terminal forms with $1$-block size $2$.

When $n =4$, formulas are still possible, but it would take a lot of {\it Sitzfleisch\/} to work out all the possible equivalence classes and their quantities.

The formulas simplify considerably when we consider only square-free choices for $d$; for example, the number of weakly terminal matrices with determinant fixed, is $\prod_{p|d}(\phi*J_2)(p) = \prod_{p|d} (p^2 + p -2) = \phi(d) d \prod_{p|d}(1+ 2/p)$. For $\pi$ a transposition, by Lemma \appcten, $\Cal S(\pi)(d) = \prod_{p|d} \Cal S(\pi)(p)  = d \prod_{p|d; \ {p\neq 2}} (1+1/p)$, and $\Cal S(123)(d) = 3^{w'(d)}$. Thus for {\it square-free\/} $d$,
$$
\Cal P \Cal H (3,d) = \frac{d\phi(d)\prod_{p|d}\( 1 + \frac 2p\) + 3d\prod_{p|d; \ {p\neq 2}} (1+1/p) + 2\cdot 3^{w'(d)}}{6}.
$$
(Recall  $w'(d)$ is the number of distinct prime divisors of $d$ that are congruent to $1$ modulo $3$.)

The middle term is $3\prod_{p|d} (1+p)$ if $d$ is odd and $2\prod_{p|d} (1+p)$ if $d$ is even. I tested the formula in Corollary \appbeig\ against Table 1 in [ALPPT] (recalling that their $I$ is our $d$) for values of $d = 30, 42, 70, 102, 105, 154, 165, 182, 186, 190, 195, 210$, as well as numerous choices of primes and products of two primes. Agreement was complete---so I am confident that the formula is  correct! [This is somewhat miraculous, as the formula is a sum of four formulas, each rather delicate.]

The first term is by far the largest, so the number is $6^{-1}\phi(d)\prod_{p|d} (p+2)\cdot (1 + \Oh{1/d})$. This is  the same as $(\phi * J_2(d))/6 $ for square-free $d$. This is also true if $ d$ is restricted to squares of primes (Proposition \appbsix).

Something rather startling occurs when we subtract from this the  number of PH-equivalence classes that contain a $1$-block size two matrix (the latter is the sum of the three numbers obtained from cases 1,2, and 3). Recall from section 7, the difference operator $\Delta$, defined by $\Delta f(x) = f(x+1) - f(x)$.

\Lem Proposition \appbsev. Let $d$ be a square-free integer. The number of PH-equivalence classes of $C \in \NS_3$ with $|\det C| = d$ and $C$ is not equivalent to a terminal form with $1$-block size two is
$$
 \frac{\phi(d) \Delta^3 f_d(-1)}{6},
$$
where $f_d (x) = \prod_{p|d} (x + p)$.

The factor $\phi(d)$ likely arises from an action of $\Z_d^*$ on the equivalence classes, presumably $(b,y) \mapsto (b,y)z$ as $z$ varies over $\Z_d^*$ (a similar phenomenon exists for the number obtained in case 3). The appearance of the third difference operator is rather mysterious. The dominant term in $\Delta^3 f_d(-1)$, at least when $\sum_{p|d} 1/p$ is large, is $\prod_{p|d} (p+2)$. We obtain that  if $d(m)$ is a sequence of square-free integers \st $\sum_{p|d(m)} 1/p \to \infty$ as $m \to \infty$, then
$$
\frac{\left|\brcs{\text{PH-equivalence classes of $C \in \NS_3$ with $|\det C| = d(m)$, no terminal form with $1$-block size two}}\right|}
{\left| \brcs{\text{PH-equivalence classes  $C \in \NS_3$, $|\det C| = d(m)$,  a terminal form  $1$-block size two}}\right |\cdot \prod_{p|d(m)} (1+ 1/p)}  \to \frac 13.
$$

If $d$ is a product of one or two primes, then $\Delta^3 f_d(-1) = 0$, consistent with Proposition \appbfiv.
If $d = pqr$, a product of three primes, then $\Delta^3 f_d(-1) = 6$, so the number of PH-equivalence classes not equivalent to a terminal form with $1$-block size two is $\phi(d)$, and in fact, the action of $\Z_d^*$ is just that of $\Z_d^*$ on itself. For example, with $d = 30$, we take
$$
C = \(\matrix 1 & 1 & 4 \\   0 & 2 & 5 \\ 0 & 0 & 15 \\
\endmatrix\); \qquad
C\op  \sim  \(\matrix 1 & 2 & 4 \\   0 & 3 & 5 \\ 0 & 0 & 10 \\
\endmatrix\):= D.
$$
Both are in terminal form, with $\Lt {J(C); J(C_{\Omega(i)}}  \iso  \Lt{\Z_{30}; \Z_3, \Z_5, \Z_2} \iso \Lt {J(D); J(D_{\Omega(i)}}$. Hence neither is PH-equivalent to a terminal form with $1$-block size two. The $8$ PH-equivalence classes of determinant $\pm30$ matrices in $\NS_3$ with no terminal form having $1$-block size two are obtained by multiplying the $(b,gy)^T$ truncated column, $(2,5)^T$, by the  integers relatively prime to $30$, that is, $1,7,11,13,17,19,23,29$ (that these are all primes is not entirely a coincidence), and then reducing modulo $15$.

More is true: $C\op$ is not PH-equivalent to $C$ (even though their invariants are identical). By calculating the ordered triples $(J(C_{\Omega(i)}))$ and $(J(D_{\Omega(i)}))$, we see that if $C$ were PH-equivalent to $D$, then the relevant permutation matrix $P$ would have to correspond to the transposition $(13)$. But a simple computation reveals that with this $P$, $D P C^{-1}$ has non-integer coefficients (specifically, the $(1,3)$ entry is $1/6$).

\SecT References

\Rf [ALPPT] {A Atanasov, C Lopez, A Perry, N Proudfoot, M Thaddeus}, Resolving toric varieties with Nash blow-ups. Experimental Mathematics 20 (2011) 288--303.

\Rf [BeH] S Bezuglyi \& D Handelman, Measures on Cantor sets\/{\rm:} the good, the ugly, the bad. Trans\ Amer Math Soc 366 (2014) 6247--6311.

\Rf [C] PM Cohn, Free ideal rings and localization in general rings. (2006) Cambridge University Press, 

\Rf [EHS] {EG Effros, David Handelman, \& Chao-Liang Shen}, Dimension groups and their affine representations. Amer J Math 102 (1980) 385--407.

\Rf [ES] {EG Effros \& Chao-Liang Shen}, Dimension groups and finite difference equations. J Operator Theory 2 (1979) 215--231.

\Rf [G] KR Goodearl, Partially ordered abelian groups with interpolation. Mathematical Surveys and Monographs, 20, American Mathematical Society, Providence RI, 1986.

\Rf [GH] KR Goodearl \& David Handelman, Metric completions of partially ordered abelian groups. Indiana Univ J Math 29 (1980) 861--895.

\Rf [Gr]  PA  Grillet, Directed colimits of free commutative semigroups. J Pure Appl Algebra 9  1 (1976) 73--87.

\Rf [GKKL]  RN Gupta{, A Khurana, D Khurana,  \& TY Lam}, Rings over which the transpose of every invertible matrix is invertible. J of Algebra  322 (2009) 1627--1636

\Rf [H] D Handelman, Free rank $n+1$ dense subgroups of $\text{\/\bf R}^{n}$ and their endomorphisms.  J Funct Anal  46  (1982), no\. 1, 1--27.

\Rf [H1] David Handelman,    Positive polynomials and product type actions of
compact groups. Mem Amer Math Soc  54  (1985),  320, xi+79 pp.

\Rf [H2] David Handelman,   Positive polynomials, convex integral polytopes, and a random walk problem. Lecture Notes in Mathematics, 1282, Springer--Verlag, Berlin, 1987, xii+136 pp.

\Rf [HW] GH   Hardy \& EM Wright, Theory of Numbers. likely a pirated edition.

\Rf [K] I Kaplansky,  Elementary divisors and modules. Trans Amer Math Soc, 66 (1949) 464Ð491

\Rf [La] TY Lam, Serre's conjecture. Lecture Notes in Mathematics  635 (1978) Berlin, New York; Springer-Verlag.

\Rf [L] G Landsberg, Uber eine Anzahlbestimmung und eine damit zusammenhangende Reihe. J Reine Angew Math 111 (1893) 87--88.

\Rf [Ma] G Maze, Natural density distribution of Hermite normal forms of integer matrices. J Number Theory 131 12 (2011) 2398--2408.

\Rf [MRW] {G Maze, J Rosenthal, \& U Wagner}, Natural density of rectangular unimodular integer matrices. Linear Algebra Appl 434 5 (2011) 1319--1324.

\Rf [M] P Moree, Counting carefree couples. http://arxiv.org/abs/math.NT/0510003 (2005).

\Rf [R] B Reznick, Lattice point simplices. Discrete Math 60 (1986), 219--242.

\Rf [R2] B Reznick, Clean lattice tetrahedra. http://de.arxiv.org/pdf/math/0606227.pdf
\Rf [TSCS] {C Torezzan, JE Strapasson, SIR Costa, RM Siquera},
Optimum commutative group codes.
Ar$\chi$iv:1205.4067v2 (2013)
\vskip 10pt

\SecT Constants' references

\noindent {\it carefree constant\/} http://oeis.org/A065463 [M]   $\prod_p (1-(2p-1)/p^3)$

\noindent {\it Landau's totient constant\/} http://oeis.org/A082695  $\zeta(2)\zeta(3)/\zeta(6) = \prod_p (1 + 1/p(p-1))$

\vskip 10pt \noindent Mathematics Department, University of Ottawa, Ottawa ON  K1N 6N5, Canada; dehsg\@uottawa.ca \& droy\@uottawa.ca

\end

Self-dual matrices, those in $\NS_n$ for which $B^{op}$ is PH-equivalent to $B$, ought to be of interest. Among other things, self-duality forces $B$ and $B^{op}$ to have the same elementary divisors. In the four-matrix example Example \thrthr\ (below), because $\det B = 30 = |I(B)|$ and $30$ is square-free, we have $\det B^{op} = 30$ and the elementary divisors are the same. If we take the leftmost matrix there, then we have
$$
B^{op} = \(\matrix  15 & 0 & 0 \\ 0 & 2 & 0 \\ -1 & -1 & 1 \\ \endmatrix\).
$$
It takes a few minutes to determine that
$$
B^{op} \text{ is PH-equivalent to both } \(\matrix  1 & 1 & 14 \\ 0 & 2 & 0 \\ 0 & 0 & 15 \\\endmatrix\) \text{ and } \(\matrix  1 & 0 & 15 \\ 0 & 1 & 14 \\ 0 & 0 & 30 \\ \endmatrix\).
$$
By Proposition \twoone, the latter is not PH-equivalent to any of the first three matrices in example \thrthr, and since the fourth is not PH-equivalent to a terminal matrix with size two $1$-block, it isn't equivalent to that one either. The list of three $(B^{op})_{\Omega(i)}$ are (using the rightmost matrix for computations)  $\Lt{\Z_{15}, \Z_2, 0}$, which is the same list as that of $B$.

As a warm-up, we find an easy example for which $I(B) \iso I(B')$ but they are distinguished by the families $(I(B_{\Omega(i)}))$ and $(I(B'_{\Omega(i)}{}))$. In this case, the absolute values of the determinants are unequal, so we already know that they cannot be PH-equivalent.

Set
$$
B = \( \matrix 1 & 0 & 1 \\ 0 & 1 & 4 \\ 0 & 0 & 8 \\\endmatrix\) \quad \text{and} \quad B' = \( \matrix 1 & 1 & 1 \\ 0 & 2 & 4 \\ 0 & 0 & 8 \\\endmatrix\).
$$
Both are in $\NS_3$ and in terminal form, with $\det B = 8$ and $\det B = 16$. It is easy to check that $I(B) \iso \Z_8 \oplus \Z_2$, and this also follows from Lemma \onetwo. To calculate $I(B')$, let $f_i$ be the rows, and $E_i$ the standard basis elements for $\Z^{1 \times 3}$. Then $E_1 = f_1 - f_2/2 + f_3/8$, so $x_1 = 8f_1 - 4f_2 + f_3$; also $E_2 = f_2 - f_3/2$, so $x_2 = 2f_2 - f_3$, and of course, $x_3 = f_3$. So it remains to find the elementary divisors (or the Smith normal form) of the matrix
$$
 \( \matrix 8 & 0 & 0\\ -4 & 2 & 0 \\ 1 & -1 & 1\\ \endmatrix \).
$$
The obvious column operations (which we are now allowed to do!) yield the Smith normal form of this matrix  as $\diag (8,2,1)$. Thus $I(B') \iso I(B)$.

For $i=1$, $I(B_{\brcs{2,3}})$ is obtained from the matrix $B$ with the first column removed, $\(\smallmatrix 0 & 1\\ 1 & 4 \\ 0 & 8 \\ \endsmallmatrix\)$; elementary row operations yield the $2 \times 2$ identity with third row consisting of zeros. So $I(B_{\brcs{2,3}}) = \brcs{0}$. The same is obviously true with $I(B_{
\brcs{1,2}})$. It is a routine calculation to see that $I(B_{\brcs{1,3}}) \iso \Z_4$. Hence the list $\Lt {I(B_{\Omega(i)})}$ is $\Lt{0,0, \Z_4}$.

On the other hand, none of the three choices for $I(B'_{\Omega(i)})$ are zero. With $i = 3$, the outcome is obviously $\Z_2$, with $i = 2$ and $i = 1$, we have the row reductions
$$
\(\matrix 1 & 1 \\ 0 & 4 \\ 0 & 8
\endmatrix\) \mapsto \(\matrix 1 & 1 \\ 0 & 4 \\ 0 & 0
\endmatrix\), \qquad \(\matrix 1 & 1 \\ 2 & 4 \\ 0 & 8
\endmatrix\) \mapsto \(\matrix 1 & 1 \\ 0 & 2 \\ 0 & 8
\endmatrix\)
\mapsto \(\matrix 1 & 1 \\ 0 & 2 \\ 0 & 0
\endmatrix\)
$$
so $I(B'_{\Omega(2)}) \iso \Z_4$ and $I(B'_{\Omega(1)}) \iso \Z_2$. Hence the list  $\Lt{I(B_{\Omega(i)})}$ is different from  $\Lt{I(B'_{\Omega(i)})} = \Lt{\Z_2, \Z_2, \Z_4}$.

Examples to distinguish them with equal determinant will appear after a few remarks.

*********************to here; then the next paragraph provides a better explanation of the earlier result (proved twice).

If when we perform this transformation, the resulting matrix has $1$-block size exceeding $r$, then it is easy to see that there exists an overset $\Omega_0 \supset \Omega$ \st $I(B_{\Omega_0}) = \brcs{0}$. There follows,

********************already done

\Lem Proposition \thrtwo. Suppose $B \in \NS_n$ and there exists $\Omega \subset \brcs{1,2,\dots, n}$ \st $I(B_{\Omega}) = \brcs{0}$, but for all supersets $\Omega_0 \supset \Omega$, $I(B_{\Omega_0}) $ is not zero. Then there exists a terminal matrix $B'$, PH-equivalent to $B$, having $1$-block size $|\Omega|$.

The converse is obvious.
*****************end already done

****the following should be streamlined or deleted

\Lem Examples \thrthr. (a) Four matrices in $\NS_3$, with the same elementary divisors ($\Lt{1,1,30}$), and in terminal form for which $I(B)$ are all isomorphic, but for which the lists, $\Lt{I(B_{\Omega(i)})}$ {\it are\/} distinct, and not all are PH-equivalent to a terminal form with $1$-block size two.  {\par}\noindent (b) A matrix $C \in \NS$ \st both $I(C)$ and $I(C\op)$ are cyclic yet  neither $C$ nor $C\op$  is PH-equivalent to a terminal form with $1$-block size $2$.

In the following, the calculations are very easy. Each of the four matrices has $\Lt{1,1,30}$ as its invariant factors.
$$
\matrix B & = & \(\matrix 1 &0& 15\\ 0 & 1 & 2\\ 0 & 0 & 30 \\
\endmatrix\) &\quad
\(\matrix 1 &0& 5\\ 0 & 1 & 6\\ 0 & 0 & 30 \\
\endmatrix\) &\quad
\(\matrix 1 &0& 3\\ 0 & 1 & 10\\ 0 & 0 & 30 \\
\endmatrix\) &\quad
\(\matrix 1 &1& 1\\ 0 & 2 & 5\\ 0 & 0 & 15 \\
\endmatrix\) &\\
I(B) & \iso & \Z_{30} & \Z_{30} &\Z_{30}&\Z_{30} &\\
\Lt{I(B_{\Omega(i)})}_{i=1}^3 & \iso & \Lt{\Z_{15}, \Z_2, 0} &
\Lt{\Z_{5}, \Z_6, 0} &
\Lt{\Z_{3}, \Z_{10}, 0} &
\Lt{\Z_{3}, \Z_{5}, \Z_2} &
\endmatrix$$
All four lists (viewed as unordered) are different (although their product groups are all $\Z_{30}$), so that there are no PH-equivalences between any of them. Moreover, the list for the fourth example contains no zero term; hence it is not PH-equivalent to a terminal matrix whose $1$-block size is two.

\noindent (b) Setting $C$ to be the opposite of the fourth example, we calculate that $\det C = 30 = \det C\op$; hence we obtain a matrix \st $I(B)$ and $I(C\op)$ are cyclic (since their orders are square-free), yet neither $B$ nor $C$ is PH-equivalent to a terminal form with $1$-block size $n-1$; compare this with Proposition \twotwo. A similar example is also discussed  in Appendix C, but with inequivalent $C$ and $C\op$ having the same lists.

\qed

Here is a more drastic family of examples, with minimal and maximal $1$-blocks.


Now for something completely different.

Let $n >k$ be integers, and let $F$ be a field (this can be modestly generalized
to rings of the form $\Z_d$; here the main application will be to $p$-element
fields, where $p$ is a prime. Define 
$$F(n,k) = \Set{M \in F^{n\times
k}}{\text{$\rk M = k$}}
$$ (we could require in addition that  all rows be nonzero
and  generate distinct subgroups of $F^{1\times k}$ instead).

Set $\Cal D_n$ to be the group of invertible diagonal $n \times n$ matrices with
entries from $F$, and let $\Cal P_n$ denote the group of $n \times n$
permutation matrices (regarded as elements of GL$(n,F)$). The group they
generate (consisting of weighted permutation matrices, will be denoted $W(n)$.
Then $W(n)$ acts from the left on $F(n,k)$ and $\gl (k,F)$ acts on the right,
yielding a $W(n) \times \gl(k,F)$ action.

\def\PP{\text{\bf P}}

As usual, define ${\PP}_{k-1}$ to be $F^k/F^\times$, that is, projective space.
Then we can interpret the $W(n) \times \gl(k,F)$ orbits as orbits of P$\gl(k,F)$
acting on lists of $n$ elements of $\PP_{k-1}$; if we had instead  defined
$F(n,k)$ so that the rows were nonzero generate distinct subgroups, then we
would have the action of P$\gl(k,F)$ acting on $n$-element subsets of
$\PP_{k-1}$.

More important for our examples is the following $W(n) \times \gl(k,F)$-stable
subset, 
$$Y(n,k) = \Set{M \in F^{n\times k}}{\text{all sets of $k$ rows of $M$
are linearly indepenent}}.$$

Denote by ${\Cal F}(n,k)$ and ${\Cal Y}(n,k)$ the respective orbit spaces of
$W(n) \times \gl (k,F)$ acting on $F(n,k)$, $Y(n,k)$. There is a type of duality
here.

We define a candidate map between orbit spaces $\Cal F(n,k) \to \Cal F(n,n-k)$,
as follows. Select $M \in F(n,k)$. Then $\Arrow M^T; F^{n\times 1}. F^{k\times
1}$ is onto (as $M$ has full rank). Thus its kernel has an ordered  basis $(v_1,
\dots, v_{n-k})$ with $n-k$ elements. Let $N$ be the $n \times (n-k)$ matrix.
Obviously $\rk N = n-k$, so $N \in F(n,n-k)$. If $(v_i')$
 is any other ordered
basis for $\ker M^T$, there exists $g \in \gl(n,k)$ \st $v_ig = v_i'$; hence all
selections of ordered bases yield $\gl(k,F)$-equivalent choices for $N$.

Now suppose that $M \sim M_0$, that is, there exist $P \in \Cal P_n$, $D \in
\Cal D_n$, and $h \in \gl (k,F)$ \st $PD Mh = M_0$. Then $M_0^T = h^T M^T D^T
P^{-1}$, and thus $\ker M_0^T = P(\Delta^T)^{-1} \ker M^T$. Then $N_0:=
(D^T)^{-1} P N$   belongs to $F(n,n-k)$, is $W(n) \times \gl(n-k,F)$-equivalent
to $N$, and $M_0^T N_0 = 0$.

It follows that the map on equivalence classes, $\Phi:[M] \mapsto [N]$ is
well-defined. Moreover, the condition $M^T N = 0$ determines $N$ up to
equivalence, and transposing, we obtain $N^T M $ is zero. The same conditions
apply, and we deduce that the map $[N] \to [M]$ is well-defined and the inverse
of $\Phi$.

The following is a little surprising, since the conditions for membership in
$Y(n,k)$ are different from those for membership in $Y(n,n-k)$ if $k \neq n-k$.
We do not assume any relation between $k$ and $n-k$.

\Lem Lemma. Suppose $M \in Y(n,k)$ and let $N \in F^{n \times (n-k)}$ satisfy
$M^T N = 0$ and $\rk N = n-k$. Then $N \in Y(n,n-k)$.

\Pf Suppose there is a set of $n-k$ rows of $N$ that is not linearly
independent. We can find a size $n$ permutation matrix $P$ so that this set
occupies the top $n-k$ rows of $PN$, which we write as $\(\smallmatrix X \\ Y \\
\endmatrix \)$, where $X$ is square of size $n-k$. Then $\rk X < n-k$, so there
exists $g \in \gl(n-k,F)$ \st the first column of $Xg$ is zero, and thus the top
$n-k$ entries of the first column of $PNg$ are all zero; the first column is $\(
\smallmatrix 0_{n-k} \\ c \\ \endmatrix \)$.

Form $M^T P^{-1} = \( \matrix Z & W\\ \endmatrix \)$  where $W$ is $k \times k$. The
rows of $W^T$ form a linearly independent set (since $M \in Y(n,k)$), hence $W$
is invertible. However,  $ \( \matrix Z & W\\ \endmatrix\) \( \smallmatrix 0_{n-k} \\
c \\ \endmatrix \) = Wc$; this is the first column of $M^T N g$, which is zero,
and thus $Wc$ is zero---forcing $c = 0$. But then $\rk PNg < n-k$, a
contradiction to $\rk N= n-k$.\qed

For the action of $W(n) \times \gl(k,F)$ on $F(n,k)$, we define the {\it
stabilizer\/} of a point, $M$, to be the subgroup of $\Cal P_n$ consisting of
$\Set{P \in \Cal P_n}{PDMg = Mj \text{ for some $D \in \Cal D_n$ and $g \in \gl
(k)$}}$. It is easy to check that this is a subgroup of $\Cal P_n$ (which we
sometimes identify with $S_n$), with the usual properties of stabilizers, e.g.,
if $M$ and $M'$ are in the same $W(n) \times \gl (k)$ orbit, then the stabilizer
of $M$ is isomorphic to that of $M'$ via an inner automorphism of $\Cal P_n$ (of
course, the word {\it inner\/} is only significant if $n=6$).

Most of the following is just what would be expected---except for part (b).

\Lem Proposition. Assume $n > k$.
\item{(a)} The map $\Arrow \Phi; \Cal F(n,k). \Cal F(n,n-k)$ is a bijection,
inducing an isomorphism between isomorphism classes of stabilizers.
\item{(b)} The restriction of $\Phi$ to  $\Cal  Y(n,k)$ is a bijection  with
$\Cal Y(n,n-k)$.

\Pf (a) We have already seen that $\Phi$ is a bijection. Select $M$ in $F(n,k)$,
and construct  $N \in F^{n,n-k}$ \st $M^T N $ is zero, as in the definition of
$\Phi$. Suppose $P $ belongs to the stabilizer of $M$; then $PDMg = M$ for
$(D,g) \in \Cal D \times \gl(k)$. Thus $M^T( D^T)^{-1} P N =  0$. Since the
columns of $N$ constitute a basis for $\ker M^T$, it follows there exists $h \in
\gl (n-k)$ \st $(D^T)^{-1}PN h = N$. Since we can write $(D^T)^{-1}P = P(P^{-1}
(D^T)^{-1} P)$ and the second factor is diagonal, we have that $P $ belongs to
the stabilizer of $N$. It follows by interchanging $M$ and $N$, that their
stabilizers are equal.

\noindent (b). This follows from xxx, and the fact that $Y(n,k)$ is stable \wrt
$W(n) \times \gl(n,k)$.
\qed

\Lem Corollary. Let $C \in \NS_{n,n-1}$, and let $B =\(\smallmatrix \I_{n-1} & a \\ 0 & d \\ \endsmallmatrix\)$ be a terminal form for $C$. The following are equivalent.
\item{(i)} $C\op \in \NS_{n,n-1}$;
\item{(ii)} $J(C\op)$ is cyclic;
\item{(iii)} $|\det C| = |\det C\op|$;
\item{(iv)} for all $1 \leq i,j \leq n-1$ and $i \neq j$, $\gcd\(d/(a_i,d), d/(a_j,d)\) = 1$.

\Pf (i) implies (ii) is trivial. (ii) implies (iii) is true  since $\Exp J(C) = \Exp J(C\op)$, and both $J(C)$ and $J(C\op)$ are cyclic, so must have the same order. (iii) implies (ii) for the same reason. 

From this point on, we may assume $C = B$.

(ii) iff (iv).  From xxx, $J(B\op) \iso \oplus_{i=1}^{n-1} \Z_{d/(a_i, d)}$. If this is cyclic, then $\gcd\(d/(a_i,d), d/(a_j,d)\) = 1$ (else the list of invariant factors will have more than one term).  The converse is trivial.

(ii,iii,iv) imply (i). Set $\Omega = S\setminus \brcs{n}$. By xxx, $J(B\op)_{\Omega})$ has order $|\prod d/(d,a_i)|/\lcm \brcs{d/(a_i,d)} =1 $. So $J(B\op)_{\Omega}) = \brcs{0}$, and since $|\Omega| =n-1$, $B\op$ has a terminal form with $1$-block of size $n-1$. 
\qed

******************all this can go 

Let $d$ be a positive integer, and let $d = \prod_{p\in U} p^{m(p)}$ be its prime factorization ($U$ is the set of its prime divisors, and $m(p)$ is the multiplicity). Pick a proper non-empty subset $T \subset U$, and define $d_1 = \prod_{p\in T} p^{m(p)}$, and $d_2 = \prod_{T^c} q^{m(q)} = d/d_1$ (elements of $T^c$ usually denoted $q$). Pick $v$ and $w$ with $(v,d_2) = (w,d_1) = 1$, and define the matrix 
$$
M(d; T,v,w) := \(\matrix  1 & 0 & d_1 v \\ 0 & 1 & d_2 w \\ 0 & 0 & d\endmatrix \).
$$
This belongs to $\NS_{3,2}$; it is in terminal form if $1 \leq v < d_2$ and $1 \leq w < d_1$. 

\Lem Proposition. Suppose $B \in \NS_{n,n-1}$ is PH-indecomposable and dual-compatible of determinant $d$. Then $n \leq 3$ and if $n=3$, then  $B$ is PH-equivalent to a matrix
$M(d; T,v,w) $ for some proper subset $T$ of the set of prime divisors of $d$.
Every such $B = M(d; T,v,w)$ is dual-compatible.

\Lem Proposition. Suppose $B \in \NS_{n,n-1}$ with $|\det B| = d$  is PH-indecomposable and $B$ is dual-conjugate, and $n \geq 3$. Then $ n = 3$ and  $B$ is PH-equivalent to the matrix of the form $M(T,v,w) $,
 $T$  is a proper nonempty  subset of $\Set{p \text{ prime}}{p|d}$,  and $v^2 d_1 + w^2 d_2 \equiv -1 \mod d$. Every such $B = M(d; T,v,w)$ is dual-conjugate.

The extra condition appearing in the characterization of dual-conjugacy, $v^2 d_1 + w^2 d_2 \equiv -1 \mod d$, limits considerably the choices for $d = |\det B|$. 

\Lem Corollary. For a positive integer $d$, and with $n \geq 3$, there exists indecomposable $B \in \NS_{n,n-1}$ with $|\det B| = d$ \st $B $ is dual-conjugate iff there exists a proper nonempty subset $T$ of the set of prime divisors of $d$ \st both
$$\eqalign{
-\prod_{p \in T; m(p) \text{ is odd}} p\quad& \text{is a square modulo $\prod_{q \in T^c} q^{m(q)}$}\cr  
-\prod_{q \in T^c; m(q) \text{ is odd}} q\quad& \text{is a square modulo $\prod_{p \in T} p^{m(p)}$}\cr 
}$$

These conditions can be substantially simplified; but we note that many integers $d$ fail to admit such a $T$ (e.g., any power of a prime; any number of the form $2^{k}$ times a square where $k \geq 1$; and many others). 

If we drop the condition that $B$ have a size $n-1$ $1$-block terminal form, then the corresponding results are conjecturally true. If $B \in \NS_n$ and $|\det B| = pq$ where $p,q$ are primes, then $B$ automatically belongs to $B_{n,n-1}$ (xxx), so all the results apply (and the condition on $d$ reduces to: (i) if $d$ is odd, both $p$ and $q$ must be congruent to $1$ modulo $4$ and $p$ is a square modulo $q$, while (ii) if $d = 2q$, then $-2$ is a square modulo $q$). 

When $n=2$, $B\in \NS_2$ is automatically dual-compatible, and it is an easy exercise to show that $B \sim \( \smallmatrix 1 & a \\ 0 & d \\ \endsmallmatrix\)$ (where $(a,d) =1$) is dual-conjugate iff $a^2 \equiv -1 \mod d$. Existence of such (for determinant $\pm d$) is equivalent to $-1$ being a square modulo $d$, which itself is equivalent to $m(2) \leq 1$ and all odd prime divisors of $d$ being congruent to $1$ module $4$. 

We first show that $M(d; T,v,w)$ has the requisite properties.
*******************up to here
******************soem of this computation below can still be used

Now we consider the question, given $d$, when there exist a pair of integers, $(v,w)$, \st $v^2 d_1 + w^2 d_2 \equiv -1 \mod d$, where $d_1 = \prod_{p \in T} p^{m(p)}$ and $d_2 = d/d_1$ for {\it some\/} proper nonempty subset $T$ of $U$. In other words, for whicd $d$ does there exist a PH-indecomposable dual conjugate matrix $B \in \NS_{n,n-1}$ of determinant $d$ for some $n \geq 3$.  (The case of $n =2$ is particularly easy: either $d$ is a product of primes all of which are congruent to $1$ modulo $4$, or $d$ is even and $d/2$ is odd with all prime divisors congruent to $1$ modulo $4$).

The condition is equivalent to $-d_1$ is a square modulo $d_2$ and $-d_2$ is a square modulo $d_1$ (examples showing that both these conditions are necessary, not just one of them, are ubiquitous). If we write $d_1 = \prod_T p^{m(p)}$ and $d_2 = \prod_{T^c} q^{m(q)}$, then necessary and sufficient are the following batch of conditions (using the Legendre symbol). 
$$\eqalign{
\Leg -1,p. \prod_{p \in T; \text{ odd } m(p) } \Leg p,q. = 1 & \text{ for all odd $q$ in $T^c$}\cr
\Leg -1,q. \prod_{q \in T^c;  \text{ odd } m(q) } \Leg q,p. = 1 & \text{ for all odd $p$ in $T$}\cr
-\prod_{p \in T; \text{ odd } m(p) } p & \text{ is a square modulo $2^{m(2)}$ if $2 \in T^c$}\cr
-\prod_{q \in T^c; \text{ odd } m(q) } q & \text{ is a square modulo $2^{m(2)}$ if $2 \in T$.}\cr
}$$
If $m(2) \leq 1$, the last two conditions are vacuous.

This yields hordes of example and non-examples. Say the  integer $d >1$ is {\it dual-admissible\/} if it satisfies these conditions. Then it is easy to check the following.

\item{(i)} If $d$ is a power of a prime, then $d$ is not dual-admissible.
\item{(ii)} If $d = 2p$ where $p$ is an odd prime, then $d$ is dual-admissible iff $\Leg-2,p. = 1$.
\item{(iii)} If $d = pq$ where $p$ and $q$ are distinct odd primes, then $d$ is dual-admissible iff both $p,q \equiv 1 \mod 4 $ and $\Leg p,q. = 1$.
\item{(iv)} If $d$ is a square, then $d$ is dual admissible iff $d$ is odd and all prime divisors are congruent to $1$ modulo $4$.
\item{(v)} If $d/2$ is a square, then $d$ is dual admissible iff $d/2$ is odd, there exists an odd prime divisor $p$ \st $\Leg -2,p. = 1$, and for all prime divisors $q$ \st $\Leg -2,q. = -1$, we must have $q \equiv 1 \mod 4$.
\item{(vi)} If $d/2^k$ is a square for some $k\geq 2$, then $d$ is not dual-admissible.
\item{} \dots and many other cases.

In some of these cases, any nonempty proper subset $T \subset U$ can be realized; in others, there may be only a single choice for $T$. The latter occurs when $d = 2\cdot 3\cdot 5$ and $3\cdot 5\cdot 7$. It is not difficult to show that if $d = pqr$ where all are odd primes, then $d$ dual-admissible implies that at least one of $p$, $q$, $r$ must be congruent to $1$ modulo $4$. 

In the case of (iii), the smallest choices for the primes are $(5,29)$ and $(13,17)$ (depending on the definition of {\it smallest\/}.

Of course, we have been assuming that $B$ has a terminal form with $1$-block  size $n-1$. It is conceivable that the results xxx hold without that assumption---that is, an indecomposable $B \in \NS_n$  that is dual-conjugate must be of size three or less, and the existence of one with determinant $d$ should be characterized in exactly the same way. (And then the size $3$ examples should be easy to deal with.) There is one special case wherein we can conclude this---if $d = pq$ is a product of two different primes, then automatically any $B\in \NS_n$ belongs to $\NS_{n,n-1}$ (xxx, later). 

Delete preceding.******************************


\Lem Lemma \thrfiv. Let $B \in \NS_n$ have a terminal form with a  $1$-block of size $n-1$, \st the entries above $d = |\det B|$ are each relatively prime to $d$. Then $J(B\op) = (\Z_d)^{n-1}$, and $J(B_{\Omega}) = 0$ for proper subsets $\Omega$ of $\brcs{1, 2, \dots, n}$.

This is elementary. More interesting is that the converse holds.

\Lem Lemma \thrsix. Suppose $B \in \NS_n$ and $J(B_{\Omega(i)})  = 0$ for all $i = 1,2,\dots, n$. Then $B$ is PH-equivalent to a terminal matrix of the form described in the previous lemma.

\Rmk Obviously $J(B_{\Omega(i)}) = 0$ for all $i$ entails $J(B_{\Omega}) = 0$ for all proper subsets $\Omega$, as the maps are onto.

\Rmk The same method yields that if $|\Omega| > 1$ and $n \in \Omega$, then $J(B_{\Omega}) = 0$ iff $\gcd\brcs{\brcs{a_j}_{ j \not\in\Omega} \cup \brcs{d}} = 1$.

\Pf Since there exists {\it one\/} $i$ \st $J(B_{\Omega(i)}) = 0$, $B$ is PH-equivalent to a terminal matrix of the form,
$$
\(\matrix \I_{n-1} & a \\ 0 & d \\ \endmatrix\)
$$
(necessarily $d = |\det B|$, since PH-equivalence preserves the absolute value of the determinant), where $a = (a_1, \dots, a_{n-1})^T$ consists of nonnegative integers less than $d$ and $\gcd\brcs{a,d} = 1$. Suppose $p$ is a prime dividing both $d$ and $a_j$ for some $j$ (this includes the possibility that $a_j = 0$). By conjugating with an obvious permutation matrix (that fixes the $n$th coordinate), we may assume $j = n-1$. Then the matrix relevant for calculating $J(B_{\Omega(n-1)})$ (delete the $n-1$st column, creating an $n \times (n-1)$ matrix) is
$$
\(\matrix \I_{n-2} & a'\\
0 & a_{n-1}\\
0 & d \\
\endmatrix\)
$$
where $a' = (a_1,\dots, a_{n-2})^T$. It is easy to check that the row space is of index $\gcd\brcs{d,a_{n-1}}$ in $\Z^{n-1}$ (alternatively: all $(n-1) \times (n-1)$ determinants are divisible by $p$). Hence $I(B_{\Omega (n-1)}) $ has order divisible by $p$, so is not zero. This contradicts the assumption, whence $\gcd\brcs{a_j,d} = 1$ for all $j$.
\qed

\Lem Proposition \thrfou. Suppose $B \in \NS_n$. 
If $B = \(\smallmatrix \I_{n-1} & a \\ 0 & d \\ \endsmallmatrix\)$ is a terminal form with $1$-block size $n-1$, where $a = (a_1, \dots, a_{n-1})^T$, then
$$
\left| \Set{j \leq n}{J(B_{\Omega(j)}) = 0}\right| = 1 + \left| \Set{i \leq n-1}{\gcd\brcs{a_i,d} = 1}\right|
$$

\Pf If $j =n$, then $B_{\Omega(n)}$ contains the $(n-1 ) \times (n-1)$ identity, so $J(B_{\Omega (n)}) = \brcs{0}$. If $ j < n$, then the same argument as in \thrsix\ yields that $J(B_{\Omega(j)}) = \brcs{0}$ iff $\gcd\brcs{a_j,d} = 1$.
\qed

\comment

There is a dual to Lemma \thrsix. If $J(B)$ and all $I(B_{\Omega})$ are isomorphic to $\Z_d$ for all subsets $\Omega$ with more than one element, then $B$ is PH-equivalent to a matrix of the form
$$
\( \matrix 1 & z_2 & z_3 & \dots & z_n \\
0 && d\I_{n-1}&&\\ \endmatrix\),
$$
where $\gcd\brcs{z_i,d} = 1$ for all $i$. The converse is true as well.

\endcomment


For each $B(X)$, define an explicit isomophism $J(B(X)) \to \Z_p^k  = \Z_p
\times \Z_p \times \dots \Z_p$ (with $k$ copies of $\Z_p$), sending (for $j = 1,
\dots, k$)  $E_{n-k+j} + r(B(X))$ to $e_j := (0,\dots,0,1,0,\dots,0)$, the $1$
appearing in the $j$th position. Since the kernel of $ J(B(X)) \to
J(B(X)_{\Omega(i)})$ is $E_i + r(B(X))$ (here $i \in \brcs{1,2,\dots, n}$), we
can identify the kernels with the following subgroups of $\Z_p^k$: for $i \leq
n-k$, $\langle -r_i (X)\rangle$ (the $i$th row of $X$, viewed as an element of
$\Z_p^k)$, since $E_i + (0^{n-k}, r_i(X)) \in r(B(X))$, and for $i > n-k$,
$e_{i-(n-k)}$.

Putting the generators into an $n \times (n-k)$ column with entries in $\Z_p$,
we obtain the matrix $M(X):= \(\smallmatrix -X \\ I_k \\ \endsmallmatrix \)$
(although the minus sign plays no role in terms of subgroups, it does play a
role when we work out the corresponding $\Phi$). The conditions on $X$ are
equivalent to $M(X) \in Y(n,k)$.

When the action of P$\gl(n,k)$ on $Y(n,k)$ is transitive, then set of subgroups
for $M(X)$ are obtainable from those of any other admissible $M(X')$---in
particular,  there is a permutation  of $\PP_{k-1}$ that induces a bijection
between the kernels associated to $M(X)_{\Omega(i)}$ and those for
$M(X')_{\Omega(i)}$. By xxx, $\Cal J(B(X)) \iso \Cal J(B(X')$. (There is no
claim that the allowable actions---on the right---preserve PH-equivalence;
simply that they preserve the $\Cal J$ isomorphism class.)

On the other hand, if the action is not transitive, then observing that every
equivalence class in $Y(n,k)$ contains a matrix of the form $\(\smallmatrix X \\
I_k \\ \endsmallmatrix\)$, we see that there are admissible choices for $X $ for which there is
no bijection between the kernels, and thus not all $\Cal J(B(X)_{\Omega}))$ are
isomorphic---even though the lists $\List{(J(B(X)_{\Omega})}_{|\Omega| = j}$ are
identical.

We have seen that if $p > 7$, then the P$\gl(k,\Z_p)$ action on $Y(5,2)$ not
transitive. we have required $Y(n,k)$ to be nonempty---this forces $p
+ 1 = |\PP_1| \geq 5$; this explains the absence of $2$ and $3$ in our
discussion.

There is still the question of whether $\Cal J(B)$ determines $\Cal J(B\op)$,
that is, if $\Cal J(B) \iso \Cal J(B')$, then $\Cal J(B\op) \iso \Cal
J(B'{}\op)$.
This is the case if $B \in \NS_{n,n-1}$ as follows from a combination of results
obtained above. It is also true here, because of the duality result.

It is easy to calculate $B(X)\op$.
$$
B(X)\op = \(\matrix pI_{n-k}& 0 \\ -X^T & I_k \endmatrix\), \quad{\text{which is
PH-equivalent to}}\quad \(\matrix I_{k}& 0 \\ -X^T &p I_{n-k} \endmatrix\).
$$
This is in the same form as $B(X)$ with $k$ and $n-k$ interchanged. The
corresponding column of generators of the kernels of $p_{\Omega(i)}$ is thus
given by $\( \smallmatrix -X^T \\ I_{n-k}\\ \endsmallmatrix\)$. Rearranging the
rows, we obtain $N(X) =\( \smallmatrix  I_{n-k} \\ -X^T \\\endsmallmatrix\) $.
As $M(X)^T = \(\smallmatrix -X^T & I_k \\ \endmatrix \)$, we have $M(X)^T N(X) =
0$. Hence \wrt the actions of $W(n)\times \gl (k,\Z_p)$ and $W(n) \times
\gl(n-k,\Z_p)$, we have $\Phi ([M(X)]) = [N(X)]$. It follows from xxxx that if
$\Cal J(B(X)) \iso \Cal J(B(X')$, then the same is true for their opposites.

We can extend the duality result somewhat, in order to cover slightly more
general situations. Insead of using primes, let $d$ be a positive integer, and
consider
$$
B_d(X) = \(\matrix I_{n-k}  & X \\ 0 & dI_k\endmatrix\),
$$
where  it is in terminal form (that is, the content of each column of $X$ is
relatively prime to $d$; but we require that the determinant of every square
submatrix of $X$ is relatively prime to $d$, in order to obtain the duality
result. This last is equivalent to $M(X)$ satisfying the (new) defining property
for $Y(n,k)$ over the ring $R = \Z_d$ (in place of a field), namely that every
subset of $k$ rows yield a determinant that is relatively prime to $d$, that is,
invertible (when we regard the entries of $X$ modulo $d$). Then the duality
result holds.

Let $R$ be a commutative ring satisfying the following two conditions.
\item{(i)} If $E,F,H$ are finitely generated $R$-modules with the last two free
and $E \oplus F = H$, then $E$ is free (stably free implies free).
\item{(ii)} If $A$ is a maximal ideal of $R$, then every unit of $R/A$ is the
image of a unit of $R$.

A lot of commutative rings satisfy (i), but (ii) is  restrictive ($\Z$ obviously
fails to satisfy it). The rings $\Z_d$ satisfy both conditions.

Then define $F(n,k) = \Set{M \in R^{n \times k}}{M^T (R^{n\times 1}) =
R^{k\times 1}}$ (this condition is implied by the existence of a set of $k$ rows
whose determinant is a unit in $R$, and is equivalent to the ideal generated by
all the determinants of size $d$ being improper), and $Y(n,k) $ consisting of
elements of $F(n,k)$ all of whose $k$-sets of rows have determinant a unit in
$R$. Then the proof of everything in and preceding xxx carries over.  This can
be used to show the opposite conjecture is satisfied for matrices in $B_d(X)$.

\Lem Corollary \onenin. Suppose that $B \in \NS_n$ and $J(B\op)$ is cyclic of prime
power order, $q^a$.
Then
$$
B \text{ is PH-equivalent to }
\( \matrix 1 && & & c(k+1) & \dots & c(n) \\
 & 1&& & 0 &  &  \\
 &  & \ddots &&&   &  \\
&&&1 &&& \\
&&&& q^{t(k+1)} &&\\
&&&&&\ddots & \\
&&&&&&q^{t(n)}
\endmatrix \),
$$
where $\gcd\brcs{c(i),q} = 1$ for each $i > k$, and $t(j)$ is increasing. In addition,
$$
B\op \text{ is PH-equivalent to }
\( \matrix \I_{n-1} & C \\
0 & q^a \\
\endmatrix \),
$$
where $C = (a_1, a_2, \dots, a_{n-1})^T$ with $k = 1 + |\Set{i}{a_i
\equiv 0 \mod q}|$.

\Rmk The condition that the order be of prime power (in addition to being cyclic) is essential. Setting $B$ to be  opposite of  example \thrthr\ (fourth) in section~3 [??where is this??xxx], we have $|J(B\op)| = 30$, which is square-free. Thus $J(B\op)$ is cyclic; but $J(B)$ is not PH-equivalent to a terminal form with $1$-block size $n-1$.

\Pf
Say $J(B\op) \iso \Z_{q^a}$ for some prime $a$ and prime $q$. We work with $B\op$, abbreviated $B'$. Then $\Z^{1\times n}/r(B') \iso J(B\op) \iso\Z_{q^a}$. Obviously  $\det B' = q^a$, so that  any terminal form for $B$ can have only powers of $q$ along the diagonal. Let $B''$ be a terminal form for $B'$; then
$$
B'' = \(\matrix \I_k & X \\ 0 & \Cal D \endmatrix\) \quad \text{where} \quad
\Cal D = \(\matrix q^{a(k+1)} \\ & q^{a(k+2)} && * \\ &&\ddots &&
\\
&&& q^{a(n)}\\ \endmatrix \),
$$
with $1 \leq a(k+1) \leq \dots \leq a(n)$ and all the strictly upper triangular entries of $\Cal D$ are either zero or divisible by a power of $p$, but that power is less than the power appearing in the diagonal entry for that column. Obviously, we may assume that $k < n-1$.

For any cyclic abelian group of prime power order, any pair of subgroups is comparable (that is, one of them is a subgroup of the other). Let $G_i$ be the subgroup of $\Z^{1\times n}/r(B'')$ generated by $E_i$. Then either $G_{n-1} \subset G_n$ or vice versa. Thus we have either
$$\eqalign{
E_n & = \sum_{j=1}^n f_j c(j) + \lambda E_{n-1} \quad \text{or}\cr
E_{n-1} & = \sum_{j=1}^n f_j c(j) + \lambda E_{n}, \cr
}$$
where $c(j), \lambda \in \Z$. Examining the first $k$ coordinates, we immediately see that in either case, $c_j = 0$ for all $1 \leq j \leq k$. In the first case, the $n$th coordinate yields $1 = \sum_{j > k} (f_j,E_n) c_j$ (where $(\ ,\ )$ denotes the usual inner product, that is, it is the $n$th coordinate of $f_j$). However, all the entries in the final column of $\Cal D$ are divisible by $q$ or zero, yielding a contradiction. In the second case, the same argument applies to the $n-1$st coordinate, again a contradiction. Hence $k = n-1$, so that the terminal form, $B''$, has $1$-block size $n-1$.

Now $B$ is PH-equivalent to ${B''}\op$, and the latter  has the form given in the statement, as is easy to check.
\qed

Here we conjecture two equalities relating the determinants of $B$ and $B\op$, and the orders of $J(B_{\Omega(i)})$ and $J((B\op)_{\Omega(i)})$. We prove them in the case that either $B$ or $B\op$
has a terminal form with $1$-block size $n-1$. (This is not such a strong restriction---for $n\geq 6$, the density of matrices $B\in \NS_n$ having a terminal form with $1$-block size $n-1$ exceeds $.80$; it is unknown what the density is for the property that either $B$ or $B\op$ has such a terminal form, but it is probably in excess of $.99$.)

\Lem Proposition \nfousev. Suppose $B \in \NS_n$ either has a terminal form of $1$-block size $n-1$ or $B\op$ does. Then 
$$\eqalign{
\frac{\prod_{i= 1}^n |J((B\op)_{\Omega(i)})| }{\prod_{i=1}^n |J(B_{\Omega(i)})|}& =\(\frac {| J(B)|}{|J(B\op)|}\)^n\cr
\(\prod_{i= 1}^n |J((B\op)_{\Omega(i)})| \){\prod_{i=1}^n |J(B_{\Omega(i)})|} & = \( {| J(B)|}{|J(B\op)|}\)^{n-2}.
}$$

\Pf Without loss of generality, we may assume that $B$ has a size $n-1$ $1$-block, and then that $B$ is already in terminal form with a $1$-block size of $n-1$, with determinant $d$; as we have seen,
$$
B= \(\matrix \I_{n-1} & a\\ 0 & d\\ \endmatrix \) \qquad B\op \sim \(\matrix 1  & a_1'   & \dots & a_{n-1}' \\ 0 && \diag (d_1,d_2, \dots, d_{n-1})\\ \endmatrix \).
$$
where $a^T = (a_1, a_2, \dots, a_{n-1})$, $d_i = d/\gcd(d,a_i)$, and $a_i'  = a_i/\gcd(d,a_i)$. Obviously $|J(B)| = d$ and $|J(B\op)| = \prod_{i=1}^{n-1} d_i$. Also obvious is that if $i \neq n$, then $J(B_{\Omega(i)}) \iso \Z/(d,a_i)$, and  $J(B_{\Omega(n)}) = \brcs{0}$. In particular, 
$$
\prod_{i=1}^n |J(B_{\Omega(i)})| = \prod_{i=1}^{n-1}(d,a_i).
$$ 

If $i \neq 1$, then the obvious column operation (permissible here, since we are calculating the Smith normal form of $(B\op)_{\Omega(i)}$) yields $|J((B\op)_{\Omega(i)})| = \prod_{j\neq i} d_i$. On the other hand, by \nfouthr(ii), we see that $(a_1', \dots, a_{n-1}')$ has image order $\lcm (d_i) = d$, and thus $J((B\op)_{\Omega(1)}) \iso \(\oplus \Z/d_i Z\)/x\Z_d$, a quotient by a free $\Z_d$-module. Hence 
$$
\prod_{i= 1}^n |J((B\op)_{\Omega(i)})|  = \frac 1 d\prod_{i=1}^{n-1} d_i^{n-2} \cdot \prod_{i=1}^n d_i = \frac{\prod_{i=1}^{n-1} d_i^{n-1}}{d}.
$$

Combining the two displayed formulas, 
$$\eqalign{
\frac{\prod_{i= 1}^n |J((B\op)_{\Omega(i)})| }{\prod_{i=1}^n |J(B_{\Omega(i)})|}& = \frac{\prod_{i=1}^{n-1} d_i^{n-1}}{d \prod_{i=1}^{n-1}(d,a_i)}\cr
& = \frac{d^{(n-1)^2 - 1}}{\prod_{i=1}^{n-1}(d,a_i)^{n-1 + 1}}\cr
& = \( \frac{d^{n-2}}{\prod (d,a_i)}\)^{n}\cr
& = \(\frac{\prod d_i}{d} \)^n  = \(\frac {| \det B\op|}{|\det B|}\)^n.\cr
}$$

Combining them differently, we obtain 
$$\eqalign{
\(\prod_{i= 1}^n |J((B\op)_{\Omega(i)})| \){\prod_{i=1}^n |J(B_{\Omega(i)})|}& = d^{-1}\prod_{i=1}^{n-1} d_i^{n-1}(d,a_i)\cr
& =\( \prod d_i^{n-2}\)d^{n-2} = \( {| \det B\op|}{|\det B|}\)^{n-2}.\cr
}$$
\qed

************Counter-example here*************